\documentclass[12pt]{amsart}

\usepackage{amsmath}
\usepackage{amsbsy}
\usepackage{amssymb}
\usepackage{amscd}
\usepackage{eq}
\usepackage[dvips]{graphicx, color}

\newcommand{\bbc}{{\mathbb C}}

\newcommand{\bbp}{{\mathbb P}}
\newcommand{\bbq}{{\mathbb Q}}
\newcommand{\bbr}{{\mathbb R}}

\newcommand{\bbz}{{\mathbb Z}}

\newcommand{\al}{{\alpha}}

\newcommand{\be}{{\beta}}

\newcommand{\gam}{{\gamma}}

\newcommand{\Gam}{{\Gamma}}

\newcommand{\Del}{{\Delta}}
\newcommand{\ep}{{\epsilon}}

\newcommand{\lam}{{\lambda}}
\newcommand{\Lam}{{\Lambda}}

\newcommand{\sig}{{\sigma}}

\newcommand{\om}{{\omega}}

\newcommand{\gB}{{\mathfrak B}}

\newcommand{\g }{{\mathfrak g}}

\newcommand{\gh}{{\mathfrak h}}

\newcommand{\gs}{{\mathfrak s}}
\newcommand{\gS}{{\mathfrak S}}
\newcommand{\gt}{{\mathfrak t}}

\newcommand{\ch}{{\operatorname{ch}}\,}
\newcommand{\aut}{{\operatorname{Aut}}\,}

\newcommand{\aff}{{\operatorname {Aff}}}
\newcommand{\kernel}{{\operatorname {Ker}}}

\newcommand{\trace}{{\operatorname {Tr}}}

\newcommand{\n}{{\operatorname {N}}}
\newcommand{\m}{{\operatorname {M}}}
\newcommand{\h}{{\operatorname {H}}}

\newcommand{\im}{{\operatorname {Im}}}
\newcommand{\re}{{\operatorname {Re}}}
\newcommand{\gal}{{\operatorname{Gal}}}
\newcommand{\gl}{{\operatorname{GL}}}

\newcommand{\spl}{{\operatorname{SL}}}

\newcommand{\strict}{{\operatorname{st}}}
\newcommand{\sst}{{\operatorname{ss}}}
\newcommand{\pfaff}{{\operatorname{Pfaff}}}
\newcommand{\sep}{{\operatorname{sep}}}

\newcommand{\Adj}{{\operatorname{Ad}}}
\newcommand{\id}{{\operatorname{id}}}

\newcommand{\sym}{{\operatorname{Sym}}}

\newcommand{\Hom}{{\operatorname{Hom}}}

\newcommand{\rep}{representation}

\newcommand{\pv}{prehomogeneous vector space}

\newcommand{\Z}{\bbz}
\newcommand{\Q}{\bbq}
\newcommand{\R}{\bbr}
\newcommand{\C}{\bbc}
\newcommand{\p}{\bbp}

\newcommand{\mk}{k^{\times}}

\newcommand{\rg}{G_{k}}

\newcommand{\rv}{V_{k}}

\newcommand{\kableadd}%
{Department of Mathematics\\ Cornell University\\
Ithaca NY 14853}

\newcommand{\sub}{\subset}

\newcommand{\ti}{\widetilde}

\newcommand{\ccd}{,\ldots,}

\newcommand{\lan}{\langle}
\newcommand{\ran}{\rangle}

\makeatletter
\def\varddots{\mathinner{
\mkern1mu%
 \raise\p@\hbox{.}\mkern2mu%
 \raise4\p@\hbox{.}\mkern2mu%
 \raise7\p@\vbox{\kern7\p@\hbox{.}}%
\mkern1mu}}
\makeatother

\theoremstyle{plain}
\newtheorem{thm}{Theorem}[section]
\newtheorem{lem}[thm]{Lemma}
\newtheorem{cor}[thm]{Corollary}
\newtheorem{prop}[thm]{Proposition}

\theoremstyle{definition}

\newtheorem{defn}[thm]{Definition}
\newtheorem{assump}[thm]{Assumption}
\newtheorem{cond}[thm]{Condition}

\theoremstyle{remark}
\newtheorem{rem}[thm]{Remark}        

\usepackage{mathrsfs,bbm}
\usepackage[all]{xy}
\usepackage{xcolor}

\newcommand{\st}{\mathrm{st}}
\newcommand{\weyl}{\mathbb W}

\newcommand{\Span}{\operatorname{Span}}

\newcommand{\coorde}{{\mathbbm e}}

\newcommand{\diag}{{\mathrm{diag}}}
\newcommand{\size}{\tiny}
\newcommand{\Ex}{\mathrm{Ex}}
\newcommand{\ot}{\otimes }
\newcommand{\op}{\oplus }

\newcommand{\Prg}{\mathrm{Prg}}
\newcommand{\Gunit}{\mathrm{Gunit}}

\newcommand{\bbma}{\mathbbm a}
\newcommand{\bbme}{\mathbbm e}

\newcommand{\bbmp}{\mathbbm p}
\newcommand{\bbmq}{\mathbbm q}

\newcommand{\soct}{\ti{\mathbbm O}}

\newcommand{\shortstrut}{\rule[-5pt]{0cm}{15pt}}
\newcommand{\sinsert}{\,\shortstrut}

\begin{document}

\address[K. Tajima]
{National Institute of Technology (KOSEN), Fukushima College, 30 Nagao, Kamiarakawa, Taira, Iwaki, Fukushima, 970-8034, Japan}
\email{tajima@fukushima-nct.ac.jp}

\address[A. Yukie]
{Mathematical Institute, Tohoku Univerrsity, 
Sendai 980-8578, Japan}
\email{akihiko.yukie.a2@tohoku.ac.jp}
\thanks{The second author was partially supported by 
Grant-in-Aid (C) (23K03052)\\}


\title{On the GIT stratification of prehomogeneous vector spaces IV}
\author{Kazuaki Tajima}
\author{Akihiko Yukie}
\maketitle

\begin{abstract}
We determine all orbits of the \pv{} 
$G=\gl_8,V=\wedge^3\aff^8 $ 
rationally over an arbitrary 
perfect field of characteristic not equal to $2$ 
in this paper.  
\end{abstract}

\section{Introduction}
\label{sec:introduction}

This is the final part of a series of four papers. 
In Part I, we determined the set $\gB$ of vectors 
which parametrizes the GIT stratification \cite{tajima-yukie} 
of the four \pv s (1)--(4) in \cite{tajima-yukie-GIT1}. 
Even tough the set $\gB$ was determined, the corresponding stratum 
$S_{\be}$ may be the empty set. In Part II 
(\cite{tajima-yukie-GIT2}) and Part III 
(\cite{tajima-yukie-GIT3}), we determined 
which stratum $S_{\be}$ is non-empty for the \pv s (1), (2) and (3) in 
\cite{tajima-yukie-GIT1}. 

In this part, we consider the following \pv:
\begin{equation}
\label{eq:PV}
G=\gl_1\times \gl_8,\;
V=\wedge^3 \aff^8
\end{equation}
where $(t,g)\in \gl_1\times \gl_8$ acts on $V$ by 
$V\ni x\mapsto t (\wedge^3 \! g) x$, which is the 
scalar multiplication by $t$ to $(\wedge^3 g) x$. 
For a general introduction to this series of papers, 
see the introduction of \cite{tajima-yukie-GIT1}.

Throughout this paper, $k$ is a fixed field. 
We shall assume that $k$ is a perfect field in 
the main theorem and in Sections  \ref{sec:non-empty}, 
\ref{sec:empty-strata}. In that case,  
the algebraic closure $\overline k$ coincides with 
the separable closure $k^{\sep}$. 
If $X,Y$ are schemes, algebraic groups, etc., over $k$ then 
$X,Y$ are said to be {\it $k$-forms} of each other if 
$X\times_k k^{\sep}\cong Y\times_k k^{\sep}$. 

Let $Q_0(v)\in \sym^2\aff^4$ be the quadratic form 
$v_1v_4-v_2v_3$. 
%
\begin{defn}
\label{defn:Ex-defn}
\begin{itemize}
\item[(1)]
$\Ex_n(k)$ is the set of conjugacy classes of 
homomorphisms from $\gal(\overline k/k)$ to $\gS_n$.
\item[(2)]
$\Prg_2(k)$ is the set of $k$-isomorphism classes of 
all the $k$-forms of $\mathrm{PGL}_2$. 
\item[(3)]
$\mathrm{QF}_4(k)$ 
is the set of $k$-isomorphism classes of 
all the algebraic groups over $k$ of the form $\text{GO}(Q)^{\circ}$ 
where $Q\in \sym^2 k^4$. Let 
$\mathrm{IQF}_4(k)\sub \mathrm{QF}_4(k)$ be the subset  
consisting of all the inner forms of $\text{GO}(Q_0)^{\circ}$. 
\item[(4)]
$\text{G}_2(k)$ is the set of $k$-isomorphism classes 
of connected simple groups of type $\text{G}_2$ over $k$ 
with trivial centers.    
\item[(5)]
Let $\mathrm{Oct}(k)$ be the set of $k$-isomorphism classes
of $k$-forms of the split octonion $\soct$.  
\item[(6)]
$\mathrm{Gunit}_2(k)$ is a set which has a surjective map 
$\gam: \mathrm{Gunit}_2(k)\to \Ex_2(k)$ such that 
$\gam^{-1}(k)$ is a single point and if $F/k$ is a quadratic
extension then $\gam^{-1}(F)$ is the set of 
equivalence classes of $2\times 2$-Hermitian matrices 
over $F$ under the action of $k^{\times}\times \gl_2(F)$. 
\end{itemize}
\end{defn}

Note that if $n=2,3$ then $\Ex_n(k)$ coincides with 
the set of $k$-isomorphism classes of 
separable extensions of $k$ of degree up to $n$. 
The set $\Gunit_2(k)$ will be described more precisely 
in Theorem \ref{thm:case4-rational-orbits}.

The following theorem is our main theorems in this part.

\begin{thm}
\label{thm:main}
Let $k$ be a perfect field.  
\begin{itemize}
\item[(1)] For the \pv{} (\ref{eq:PV}), 
there are $21$ non-empty strata $S_{\be}$. 
\item[(2)]
Suppose that $\ch(k)\not=2$. 
If $S_{\be}\not=\emptyset$ then $G_k\backslash S_{\be\,k}$ 
is either (i) a single point (abbreviated as SP from now on) 
(ii) $\Ex_2(k)$ (iii) $\Ex_3(k)$ 
(iv) $\Prg_2(k)$ (v) $\mathrm{IQF}_4(k)$ (vi) $\mathrm{G}_2(k)$ 
or (vii) $\Gunit_2(k)$. 
Moreover the number of $S_{\be}$'s 
for (i)--(vii) are as follows.

\vskip 15pt
\hskip -15pt
\begin{tabular}{|c|c|c|c|c|c|c|c|}
\hline
\rule[-8pt]{-0.2cm}{22pt}
\upshape Type & $\;$ \hskip 5pt \upshape SP \hskip 5pt $\;$ 
& $\Ex_2(k)$ & $\Ex_3(k)$ & $\Prg_2(k)$ & $\mathrm{IQF}_4(k)$ &
$\mathrm{G}_2(k)$ & $\mathrm{Gunit}_2(k)$\\
\hline
\rule[-8pt]{-0.2cm}{22pt}
\upshape $\# S_{\be}$ & \upshape $12$ & \upshape $3$ 
& \upshape $2$ & \upshape $1$ & \upshape $1$ & \upshape $1$ & $1$ \\
\hline
\end{tabular}

\end{itemize}
\end{thm}
\vskip 10pt

As we pointed out in Part I, the orbit decomposition of 
the present \pv{} is known over $\C$ (see \cite{ozekib}). 
Our approach answers rationality questions and provides 
the inductive structure of orbits rationally over $k$.  
Sato--Kimura (\cite{saki} over $\C$) and Igusa (\cite[p.275]{igusac}) 
proved that the set of generic rational orbits of the \pv{} (\ref{eq:PV})
is in bijective correspondence with the set of 
$k$-isomorphism classes of 
$k$-forms of $\spl_3$ (also see \cite[p.127]{yukiej}). 

In the process of proving Theorem \ref{thm:main}, we are forced to
determine the set of generic rational orbits of the following 
vector space $(G,V)$. 
\begin{equation*}
G = \gl_4\times \gl_3\times \gl_1, \; 
V= \wedge^2\aff^3\oplus \wedge^2 \aff^4\otimes \aff^3.
\end{equation*}
This case is very close to the case (4b) of 
\cite[p.453]{kimura-2simple-typeI-inv}.  
We also determine the set of rational orbits for some 
easy cases and point out that the assumptions on $\ch(k)$ 
for some cases in \cite{wryu}, \cite{yukiej}, \cite{yukiek}
can be weakened.  

The organization of this paper is as follows.
In Section \ref{sec:notation}, we discuss the notation used 
in this Part IV. We review the notation related to the 
GIT stratification in Section \ref{sec:notation-related-git}. 
The set of rational orbits of some \pv s will be reviewed in 
Section \ref{sec:orbits1} and the restriction on $\ch(k)$ is removed
for some cases.  The consideration of such orbits will be 
necessary in the process of the consideration of rational
orbits in unstable strata $S_{\be}$.  
Section \ref{sec:orbits2} concerns a certain \pv{} which 
requires some work to determine the set of rational orbits. 
We prove that there are precisely $21$ strata $S_{\be}$ 
which are non-empty in Section \ref{sec:non-empty}. 
It will be proved in Section \ref{sec:empty-strata} that 
the strata $S_{\be}$ other than those in Section 
\ref{sec:non-empty} are the empty set.

\section{Notation}
\label{sec:notation}

We discuss notations used in this part.   

When we quote statements from previous parts, we write
them as Theorem I--2.1, Lemma II--4.4.  

We often have to refer to a set consisting of a single point. 
We write SP for such a set. If $V$ is a vector space then 
we denote the dual space of $V$ by $V^*$. 

Let $\gl_n$ (resp. $\spl_n$) be the general linear group 
(resp. special linear group), 
$\m_{n,m}$ the space of $n\times m$ matrices and $\m_n=\m_{n,n}$. 
If $A=(a_{ij})\in \m_n$ then $\trace(A)=a_{11}+\cdots+a_{nn}$. 
For $1\leqq i\leqq n,1\leqq j\leqq m$, 
$E_{ij}$ is the matrix whose only non-zero entry 
is the $(i,j)$-entry $1$.  $E_{ij}$ is called a matrix unit.  
We express the set of $k$-rational points by 
$\gl_n(k)$, etc. We sometimes use 
the notation $[x_1\ccd x_n]$ to express column vectors 
to save space.  We denote the unit matrix of dimension $n$ 
by $I_n$.  We use the notation $\diag(g_1\ccd g_m)$ 
for the block diagonal matrix whose diagonal blocks are 
$g_1\ccd g_m$.  If $t_1\ccd t_n\in \gl_1$, we may also use the 
notation 
\begin{equation*}
a_n(t_1\ccd t_n) = \diag(t_1\ccd t_n).  
\end{equation*}
For $u=(u_{ij})\in \aff^{n(n-1)/2}$ ($1\leqq j<i\leqq n$), 
let $n_n(u)$ be the lower triangular matrix whose diagonal entries are $1$ 
and the $(i,j)$-entry is $u_{ij}$ for $i>j$. 
If $v_1\ccd v_m$ are elements of a vector space $V$ then 
$\lan v_1\ccd v_m\ran$ is the subspace spanned by 
$v_1\ccd v_m$. We may use the notation 
$\text{Span}\{v_1\ccd v_m\}$ also.  
For the rest of this paper, 
tensor products are always over $k$.

We use parabolic subgroups of $\gl_{8}$ which consist 
of lower triangular  blocks. 
Let $j_0=0<j_1<\cdots<j_{N}=8$. 
We use the notation $P_{[j_{1}\ccd j_{N-1}]}$ 
(resp. $M_{[j_{1}\ccd j_{N-1}]}$)
for the parabolic subgroup (resp. reductive subgroup) of $\gl_{8}$ 
in the form 
\begin{equation*}
\begin{pmatrix} 
P_{11} & 0 \cdots 0 & 0 \\
\vdots & \vdots & 
\begin{matrix}
0 \\ \vdots \\ 0
\end{matrix}  \\ 
P_{N1} & \cdots & P_{NN} 
\end{pmatrix}, \quad
\begin{pmatrix} 
M_{11} & 0 \cdots 0 & 0 \\
\begin{matrix}
0 \\ \vdots \\ 0
\end{matrix}
& \ddots & 
\begin{matrix}
0 \\ \vdots \\ 0
\end{matrix}  \\ 
0 & 0 \cdots 0 & M_{NN} 
\end{pmatrix}
\end{equation*}
where the size of $P_{kl},M_{kl}$ is 
$(j_{k}-j_{k-1})\times (j_{l}-j_{l-1})$.

If $N=1$ then we replace $[j_{1}\ccd j_{N-1}]$ by 
$\emptyset$. 
Let 
\begin{align*}
& M^{\st}_{[j_{1}\ccd j_{N-1}]} = \spl_{8}\cap
M_{[j_{1}\ccd j_{N-1}]}
\end{align*}
and $M^s_{[j_{1}\ccd j_{N-1}]}$ 
be the semi-simple part of 
$M_{[j_{1}\ccd j_{N-1}]}$.

We consider many \rep s of groups of the 
form $M_{[j_{1}\ccd j_{N-1}]}$
in later sections. 
We use notations such as 
\begin{equation}
\label{eq:Lam-defn}
\Lam^{m,i}_{[c,d]}.
\end{equation}
The meaning of this notation is that this is 
$\wedge^i \aff^m$ as a vector space where 
$\aff^m$ is the standard \rep{} of $\gl_m$ 
and the index $[c,d]$ means that 
the block from the $(c,c)$-entry to 
the $(d,d)$-entry of the $\gl_{8}$ 
of $M_{[j_{1}\ccd j_{N-1}]}$
is acting on this vector space. 
For example, if $N=3$ then 
$M_{[3,4]}$ consists of elements of the form 
$\diag(g_1,t,g_2)$ where $g_1\in\gl_3,g_2\in\gl_4,t\in \gl_1$. 
Then $\Lam^{4,1}_{[5,8]}$ is the standard \rep{} of $g_2\in \gl_4$
identified with the element $\diag (I_3,1,g_4)$.   
The trivial \rep{} of $M^s_{[j_{1}\ccd j_{N-1}]}$ is denoted by $1$. 

If $G$ is an algebraic group over $k$ with unit element $e$, 
the identity component of $G$ is denoted by $G^{\circ}$.  
If $F/k$ is an extension field, the group of $F$-rational 
points of $G$ is denoted by $G(F)$ or $G_F$. 
If $V$ is a vector space over $k$, the set of 
$F$-rational points of $V$ is denoted by $V_F$. 
Sometimes $V_F$ is identified with $V\otimes_k F$. 

If $\chi_1,\chi_2:G\to \gl_1$ 
are characters of $G$ then we say that $\chi_1,\chi_2$ are
{\it proportional} if there exist integers $a,b>0$ such that 
$\chi_1^a=\chi_2^b$.  
When we describe elements of the tangent space 
$\mathrm{T}_e(G)$, we use the ring of dual numbers $k[\ep]/(\ep^2)$. 
We abuse the notation and write $\ep$ for the class of $\ep$ 
modulo $\ep^2$ also.

\section{Notation related to the GIT stratification}
\label{sec:notation-related-git}

We explain notations related to the GIT stratification 
in this section.  
We basically follow the notation of \cite{tajima-yukie-GIT3} 
which is similar to that of \cite{tajima-yukie-GIT1}, except 
that we use the notations such as $G_{\st}$ instead of 
$G_1$ in \cite{tajima-yukie-GIT1}, \cite{tajima-yukie} 
because of many indices in this paper. 

Let $(G,V)$ be the \pv{} in (\ref{eq:PV}). We put  
\begin{equation}
\label{eq:T0-defn}
\begin{aligned}
G_{\st} & = \spl_8, \; 
T_{\st} = \{a_n(t_1\ccd t_8)\mid t_1\ccd t_8\in \gl_1,t_1\cdots t_8=1\}, \\
T_0 & =\{(t_0,t_1I_8)\mid t_0,t_1\in\gl_1\}.  
\end{aligned}
\end{equation}
$T_{\st}\sub G_{\st}$ is a maximal split torus. 
$T_0$ is the center of $G$, $G=T_0\cdot G_{\st}$ and 
$T_0\cap G_{\st}$ is finite.  
Let $X_*(T_{\st}),X^*(T_{\st})$ be the group of one parameter subgroups
and the character group respectively. We abbreviate 
`one parameter subgroup' as `1PS'.  

We put 
\begin{equation*}
\mathfrak {t}=X_*(T_{\text{st}})\otimes \R, \;
\mathfrak {t}_{\Q}=X_*(T_{\text{st}})\otimes \Q, \;
\mathfrak t^*=X^*(T_{\text{st}})\otimes \R, \;
\mathfrak t^*_{\Q}=X^*(T_{\text{st}})\otimes \Q.
\end{equation*}
We can express $\gt^*$ as follows:
\begin{equation*}
\gt^* = \left\{(a_1\ccd a_8)\in\R^8 
\,\vrule\, \sum_{i=1}^8 a_{i}=0 \right\}. 
\end{equation*}
Let 
\begin{equation*}
\gt^*_+ = \{(a_1\ccd a_8)\in \gt^*\mid a_1\leqq \cdots \leqq a_8\}.
\end{equation*}
This is the Weyl chamber with respect to the Weyl group 
$\weyl=\gS_8$. 

There is a natural pairing 
\begin{equation*}
\langle \;  ,\;  \rangle :X^*(T_{\text{st}})\times 
X_*(T_{\text{st}})\rightarrow \Z
\end{equation*}
defined by 
$t^{\langle \chi ,\lambda \rangle} =\chi (\lambda (t))$ for 
$\chi \in X^*(T_{\text{st}}),\lambda \in X_*(T_{\text{st}})$. 
This is a perfect paring (\cite[pp.113--115]{borelb}). 
There exists an inner product 
$(\;  , \; )$
on $\mathfrak{t}$ which is invariant under the actions of $\weyl$.
We may assume that this inner product is rational, i.e., 
$(\lambda,\nu)\in \Q$ for all 
$\lambda,\nu\in \mathfrak t_{\Q}$. 
Let $\|\; \|$ be the norm on $\mathfrak t$ defined 
by $(\;  , \; )$. 

For $\lambda\in\mathfrak t$, 
let $\beta =\beta(\lambda)$ be the element of
$\mathfrak t^*$ such that 
$\langle \beta, \nu\rangle = (\lambda, \nu)$
for all $\nu\in \mathfrak t$. 
The map $\lambda\mapsto \beta(\lambda)$ is a bijection 
and we denote the inverse map by 
$\lambda=\lambda(\beta)$. 
There is a unique positive rational
number $a$ such that $a\lambda(\beta)\in X_*(T_{\text{st}})$
and is indivisible. We use the notation 
$\lambda_{\beta}$ for $a\lambda(\beta)$. 
Identifying $\mathfrak t$ with $\mathfrak t^*$ by the map $\lam\mapsto \be(\lam)$, 
we have a $\weyl$-invariant inner product 
$(\;  , \;  )_{*}$ on $\mathfrak{t^*}$,

Let $\bbme_i$ be the coordinate vector of 
$\aff^8$ with respect to the $i$-th coordinate.
We put $e_{i_1i_2,i_3}=\coorde_{i_1}\wedge \coorde_{i_2}\wedge \coorde_{i_3}$
for $i_1,i_2,i_3=1\ccd 8$. The numbering 
used in \cite{tajima-yukie-GIT1} 
for (\ref{eq:PV}) is as follows. 

\vskip 10pt

\small 

\begin{center}
 
\begin{tabular}{|c|c|c|c|c|c|c|c|c|c|c|c|c|c|}
\hline
\rule[-3pt]{0cm}{15pt}
1 & 2 & 3 & 4 & 5 & 6 & 7 & 8 & 9 & 10& 11 & 12 & 13 & 14 \\
\hline
\sinsert $e_{123}$ & $e_{124}$ & $e_{125}$ & $e_{126}$ & $e_{127}$ 
& $e_{128}$ & $e_{134}$ & $e_{135}$ & $e_{136}$ 
& $e_{137}$ & $e_{138}$ & $e_{145}$ & $e_{146} $ & $e_{147}$  \\
\hline
\rule[-3pt]{0cm}{15pt}
15 & 16 & 17 & 18 & 19 & 20 & 21 & 22 & 23 & 24 & 25 & 26 & 27 &28  \\
\hline
\sinsert $e_{148}$ & $e_{156}$ & $e_{157}$ & $e_{158}$ & $e_{167}$ 
& $e_{168}$ & $e_{178}$ & $e_{234}$ & $e_{235}$ & $e_{236}$ 
& $e_{237}$ & $e_{238}$ & $e_{245}$ & $e_{246}$ \\
\hline
\rule[-3pt]{0cm}{15pt}
29 & 30 & 31 & 32 & 33 & 34 & 35 & 36 
& 37 & 38 & 39 & 40  & 41 & 42 \\
\hline
\sinsert $e_{247}$ & $e_{248}$ & $e_{256}$ & $e_{257}$ & $e_{258}$ 
& $e_{267}$ & $e_{268}$ & $e_{278}$ & $e_{345}$ & $e_{346}$ 
& $e_{347}$ & $e_{348}$ & $e_{356}$ & $e_{357}$ \\
\hline
\rule[-3pt]{0cm}{15pt}
43 & 44 & 45 & 46 & 47 & 48 & 49 & 50 
& 51 & 52 & 53 & 54 & 55 & 56 \\
\hline
\sinsert $e_{358}$ & $e_{367}$ & $e_{368}$ & $e_{378}$ & $e_{456}$ 
& $e_{457}$ & $e_{458}$ & $e_{467}$ & $e_{468}$ & $e_{478}$ 
& $e_{567}$ & $e_{568}$ & $e_{578}$ & $e_{678}$ 
\\
\hline
\end{tabular}

\end{center}

\normalsize

\vskip 10pt

We denote the $i$-the coordinate vector by $\bbma_i$. For example, 
$\bbma_1=e_{123},\bbma_{56}=e_{678}$.
Let $\gam_i\in\gt^*$ be the weight of $\bbma_i$ with respect to 
the action of $T_{\st}$. For example, 
$\gam_1= \tfrac 18(5,5,5,-3,-3,-3,-3,-3)$.  
If $I\sub \{\gam_1\ccd \gam_{56}\}$ is a finite set, 
we denote its convex hull by $\text{Conv}(I)$. If 
$0\notin \text{Conv}(I)$, let $\be(I)\in \gt^*_{\Q}$ be the 
closest point to the origin.  Let $\gB$ be the set of 
all such $\be(I)$ which belongs to $\gt^*_+$. 
This is the parametrizing set
for the GIT stratification.  

Suppose that $\be\in\gB$.  
$G$ acts on $T_{\st}$ by conjugation and so acts on 
$\gt^*$ also. Let $M_{\be}$ be the stabilizer of $\be$.  
For example, if $\be=\tfrac 1{56}(-5,-1,-1,-1,-1,3,3,3)$, 
then $M_{\be}=\gl_1\times M_{[1,5]}$. 
There is a minimum positive integer $a>0$ 
such that $a\be\in X^*(T_{\st})$ is induced from a character  
of $M_{\be}\cap \gl_8$. 
We denote this character by $\chi_{\be}$.  
We define
\begin{align*}
M^{\st}_{\be} & = (M_{\be}\cap G_{\st})^{\circ}, \;
G_{\st,\be} = \{g\in M^{\st}_{\be}\mid \chi_{\be}(g)=1\}^{\circ}, \\
P_{\be} & = \left  \{p\in G\; \Big |\; 
\lim_{t\rightarrow 0}\lambda_{\be} (t)p\lambda_{\be} (t)^{-1} 
\; \textrm{exists} \right \}, \\
U_{\be} & = \left  \{p\in G\; \Big |\; 
\lim_{t\rightarrow 0}\lambda_{\be} (t)p\lambda_{\be}(t)^{-1}= 1 \right \}, \\
& Y_{\beta }= \Span \{\bbma_i\,|\, (\gamma_i,\beta )_{*}
\geq (\beta ,\beta )_{*}\}, \\ 
& Z_{\beta }= \Span \{\bbma_i\,|\, (\gamma_i,\beta )_{*}
=(\beta ,\beta )_{*}\}, \\
& W_{\beta }=\Span \{ \bbma_i\,|\, (\gamma_i,\beta )_{*}
>(\beta ,\beta )_{*}\}.
\end{align*}
Clearly $Y_{\beta}=Z_{\beta}\oplus W_{\beta}$. 
$P_{\be}$ is a parabolic subgroup of 
$G$ (\cite[p.148]{Springer-LAG}) with Levi part $M_{\be}$ 
and unipotent radical $U_{\be}$.

Note that $T_0$ is a split torus acting on 
$V$ by scalar multiplication. 
We measure the stability of $V$ by $G_{\text{st}}$ 
ignoring the scalar multiplication. 
We ignore the direction of $\be$ 
from $M^{\st}_{\be}$ when we consider the stability 
of points in $Z_{\be}$. However, 
we consider the action of $M_{\be}$ when we consider 
the set of rational orbits.  

The group $G_{\st,\beta}$ acts on $Z_{\beta}$. 
Let $\p(V)^{\sst}$ (resp. $\mathbb P(Z_{\beta })^{\sst}$)
be the set of semi-stable points of $\p(V)$ (resp. 
$\mathbb P(Z_{\beta })$) with respect to 
the action of $G_{\text{st}}$ (resp. $G_{\text{st},\be}$). 
Since there is a difference between $V$ and $\p(V)$
(resp. $Z_{\beta}$ and $\mathbb P(Z_{\beta })$), 
we removed some scalar directions from 
$G$ (resp. $M^{\st}_{\beta}$) and considered 
stability with respect to 
$G_{\text{st}}$ (resp. $G_{\text{st},\be}$).  
For the notion of semi-stable points, see \cite{mufoki}. 
We regard $\mathbb P(Z_{\beta })^{\sst}$ as a subset of $\p(V)$. 

Put 
\begin{align*}
& V^{\sst} = \pi_V^{-1}(\p(V)^{\sst}),\; 
Z_{\beta }^{\sst}=\pi_V^{-1}(\mathbb P(Z_{\beta })^{\sst }), \;
Y_{\beta }^{\sst}=\{(z,w)\,|\, z\in Z_{\beta}^{\sst},w\in W_{\beta}\}. 
\end{align*}
We define $S_{\beta }=GY_{\beta }^{\sst }$. 
Note that $S_{\beta }$ can be the empty set. 
We denote the set of $k$-rational points of 
$S_{\beta}$, etc., by $S_{\beta \,k}$, etc.

The following theorem is 
\cite[p.264, Corollary 1.4]{tajima-yukie}. 
\begin{thm} 
\label{KKN} 
Suppose that $k$ is a perfect field. Then  we have 
\begin{align*} 
V_k\setminus \{0\} = V^{\sst } _k
\coprod \coprod_{\beta \in \mathfrak{B}} S_{\beta \, k}. 
\end{align*}
Moreover, 
$S_{\beta \,k}\cong G_{k}\times_{P_{\beta \,k}} Y_{\beta \,k}^{\sst }$. 
\end{thm}  

The set $\gB$ for the present case was computed 
by computer computations in Section I--9.

\section{Rational orbits (1) }
\label{sec:orbits1}

In this section and the next section, we consider
rational orbits of some \pv s which appear as 
$(M_{\be_i},Z_{\be_i})$. We consider
the action without the $\gl_1$-factor in some cases.  
We describe known cases of rational orbits of 
some \pv s in this section. 
Even though the results are known for all cases 
(see \cite{wryu}, \cite{yukiek}), we have to describe
typical points in $Z_{\be_i}^{\sst}$.  
Also we would like to remove the restriction on $\ch(k)$ in some cases 
as much as possible. For that purpose, 
we reprove some known results in 
\cite{wryu}, \cite{yukiek} in a slightly different manner using the
notion of regularity.

For the rest of this paper, 
$\{\bbmp_{n,1}\ccd \bbmp_{n,n}\}$ is the standard basis of $\aff^n$. 
We put $p_{n,i_1\cdots i_m}=\bbmp_{n, i_1}\wedge \cdots\wedge \bbmp_{n, i_m}$
for $1\leqq i_1\ccd i_m\leqq n$. 
We use these notations to distinguish such bases in the consideration 
of individual strata 
from the notations $\{\bbme_1\ccd \bbme_8\}$, 
$e_{i_1i_2i_3}$ we use for bases of $\aff^8$ and $V$.

When we have to consider more than one standard \rep s of 
$\gl_n$'s, we may use different letters such as 
$\mathbbm q_{n,i},\mathbbm r_{n,i}$ to avoid confusion. 
In that case, we use the notation such as  
$q_{n,i_1\cdots i_m}$, $r_{n,i_1\cdots i_m}$. 
We identify $\wedge^n \aff^n$ with $\aff^1$ 
so that $p_{n,12\cdots n},q_{n,12\cdots n}$, etc., correspond to $1$.

We frequently use the following proposition to show that 
some \rep s are regular \pv s. For the definition of 
the regularity of \pv s, see Definition III--4.5. 
If $(G,V)$ is a regular \pv, $w\in V$ and $U=G\cdot w\sub V$ 
is Zariski open then $U_{k^{\sep}}=G_{k^{\sep}} w$ 
by Proposition III--4.4. 
Let $G$ be an algebraic group, $V$ a finite dimensional 
\rep{} of $G$ over $k$, $w\in V$, $e\in G$ the unit element and 
$\mathrm{T}_e(G_w)$ the tangent space of $G_w$ at $e$.  

\begin{prop}
\label{prop:dimensions-compatible}
Let $G$ be an algebraic group, $V$ a finite dimensional \rep{} 
both defined over $k$ and $e\in G$ the unit element. 
If $\dim \mathrm{T}_e(G_w)=\dim G-\dim V$ then 
$G\cdot w\sub V$ is Zariski open and the group scheme 
$G_w$ is smooth over $k$. 
\end{prop}

We briefly review the cohomological argument 
to determine the set of rational orbits. 

Suppose that $G$ is a product of general linear groups, 
$(G,V)$ is a regular \pv, $w\in V$, $G\cdot w\sub V$ 
is Zariski open, $n>0$, $H\sub G_{w\,k}$ is a subgroup isomorphic
to $\gS_n$ and $G_w/G_w^{\circ}\cong \gS_n$ where $H$ maps to 
$\gS_n$ isomorphically. Since $\Gam=\gal(k^{\sep}/k)$ acts on $\gS_n$ 
trivially, $\h^1(k,\gS_n)$ can be identified with the set of 
conjugacy classes of all the anti-homomorphisms from $\Gam$ to $\gS_n$. 
By considering the extensions of $k$ corresponding to the kernels of 
anti-homomorphisms, we have a map 
\begin{equation}
\h^1(k,\gS_n) \to \Ex_n(k)
\end{equation}
and this map is bijective if $n=2,3$. 
Therefore, if $n=2,3$, the map 
$\h^1(k,G_w)\to \h^1(k,\gS_n)$
can be identified with a map 
\begin{equation*}
\rg \backslash V^{\sst}_k \to \Ex_n(k). 
\end{equation*}
We denote this map by $\gam_V$.  

Let $U=G\cdot w$ and $x\in U_k$. Then $x\in G_{k^{\sep}} w$
since $(G,V)$ is a regular \pv.   
\begin{equation*}
\xymatrix{
1\ar[r] & G_x^{\circ} \ar[r] & G_x \ar[r] & G_x/G_x^{\circ}\ar[r] & 1
}
\end{equation*}
is an exact sequence, i.e., $G_x^{\circ}\triangleleft G_x$,  
$G_x\to G_x/G_x^{\circ}$ is the natural homomorphism and 
$G_x/G_x^{\circ}$ is a $k$-form of $\gS_n$.  

The following lemma is \cite[p.120, Lemma (1.8)]{yukiel}. 

\begin{lem}
\label{lem:Gx-exact}
{\upshape (1)} $(G_x/G_x^{\circ})_k$ acts on $\h^1(k,G_x^{\circ})$ 
and 
\begin{equation*}
\xymatrix{
1 \ar[r] & (G_x/G_x^{\circ})_k\backslash \h^1(k,G_x^{\circ})
\ar[r] & \h^1(k,G_x) \to \h^1(k,G_x/G_x^{\circ}) \ar[r] & 1
}
\end{equation*}
is exact, i.e., the map 
$(G_x/G_x^{\circ})_k\backslash \h^1(k,G_x^{\circ})\to \h^1(k,G_x^{\circ})$
is injective and the inverse image of the trivial element of 
$\h^1(k,G_x/G_x^{\circ})$ is the image of 
$(G_x/G_x^{\circ})_k\backslash \h^1(k,G_x^{\circ})$. 

{\upshape (2)}
The trivial element of 
$\h^1(k,G_x)$ corresponds to 
$\gam_V(G_k\cdot x)\in \Ex_n(k)$ and 
\begin{equation*}
\gam_V^{-1}(\gam_V(G_k \cdot x)) 
\cong (G_x/G_x^{\circ})_k\backslash \h^1(k,G_x^{\circ})
\end{equation*}

{\upshape (3)}
If $c\in \h^1(k,G_x^{\circ})$, $g\in G_{k^{\sep}}$ 
and $c$ is represented by a 1-cocycle 
of the forms $\{g^{-1}g_{\sig}\}_{\sig\in\gal(k^{\sep}/k)}$, 
then the orbit corresponding to $c$ is $G_k\cdot gx$.  
\end{lem}

We use the following Witt's theorem often.

\begin{thm}
\label{thm:alternating-matrix}
Let $n>0$ be an integer, $G=\spl_n$ and $V=\wedge^2 \aff^n$. 
We identify $V$ with the space of alternating matrices with diagonal 
entries $0$ (this assumption is necessary if $\ch(k)=2$).  
\begin{itemize}
\item[(1)]
\setlength{\itemsep=2pt}
If the rank of $A\in V_k$ is $m$ then $m$ is even.  
\item[(2)]
Suppose that $m=2l$. There exists $g\in G_k$ such that 
if $B=(b_{ij})=gA{}^tg$ then $b_{ij}=0$ unless 
$(i,j)=(n-2l+1,n-2l+2),(n-2l+2,n-2l+1)\ccd (n-1,n),(n,n-1)$. 
In particular, if $n$ is odd then the first row and the first column 
are zero.  
\end{itemize}
\end{thm}

We now start the consideration of individual cases.

\subsection{Case I}

We consider the \rep s 
\begin{align}
\label{eq:332-tensor-standard}
G & = \gl_3^2\times \gl_2, \; V = \m_3 \otimes \aff^2, \\
\label{eq:62-wedge-tensor-standard}
G & = \gl_6\times \gl_2, \; V = \wedge^2 \aff^6\otimes \aff^2. 
\end{align}

These cases have been considered in \cite[\S\S 3,4]{wryu}.  
However, the assumption $\ch(k)\not=2,3$ were made in \cite[\S\S 3,4]{wryu}
and the purpose of this subsection is to point out that 
this assumption is not necessary.

The case (\ref{eq:332-tensor-standard}) has been 
considered in Section III--5 and so 
we will be brief.  We identify $V$ with the space of 
pairs $x=(x_1,x_2)$ of $3\times 3$ matrices. 
The element $(g_1,g_2)\in\gl_3^2$ acts on $\m_3$ 
by $\m_3\ni X \mapsto g_1 X {}^tg_2 \in \m_3$. 
If we express $V=\m_3\otimes \aff^2$, 
The element $(g_1,g_2,g_3)\in G$ acts on $V$ by 
$(g_1,g_2)\otimes g_3$.

In the case (\ref{eq:62-wedge-tensor-standard}), we identify $V$ 
with the space of pairs $x=(x_1,x_2)$ 
of $6\times 6$ alternating matrices with diagonal entries $0$.   
The element $(g_1,g_2)\in G$ acts on $V$ by 
$\wedge^2 g_1\otimes g_2$. If we identify elements of $\wedge^2 \aff^6$ 
with alternating matrices (with diagonal entries $0$), 
the action of $g_1$ is by 
$\wedge^2 \aff^6 \ni X \mapsto g_1 X \,{}^t\! g_1\in\wedge^2 \aff^6$.

Let 
\begin{equation}
\label{eq:case7-R-defn}
J = \begin{pmatrix}
0 & 1 \\ 
-1 & 0
\end{pmatrix}, \quad 
w = \begin{cases}
\left(
\begin{pmatrix}
1 & 0 & 0 \\ 0 & -1 & 0 \\ 0 & 0 & 0
\end{pmatrix}
+ \begin{pmatrix}
0 & 0 & 0 \\ 0 & 1 & 0 \\ 0 & 0 & -1
\end{pmatrix}\right) & (\ref{eq:332-tensor-standard}), \\[25pt]
\left(
\begin{pmatrix}
J & 0 & 0 \\ 0 & -J & 0 \\ 0 & 0 & 0 
\end{pmatrix}, 
\begin{pmatrix}
0 & 0 & 0 \\ 0 & J & 0 \\ 0 & 0 & -J 
\end{pmatrix}
\right)
& (\ref{eq:62-wedge-tensor-standard}).   
\end{cases}
\end{equation}
\begin{prop}
\label{prop:tangent-cubic-h} 
$\dim \mathrm{T}_e(G_w) = 4,10$ for the cases 
(\ref{eq:332-tensor-standard}), (\ref{eq:62-wedge-tensor-standard})
respectively.  
\end{prop}
\begin{proof}
The statement for the case (\ref{eq:332-tensor-standard}) 
has been verified in Lemma III--5.10. 
So we only verify this proposition for the case 
(\ref{eq:62-wedge-tensor-standard}).   
Let $A\in \m_6(k),B=(b_{ij})\in\m_2(k)$ 
and $X=(I_6+\ep A,I_2 + \ep B)$. We express $A$ 
in block form as $A=(A_{ij})$ where $A_{ij}\in\m_2(k)$ 
for $i,j=1,2,3$. 
Suppose that $Xw=w$. For $x=(x_1,x_2)\in V$ and variables 
$v_1,v_2$, let $F_x(v)=\pfaff(v_1x_1+v_2x_2)$ 
be the Pfaffian. There are two choices for the Pfaffian, 
but we are only interested in whether or not the value is non-zero.  
So the choice of the Pfaffian does not matter.

$F_x(v)$ is a binary cubic form of $v=(v_1,v_2)$ 
and $F_{(g_1,g_2)x}(v) = (\det g_1) F_x (vg_2)$ 
regarding $v$ as a row vector. Let $P(x)$ be the 
discriminant of the binary cubic form $F_x$. 
Then $P(x)$ is a homogeneous polynomial of degree
$12$ of $x\in V$ and 
\begin{equation}
\label{eq:62-relative-invariant}
P(gx) = (\det g_1)^4 (\det g_2)^6 P(x).  
\end{equation}

Since $F_w(v) = \pm v_1v_2(v_1-v_2)$, 
\begin{equation*}
(1+\ep \trace(A))(I_2+\ep B) v_1v_2(v_1-v_2)=v_1v_2(v_1-v_2). 
\end{equation*}
Explicitly, 
\begin{align*}
& \trace(A) v_1v_2(v_1-v_2) + 2(b_{11}v_1+b_{21}v_2)v_1v_2 
+ v_1^2(b_{12}v_1+b_{22}v_2) \\
& - (b_{11}v_1+b_{21})v_2^2 - 2v_1v_2(b_{12}v_1+b_{22}v_2) \\ 
& = b_{12}v_1^3 + * v_1^2v_2+*v_1v_2^2 - b_{21}v_2^3 =0.  
\end{align*}
Therefore, $b_{12}=b_{21}=0$. 

This implies that the condition $Xw=w$ is explicitly,  
\begin{align*}
& A\begin{pmatrix}
J & 0 & 0 \\ 0 & -J & 0 \\ 0 & 0 & 0 
\end{pmatrix}
+ \begin{pmatrix}
J & 0 & 0 \\ 0 & -J & 0 \\ 0 & 0 & 0 
\end{pmatrix} 
{}^t A
+ b_{11}\begin{pmatrix}
J & 0 & 0 \\ 0 & -J & 0 \\ 0 & 0 & 0 
\end{pmatrix}
= 0, \\
& A \begin{pmatrix}
0 & 0 & 0 \\ 0 & J & 0 \\ 0 & 0 & -J 
\end{pmatrix}
+ \begin{pmatrix}
0 & 0 & 0 \\ 0 & J & 0 \\ 0 & 0 & -J 
\end{pmatrix}
{}^t A
+ b_{22} \begin{pmatrix}
0 & 0 & 0 \\ 0 & J & 0 \\ 0 & 0 & -J 
\end{pmatrix}
=0, 
\end{align*}
which means that 
\begin{align*}
& \begin{pmatrix}
A_{11}J & -A_{12}J & 0 \\
A_{21}J & -A_{22}J & 0 \\
A_{31}J & -A_{32}J & 0 
\end{pmatrix}
+ \begin{pmatrix}
J \,{}^t\! A_{11} & J \,{}^t\! A_{12} & J \,{}^t\! A_{13} \\
-J \,{}^t\! A_{21} & -J \,{}^t\! A_{22} & -J \,{}^t\! A_{23} \\
0 & 0 & 0 
\end{pmatrix}
+ b_{11}\begin{pmatrix}
J & 0 & 0 \\ 0 & -J & 0 \\ 0 & 0 & 0
\end{pmatrix}
= 0, \\
& \begin{pmatrix}
0 & A_{12}J & -A_{13}J \\
0 & A_{22}J & -A_{23}J \\
0 & A_{32}J & -A_{33}J 
\end{pmatrix}
+ \begin{pmatrix}
0 & 0 & 0 \\
J \,{}^t\! A_{21} & J \,{}^t\! A_{22} & J \,{}^t\! A_{23} \\
-J \,{}^t\! A_{31} & -J \,{}^t\! A_{32} & -J \,{}^t\! A_{33} 
\end{pmatrix}
+ b_{22} \begin{pmatrix}
0 & 0 & 0 \\ 0 & J & 0 \\ 0 & 0 & -J 
\end{pmatrix}=0.   
\end{align*}
So $A_{ij}=0$ for $i\not=j$ and 
\begin{equation*}
\trace(A_{11}) + b_{11}
= \trace(A_{22}) + b_{11} = 0, \quad 
\trace(A_{22}) + b_{22}
= \trace(A_{33}) + b_{22} = 0.  
\end{equation*}
It is easy to show that the set of solutions to this 
system of linear equations is a vector space of dimension 
$10$.  
\end{proof}

Let $H\sub G$ be the following subgroups 
\begin{align*}
& \left\{
\left(
a_3(t_{11},t_{12},t_{13}), 
a_3(t_{21},t_{22},t_{23}),
a_2(t_{31},t_{32}) 
\right)
\, \vrule \, 
\begin{matrix}
{}^{\forall} i,j, t_{ij}\in\gl_1, \\
t_{11}t_{21}t_{31} =  t_{12}t_{22}t_{31}=1, \\
t_{12}t_{22}t_{32} =  t_{13}t_{23}t_{33}=1
\end{matrix}
\right\}, \\
& \left\{
\left(
\diag(g_{11},g_{12},g_{13}), 
t_2 I_2 
\right)
\, \vrule \, g_{1i}\in \gl_2,t_2\in \gl_1, 
\det(g_{1i})t_2 =1 \; (i=1,2,3)
\right\} 
\end{align*}
for the cases (\ref{eq:332-tensor-standard}), 
(\ref{eq:62-wedge-tensor-standard}) respectively.

It is easy to verify that $H\sub G_w^{\circ}$. 
Since the dimensions are the same $H= G_w^{\circ}$. 
Since $G_w^{\circ}$ is smooth and reductive, 
we obtain the following proposition. 

\begin{prop}
\label{prop:case7-regular}
$(G,V)$ is a regular \pv.   
\end{prop}

The above proposition implies that 
$V^{\sst}_{k^{\sep}} = G_{\sep}\cdot w$. 
Therefore, it is possible to use the 
usual cohomological argument and there is a 
surjective map $\gam_V:\rg\backslash \rv^{\sst}\to \Ex_3(k)$.  
Suppose that $x\in \rv^{\sst}$ and $\gam_V(x)$ is the splitting 
field of $k'/k$ which is a separable quadratic or cubic extension of $k$.  
Similarly as in \cite[pp.306,309, Propositions 3.7,4.7]{wryu}, 
\begin{equation*}
G_x^{\circ} \cong 
\gl_1^2 \times \text{R}_{k'/k} \gl_1, \; 
\gl_1 \times \text{R}_{k'/k} \gl_1 
\end{equation*}
for the case (\ref{eq:332-tensor-standard}) 
according as $[k':k]=2,3$ and 
\begin{equation*}
G_x^{\circ} \cong 
\begin{cases}
\{(g_1,g_2)\in \gl_2\times \text{R}_{k'/k} \gl_2 
\mid \det g_1=\det g_2\} & [k':k]=2, \\ 
\{g\in \text{R}_{k'/k} \gl_2 \mid \det g\in \gl_1 \} 
& [k':k]=3 
\end{cases}
\end{equation*}
for the case (\ref{eq:62-wedge-tensor-standard}).  
So $\h^1(k,G_x^{\circ})=\{1\}$ for all the cases. 
Therefore, we obtain the following theorem 
without any restriction on $\ch(k)$. 

\begin{thm}
\label{thm:case7-rational-orbits}
$\gam_V$ induces a bijective map 
$\rg\backslash \rv^{\sst}\to \Ex_3(k)$. 
\end{thm}

\subsection{Case II}

We consider the \rep{} 
\begin{align}
\label{eq:42-wedge-tensor-standard}
G & = \gl_4\times \gl_2, \; V = \wedge^2 \aff^4\otimes \aff^2.
\end{align}

This case has been considered in \cite[\S 4]{wryu} 
and the assumption on $\ch(k)$ was removed in 
Subsection III--5.3. 
We summarize known properties of this case as follows.

Elements of $V$ can be expressed as 
pairs of $4\times 4$ alternating matrices 
$x=(x_1,x_2)$ with diagonal entries $0$.  
For variables $v=(v_1,v_2)$, 
the Pfaffian $F_x(v) = \pfaff(v_1x_1+v_2x_2)$ 
is a binary quadratic form. The 
discriminant $P(x)$ of $F_x(v)$ 
is a homogeneous polynomial of degree $4$ 
such that $P(gx)=(\det g_1)^2(\det g_2)^2P(x)$. 

Let $J$ be the matrix in (\ref{eq:case7-R-defn}) and 
\begin{equation*}
w = \left(
\begin{pmatrix}
J & 0 \\
0 & 0 \\
\end{pmatrix},
\begin{pmatrix}
0 & 0 \\
0 & J 
\end{pmatrix}
\right). 
\end{equation*}
\begin{prop}
\label{prop:42-regular}
\begin{itemize}
\item[(1)]
$(G,V)$ is a regular \pv{} and the orbit $G\cdot w\sub V^{\sst}$ 
is Zariski open. 
\item[(2)]
$\rg\backslash \rv^{\sst}$ is in bijective correspondence with 
$\Ex_2(k)$ by associating the field generated by a root of $F_x$ 
to $\rg\cdot x$.  
\end{itemize}
\end{prop}

\subsection{Case III}

We consider the \rep{} 
\begin{align}
\label{eq:7-wedge3}
G & = \gl_1\times \gl_7, \; V = \wedge^3 \aff^7.
\end{align}
The action of $\gl_7$ on $V$ is induced by the standard \rep{} of $\gl_7$ 
and $t\in \gl_1$ acts by the scalar multiplication by $t$.  
This case has been considered in \cite{igusac} and \cite{yukiek}. 
However, it was assumed that $\ch(k)=0$ in \cite{igusac} and \cite{yukiek}. 
The purpose of this section is to point out that 
the argument in \cite{yukiek} works for most part 
regardless of $\ch(k)$ and for the correspondence 
with isomorphism classes of $k$-forms of octonions 
if $\ch(k)\not=2$.  
Let $\text{G}_2(k),\mathrm{Oct}(k)$ 
be as in Definition \ref{defn:Ex-defn} (4), (5).

We denote $\bbmp_{7,i}$, $p_{7,i_1i_2i_3}$ by $\bbmp_i$, 
$p_{i_1i_2i_3}$ in this subsection.  
Let 
\begin{equation}
\label{eq:w-defn-case9}
w = p_{234} + p_{567} + p_{125} + p_{136} + p_{147} \in V. 
\end{equation}
We remind the reader that 
$p_{234}=\bbmp_{2}\wedge \bbmp_{3}\wedge \bbmp_{4}$, etc.  

We determine the tangent space $\mathrm{T}_e(G_w)$ of the stabilizer 
where $e=(1,I_7)$.  Let 
\begin{equation*}
A(a\ccd f,X) = 
\begin{pmatrix}
0  & 2d & 2e & 2f & 2a & 2b & 2c  \\
a &&&&  0  &  f  &  -e  \\
b  &&  X  &&  -f  &  0  &  d  \\ 
c  &&&&  e  &  -d  &  0  \\ 
d  &  0  &  -c  &  b  &&&&  \\ 
e  &  c  &  0  &  -a  
&&  -{}^t\! X  &&  \\
f  &  -b  &  a  &  0  &&&& 
\end{pmatrix} 
\end{equation*}
for $a\ccd f\in k,\; X\in \m_3(k)$ where $\trace(X)=0$. 
Let 
\begin{align*}
L_1 & = \{(0,A(a\ccd f,X))\mid a\ccd f\in k,X\in\m_3(k),\trace(X)=0\}, \\
L_2 & = \{(-3t,tI_7)\mid t\in k\}, \\
L & = \{(t,A)\mid t \in k,A\in\m_7(k),\; (e+\ep (t,A))w=w\}. 
\end{align*}

$\mathrm{T}_e(G_w)$ can be identified with $L$. 
By long but straightforward computations, one can verify that
$L_1\sub L$.  
We define a character $\chi_0:G\to \gl_1$ by 
\begin{equation}
\label{eq:chi0-defn}
\chi_0(t_0,g_1)=t_0^2(\det g_1)
\end{equation}
for $t_0\in\gl_1,g_1\in\gl_7$.  

\begin{prop}
\label{prop:case9-tangent-space}
\begin{itemize}
\item[(1)]
$\mathrm{T}_e(G_w)=L_1\oplus L_2$, $\dim L_1=14$, 
$\dim \mathrm{T}_e(G_w)=15$, $G_w$ is smooth over $k$, 
the dimension of $G_w$ as an algebraic group is $15$ 
and $G\cdot w\sub V$ is Zariski open. 
\item[(2)]
$\mathrm{T}_e(G_w\cap \spl_7) = \mathrm{T}_e(G_w \cap \kernel(\chi_0))=L_1$  
and 
$(G_w \cap \kernel(\chi_0))^{\circ}=(G_w\cap \spl_7)^{\circ}$.  
\end{itemize}
\end{prop}
\begin{proof}
(1) $\dim L_1=14$ is obvious.  
We show that $\dim L/L_1=1$. 
Suppose that 
$e+\ep(t,A)\in \mathrm{T}_e(G_w)$. 
By subtracting an element of the form  
$e+\ep (0,A(a\ccd f,X))$ where $a\ccd f\in k,X\in\m_3(k)$ 
and $\trace(X)=0$, we may assume 
\begin{equation}
\label{eq:dim15-assumption}
a_{21}=\cdots = a_{71}=a_{22}=a_{23}=a_{24}
=a_{32} =a_{33} =a_{34}=a_{42}=a_{43}=0. 
\end{equation}

By $I_7+\ep A$,   
\begin{align*}
& \bbmp_{1} \mapsto \bbmp_{1} + \ep a_{11} \bbmp_{1}, \\
& \bbmp_{2} \mapsto \bbmp_{2} 
+ \ep (a_{12} \bbmp_{1}+a_{52} \bbmp_{5}
+a_{62} \bbmp_{6}+a_{72} \bbmp_{7}), \\
& \bbmp_{3} \mapsto \bbmp_{3} + \ep (a_{13} \bbmp_{1}+a_{53} \bbmp_{5}
+a_{63} \bbmp_{6}+a_{73} \bbmp_{7}), \\
& \bbmp_{4} \mapsto \bbmp_{4} 
+ \ep (a_{14} \bbmp_{1}+a_{44}\bbmp_{4}+a_{54}\bbmp_{5}
+a_{64} \bbmp_{6}+a_{74}\bbmp_{7})
\end{align*}
and 
\begin{equation*}
\bbmp_{j} \mapsto \bbmp_{j} + \ep \sum_{i=1}^7 
a_{ij} \bbmp_{i} 
\end{equation*}
for $j=5,6,7$. 

Then $(e+\ep (t,A))w -w$ is $\ep$ times the following. 
\begin{align*}
& \quad (a_{14}-a_{26}+a_{35})p_{123}
+(-a_{13}-a_{27}+a_{45})p_{124}
+(t+a_{11}+a_{55})p_{125} \\
& \quad + a_{65} p_{126} 
+ a_{75}p_{127} + (a_{12}-a_{37}+a_{46})p_{134} 
+ a_{56}p_{135}+(t+a_{11}+a_{66})p_{136} \\
& \quad + a_{76}p_{137} + a_{57}p_{145}+ a_{67}p_{146}
+ (t+a_{11}+a_{44}+a_{77})p_{147} + (a_{17}+a_{53}-a_{62})p_{156} \\
& \quad + (-a_{16}+a_{54}-a_{72})p_{157} 
+ (a_{15}+a_{64}-a_{73})p_{167} 
+ (t+a_{44}) p_{234} + a_{54} p_{235} \\
& \quad + (a_{21}+a_{64})p_{236} + a_{74}p_{237}
- a_{53} p_{245} - a_{63}p_{246} -a_{73} p_{247} + a_{27} p_{256} \\
& \quad -a_{26}p_{257} + a_{25}p_{267} 
+ a_{52}p_{345} + a_{62}p_{346} + a_{72}p_{347}
+ a_{37} p_{356} -a_{36}p_{357} \\
& \quad  + a_{35}p_{367}
+ a_{47}p_{456} -a_{46} p_{457} 
+ a_{45}p_{467}+(t+a_{55}+a_{66}+a_{77})p_{567}. 
\end{align*}
If this is $0$, $a_{ij}=0$ for $i\not=j$ and 
\begin{align*}
& a_{44}+t=0, \; t+a_{55} + a_{66} + a_{77}=0, \;
t+ a_{11}+a_{55}=0, \\
& t+ a_{11}+a_{66}=0, \; 
t+a_{11}+a_{44}+a_{77}=0. 
\end{align*}
This implies that 
\begin{equation}
\label{eq:dim14-step}
t = 3a_{77}, \; a_{11} = -a_{77}, \; a_{44} = -3a_{77}, \;
a_{55} = -2a_{77}, \; a_{66} = -2a_{77}. 
\end{equation}
Therefore, $\dim L = 14 + 1 = 15$. 

Since elements of $L$ satisfy (\ref{eq:dim14-step}) 
modulo elements of the form $(0,A(a\ccd f,X))$, 
if $(t,A)\in L$ then $t=3s$ for some $s\in k$. 
Then $(t,A)-(3s,sI_7)\in L_1$, which implies 
that $L=L_1+L_2$.  If $(-3t,tI_7)\in L_1$, 
then $t=0$ by considering the $(1,1)$-entry. 
Therefore, $L=L_1\oplus L_2$. 

Since 
\begin{equation*}
\dim L= 15 = \dim (\gl_1\times \gl_7)-\dim V,
\end{equation*}
Proposition \ref{prop:dimensions-compatible} implies 
that $G\cdot w\sub V$ is Zariski open and $G_w$ is smooth 
over $k$. Therefore, the dimension of $G_w$ as an 
algebraic group is $15$. 

(2) 
Obviously, $G_w\cap \spl_7\sub G_w\cap \kernel(\chi_0)$, 
which implies that 
\begin{math}
\mathrm{T}_e(G_w\cap \spl_7) \sub \mathrm{T}_e(G_w\cap \kernel(\chi_0)). 
\end{math}
We prove  that $\mathrm{T}_e(G_w\cap \kernel(\chi_0))=L_1$. 
Suppose that that $e+\ep(t,A)\in \mathrm{T}_e(G_w\cap \kernel(\chi_0))=L_1$.
We assume (\ref{eq:dim15-assumption}) and prove that 
$t=0,A=0$. 

The equations (\ref{eq:dim14-step}) are still satisfied. 
By assumption, 
\begin{equation*}
0 = 2t+\trace(A) = 6a_{77}-a_{77}-3a_{77}-2a_{77}-2a_{77}+a_{77}
= -a_{77}. 
\end{equation*}
Therefore, $a_{77}=0$, which implies that $t=0$ and $A=0$. 
This implies that 
\begin{equation*}
\mathrm{T}_e(G_w\cap \kernel(\chi_0)) = L_1
\end{equation*}
and $\dim \mathrm{T}_e(G_w\cap \kernel(\chi_0)) =14$. 
Since $L_1\sub \mathrm{T}_{I_7}(\spl_7)$, 
$L_1\sub \mathrm{T}_e(G_w\cap \spl_7)$, which implies 
that $\mathrm{T}_e(G_w\cap \spl_7)=L_1$ also.  

Let $W=V\oplus \aff^1$ and define an action of 
$g=(t_0,g_1)\in G$ on $W$ by 
$W\ni (x,a)\mapsto (g\cdot x,(\det g_1)a)\in W$. 
Then $G_{(w,1)}=G_w \cap \spl_7$ and 
\begin{equation*}
\mathrm{T}_e(G_{(w,1)})=\mathrm{T}_e(G_w)\cap \mathrm{T}_{I_7}(\spl_7)
= \mathrm{T}_e(G_w\cap \spl_7)=L_1. 
\end{equation*}
Since $\dim L_1 = 14 = \dim (\gl_1\times \gl_7)-\dim W$, 
$G_w\cap \spl_7$ is smooth over $k$ and its dimension 
as an algebraic group is $14$. 
Since $G_w\cap \kernel(\chi_0)$ is defined by a single 
equation on $G_w$, the dimension of $G_w\cap \kernel(\chi_0)$ 
is $14$.  Therefore, 
$(G_w\cap \kernel(\chi_0))^{\circ}=(G_w\cap \spl_7)^{\circ}$. 
\end{proof}

We show that $G_w$ is reductive and $(G_w\cap \spl_7)^{\circ}$ 
is a simple group of type $\text{G}_2$. 
For that purpose, we remind the reader a criterion 
for a connected group $G$
to be reductive or semi-simple. Since $G$ is smooth over $k$, 
$G$ is connected if and only if $G\times_k \overline k$ is connected.  
So we assume that $k=\overline k$. 

Let $G$ be a connected algebraic group,
and $T\sub G$ a subtorus.  We denote the Lie algebras of
$G,T$ by $\g,\gt$.  
Consider a weight decomposition of
$\g$ with respect to $T$ such as 
\begin{equation}  
\label{weight-decomposition}
\g = \gt\oplus \oplus_{\al\in \Phi}\g_{\al}
\end{equation}
where $\g_{\al}$ is the weight space with weight $\al$.  
If $r\in N_G(T)$ (the normalizer) and $\al\in\Phi$, 
then the element of $X^*(T)$ defined by 
$t\to \al(r^{-1}tr)$ is denoted by 
$r(\al)$.  Note that $\Adj(r)\g_{\al}
=\g_{r(\al)}$ and so $r(\al)\in\Phi$.
Suppose $G$ is imbedded in 
$\gl_n$ for some $n$ as a closed subgroup.
Let $N^+\sub \gl_n$ (resp. $N^-\sub \gl_n$)
the subgroup of upper triangular matrices 
(resp. lower triangular matrices) 
with diagonal entries $1$, and  
$\overline T\sub \gl_n$ the subgroup of diagonal
matrices.  Then  $B^+=\overline T N^+,\; 
B^-=\overline T N^-$ are Borel subgroups of $\gl_n$.

We consider the following assumption.
\begin{assump}  
\label{semi-simple-condition}
\begin{itemize}
\item[(1)] For all $\al\in \Phi$, $\dim \g_{\al}=1$.
\item[(2)] The set $\Phi$ does not contain $0$.  
\item[(3)] For any $\al\in\Phi$,
$[\g_{\al},\g_{-\al}]\not=0$.  
\item[(4)] For any $\al\in\Phi$,
there exists $r_{\al}\in N_G(T)$ such that
$r_{\al}(\al)=-\al$.  
\item[(5)] The relation
\begin{equation*}  
(\cap_{\al\in\Phi} \kernel(\al))^{\circ}
=\{e\}
\end{equation*}
holds.
\end{itemize}
\end{assump}
\begin{prop}  
\label{borel-condition}
If Assumption \ref{semi-simple-condition}(1)--(4)
are satisfied then $T$ is a maximal torus and
$G$ is reductive.  
If furthermore Assumption \ref{semi-simple-condition}(5)
is satisfied then $G$ is semi-simple. 
\end{prop}
\begin{proof}  Suppose that 
Assumption \ref{semi-simple-condition} (1)--(4)
are satisfied.  
Let $T\sub S$ be a maximal torus.  If 
$T\not=S$, the Lie algebra $\gs$ of $S$
strictly contains $\gt$.  Since 
$S\sub Z_G(T)$, the adjoint \rep\ of
$T$ on $\gs$ is trivial.  Therefore, 
$\gs$ has no non-zero weights of $T$.
Since $\Phi$ does not contain $0$, 
this is a contradiction.  Hence, $T$ is a maximal
torus.  

Suppose that the unipotent radical $R_u(G)$ is not trivial.  
This is a connected, solvable, unipotent, and normal subgroup
of $G$.  Let $D_1(G)\supset D_2(G)\supset \cdots$ be the
derived series of $G$.  
Choose $n>0$ so that $D_n(R_u(G))\not=\{e\}$
and $D_{n+1}(R_u(G))=\{e\}$.  Then 
$H = D_{n+1}(R_u(G))$ is a 
connected, abelian, unipotent, and normal subgroup
of $G$.  Let $\gh$ be the Lie algebra of $H$.  
Since $T$ is a torus, $T$ is conjugate to
a subgroup of $\overline T$.  So $\gt$ consists 
of semi-simple elements.  Since 
$H$ is unipotent and abelian, $H$ is conjugate
to a subgroup of $N^+$.  So $\gh$ consists
of nilpotent elements.  Therefore, $\gt\cap \gh=\{0\}$.

Since $\gh\not=\{0\}$, there exists $\al\in\Phi$
such that $\g_{\al}\sub\gh$.  Since $H$ is a 
normal subgroup of $G$, $\Adj(r_{\al})(\gh)=\gh$.
Since $r_{\al}(\al)=-\al$, this implies that 
$\g_{-\al}\sub\gh$.  Since $[\g_{\al},\g_{-\al}]\not=0$,
$\gh$ is not abelian, which is a contradiction.  
Hence, $R_u(G)=\{e\}$ and so $G$ is reductive.  

Suppose furthermore, 
Assumption \ref{semi-simple-condition} (5)
is satisfied.  Since $G$ is reductive, 
$R(G)=Z(G)^{\circ}$.  Since $Z(G)^{\circ}$
acts trivially on $\g$ by the adjoint \rep, 
\begin{equation*}
Z(G)^{\circ}\sub (\cap_{\al\in\Phi} \kernel(\al))^{\circ}
=\{e\}
.\end{equation*}
Hence, $R(G)=\{e\}$ and so $G$ is semi-simple. 
\end{proof}

If a group is reductive, its Lie algebra has a root system. 
The above proposition says showing the existence of a root
system proves that the group is reductive.  

Let 
\begin{equation*}
\ti T=\{(t^{-3},tI_7)\mid t\in\gl_1\} \cong \gl_1. 
\end{equation*}
Obviously $\ti T\sub G_w^{\circ}$ ($\ti T$ is connected)
and $\mathrm{T}_e(\ti T)=L_2$.  

\begin{prop}
\label{prop:case9-regular-pv}
\begin{itemize}
\item[(1)]
$G_w\cong (G_w\cap \spl_7)\times \ti T$ 
\item[(2)]
$G_w\cap \spl_7$ is a connected simple group of type $\text{G}_2$ 
with trivial center. 
\item[(3)]
$G_w$ is connected and reductive.  
\end{itemize}
\end{prop}
\begin{proof}
(1) Let $\phi:(G_w\cap \spl_7)\times \ti T\to G_w$ 
be the homomorphism which takes the product. 
If $g\in G_w$, put $\chi_0(g)=t$.  Since 
$\chi_0(t^{-3},tI_7)=t$ also, 
$\chi_0((t^{-3},tI_7)^{-1}g)=1$. 
This implies that 
$(t^{-3},tI_7)^{-1}g\in \kernel(\chi_0)=G_w\cap \spl_7$. 
Therefore, $\phi$ is surjective.  
If $(t^{-3},tI_7)\cap (G_w\cap \spl_7)$, 
then $\chi_0(t^{-3},tI_7)=t=1$. 
Therefore, $\phi$ is injective. 
Since the map of the tangent spaces is an 
isomorphism by Proposition \ref{prop:case9-tangent-space}, 
$\phi$ is an isomorphism of algebraic groups.

(2) We first prove that 
$(G_w\cap \spl_7)^{\circ}$ is a simple group 
of type $\text{G}_2$.    

Let
\begin{equation}
\label{eq:dt1t2-defn}
\begin{aligned}
d(t_1,t_2) & = 
\diag(1,t_1,t_2,(t_1t_2)^{-1},t_1^{-1},t_2^{-1},t_1t_2), \\
H(s_1,s_2) & = \diag(0,a_1,a_2,-a_1-a_2,-a_1,-a_2,a_1+a_2), \\
X & = \begin{pmatrix}
s_1 & u_1 & u_2 \\
v_1 & s_2 & u_3 \\
v_2 & v_3 & -s_1-s_2
\end{pmatrix}
\end{aligned}
\end{equation}
for $t_1,t_2\in\gl_2,s_1,s_2,a_1,a_2,u_1\ccd v_3\in k$ 
and $T=\{d(t_1,t_2)\mid t_1,t_2\in\gl_1\}$.  
Define characters $e_1,e_2$ of $T$ by 
\begin{equation*}
e_1(d(t_1,t_2))=t_1,\; e_2(a(t_1,t_2))=t_2.  
\end{equation*}

Let $X_a=A(1,0\ccd 0)$. We define 
$X_b\ccd X_f,X_{u_1}\ccd X_{v_3}$ similarly.  
Then these are weight vectors with respect to 
the adjoint action by $d(t_1,t_2)$.  
Let $\al_1=e_1-e_2,\al_2=e_2$. 
The set of weights are 
\begin{equation}
\label{eq:root-system-G2}
\Phi= \{\pm \al_1,\pm \al_2,\pm(\al_1+\al_2),\pm (\al_1+2\al_2),
\pm (\al_1+3\al_2),\pm (2\al_1+3\al_2)\}. 
\end{equation}
One can verify by computations that 
\begin{align*}
& [X_a,X_d]=H(2,-1),\; [X_b,X_e]=H(-1,2), \;
[X_c,X_f]=H(-1,2), \\
& [X_{u_1},X_{v_1}]=H(1,-1),\; [X_{u_2},X_{v_2}]=H(1,0), \;
[X_{u_3},X_{v_3}]=H(0,1). 
\end{align*}
So Assumption \ref{semi-simple-condition} 
(1)--(3) are satisfied.  

Let 
\begin{equation*}
\Lam_1 = 
\begin{pmatrix}
0 & 1 & 0 \\
1 & 0 & 0 \\
0 & 0 & -1
\end{pmatrix}, \;
\Lam_2 = 
\begin{pmatrix}
0 & 0 & 1 \\
0 & 1 & 0 \\
1 & 0 & 0 
\end{pmatrix}, \;
r_1 = 
\begin{pmatrix}
1 & 0 & 0 \\
0 & 0 & \Lam_1 \\
0 & \Lam_1 & 0 
\end{pmatrix}, \;
r_1 = 
- \begin{pmatrix}
1 & 0 & 0 \\
0 & 0 & \Lam_2 \\
0 & \Lam_2 & 0 
\end{pmatrix}. 
\end{equation*}
Then one can verify Assumption \ref{semi-simple-condition} (4) 
by elements of the subgroup generated by $r_1,r_2$.

If $d(t_1,t_2)\in \kernel(\al_1)\cap \kernel(\al_2)$, 
then $t_1t_2^{-1}=t_2=1$. This implies that 
\begin{math}
\bigcap_{\al\in\Phi} \kernel(\al)=\{e\}.  
\end{math}
Proposition \ref{borel-condition} implies that 
$(G_w\cap \gl_7)^{\circ}$ is a semi-simple group.  

Let 
\begin{equation*}
\exp(u X_a) = I_7 + u X_a + \frac 12 u^2 X_a^2. 
\end{equation*}
Then despite the denominator $2$, the 
entries of the resulting matrix are integral 
polynomials of $u$. Since $X_a^3=0$,  
if we put 
\begin{equation*}
U_a = \{\exp(u X_a)\mid u\in \aff^1\}, 
\end{equation*}
then this is a subgroup of $\spl_7$ 
with Lie algebra $\{u X_a\mid u\in\aff^1\}$. 
We can define $U_b$, etc., similarly.  
Then $U_a,U_b,U_f,U_{u_1},U_{u_2},U_{u_3}$ and $T$
generate a Borel subgroup of $(G_w\cap \spl_7)^{\circ}$.  
The set of roots coincides with $\Phi$ in 
(\ref{eq:root-system-G2}) and $(G_w\cap \spl_7)^{\circ}$ 
becomes a simple group of type $\text{G}_2$.

We next show that if a diagonal matrix 
$g=\diag(t_1\ccd t_7)$ with $t_1\ccd t_7\in\gl_1$
commutes with all elements of 
$(G_w\cap \spl_7)^{\circ}$ then $g=t_1 I_7$.  
Since $g$ commutes with elements of $U_a,U_b,U_c$, 
$t_2t_1^{-1}=t_3t_1^{-1}=t_4t_1^{-1}=1$. 
So $t_2=t_3=t_4=t_1$. Similarly, $t_5=t_6=t_7=t_1$. 
This implies that $g=t_1 I_7$.

Suppose that $g\in G_w\cap \spl_7$.  
Since all automorphisms of simple groups of 
type $\text{G}_2$ are inner, there exists 
$g_1\in (G_w\cap \spl_7)^{\circ}$ such that 
$ghg^{-1}=g_1hg_1^{-1}$ for all 
$h\in (G_w\cap \spl_7)^{\circ}$.
Replacing $g$ by $gg_1^{-1}$, we may assume 
that $g$ commutes with all elements of 
$(G_w\cap \spl_7)^{\circ}$.  Then 
there exists $t\in \gl_1$ such that $g=tI_7$. 
Since $g\in G_w,\spl_7$,  $t^3=t^7=1$.  
Therefore, $t=1$ and $g=I_7$. This shows
that $G_w\cap \spl_7$ is connected.  
The center of $G_w\cap \spl_7$ 
is trivial by the same argument. 

(3) follows from (1), (2) since $\ti T$ is a torus.  
\end{proof}

The following corollary follows from 
Propositions \ref{prop:case9-tangent-space}, 
\ref{prop:case9-regular-pv}. 

\begin{cor}
\label{cor:case3-G2-regular}
$(G,V)$ is a regular \pv. 
\end{cor}

\begin{prop}
\label{prop:r-orbit-G2}
$\rg\backslash \rv^{\sst}$ is in bijective 
correspondence with the set of $k$-isomorphism
classes of $k$-forms of connected simple group 
of type $\text{G}_2$ with trivial center 
by associating $G_x\cap \spl_7$ to $x\in \rv^{\sst}$.  
\end{prop}
\begin{proof}
It follows from the above corollary that 
$\rg\backslash \rv^{\sst}$ is in bijective 
correspondence with 
\begin{equation*}
\h^1(k,G_w) \cong \h^1(k,G_w\cap \spl_7)
\times \h^1(k,\ti T) 
\cong \h^1(k,G_w\cap \spl_7). 
\end{equation*}
Since $\h^1(k,\gl_1\times \gl_7)=\{1\}$, 
the orbit in $\rg\backslash \rv^{\sst}$
which corresponds to a 1-cocycle 
$\{h_{\sig}\}_{\sig\in \gal(k^{\sep}/k)}$
with coefficients in $G_w\cap \spl_7$ 
is $\rg\cdot gw$ where $g\in G_{k^{\sep}}$ 
such that $g^{-1}g^{\sig}=h_{\sig}$ 
for all $\sig\in \gal(k^{\sep}/k)$. 

Let $x=gw$. We describe $G_x\cap \spl_7$. 
It is easy to see that 
\begin{math}
G_{x\,k^{\sep}} = g G_{w\, k^{\sep}} g^{-1}. 
\end{math}
Since $G_{x\, k^{\sep}}$ modulo its center 
is isomorphic to $G_{x\, k^{\sep}}\cap \spl_7(k^{\sep})$, 
\begin{equation*}
G_{x\,k^{\sep}}\cap \spl_7(k^{\sep}) 
= g (G_{w\, k^{\sep}} \cap \spl_7(k^{\sep})) g^{-1}. 
\end{equation*}
If $a\in G_{w\,k^{\sep}}$ then $gag^{-1}\in G_{x\, k}$ 
if and only if 
\begin{math}
g^{\sig}a^{\sig}(g^{\sig})^{-1}
= gag^{-1}
\end{math}
for all $\sig\in \gal(k^{\sep}/k)$. 
This is equivalent to 
\begin{equation*}
g^{-1}g^{\sig} a^{\sig} (g^{-1}g^{\sig})^{-1}
= h_{\sig} a^{\sig} h_{\sig}^{-1} =a.   
\end{equation*}

The element of $\aut G_{w\, k^{\sep}}\cap \spl_7(k^{\sep})$
which corresponds to $h_{\sig}$ is 
\begin{equation*}
G_{w\, k^{\sep}}\cap \spl_7(k^{\sep}) \ni a 
\mapsto h_{\sig} a h_{\sig}^{-1} \in   
G_{w\, k^{\sep}}\cap \spl_7(k^{\sep}). 
\end{equation*}
Therefore, the $k$-form of $G_w\cap \spl_7$ 
which corresponds to $\{h_{\sig}\}$ is 
$G_{gw}\cap \spl_7$.  
\end{proof}

An equivariant map $V\ni x \mapsto \overline {Q}_x\in \sym^2 \aff^7$ 
was constructed in \cite[p.1694, Definition (2.13)]{yukiem} 
and 
\begin{equation*}
\overline Q_w = 6(-p_1^2+p_2p_5+p_3p_6+p_4p_7)
\end{equation*}
(we used the notation $p_i$ for $\bbmp_i$ 
since we regard it as a polynomial rather than a basis element).  
$Q_x$ is defined by homogeneous cubic polynomials of $x\in V$
and Proposition \ref{prop:case9-tangent-space} (1) 
implies that $w$ is universally generic 
(see \cite[p.279]{kable-yukie-quinary}, 
Definition III--4.9). 
So the above map is actually divisible by $6$
by \cite[p.279, Proposition 1]{kable-yukie-quinary}, 
Proposition III--4.10. 
This implies that 
$Q_x\stackrel{\mathrm{def}}=\tfrac 16 \overline Q_x$ 
is defined over $\Z$ and 
\begin{equation*}
Q_w = -p_1^2+p_2p_5+p_3p_6+p_4p_7.  
\end{equation*}
Let $P(x)$ be the discriminant of $Q_x$. 
Then $P(w)= - 1$. 

Even tough it is possible to associate octonions 
to elements of $\rg\backslash \rv^{\sst}$, 
we have to assume $\ch(k)\not=2$. It is probably possible
to remove the condition $\ch(k)\not=2$, but we consider 
this issue elsewhere.  We briefly mention the outline 
of the method to associate octonions to 
elements of $\rg\backslash \rv^{\sst}$. 

Let $x\in\rv^{\sst}$ and $W=(\aff^7)^*$. 
$Q_x$ can be regarded as a function on $W$.  
Let $B_x$ be the symmetric bilinear form 
on $W$ defined by 
\begin{equation*}
B_x(v_1,v_2) = Q_x(v_1+v_2)-Q_x(v_1)-Q_x(v_2)  
\end{equation*}
for $v_1,v_2\in W$

For elements of $\m_2$, define their conjugates as follows:
\begin{equation*}
\overline{
\begin{pmatrix}
a & b \\ c & d
\end{pmatrix}
}
= \begin{pmatrix}
d & -b \\
-c & a
\end{pmatrix}. 
\end{equation*}
Then the split octonion $\soct$ is defined as 
$\m_2(+)$ (see \cite[1692, Definition (2.1)]{yukiek}). 
It is assumed in \cite{yukiek} that $\ch(k)=0$, 
but the definition of 
$\soct$ works regardless of $\ch(k)$. Since we defined
the conjugate directly, the assumption $\ch(k)\not=2$ 
is not necessary for this part.  $B_x$ is $2$ times 
the bilinear form in \cite{yukiek}.  

$x$ can be regarded as an alternating 
multilinear map $W^3\to \aff^1$. 
Let $\im(\soct)$ be the orthogonal complement 
to $k\cdot 1\sub \soct$ with respect to $B_x$. 
We have to assume $\ch(k)\not=2$ for this part.  
If we choose a basis as \cite[1694, (2.9)]{yukiek}), 
then 
\begin{equation*}
w(v_1,v_2,v_3)=B_w(v_1,v_2\cdot v_3). 
\end{equation*}
Note that $Q_w$ is $1/6$ times $Q_w$ in \cite{yukiek}, 
but $1/2$ is not multiplied to make the bilinear form 
$B_w$. So the above equation is compatible with 
\cite[1695, (2.18)]{yukiek}). 
There is a relative invariant polynomial $\Del(x)$ on $V$ 
of degree $7$ (\cite[p.86]{saki}). Since $w$ 
is universally generic, we may assume that $\Del(w)=1$. 
For $w=(a,v)\in k\oplus W$, we denote 
$a,v$ by $\re(w),\im(w)$ respectively. 
For $v_1,v_2\in W$, define 
$(v_1\cdot v_2)_x\in W$ as the unique element 
such that 
\begin{equation*}
x(v_3,v_1,v_2)=B_x(v_3,(v_1\cdot v_2)_x)
\end{equation*}
for all $v_3\in W$. 
We define a non-associative $k$-algebra structure 
on $k\oplus W$ by 
\begin{equation*}
\re(v_1v_2) = -\Del(x)^{-1}B_x(v_1,v_2), \;
\im(v_1v_2) = v_2
\end{equation*}
Let $\mathbbm O_x$ be this $k$-algebra. 
Note that we need the assumption $\ch(k)\not=2$ 
to define the projection $\mathbbm O_x\to k$. 
Then the argument in \cite{yukiek} works and 
$\rg\backslash \rv^{\sst}$ is in bijective 
correspondence with $k$-isomorphism classes 
of $k$-forms of $\soct$.

\subsection{Case IV}

Let (i) $G=\gl_6$ or (ii) $\gl_1\times \gl_6$ 
and $V=\wedge^3 \aff^6$.
This case was considered in Section 1 of \cite{yukiek}. 
However, since $\ch(k)=0$ was assumed in \cite{yukiek}, 
we verify that $(G,V)$ is a regular \pv. 
We use the notation such as 
$\bbmp_{6,1}\ccd \bbmp_{6,6}$ and 
$p_{6,i_1i_2i_3}=\bbmp_{6,1}\wedge \bbmp_{6,2}\wedge \bbmp_{6,3}$, etc. 

Let 
\begin{equation*}
w = p_{6,123} + p_{6,456}.  
\end{equation*}
\begin{prop}  
\label{regular-wedge36}
The \rep\ $(G,V)$ is a regular \pv. 
Moreover, 
$G_w^{\circ}\cong \spl_3\times \spl_3$ or 
$\{(g_1,g_2)\in \gl_3\times \gl_3\mid \det g_1=\det g_2\}$ 
for the cases (i), (ii) respectively.    
\end{prop} 
\begin{proof}
We only prove this proposition for the case (ii).  
The consideration for the case (i) is similar.  
Let $t\in k,A=(a_{ij})\in\m_6(k)$. Then 
$(1+\ep t)(I_6+\ep A) w = w$ 
if and only if 
\begin{equation}
\label{eq:case98-regular}
\begin{aligned}
& t p_{6,123} + 
(a_{11}\bbmp_{6,1} + \cdots + a_{61}\bbmp_{6,6}) \wedge p_{23}
- (a_{12}\bbmp_{6,1} + \cdots + a_{62}\bbmp_{6,6}) \wedge p_{13} \\
& + (a_{13}\bbmp_{6,1} + \cdots + a_{63}\bbmp_{6,6}) \wedge p_{12} 
+ tp_{6,456} + 
+ (a_{14}\bbmp_{6,1} + \cdots + a_{64}\bbmp_{6,6}) \wedge p_{56} \\
& - (a_{15}\bbmp_{6,1} + \cdots + a_{65}\bbmp_{6,6}) \wedge p_{46}
+ (a_{16}\bbmp_{6,1} + \cdots + a_{66}\bbmp_{6,6}) \wedge p_{45} \\
& = 0.  
\end{aligned}
\end{equation}
Since terms like $p_{6,23i}$ for $i=4,5,6$ appear only 
in the second term, 
$a_{14},a_{15},a_{16}=0$. The consideration is similar for other 
cases and $A$ is in the following form 
\begin{equation*}
\begin{pmatrix}
X & 0 \\
0 & Y
\end{pmatrix}
\end{equation*}
where $X,Y\in\m_3$ and $t+\trace(X)=t+\trace(Y)=0$.  
Conversely, matrices of the above form satisfies 
(\ref{eq:case98-regular}). 

Elements of the form
\begin{equation*}
\left(
t,
\begin{pmatrix}
A & 0 \\
0 & B 
\end{pmatrix}
\right)
\end{equation*}
where $t\in\gl_1$, $A,B\in\gl_3$ and $\det A=\det B=t^{-1}$ fix $w$. Since 
$\dim \mathrm{T}_{I_6}(G_w)=17 = \dim G-\dim V = 37-20$, 
$G_w^{\circ}$ consists of all the elements of the above form 
and is smooth over $k$. Therefore, 
$(G,V)$ is a regular \pv{} and $G\cdot w\sub V$ is Zariski 
open.  
\end{proof}

Since $(G,V)$ is a regular \pv, $V^{\sst}_{k^{\sep}}=G_{k^{\sep}} \cdot w$. 
The above proposition implies that $w$ is universally generic 
(\cite[p.279]{kable-yukie-quinary}, 
Definition III--4.9). 
Let $W=\aff^6$.  There exists an equivariant map 
\begin{equation*}
S:V\ni x\mapsto S_x\in \Hom(W,W)
\end{equation*}
(see \cite[pp.80,81]{saki}, \cite[p.1691]{yukiek})
such that 
\begin{equation*}
S_w(\bbmp_{6,i})= 
\begin{cases}
\bbmp_{6,i} & i=1,2,3, \\
-\bbmp_{6,i} & i=4,5,6 
\end{cases}
\end{equation*}
where the action of $g=(t_0,g_1)\in G$ on $f\in \Hom(W,W)$ is 
$t_0^2(\det g) gf(g^{-1}x)$ (without the factor $t_0^2$ in the case (i)).  
Note that 
$k=\C$ in \cite{saki} and $\ch(k)=0$ in \cite{yukiek}. 
However, the argument works for all $k$ for this part.  
There exists a polynomial $P(x)$ on $V$ such that 
$S_x^2=t_0^4 P(x)I_6$ (without the factor $t_0^4$ in the case (i)).  
This $P(x)$ is a homogeneous invariant polynomial
of degree $4$ satisfying $P(w)=1$ and $P(gx)=t_0^4 (\det g_1)^2P(x)$ 
(without the factor $t_0^2$ in the case (i)).

Suppose that $x\in V^{\sst}$. Let 
$\overline T_x = \wedge^2 S_x - P(x) \id_{\wedge^2 W}
\in \Hom(\wedge^2 W,\wedge^2 W)$.  
Then 
\begin{equation}
\label{eq:barT-equivariant}
\overline T_{gx}(\om) = t_0^4(\det g_1)^2 \wedge^2 \!g \, 
\overline T_x((\wedge^2 \!g)^{-1} \om)
\end{equation}
(without the factor $t_0^2$ in the case (i)). 
Note that 
\begin{equation*}
\wedge^2 \!g \, \id_{\wedge^2 W}((\wedge^2 \!g)^{-1} \om) = \om.  
\end{equation*}
\begin{prop}
\label{prop:Txhalf}
There exists a map 
\begin{math}
V\ni x \mapsto T_x\in \Hom(\wedge^2 W,\wedge^2 W)
\end{math}
with the property 
\begin{equation*}
T_{gx}(\om) = t_0^4(\det g_1)^2 \wedge^2 \!g \, 
T_x((\wedge^2 \!g)^{-1} \om)
\end{equation*}
(without the factor $t_0^2$ in the case (i)) 
such that $\overline T_x=2T_x$ for all $x\in V$, 
$T_w(p_{6,ij})=0$ if $i,j\leqq 3$ or $i,j\geqq 4$ and 
$T_w(p_{6,ij})=-p_{6,ij}$ if $i\leqq 3$ and $j\geqq 4$.  
\end{prop}
\begin{proof}
If $i,j\leqq 3$ then 
\begin{math}
S_w(\bbmp_{6,i}) = \bbmp_{6,i}, \;
S_w(\bbmp_{6,j}) = \bbmp_{6,j}.
\end{math}
So $\overline T_w(p_{6,ij})=p_{6,ij}-p_{6,ij}=0$. 
Similarly, $\overline T_w(p_{6,ij})=0$ if $i,j\geqq 4$. 
If $i\leqq 3,j\geqq 4$ then 
\begin{equation*}
\overline T_w(p_{6,ij})= \bbmp_{6,i}\wedge (-\bbmp_{6,j})-p_{6,ij}
= -2p_{6,ij}.   
\end{equation*}

Note that $\overline T_x$ is defined over $\Z$.  
As we pointed out above, $w$ is universally generic and 
$\overline T_w$ is divisible by $2$ over $\Z$. Therefore, 
by \cite[p.279, Proposition 1]{kable-yukie-quinary}, 
Proposition III--4.10, 
there exists an equivariant map over $\Z$ 
\begin{equation*}
V\ni x \mapsto T_x\in \Hom(\wedge^2 W,\wedge^2 W)
\end{equation*}
such that $\overline T_x=2T_x$. If we consider 
the natural map $\Z\to k$, the induced map   
$T_x$ satisfies the statement of the proposition. 
\end{proof}

Let 
\begin{equation}
\label{eq:Tx-image}
U(x) = \im(T_x).
\end{equation}
This is a subspace of $\wedge^2 W$. Any element $\om$ of $U(x)$ 
defines an alternating bilinear form $B_{\om}$ on $W^*$. 
A subspace $E\sub W^*$ is said to be totally isotropic 
with respect to $U(x)$ if $B_{\om}(v,w)=0$ for all 
$v,w\in E$ and $\om\in U(x)$.  

Let $\{\bbmp_{6,1}\ccd \bbmp_{6,6}\}$ be the 
standard basis of $W$ and $\{\bbmq_{6,1}\ccd \bbmq_{6,6}\}$ 
the dual basis. We define 
\begin{equation}
\label{eq:X1X2-defn}
X_1 = \lan \bbmq_{6,1},\bbmq_{6,2},\bbmq_{6,3}\ran, \;
X_2 = \lan \bbmq_{6,4},\bbmq_{6,5},\bbmq_{6,6}\ran.  
\end{equation}
\begin{prop}
\label{prop:wegde63-isotropic}
The only totally isotropic subspaces 
of dimension $3$ with respect to $U(w)$ are $X_1$ and $X_2$.  
\end{prop}
\begin{proof}
Suppose that $i,j\leqq 3$ or $i,j\geqq 4$.  
Then Proposition \ref{prop:Txhalf} implies that 
$X_1=\lan \bbmq_{6,1},\bbmq_{6,2},\bbmq_{6,3}\ran$ and 
$X_2=\lan \bbmq_{6,4},\bbmq_{6,5},\bbmq_{6,6}\ran$ 
are totally isotropic subspaces.  

Let $E\sub W^*$ be a totally isotropic 
subspace of dimension $3$ with respect to $U(w)$.  
Suppose that there exist 
\begin{equation*}
v = \sum_{i=1}^3 a_i \bbmq_{6,i},\; 
w = \sum_{j=4}^6 b_j \bbmq_{6,j} \in E 
\end{equation*}
and $a_{i_0},b_{j_0}\not=0$.  
Let $\om=p_{6,i_0j_0}$. Then $\om \in U(w)$ and 
$B_{\om}(v,w)\not=0$, which is a contradiction.  
Therefore, $E$ cannot contain non-zero elements 
of $X_1$ and $X_2$ simultaneously.  

Suppose that there exists 
\begin{equation*}
v = \sum_{i=1}^3 a_i \bbmq_{6,i} 
+ \sum_{j=4}^6 b_j \bbmq_{6,j} \in E  
\end{equation*}
where $a_i,b_j\in k$, $i_0\leqq 3,j_0\geqq 4$
and $a_{i_0},b_{j_0}\not=0$.  Let
\begin{equation*}
w = \sum_{i=1}^3 c_i \bbmq_{6,i} 
+ \sum_{j=4}^6 d_j \bbmq_{6,j} \in E  
\end{equation*}
where $c_i,d_j\in k$. 

If $i_1\leqq 3$ and $c_{i_1}\not=0$ then 
\begin{equation*}
B_{\om}(v,w) = a_{i_1}d_{j_0}-b_{j_0}c_{i_1} = 0.   
\end{equation*}
for $\om=p_{6,i_1j_0}$.  
This implies that $a_{i_1},d_{j_0}\not=0$ 
and $c_{i_1}= a_{i_1}d_{j_0}/b_{j_0}$. 
The above equation holds without the condition 
$c_{i_1}\not=0$.  So 
$c_i= a_id_{j_0}/b_{j_0}$ for $i=1,2,3$. 
In particular, $c_{i_0}\not=0$ and 
$c_{i_0}/a_{i_0}= d_{j_0}/b_{j_0}$.

Let $j\geqq 4$ and $\om=p_{6,i_0,j}$. Then 
\begin{equation*}
B_{\om}(v,w) = a_{i_0}d_j- b_j c_{i_0} = 0.   
\end{equation*}
This implies that $d_j = b_j c_{i_0}/a_{i_0}=b_jd_{j_0}/b_{j_0}$. 
So $w$ is a scalar multiple of $v$.  The case 
$d_j\not=0$ for some $j\geqq 4$ is similar. 
Therefore, $\dim E=1$, which is a contradiction.   
This implies that elements of $E$ belong to $X_1$ or $X_2$. 
Since $E$ cannot contain non-zero elements of $X_1$ and $X_2$ 
simultaneously, $E$ has to be $X_1$ or $X_2$.  
\end{proof}

Let $F=k(\al_1)$ be a quadratic Galois extension of $k$  
and $\al_2\not=\al_1$ the conjugate of $\al_1$ over $k$.  
Then the argument in Section 1 of \cite{yukiek} works 
as long as we replace the matrix 
\cite[p.1690, (1.8)]{yukiek} by 
\begin{equation*}
g_{\al} = 
\begin{pmatrix}
1 & 0 & 0 & 1 & 0 & 0 \\
\al_1 & 0 & 0 & \al_2 & 0 & 0 \\
0 & 1 & 0 & 0 & 1 & 0 \\
0 & \al_1 & 0 & 0 & \al_2 & 0 \\
0 & 0 & 1 & 0 & 0 & 1 \\
0 & 0 & \al_1 & 0 & 0 & \al_2 
\end{pmatrix}.  
\end{equation*}
Therefore, we obtain the following proposition. 

\begin{thm}
\label{thm:case98-interpretation}
$G_k \backslash \rv^{\sst}$ 
is in bijective correspondence with $\Ex_2(k)$. 
\end{thm}
\begin{prop}
\label{prop:Ux-totally-isotropic}
Let $g_x\in \gl_6({k^{\sep}})$ and 
$x=g_x w\in V_{k^{\sep}}^{\sst}$.  Then 
the only totally isotropic subspaces over $k^{\sep}$ 
of dimension $3$ with respect to $U(x)$ are
$\,{}^t\!g_x^{-1}(X_1),{}^t\!g_x^{-1}(X_2)$.   
\end{prop}
\begin{proof}
Let $a\in W$, $v\in W^*$ and $g\in\gl_6$. 
We regard $a$ as a map from $W^*$ to $\aff^1$. 
Then $(ga)(v) = a(\,{}^t\!g v)$.  
If $\om\in U(w)$, $\om'= (\wedge^2 g_x)\, \om$
and $v,w\in W^*_{k^{\sep}}$ then 
\begin{equation*}
B_{\om'}(v,w) = B_{\om}(\,{}^t\!g_x v,{}^t\!g_x w).  
\end{equation*}
Since $(\wedge^2 g_x)$ induces an isomorphism 
$U(w)\to U(x)$ over $k^{\sep}$, 
the subspace $E\sub W^*_{k^{\sep}}$ is 
totally isotropic with respect to $U(x)$ 
if and only if ${}^t\!g_x(E)$ is totally isotropic 
with respect to $U(w)$.  Therefore,   
if such $E$ is of dimension $3$ then 
$E = {}^t\!g_x^{-1}(X_1),{}^t\!g_x^{-1}(X_2)$.  
\end{proof}

The above proposition implies that 
${}^t\!g_x^{-1}(X_1),{}^t\!g_x^{-1}(X_2)$ 
are the two totally isotropic subspaces 
of $U(x)$.  Let $F=k(\al_1)\in\Ex_2(k)$.  
Then $x_{\al}=g_{\al}w$ corresponds to $F$.  
If $x\in \rv^{\sst}$ corresponds to $F$, 
there exists $g\in\gl_6(k)$ such that 
$x=gg_{\al}w$.  Then 
the only totally isotropic subspace over $k^{\sep}$ 
of dimension $3$ with respect to $U(x)$ are
\begin{equation*}
{}^t\!g^{-1}{}^t\!g_{\al}^{-1}(X_1), \; 
{}^t\!g^{-1}{}^t\!g_{\al}^{-1}(X_2).  
\end{equation*}

These subspaces are defined over $F$. 
Let $\sig\in\gal(F/k)$ be the non-trivial element and 
\begin{equation*}
\tau = \begin{pmatrix}
0 & I_3 \\
I_3 & 0 
\end{pmatrix}.
\end{equation*}
Then $g_{\al}^{\sig}=g_{\al}\tau$. So 
\begin{equation*}
({}^tg^{-1}{}^tg_{\al}^{-1})^{\sig}
= {}^tg^{-1}{}^tg_{\al}^{-1}\tau.
\end{equation*}

Since 
\begin{math}
{}^tg^{-1}{}^tg_{\al}^{-1}(X_1)
\end{math}
is spanned by the first three columns of 
${}^tg^{-1}{}^tg_{\al}^{-1}$ and multiplying $\tau$ 
from the right exchanges the first three columns 
and the last three columns, 
\begin{equation*}
({}^tg^{-1}{}^tg_{\al}^{-1}(X_1))^{\sig}
= {}^tg^{-1}{}^tg_{\al}^{-1}(X_2).
\end{equation*}
Therefore, these two subspaces are not defined over $k$ 
and we obtain the following proposition.  

\begin{prop}
If $x\in \rv^{\sst}$ corresponds to a quadratic field $F$, 
then $F$ is the field of definition of the two 
totally isotropic subspaces over $k^{\sep}$ of dimension $3$ 
of $W^*_{k^{\sep}}$ with respect to $U(x)$.  
\end{prop}

\subsection{Case V}

Let 
\begin{equation}
\label{eq:322-castling}
G=\gl_3\times \gl_2^2, \; W_1=\aff^3, \; W_2=W_3=\aff^2
\end{equation}
and $V=W_1\otimes W_2\otimes W_3$. 
We express elements of $V$ as 
$x=(x_1,x_2,x_3)$ where $x_1,x_2,x_3\in\m_2$, 
which corresponds to 
\begin{equation*}
\bbmp_{3,1}\otimes x_1
+\bbmp_{3,2}\otimes x_2
+\bbmp_{3,3}\otimes x_3.
\end{equation*}
The action of $g=(g_1,g_2,g_2)$ on $x$ is defined by 
\begin{equation*}
g \begin{pmatrix}
x_1 \\
x_2 \\
x_3
\end{pmatrix}
= g_1 \begin{pmatrix}
g_2 x_1 \,{}^t\! g_3 \\ 
g_2 x_2 \,{}^t\! g_3 \\ 
g_2 x_3 \,{}^t\! g_3 
\end{pmatrix}.  
\end{equation*}

The \rep{} $(G,V)$ is related to two of the strata.  
Let $H=\gl_2^2$ and $W=(W_2\otimes W_3)^*$ (the dual space). 
Then $(G,V)$ and $(H,W)$ 
are the Castling transform of each other.  
Let 
\begin{equation*}
\Phi:C\ni (x_1,x_2,x_3) \mapsto x_1\wedge x_2\wedge x_3
\in \wedge^3 (W_2\otimes W_3) \cong W.
\end{equation*}

The following proposition is pointed out 
in \cite[p.8]{tajima-yukie-GIT2}.  

\begin{prop}
\label{prop:322-orbit-SP}
\begin{itemize}
\item[(1)]
$V^{\sst}\not=\emptyset$. 
\item[(2)]
$\Phi$ induces a bijection 
$\rg\backslash \rv^{\sst}\to H_k\backslash W_k^{\sst}$. 
\item[(3)]
Let 
\begin{equation*}
w = \bbmp_{3,1}\otimes (-E_{11}+E_{22})
+ \bbmp_{3,2}\otimes E_{12}
+ \bbmp_{3,3}\otimes E_{21} \in V
\end{equation*}
(see Notation for $E_{11}$, etc.).  
Then $\Phi(w)=I_2$ and 
$\rv^{\sst}=\rg \cdot w$. 
\end{itemize}
\end{prop}

Let $P(x)=\det \Phi(x)$.  Then $P(x)$ is a homogeneous 
polynomial of degree $6$ of $x\in V$ and 
\begin{equation}
\label{eq:223-relative-inv-poly}
P(gx) = (\det g_1)^2(\det g_2)^3(\det g_3)^3P(x).
\end{equation}

\section{Rational orbits (2) }
\label{sec:orbits2}

Let $G=\gl_4\times \gl_3\times \gl_1$ 
and $V=\wedge^2\aff^4\otimes \aff^3\oplus \wedge^2\aff^3$.  
We define an action of $g=(g_1,g_2,t_3)\in G$ on $V$ by 
\begin{math}
gx = ((\wedge^2 \hskip -1pt g_1\otimes g_2) \, x_1, 
t_3 \wedge^2 \hskip -3pt g_2 \, x_2)
\end{math}
for $x=(x_1,x_2)\in V$. 

Let $\Phi_1:V\to \wedge^2 \aff^4$ be the map
obtained by taking the wedge product of the 
second component of $x_1$ and $x_2$. 
We identify 
$\bbmp_{3,1}\wedge \bbmp_{3,2}\wedge \bbmp_{3,3}$
with $1\in k$. Then 
\begin{equation*}
\Phi_1(gx) = (\det g_2)\, t_3 \wedge^2 \! g_1 \, \Phi_1(x)
\end{equation*}
for the above $g$.  

For 
\begin{equation*}
M = \begin{pmatrix}
0 & a & b & c \\
-a & 0 & d & e \\
-b & -d & 0 & f \\ 
-c & -e & -f & 0
\end{pmatrix}
\end{equation*}
we choose the Pfaffian of $M$ as $\pfaff(M)=af-be+cd$. 
Let $P_1(x) = \pfaff(\Phi_1(x))$.  
$P_1(x)$ is a homogeneous polynomial of degree $4$ of 
$x\in V$ and 
\begin{equation}
\label{eq:be4-phi1-inv}
P_1((g_1,g_2)x) = (\det g_1)(\det g_2)^2t_3^2 P_1(x). 
\end{equation}

We express $x_1$ as $x_1=(x_{11},x_{12},x_{13})$ where 
$x_{11},x_{12},x_{13}$ are $4\times 4$ alternating matrices
with diagonal entries $0$. 
Let $v_1,v_2,v_3$ be variables, $v=(v_1,v_2,v_3)$ and 
$\Phi_2(x_1,x_2)(v)$ be the Pfaffian of 
$v_1x_{11}+v_2 x_{12}+v_3x_{13}$.  
This is a quadratic form on $k^3$. $\Phi_2$ is equivariant 
with respect to the action of $G$ as follows: 
\begin{equation*}
\Phi_2(gx)(v) = (\det g_1) g_2\Phi_2(x)(v).  
= (\det g_1) \Phi_2(x)(vg_2)
\end{equation*}
where $v$ is regarded as a row vector.  
Let $P_2(x)$ be the discriminant of 
$\Phi_2(x)(v)$.  Then $P_2(x)$ is a 
homogeneous polynomial of degree $6$ of $x\in V$
and 
\begin{equation}
\label{eq:be4-phi2-inv}
P_2((g_1,g_2)x) = (\det g_1)^3 (\det g_2)^2 P_2(x).  
\end{equation}

Let 
\begin{align*}
& A_1 = 
\begin{pmatrix}
J & 0 \\
0 & 0 
\end{pmatrix}, \;
A_2 = \begin{pmatrix}
0 & I_2 \\
-I_2 & 0 
\end{pmatrix}, \;
A_3 = \begin{pmatrix}
0 & 0 \\
0 & J 
\end{pmatrix},\; \\
& A = A_1\otimes \bbmp_{3,1}
+A_2\otimes \bbmp_{3,2}
+A_3\otimes \bbmp_{3,3}, \\
& B = \bbmp_{3,1} \wedge \bbmp_{3,3}
\end{align*}
where $J$ is the matrix in (\ref{eq:case7-R-defn}). 
Then 
\begin{equation}
\label{eq:be4-invariants}
\begin{aligned}
& \Phi_1(A,B) = -A_2, \; 
P_1(A,B) = 1, \\
& \Phi_2(A,B)(v) = v_1v_3-v_2^2, \; 
P_2(A,B) = 1.  
\end{aligned}
\end{equation}

Let $a,b>0$ be integers and 
\begin{equation*}
\chi(g)=((\det g_1)(\det g_2)^2t_3^2)^a 
((\det g_1)^3(\det g_2)^2)^b
\end{equation*}
for $g\in G$.  
If we define the stability with respect to 
this character and the action of 
$\{g=(g_1,g_2,t_3)\in G\mid (\det g_1)(\det g_2)t_3=1\}$, 
then $V^{\sst}=\{x\in V\mid P_1(x),P_2(x)\not=0\}$.  
Let 
\begin{equation}
\label{eq:case4-R-defn}
w=(A,B). 
\end{equation}
(\ref{eq:be4-invariants}) implies the following proposition. 

\begin{prop}
\label{prop:case4-vss-nonempty}
$w\in V^{\sst}$ and so $V^{\sst} \not=\emptyset$.  
\end{prop}

Let 
\begin{equation}
\label{eq:ght-defn}
\phi(t,h) = \left(
\begin{pmatrix}
h & 0 \\
0 & {}^t h^{-1}t^{-1}
\end{pmatrix}, 
\begin{pmatrix}
(\det h)^{-1} & 0 & 0 \\
0 & t & 0 \\
0 & 0 & (\det h)t^2
\end{pmatrix},
t^{-2}\right)
\end{equation}
for $h\in\gl_2$ and $t\in\gl_1$.  
We define 
\begin{align*}
H & = 
\{\phi(t,h)\mid h\in\gl_2,t\in\gl_1\}\sub G, \\
\tau & = \left(
\begin{pmatrix}
0 & I_2 \\
I_2 & 0 
\end{pmatrix}, 
\begin{pmatrix}
0 & 0 & 1 \\
0 & -1 & 0 \\
1 & 0 & 0 
\end{pmatrix}, -1
\right)\in G.
\end{align*}
Then $H\sub G_w$ and $\tau\in G_w$.

Since $H$ is the image of the homomorphism
\begin{equation*}
\gl_2\times \gl_1\ni (h,t)
\mapsto \phi(t,h)\in G, 
\end{equation*}
$H$ is a closed subgroup of $G$. 
Note that the image of a homomorphism 
of algebraic groups is a closed subgroup.

\begin{prop}
\label{prop:be4-stabilizer}
Suppose that $\ch(k)\not=2$. Then $(G,V)$ is a 
regular \pv, 
$G_w^{\circ} = H$ and $G_w/G_w^{\circ}\cong \Z/2\Z$ 
is represented by $\tau$.    
\end{prop}
\begin{proof}
Let $k[\ep]/(\ep^2)$ be the ring of dual numbers and 
$e=(I_4,I_3,1)\in G$ the unit element.  
We first determine the tangent space of $G_w$ at $e$.  
Let 
\begin{equation*}
X=(I_4+\ep X_1,I_3+\ep X_2,1+\ep a)
\end{equation*}
where 
$X_1\in\m_4(k),X_2=\m_3(k)$ and $a\in \mk$.

Suppose that $X\in \text{T}_e(G_w)$. 
Then $\det (I_4 + \ep X_1) (I_3+ \ep X_2)$ fixes 
$\Phi_2(w)(v)=v_1v_3-v_2^2$. 
Let $Q(v)=v_1v_3-v_2^2$. The action of 
$I_3+ \ep X_2$ on $Q(v)=v_1v_3-v_2^2$ is 
$Q(v (I_3+ \ep X_2))$. Let $X_2=(x_{2,ij})$. 
Since the entries of $v (I_3+ \ep X_2)$ are 
\begin{align*}
v_1 + \ep (x_{2,11}v_1+x_{2,21}v_2+x_{2,31}v_3), \\ 
v_2 + \ep (x_{2,12}v_1+x_{2,22}v_2+x_{2,32}v_3), \\ 
v_3 + \ep (x_{2,13}v_1+x_{2,23}v_2+x_{2,33}v_3), 
\end{align*}
$(1+\ep \trace(X_1))Q(v(I_3+\ep X_2))-(v_1v_3-v_2^2)$ 
is $\ep$ times 
\begin{align*}
& \trace(X_1)(v_1v_3-v_2^2) 
+ v_1(x_{2,13}v_1+x_{2,23}v_2+x_{2,33}v_3) \\
& \quad + v_3 (x_{2,11}v_1+x_{2,21}v_2+x_{2,31}v_3) 
- 2v_2 (x_{2,12}v_1+x_{2,22}v_2+x_{2,32}v_3) \\
& = x_{2,13}v_1^2 
+ (x_{2,23}-2x_{2,12})v_1v_2
+ (\trace(X_1) + x_{2,11}+x_{2,33})v_1v_3  \\
& \quad -(\trace(X_1) + 2x_{2,22})v_2^2
+ (x_{2,21}-2x_{2,32})v_2v_3
+ x_{2,31}v_3^2 = 0.  
\end{align*}

This implies that 
\begin{equation}
\begin{aligned}
& x_{2,13},
x_{2,23}-2x_{2,12},
\trace(X_1) + x_{2,11}+x_{2,33}, \\
& \trace(X_1) + 2x_{2,22},
x_{2,21}-2x_{2,32},
x_{2,31}=0. 
\end{aligned}
\end{equation}
So $X_2$ is in the following form:
\begin{equation*}
X_2 = \begin{pmatrix}
x_{2,11} & x_{2,12} & 0 \\
2x_{2,32} & x_{2,22} & 2x_{2,12} \\
0 & x_{2,32} & 2x_{2,22}-x_{2,11}
\end{pmatrix}.
\end{equation*}
Note that $2x_{2,22}=x_{2,11}+x_{2,33}$. 

Now we consider the direct action of $X$ on $(A,B)$. 
The second component of $X(A,B)$ is 
\begin{equation*}
(1+\ep a) ((1+\ep x_{2,11})\bbmp_{3,1}+2\ep x_{2,32}\bbmp_{3,2})
\wedge (2\ep x_{2,12}\bbmp_{3,2}+(1+\ep(2x_{2,22}-x_{2,11}))\bbmp_{3,3})
\end{equation*}
and this is $\bbmp_{3,1}\wedge \bbmp_{3,3}$ if $X\in \text{T}_e(G_w)$. 
So 
\begin{equation*}
(a+x_{2,11}+2x_{2,22}-x_{2,11})\bbmp_{3,1}\wedge \bbmp_{3,3}
+ 2x_{2,32}\bbmp_{3,2}\wedge \bbmp_{3,3} 
+ 2x_{2,12}\bbmp_{3,1}\wedge \bbmp_{3,2}=0.
\end{equation*}
Therefore, 
\begin{equation*}
a+2x_{2,22}=2x_{2,32}=2x_{2,12}=0. 
\end{equation*}
Since $\ch(k)\not=2$ by assumption, 
$x_{2,32}=x_{2,12}=0$, $a=-2x_{2,22}$.  
This implies that $X_2$ is in the form:
\begin{equation*}
X_2 = \begin{pmatrix}
x_{2,11} & 0 & 0 \\
0 & x_{2,22} & 0 \\
0 & 0 & 2x_{2,22}-x_{2,11}
\end{pmatrix}.
\end{equation*}

We express 
\begin{equation*}
X_1 = \begin{pmatrix}
Y_1 & Y_2 \\
Y_3 & Y_4
\end{pmatrix}
\end{equation*}
where $Y_1\ccd Y_4\in\m_2$.
The first component of $X(A,B)$ is 
$A+\ep(W_1,W_2,W_3)$ where
\begin{align*}
W_1 & = X_1 A_1 + A_1 \, {}^t \! X_1 + x_{2,11} A_1 =0, \\
W_2 & = X_1 A_2 + A_2 \, {}^t \! X_1 + x_{2,22} A_2 =0, \\
W_3 & = X_1 A_3 + A_3 \, {}^t \! X_1 + (2x_{2,22}-x_{2,11}) A_3.
\end{align*}
If $M\in\m_2$ then $MJ+J\,{}^t\!M=\trace(M)J$.

Explicitly, the condition $W_1=W_2=W_3=0$ is as follows:
\begin{align*}
& \begin{pmatrix}
Y_1 J & 0 \\
Y_3 J & 0 
\end{pmatrix}
+ \begin{pmatrix}
J \,{}^tY_1 & J \,{}^tY_3 \\
0 & 0 
\end{pmatrix}
+ \begin{pmatrix}
x_{2,11} J & 0 \\
0 & 0 
\end{pmatrix} = 0, \\
& \begin{pmatrix}
-Y_2 & Y_1 \\
-Y_4 & Y_3 
\end{pmatrix}
+ \begin{pmatrix}
{}^tY_2 & {}^tY_4 \\
-\,{}^tY_1 & -\,{}^tY_3
\end{pmatrix} 
+ \begin{pmatrix}
0 & x_{2,22} I_2 \\
-x_{2,22}I_2 & 0 
\end{pmatrix}=0, \\
& \begin{pmatrix}
0 & Y_2 J \\
0 & Y_4 J 
\end{pmatrix}
+ \begin{pmatrix}
0 & 0 \\
J\,{}^t Y_2 & J\,{}^tY_4
\end{pmatrix} 
+ \begin{pmatrix}
0 & 0 \\
0 & (2x_{2,22}-x_{2,11})J
\end{pmatrix} = 0.  
\end{align*}
By these equations, 
\begin{equation*}
Y_2=Y_3=0, \;
Y_1+{}^tY_4+x_{2,22}I_2=0, \;
Y_4J + J{}^tY_4 + (2x_{2,22}-x_{2,11})J=0.
\end{equation*}
So 
\begin{math}
Y_4 = -{}^tY_1-x_{2,22}I_2.  
\end{math}
Since $Y_4J + J{}^tY_4=\trace(Y_4)J$, 
\begin{equation*}
\trace(Y_4)+2x_{2,22}-x_{2,11}=-\trace(Y_1)-2x_{2,22}+2x_{2,22}-x_{2,11}
= -\trace(Y_1)-x_{2,11}=0.    
\end{equation*}
This implies that $x_{2,11}=-\trace(Y_1)$.  
Therefore, 
\begin{equation*}
X_1 = \begin{pmatrix}
Y_1 & 0 \\
0 & -{}^t Y_1 -x_{2,22}I_2
\end{pmatrix}, \; 
X_2 = \begin{pmatrix}
-\trace (Y_1) & 0 & 0 \\
0 & x_{2,22} & 0 \\
0 & 0 & 2x_{2,22}+\trace (X_1)
\end{pmatrix}, \;
a = -2x_{2,22}. 
\end{equation*}

Since $\dim G =26$, $\dim V=21$ and 
\begin{equation*}
\dim \text{T}_e(G_w) = 5 = \dim G - \dim V, 
\end{equation*}
$G\cdot w\sub V$ is Zariski open, 
$G_w$ is smooth over $k$ and $\dim G_w = 5$.  
Since $H\sub G_w$ is a connected closed subgroup and $\dim H=5$, 
$H=G_w^{\circ}$. 
Since $H\cong \gl_2\times \gl_1$ which is reductive, 
$(G,V)$ is a regular \pv{} and so 
if $U=G\cdot w$ then $U_{k^{\sep}}$ is a single 
$G_{k^{\sep}}$-orbit.

Let $g\in G_w$. The map 
\begin{equation*}
I(g):G_w^{\circ} \ni h \mapsto ghg^{-1} \in G_w^{\circ}
\end{equation*}
is an automorphism of $G_w^{\circ}$. 
Since $\spl_2\sub G_w^{\circ}$ does not have a non-trivial 
character, $I(g)$ induces an automorphism of $\spl_2$. 
Since automorphisms of $\spl_2$ are inner, 
we may assume that $g$ commutes with 
elements of the form $\phi(1,\diag(t_0,t_0^{-1}))$ 
by multiplying an element of the form 
$\phi(1,h)$ ($h\in\spl_2$) if necessary.

Also $I(g)$ maps the center of $G_w^{\circ}$ 
to itself. Therefore, $I(g)$ induces an automorphism 
of the center of $G_w^{\circ}$, which  
consists of $\phi(t_1,t_2I_2)$ in (\ref{eq:ght-defn}) 
for all $t_1,t_2\in\gl_1$.  
$G_w^{\circ}$ naturally acts on $\aff^4$ and $\aff^3$. 
The weights of the coordinates regarding the 
action of $\phi(t_1,\diag(t_0t_2,t_0^{-1}t_2))$ 
are 
\begin{equation*}
t_0t_2,
t_0^{-1}t_2, 
t_0^{-1}t_1^{-1}t_2^{-1},
t_0t_1^{-1}t_2^{-1},
t_0^{-2},
t_1, 
t_0^2t_1^2 
\end{equation*}
respectively.

Since these are distinct, 
the action of $g$ induces a permutation
of weight spaces. Moreover, the 
weights with respect to $t_0$ do not change. 
Therefore, the second component of $g$ 
is a diagonal matrix and 
the first component of $g$ 
either fixes $\bbmp_{4,1},\bbmp_{4,4}$ 
or exchange them. Similarly, 
the first component of $g$ 
either fixes $\bbmp_{4,2},\bbmp_{4,3}$ 
or exchange them.

By taking the conjugation by $\tau$, 
the second component of $g$ stays diagonal
and we may assume that 
the first component of $g$ is in one of the following 
forms:
\begin{equation*}
\begin{pmatrix}
* &&& \\
& * && \\
&& * & \\
&&& *
\end{pmatrix}, \;
\begin{pmatrix}
* &&& \\
& & * & \\
& * && \\
&&& *
\end{pmatrix}. 
\end{equation*}
By the condition $g\in G_w$, 
the only possibility is that 
the first component is a diagonal matrix. 
Then it is easy to show that 
$g\in G_w^{\circ}$. 
\end{proof}

We assume $\ch(k)\not=2$ for the rest of this section. 
Since $V$ has two irreducible factors and there are 
two relative invariant polynomials (\ref{eq:be4-phi1-inv}), 
(\ref{eq:be4-phi2-inv}), the following corollary 
follows from Corollary III--4.7.

\begin{cor}
\label{cor:case4-single-orbit}
$U=\{x\in V\mid P_1(x),P_2(x)\not=0\}$.
\end{cor}

Note that $\lan \tau \ran\cong \Z/2\Z$ is an abelian group 
and $\gal(k^{\sep}/k)$ acts on this group trivially.  
Therefore, $\h^1(k,\Z/2\Z)$ is in bijective correspondence with 
the set of homomorphisms from $\gal(k^{\sep}/k)$ to 
$\Z/2\Z$.  Suppose that $\chi:\gal(k^{\sep}/k)\to \Z/2\Z$ 
is the non-trivial homomorphism. The kernel of $\chi$ 
is a closed subgroup of 
$\gal(k^{\sep}/k)$ of index $2$.  Let $F/k$ be the 
corresponding quadratic Galois extension.  
If the element of $\h^1(k,\Z/2\Z)$ which corresponds to 
$\chi$ is represented by 
a 1-cocycle of the form $\{g^{-1}g^{\sig}\}_{\sig\in}$ 
where $g\in G_{k^{\sep}}$ then $g^{-1}g^{\sig}$ has to be $1$ 
for $\sig\in \kernel(\chi)$.  Therefore, $g\in G_F$.  

We look for $g\in G_F$ such that $g^{-1}g^{\sig}=\tau$ 
for the non-trivial element $\sig$ of $\gal(F/k)$. 
Since $\ch(k)\not=2$, there exists $a\in k\setminus k^2$ 
such that $F=k(\sqrt a\,)$.  
Let 
\begin{align*}
g_F(a) & = \left(
\begin{pmatrix}
I_2 & I_2 \\
\sqrt a I_2 & -\sqrt a I_2 
\end{pmatrix}, 
\frac 12\begin{pmatrix}
1 & 0 & 1 \\
0 & \sqrt a & 0 \\
\sqrt a  & 0 & -\sqrt a 
\end{pmatrix}, 
-\frac 12 \sqrt a
\right), \\
x_F(a) & =(x_{F,1}(a),x_{F,2}(a))\stackrel{\rm{def}}=g_F(a)(A,B). 
\end{align*}
Then 
\begin{align*}
x_{F,1}(a)
& = \left(
\begin{pmatrix}
J & 0 \\
0 & J
\end{pmatrix},
\begin{pmatrix}
0 & -I_2 \\
a I_2 & 0 
\end{pmatrix}, 
\begin{pmatrix}
0 & J \\
a J & 0
\end{pmatrix}
\right), \\ 
x_{F,2}(a) & = a \, \bbmp_{3,1}\wedge \bbmp_{3,3}.
\end{align*}
\begin{prop}
\label{prop:F-characterize}
In the above situation, $F=k(\sqrt{\Phi_2(x_F(a))}\,)$. 
\end{prop}
\begin{proof}
Since $\Phi_2(x_F(a))=(4a)^3\cdot (-2a)^2\Phi_2(A,B)=2^8 a^5$, 
\begin{equation*}
k(\sqrt{\Phi_2(x_F(a))}\,)=k(2^4a^2\sqrt a\,)=k(\sqrt a\,)=F. 
\end{equation*}
\end{proof}

Let $W\sub \m_2(F)$ be the $k$-vector space of 
$2\times 2$ Hermitian matrices over $F$ 
and $H=\gl_1\times \text{R}_{F/k} \gl_2$ where
$\text{R}_{F/k} \gl_2$ is the restriction of the scalar. 
For $h\in\gl_2(F)$, we put $h^*={}^t h^{\sig}$. 
$(t,h)\in H$ acts on $M\in W$ by 
$(t,h)\cdot M= t hM h^*$. 

Let 
\begin{equation*}
\text{GU}_2 = \{h\in \text{R}_{F/k} \gl_2\mid 
{}^{\exists} t\in\gl_1,\; thh^*= I_2\}.  
\end{equation*}
We regard that 
\begin{equation}
\label{eq:GU2-expression}
\text{GU}_2 = \{(t,h)\in \gl_1\times \text{R}_{F/k} \gl_2\mid 
t h h^* = I_2\}.  
\end{equation}
Simple Lie algebra computations show that the stabilizer 
$H_{I_2}$ of $I_2$ is smooth and isomorphic to $\text{GU}_2$ 
(since $\ch(k)\not=2$).  
Therefore, $(H,W)$ is a  regular \pv.

Since $\h^1(k,H)=\{1\}$, $\h^1(k,H_{I_2})\cong \h^1(k,\text{GU}_2)$ 
is in bijective correspondence with 
$H_k\backslash W^{\sst}_k$. If $M_1,M_2\in W_k$ 
are in the same $H_k$-orbit, we say that 
$M_1,M_2$ are equivalent.  

The Gram-Schmidt process shows that any element of $H_k\backslash W_k^{\sst}$
is represented by a diagonal matrix. By the action of the $\gl_1$-factor, 
any element of $H_k\backslash W_k^{\sst}$ is represented by 
a matrix of the form 
\begin{equation*}
\diag(s,1) = \begin{pmatrix}
s & 0 \\
0 & 1 
\end{pmatrix}
\end{equation*}
where $s\in k^{\times}$.  

\begin{prop}
\label{prop:s1s20unit-equivalent}
Suppose that $s_1,s_2\in \mk$ and 
$\diag(s_1,1)$, $\diag(s_2,1)$ are equivalent. 
Then there exists $t\in F^{\times}$ such that 
$s_1=s_2\n_{F/k}(t)$ or $s_1=s_2^{-1}\n_{F/k}(t)$.
\end{prop}
\begin{proof}
By assumption, there exists 
\begin{equation*}
\begin{pmatrix}
a & b \\
c & d
\end{pmatrix}\in \gl_2(F)
\end{equation*}
such that 
\begin{align*}
\begin{pmatrix}
s_1 & 0 \\
0 & 1 
\end{pmatrix}
& = \begin{pmatrix}
a & b \\
c & d
\end{pmatrix}
\begin{pmatrix}
s_2 & 0 \\
0 & 1 
\end{pmatrix}
\begin{pmatrix}
a^{\sig} & c^{\sig} \\
b^{\sig} & d^{\sig}
\end{pmatrix} \\
& = 
\begin{pmatrix}
s_2\n_{F/k}(a) + \n_{F/k}(b) & s_2ac^{\sig} + bd^{\sig} \\
s_2a^{\sig}c + b^{\sig}d & s_2\n_{F/k}(c)+\n_{F/k}(d)
\end{pmatrix}. 
\end{align*}

This implies that $s_2ac^{\sig} =- bd^{\sig}$. 
So 
\begin{equation*}
s_2^2\n_{F/k}(a)\n_{F/k}(c)=\n_{F/k}(b)\n_{F/k}(d). 
\end{equation*}
If $b=0$ then $a,d\not=0$ and $c=0$. Therefore, 
\begin{equation*}
s_1= \frac {s_2\n_{F/k}(a)+\n_{F/k}(b)} {s_2\n_{F/k}(c)+\n_{F/k}(d)}
= s_2\n_{F/k}(a/d).  
\end{equation*}
If $a=0$ then $b,c\not=0$ and $d=0$. Therefore,   
\begin{equation*}
s_1=\frac {s_2\n_{F/k}(a)+\n_{F/k}(b)} {s_2\n_{F/k}(c)+\n_{F/k}(d)}
= \frac {\n_{F/k}(b)} {s_2\n_{F/k}(c)} 
= {s_2^{-1}\n_{F/k}(b/c)}.   
\end{equation*}

Suppose that $a,b,c,d\not=0$. Since 
\begin{equation*}
s_2\n_{F/k}(a)=s_2^{-1}\n_{F/k}(c)^{-1}\n_{F/k}(b)\n_{F/k}(d),  
\end{equation*}
we have 
\begin{align*}
s_1 & = \frac {s_2\n_{F/k}(a)+\n_{F/k}(b)} {s_2\n_{F/k}(c)+\n_{F/k}(d)} 
= \frac {\n_{F/k}(b/c)(\n_{F/k}(c)+s_2^{-1}\n_{F/k}(d))} 
{s_2\n_{F/k}(c)+\n_{F/k}(d)} \\
& = {s_2^{-1}\n_{F/k}(b/c)}.   
\end{align*}
This exhausts all the cases.  
\end{proof}

Conversely, $\diag(s,1)$ and $\diag(s^{-1},1)$ 
are equivalent since 
\begin{equation*}
\begin{pmatrix}
s & 0 \\
0 & 1 
\end{pmatrix}
= s 
\begin{pmatrix}
0 & 1 \\
1 & 0 
\end{pmatrix}
\begin{pmatrix}
s^{-1} & 0 \\
0 & 1 
\end{pmatrix}
\begin{pmatrix}
0 & 1 \\
1 & 0 
\end{pmatrix}. 
\end{equation*}
\begin{lem}
\label{lem:case4-tau-conjugate}
$G_{x_F(a)}/G_{x_F(a)}^{\circ}\cong \Z/2\Z$ 
with the trivial action of the Galois group.  
The non-trivial element of 
$G_{x_F(a)}/G_{x_F(a)}^{\circ}$ 
is represented by 
\begin{equation*}
g_F(a) \tau g_F(a)^{-1}
= \left(
\begin{pmatrix}
I_2 & 0 \\
0 & -I_2 
\end{pmatrix}, 
\begin{pmatrix}
1 & 0 & 0 \\
0 & -1 & 0 \\
0 & 0 & -1
\end{pmatrix}, 
-1 \right). 
\end{equation*}
\end{lem}
\begin{proof}
Since 
\begin{align*}
& \begin{pmatrix}
I_2 & I_2 \\
\sqrt a I_2 & -\sqrt a I_2
\end{pmatrix}
\begin{pmatrix}
0 & I_2 \\
I_2 & 0 
\end{pmatrix}
\begin{pmatrix}
I_2 & I_2 \\
\sqrt a I_2 & -\sqrt a I_2
\end{pmatrix}^{-1} 
= \begin{pmatrix}
I_2 & 0 \\
0 & -I_2 
\end{pmatrix}, \\
& \begin{pmatrix}
1 & 0 & 1 \\
0 & \sqrt a & 0 \\
\sqrt a  & 0 & -\sqrt a 
\end{pmatrix}
\begin{pmatrix}
0 & 0 & 1 \\
0 & -1 & 0 \\
1 & 0 & 0 
\end{pmatrix}
\begin{pmatrix}
1 & 0 & 1 \\
0 & \sqrt a & 0 \\
\sqrt a  & 0 & -\sqrt a 
\end{pmatrix}^{-1} 
= \begin{pmatrix}
1 & 0 & 0 \\
0 & -1 & 0 \\
0 & 0 & -1
\end{pmatrix}, 
\end{align*}
the lemma follows.
\end{proof}

Suppose that $h\in\gl_2(F),t\in F^{\times}$. 
Since $g_F(a)^{\sig}=g_F(a)\tau$, 
\begin{equation*}
g_F(a) \phi(t,h)g_F(a)^{-1} \in G_{x_F(a)\, k}
\end{equation*}
if and only if 
\begin{equation*}
h^{\sig} = {}^t h^{-1}t^{-1}, \; 
t^{\sig} = t.  
\end{equation*}
So $t\in k^{\times}$ and $t h h^* = I_2$. 
Note that $t h h^* = I_2$ if and only if 
${}^th h^{\sig} = t^{-1}I_2$.  
Therefore, $G_{x_F(a)}^{\circ}\cong \text{GU}_2$. 
Even though we only considered $k$-rational points, 
a similar consideration involving 
\cite[p.17, Theorem]{mum} as in 
\cite[p.124, Proposition (2.14)]{yukiel}
shows that it is an isomorphism of algebraic groups. 
Note that if $t\in F^{\times}$ and $h\in\gl_2(F)$, 
then the two diagonal blocks of the first 
component of $\phi(t,h)$ corresponds to 
the two $\gl_2(F)$ components of 
$\text{GU}_2(F)\cong F^{\times}\times \gl_2(F)\times \gl_2(F)$.

%

We would like to find the orbit 
$G_k\backslash \rv^{\sst}$ which 
corresponds to $A(s) = \diag(s,1)$. 

\begin{itemize}
\item[(i)]
We find an element $B(s)\in H_F$ 
such that $B(s) \cdot I_2 = \diag(s,1)\in W$. 
Then $\diag(s,1)$ corresponds to the 
cohomology class represented by 
the map $\Z/2\Z\to \operatorname{GU}_2(F)$ 
whose value at $\sig$ is $\{B(s)^{-1}B(s)^{\sig}\}$. 
\item[(ii)]
Then we find $C(s)\in G_F$ such that 
\begin{math}
C(s)^{-1}C(s)^{\sig}
\end{math}
corresponds to $\{B(s)^{-1}B(s)^{\sig}\}$. 
\item[(iii)]
Lemma \ref{lem:Gx-exact} (3) implies that 
$C(s)x_F(a)$ is the element which corresponds 
to $\diag(s,1)$. 
\end{itemize}

\vskip 5pt
\noindent
{\bf Step (i)} 

Note that 
\begin{equation*}
H_{F} \cong F^{\times} \times \gl_2(F) \times \gl_2(F), \; 
H_{I_2\, F} \cong \operatorname{GU}_2(F) \cong F^{\times} \times \gl_2(F)
\end{equation*}
where $H_{I_2\, F}$ is the group of $F$-rational points 
of the stabilizer of $I_2$. 
\begin{equation*}
\operatorname{GU}_2(k)\sub \operatorname{GU}_2(F)\cong 
F^{\times}\times \gl_2(F)
\end{equation*}
and the inclusion map 
\begin{math}
\operatorname{GU}_2(F) \to F^{\times} \times \operatorname{R}_{F/k}(F) 
\cong F^{\times} \times \gl_2(F)\times \gl_2(F) 
\end{math}
is given by 
\begin{equation*}
F^{\times} \times \gl_2(F) \ni (t,h) \mapsto 
(t,h,t^{-1} {}^t h^{-1}) \in 
F^{\times} \times \gl_2(F)\times \gl_2(F).  
\end{equation*}

The action of the non-trivial element $\sig\in\gal(F/k)$ 
on $H_F$ is 
\begin{equation*}
F^{\times}\times \gl_2(F)\times \gl_2(F)\ni 
(t,h_1,h_2) \mapsto (t^{\sig},h_2^{\sig},h_1^{\sig})
\in F^{\times}\times \gl_2(F)\times \gl_2(F).
\end{equation*}
It induces an action of $\sig$ on 
$\text{GU}_2(F)$ as follows:
\begin{equation*}
F^{\times} \times \gl_2(F)\ni (t,h) \mapsto 
(t^{\sig},(t^{\sig}h^*)^{-1}).  
\end{equation*}
Therefore, $\text{GU}_2(k)=\{(t,h)\in \mk\times \gl_2(F)\mid thh^*=I_2\}$
(see (\ref{eq:GU2-expression})).

We can identify $W\otimes_k F \cong \m_2(F)$.  
With this identification, 
the action of $(t,h_1,h_2)\in H_F$ 
on $W\otimes_k F \cong \m_2(F)$ is 
\begin{equation*}
\m_2(F) \ni M \mapsto  t h_1 M {}^t h_2 \in \m_2(F). 
\end{equation*}
The action of $\sig$ on $W\otimes_k F \cong \m_2(F)$ 
is $M\mapsto {}^t \! M^{\sig}$.

Let $s\in \mk$, $A(s)=\diag(s,1)$ and 
$B(s)=(1,A(s),I_2)\in H_F$. Then 
$B(s)\cdot I_2 = A(s)$.  Therefore, 
this $B(s)$ is what we were looking for.

\vskip 5pt
\noindent
{\bf Step (ii)} 

The group $G_{x_F(a)\, F}^{\circ}$ of $F$-rational points 
of $G_{x_F(a)}^{\circ}$ consists of elements of the form 
\begin{equation*}
g_F(a) \phi(t,h) g_F(a)^{-1}
\end{equation*}
for $h\in \gl_2(F),t\in F^{\times}$.  
Define a map $\Psi:\operatorname{GU}_2(F)\to H_F$
by 
\begin{equation*}
\Psi(t,h) = g_F(a) \phi(t,h) g_F(a)^{-1}
\end{equation*}
for $t\in F^{\times},h\in\gl_2(F)$. 
Then $\Psi$ is obviously a homomorphism. 

Since 
\begin{align*}
\Psi(t,h)^{\sig} 
& = g_F(a) \tau \phi(t^{\sig},h^{\sig}) \tau g_F(a)^{-1}
= g_F(a) \tau \phi(t^{\sig},h^{\sig}) \tau g_F(a)^{-1} \\
& = g_F(a) \phi(t^{\sig},(t^{\sig}h^*)^{-1}) g_F(a)^{-1}
= \Psi(t^{\sig},(t^{\sig}h^*)^{-1}), 
\end{align*}
$\Psi$ is equivariant with respect to the action of 
$\gal(F/k)$.  Since $s^{\sig}=s$, 
\begin{equation*}
B(s)^{-1}B(s)^{\sig} = (1,A(s^{-1}),A(s))
\end{equation*}
and this is the first component of 
$\phi(1,A(s^{-1}))$.  Therefore, 
$\Psi(1,A(s^{-1}))$ corresponds to the 
1-cocycle which is associated with  
the Hermitian matrix $A(s)$.  

\vskip 5pt
\noindent
{\bf Step (iii)} 

Let 
\begin{equation*}
C(s) = g_F(a) 
\left(
\begin{pmatrix}
A(s) & \\
& I_2
\end{pmatrix}, 
\begin{pmatrix}
s^{-1} & 0 & 0 \\
0 & 1 & 0 \\
0 & 0 & 1 
\end{pmatrix},
1 \right)
g_F(a)^{-1}. 
\end{equation*}
Since $s^{\sig}=s$, 
\begin{align*}
C(s)^{-1}C(s)^{\sig}
& = g_F(a) 
\left(
\begin{pmatrix}
A(s^{-1}) & \\
& I_2 
\end{pmatrix}, 
\begin{pmatrix}
s & 0 & 0 \\
0 & 1 & 0 \\
0 & 0 & 1 
\end{pmatrix}, 
1 \right)
g_F(a)^{-1} \\
& \quad 
\times 
g_F(a) \tau 
\left(
\begin{pmatrix}
A(s) & \\
& I_2 
\end{pmatrix}, 
\begin{pmatrix}
s^{-1} & 0 & 0 \\
0 & 1 & 0 \\
0 & 0 & 1 
\end{pmatrix}, 
1 \right)
\tau 
g_F(a)^{-1} \\
& = g_F(a) 
\left(
\begin{pmatrix}
A(s^{-1}) & \\
& A(s)
\end{pmatrix}, 
\begin{pmatrix}
s & 0 & 0 \\
0 & 1 & 0 \\
0 & 0 & s^{-1} 
\end{pmatrix}, 
1 \right)
g_F(a)^{-1} \\
& = g_F(a) \phi(1,A(s^{-1})) g_F(a)^{-1} 
= \Psi(1,A(s^{-1})).  
\end{align*}
Therefore, $\rg\cdot C(s)x_F(a)$ is the orbit which 
corresponds to the Hermitian matrix $A(s)$. 
Since 
\begin{equation*}
g_F(a)\tau g_F(a)^{-1} 
C(s^{-1})C(s)^{\sig} (g_F(a)\tau g_F(a)^{-1})^{-1}
= C(s)C(s^{-1})^{\sig}, 
\end{equation*}
The action of $g_F(a)\tau g_F(a)^{-1}$ is 
$C(s^{-1})C(s)^{\sig}\mapsto C(s)C(s^{-1})^{\sig}$. 
So the representatives of $\gam_V^{-1}(F)$ are 
\begin{equation*}
x_F(a,s) = C(s^{-1})x_F(a)
\end{equation*}
where $x_F(a,s),x_F(a,s^{-1})$ are identified.  
Explicitly, 
$x_F(a,s)=(x_{F,1}(a,s),x_{F,2}(a,s))$ where 
\begin{align*}
x_{F,1}(a,s) & = 
\left(
\begin{pmatrix}
J & 0 \\
0 & J 
\end{pmatrix},
\begin{pmatrix}
0 & -A(s) \\
aA(s) & 0 
\end{pmatrix}
\begin{pmatrix}
0 & J \\
aJ & 0 
\end{pmatrix}
\right), \\
x_{F,2}(a,s) & = as \bbmp_{3,1}\wedge \bbmp_{3,3}.    
\end{align*}

Let 
\begin{equation*}
\Gunit_2(k) = \{k\} \coprod \bigcup_F \left( 
(\Z/2\Z)\backslash \mk/\n_{F/k}(F^{\times})
\right)
\end{equation*}
where $k$ is the trivial extension $k/k$ and 
$F$ runs through all the quadratic 
extensions of $k$ and $\Z/2\Z$ acts on $\mk$ 
by $s\mapsto s^{-1}$. 
Note that since $\ch(k)\not=2$, 
any quadratic extension is a Galois extension.  
Let $\gam:\Gunit_2(k)\to \Ex_2(k)$ be the map such that 
$\gam(k)=k$ and $\gam((\Z/2\Z)\cdot s \cdot \n_{F/k}(F^{\times}))=F$. 

We obtain the following theorem 
by the above consideration. 

\begin{thm}
\label{thm:case4-rational-orbits}
$\rg \backslash V^{\sst}$ is in bijective correspondence 
with $\Gunit_2(k)$ where $\rg\cdot w$ corresponds to 
the trivial extension $k/k$ and if $F=k(\sqrt a)$ 
is a quadratic extension of $k$ then
$\rg\cdot x_{F}(a,s)$ corresponds to 
$s\in\mk$ in $(\Z/2\Z)\backslash \mk/\n_{F/k}(F^{\times})$. 
\end{thm}

\section{Non-empty strata}
\label{sec:non-empty}

In this section and the next, $G=\gl_1\times \gl_8$.  
The set $\gB$ which parametrizes the GIT stratification 
for (\ref{eq:PV}) consists of $\be_1\ccd \be_{183}$
of the table in Section I--8. 
We shall prove that $S_{\be_i}\not=\emptyset$ for 
\begin{equation}
\label{eq:list-non-empty-tri8}
\begin{aligned}
i & = 1,4,7,9,13,32,33,37,63,74,81,98,108, \\
& \quad \hskip 4pt 152, 154,172,173,179,180,182,183 
\end{aligned}
\end{equation}
for the \pv{} (\ref{eq:PV}) in this section. We shall 
prove that $S_{\be}=\emptyset$ for other $\be_i$'s in the next section

We explain how we proceed to show that $S_{\be_i}\not=\emptyset$ 
and to describe $\rg\backslash S_{\be_i\, k}$.

When we consider individual cases, we first 
prove that $S_{\be_i}\not=\emptyset$ regardless of $\ch(k)$. 
It is guaranteed by the general theory that 
$S_{\be_i\, k}$ has an inductive structure 
\begin{equation}
\label{eq:S-beta-ind-structure}
S_{\be_i\, k}\cong \rg\times_{P_{\be_i}} Y_{\be_i\, k}^{\sst} 
\end{equation}
and 
\begin{equation*}
Y_{\be_i\, k}^{\sst} =\{(z,w)\in Z_{\be_i\, k}\oplus W_{\be_i\, k}\mid 
z\in Z_{\be_i\, k}^{\sst}\}. 
\end{equation*}
So it is enough to show that $Z_{\be_i\, k}^{\sst}\not=\emptyset$.  

There are two types of definitions of \pv s. 
One only requires the existence of a Zariski open orbit 
and a relative invariant polynomial.  
Definition III--4.1 is such a definition. 
The other definition fixes the character which is related 
to the relative invariant polynomial.  
Definition III--4.8  is such a definition. 
The formulation of the GIT stratification uses the latter definition 
and the representation does not have to be prehomogeneous.  
If the \rep{} $(M_{\be_i},Z_{\be_i})$ is irreducible, 
the two definitions are essentially the same.  However, 
If the \rep{} $(M_{\be_i},Z_{\be_i})$ is reducible, 
one has to find a relative invariant polynomial 
$P(x)$ on $Z_{\be_i}$ such that $P(gx)=\chi(g)P(x)$ 
and $\chi(g)$ is proportional (see Notation for the definition) 
to the character $\chi_{\be_i}(g)$ 
(see Section \ref{sec:notation-related-git} for the definition).    
We remind the reader that the stability of points in 
$Z_{\be_i}$ is defined by the group action of 
$G_{\be_i}^{\st}$, ignoring the direction of $\be_i$ 
from $M_{\be_i}^{\st}$.  So if $P(x)$ is as above and $\chi(g)$
is proportional to $\chi_{\be_i}(g)$ then 
$P(x)$ is an invariant polynomial with respect to 
the action of $G_{\be_i}^{\st}$.  

The next step is to describe the set of orbits 
$M_{\be_i\, k}\backslash Z_{\be_i\, k}^{\sst}$.  
In this step, results of Sections \ref{sec:orbits1}, 
\ref{sec:orbits2} may be used.  
The set $M_{\be_i\, k}\backslash Z_{\be_i\, k}^{\sst}$ 
will be described in the individual consideration.

Because of the inductive structure (\ref{eq:S-beta-ind-structure}), 
\begin{math}
\rg\backslash S_{\be_i\, k}
\cong P_{\be_i\, k}\backslash Y_{\be_i\, k}^{\sst}. 
\end{math}
It turns out that 
\begin{math}
P_{\be_i\, k}\backslash Y_{\be_i\, k}^{\sst}
\cong M_{\be_i\, k}\backslash Z_{\be_i\, k}^{\sst}
\end{math}
for all the cases.  To show this, we use the argument of 
\cite[pp.264,265]{tajima-yukie-GIT3}. 
We have to choose $R_i\in Z_{\be_i\, k}^{\sst}$ 
so that the following condition is satisfied 
(see Condition III--12.6). 

\begin{cond}
\label{cond:non-empty-strata-sorbit}
\begin{itemize}
\item[(1)]
$Z_{i\, k^{\sep}}= M_{\be_i \, k^{\sep}}R_i$.  
\item[(2)]
If $y\in W_{i\, k^{\sep}}$ then there exists 
$n\in U_{\be_i\, k^{\sep}}$ such that 
$nR_i = (R_i,y)$. 
\item[(3)]
$G_{R_i}\cap U_{\be_i}$ is connected.  
\end{itemize}
\end{cond}

We remind the reader that $k$ is a perfect field, 
which implies that connected unipotent groups are 
split (see \cite[p.205, 15.5 Corollary (ii)]{borelb}).  
The following proposition 
was proved in Proposition III--12.7 
using Hilbert's theorem 90.  

\begin{prop}
\label{prop:Y2Z}
If Condition \ref{cond:non-empty-strata-sorbit} 
is satisfied, for any $x\in Z_{i\, k}^{\sst}$ and 
$y\in W_{i\, k}$, there exists $n\in U_{\be_i\, k}$ 
such that $nx = (x,y)$.  
\end{prop}

The above proposition implies that 
the map 
\begin{equation*}
M_{\be_i\, k}\backslash Z_{\be_i\, k}^{\sst} 
\to P_{\be_i\, k}\backslash Y_{\be_i\, k}^{\sst} 
\end{equation*}
is surjective.  Since $U_{\be_i}$ does not change 
the $Z_i$ component of $Z_{\be_i}\oplus W_{\be_i}$, 
the above map is injective also. Therefore, 
\begin{equation*}
\rg\backslash S_{\be_i\, k}\cong 
M_{\be_i\, k}\backslash Z_{\be_i\, k}^{\sst}.  
\end{equation*}

If $(M_{\be_i},Z_{\be_i})$ is a regular \pv{},  
$Z_{\be_i}$ is a direct sum $W_1\oplus \cdots \oplus W_N$
of subrepresentations of $M_{\be_i}$ and 
there are $N$ independent relative invariant 
polynomials $\Del_1(x)\ccd \Del_N(x)$ of $x\in Z_{\be_i}$,  
then Condition \ref{cond:non-empty-strata-sorbit} (1) 
is satisfied if $\Del_1(R_i)\ccd \Del_N(R_i)\not=0$ 
(Corollary 4.7).  
In the case where we cannot find $N$ independent 
relative invariant polynomials, we show 
Condition \ref{cond:non-empty-strata-sorbit} (1) 
explicitly.

Let $m=\dim W_{\be_i}$. 
It turns out that $\dim U_{\be_i}\geqq m$ for all the cases.  
Put $\dim U_{\be_i}=n+m$.   
We express elements of $U_{\be_i}$ as $n(u)$ 
where $u=(u_{ij})$. 
We number $u_{ij}$'s as $u^{(1)}\ccd u^{(n+m)}$. 
$W_i$ has a basis consisting of 
elements of the form $e_{i_1i_2i_3}$. 
We fix such a basis and consider the coordinate
with respect to the basis chosen.  
To show Condition \ref{cond:non-empty-strata-sorbit} (2), (3), 
we show the following condition.  

\begin{cond}
\label{cond:unipotent-connected}
Let $n(u)R_i=(R_i,y_1(u)\ccd y_m(u))$. 
By changing the order of $y_1\ccd y_m$ 
and $u^{(1)}\ccd u^{(n+m)}$ if necessary, 
\begin{equation*}
y_i(u) = \pm u^{(n+i)} + F_i(u)
\end{equation*}
where $F_i(u)$ is a polynomial of 
$u^{(1)}\ccd u^{(n+i-1)}$.  
\end{cond}
 
If Condition \ref{cond:unipotent-connected} is satisfied, 
Condition \ref{cond:non-empty-strata-sorbit} (2) is obviously 
satisfied. Since $n(u)\in G_{R_i}\cap U_{\be_i}$ 
if and only if $F_1(u)=\cdots=F_n(u)=0$, 
Condition \ref{cond:non-empty-strata-sorbit} (3) 
is also satisfied (see Lemma III--12.9). 
Therefore, our strategy is to show
that Condition \ref{cond:unipotent-connected} is satisfied 
for individual cases.

The following table describes $M_{\be_i}$, $Z_{\be_i}$ as a 
\rep{} of $M^s_{\be_i}$ (the semi-simple part of $M_{\be_i}$), 
the coordinates of $Z_{\be_i},W_{\be_i}$ 
and 
\begin{math}
G_k\backslash S_{\be_i\, k}\cong P_{\be_i\,k}\backslash Y^{\sst}_{\be_i\,k}
\end{math}
for all of the above $i$'s.   
We carry out the above plan for individual cases 
after the table.

\vskip 10pt

\begin{center} 

\begin{tabular}{|c|l|}
\hline 
\size $i$ & \hskip 1.2in 
\size $M_{\be_i}$, \quad $Z_{\be_i}$ 
as a \rep{} of the semi-simple part of $M_{\be_i}$ \\
\hline
\size $G_k\backslash S_{\be_i\, k}$  
&  \hskip 1.2in 
\size coordinates of $Z_{\be_i}$,
coordinates of $W_{\be_i}$ \\
\hline \hline 

\size $1$ & 
\size $M_{[3]}\cong \gl_5\times \gl_3$
,\quad  \size $\Lam^{3,1}_{[1,3]} \otimes \Lam^{5,2}_{[4,8]}$ \\
\hline 
\size  \rule[-24pt]{0cm}{52pt}
$\mathrm{Prg}_2(k)\atop \ch(k)\not=2$ & 
\size 
\begin{tabular}{|c|c|c|c|c|c|c|c|c|c|c|c|c|c|c|c|c|c|c|c|c|c|c|c|}
\hline
 $12,13,14,15,16,17,18,19,20,21$  \\
\hline  
$x_{145},x_{146},x_{147},x_{148},x_{156},x_{157},x_{158},x_{167},
x_{168},x_{178}$  \\
\hline 
$27,\dots, 46$ \\
\hline 
$x_{245},\dots ,x_{378}$ \\
\hline
\end{tabular}
\begin{tabular}{|c|c|c|c|c|c|c|c|c|c|c|c|c|c|c|c|c|c|c|c|c|c|c|c|}
\hline
$47,48,49,50,51$ \\  
\hline
$x_{456},x_{457},x_{458},x_{467},x_{468}$ \\
\hline 
$52,53,54,55,56$ \\
\hline 
$x_{478},x_{567},x_{568},x_{578},x_{678}$ \\
\hline
\end{tabular}

\\
\hline \hline

\size $4$ & 
 \size $M_{[1,5]}\cong \gl_4\times \gl_3 \times \gl_1$
,\quad 
 \size $\Lam^{3,2}_{[6,8]}\oplus \Lam^{4,2}_{[2,5]}\otimes \Lam^{3,1}_{[6,8]}$ \\

\hline 
\size $\mathrm{Gunit}_2(k)\atop \ch(k)\not=2$ 
\rule[-24pt]{0cm}{52pt} & 
\size \begin{tabular}{|c|c|c|c|c|c|c|c|c|c|c|c|c|c|c|c|c|c|c|c|c|c|c|c|}
\hline
$19,20,21$ & $24,25,26,28,29,30$ \\ 
\hline
$x_{167},x_{168},x_{178}$ & 
$x_{236},x_{237},x_{238},x_{246},x_{247},x_{248}$ \\
\hline 
& $31,\dots ,49$ \\
\hline
& $x_{256},\dots ,x_{458}$ \\
\hline 
\end{tabular}

\begin{tabular}{|c|c|c|c|c|c|c|c|c|c|c|c|c|c|c|c|c|c|c|c|c|c|c|c|}
\hline
$34,35,36,44,45,46$ & $56 $  \\
\hline
$x_{267},x_{268},x_{278},x_{367},x_{368},x_{378}$ & $x_{678}$\\
\hline 
$50,51,52,53,54,55$ & \\
\hline 
$x_{467},x_{468},x_{478},x_{567},x_{568},x_{578}$ & \\
\hline
\end{tabular}
\\
\hline \hline

\size $7$ & 
\size $M_{[6]}\cong \gl_6\times \gl_2$
,\quad 
 \size $\Lam^{6,2}_{[1,6]}\otimes \Lam^{2,1}_{[7,8]}$ \\

\hline 
\size $\Ex_3$ \rule[-24pt]{0cm}{52pt} & 
\size 
\begin{tabular}{|c|c|c|c|c|c|c|c|c|c|c|c|c|c|c|c|c|c|c|c|c|c|c|c|}
\hline
$5,6,10,11,14,15,17,18,19,20$  \\
\hline
$x_{127},x_{128},x_{137},x_{138},x_{147},x_{148},x_{157},x_{158},
x_{167},x_{168}$ \\
\hline
$25,26,\dots,54$ \\ 
\hline
$x_{237},x_{238},\dots ,x_{568}$ \\
\hline
\end{tabular}

\begin{tabular}{|c|c|c|c|c|c|c|c|c|c|c|c|c|c|c|c|c|c|c|c|c|c|c|c|}
\hline
$21,36,46,52,55,56$  \\
\hline
$x_{178},x_{278},x_{378},x_{478},x_{578},x_{678}$ \\
\hline
\end{tabular}

\\
\hline \hline 

\size $9$ & 

\size $M_{[1]}\cong \gl_7\times \gl_1$,
\quad  \size $\Lam^{7,3}_{[2,8]}$ \\

\hline 
\size  $\text{G}_2(k)$ \rule[-14pt]{0cm}{32pt} & 
\size 
\begin{tabular}{|c|c|c|c|c|c|c|c|c|c|c|c|c|c|c|c|c|c|c|c|c|c|c|c|}
\hline
$22,23,24,25,26,27,28,29,30,31,\dots ,56 $ \\
\hline
$x_{234},x_{235},x_{236},x_{237},x_{238},x_{245},x_{246},x_{247},
x_{248},x_{256},\dots ,x_{678} $\\
\hline
\end{tabular}
\quad 
none

\\
\hline 
\end{tabular}

\begin{tabular}{|c|l|}
\hline 

\size $13$ & 
 \size $M_{[1,4,7]}\cong \gl_3^2\times \gl_1^2$
,\quad
 \size $\Lam^{3,1}_{[5,7]}\oplus \Lam^{3,2}_{[2,4]}\oplus 
 \Lam^{3,1}_{[2,4]}\otimes \Lam^{3,2}_{[5,7]}$ \\

\hline 
\size SP \rule[-24pt]{0cm}{52pt} & 
\size 
\begin{tabular}{|c|c|c|c|c|c|c|c|c|c|c|c|c|c|c|c|c|c|c|c|c|c|c|c|}
\hline
$18,20,21$ & $26,30, 40$  \\
\hline
$x_{158},x_{168},x_{178}$ & $x_{238},x_{248}, x_{348}$ \\
\hline
\multicolumn{2}{|c|}{
$31,32,34,41,\dots,50$} \\
\hline 
\multicolumn{2}{|c|}{
$x_{256},x_{257},x_{267},x_{356},\dots, x_{467}$} \\
\hline
\end{tabular}

\begin{tabular}{|c|c|c|c|c|c|c|c|c|c|c|c|c|c|c|c|c|c|c|c|c|c|c|c|}
\hline
$33,35,36,43,\dots ,52$ & $53$ & $54,55,56 $  \\
\hline
$x_{258},x_{268},x_{278},x_{358},\dots ,x_{478}
$ & $x_{567}$ & $x_{568},x_{578},x_{678}$\\
\hline
\end{tabular}

\\
\hline \hline

\size $32$ & 
 \size $M_{[4]}\cong \gl_4^2$
,\quad
 \size $\Lam^{4,1}_{[1,4]}\otimes \Lam^{4,2}_{[5,8]}$ \\

\hline 
\size $\Ex_2(k)$ \rule[-14pt]{0cm}{32pt} & 
\size 
\begin{tabular}{|c|c|c|c|c|c|c|c|c|c|c|c|c|c|c|c|c|c|c|c|c|c|c|c|}
\hline
$16,17,18,19,20,21,31,32,33,\dots ,52 $  \\
\hline
$x_{156},x_{157},x_{158},x_{167},x_{168},x_{178},x_{256},x_{257},x_{258},
\dots ,x_{478}  $\\
\hline
\end{tabular}

\begin{tabular}{|c|c|c|c|c|c|c|c|c|c|c|c|c|c|c|c|c|c|c|c|c|c|c|c|}
\hline
$53,54,55,56 $  \\
\hline
$x_{567},x_{568},x_{578},x_{678} $\\
\hline
\end{tabular}

\\
\hline \hline 

\size $33$ & 
\size $M_{[1,3,6]}\cong \gl_3\times \gl_2^2\times \gl_1$
,\quad
 \size $1\oplus \Lam^{2,1}_{[2,3]}\otimes \Lam^{3,1}_{[4,6]}
 \otimes \Lam^{2,1}_{[7,8]}\oplus 1$ \\

\hline 
\size SP \rule[-24pt]{0cm}{52pt} & 
\size 
\begin{tabular}{|c|c|c|c|c|c|c|c|c|c|c|c|c|c|c|c|c|c|c|c|c|c|c|c|}
\hline
$21$ & $29,30,32,33,34,35$ & $47$ \\
\hline
$x_{178}$ & $x_{247},x_{248},x_{257},x_{258},x_{267},x_{268}$ 
& $x_{456} $\\
\hline
& $39,40,42,43,44,45$ & \\
\hline
& $x_{347},x_{348},x_{357},x_{358},x_{367},x_{368}$ & \\
\hline
\end{tabular}

\begin{tabular}{|c|c|c|c|c|c|c|c|c|c|c|c|c|c|c|c|c|c|c|c|c|c|c|c|}
\hline
$36,46$ & $54,55,56$ \\
\hline
$x_{278},x_{378}$ & $x_{478},x_{578},x_{678}$ \\
\hline
\multicolumn{2}{|c|}{$48,49,50,51,52,53$} \\
\hline
\multicolumn{2}{|c|}{
$x_{457},x_{458},x_{467},x_{468},x_{567},x_{568}$} \\
\hline
\end{tabular}

\\

\hline \hline

\size $37$ & 
\size $M_{[3,6]}\cong \gl_3^2\times \gl_2$
,\quad \size $\Lam^{3,1}_{[1,3]}\otimes \Lam^{3,1}_{[4,6]}\otimes \Lam^{2,1}_{[7,8]}\oplus 1$ \\

\hline 
\size $\Ex_3(k)$ \rule[-24pt]{0cm}{52pt} & 
\size 
\begin{tabular}{|c|c|c|c|c|c|c|c|c|c|c|c|c|c|c|c|c|c|c|c|c|c|c|c|}
\hline
$14,15,17,18,19,20,29,30,\dots ,45 $ & $47 $ \\
\hline
$x_{147},x_{148},x_{157},x_{158},x_{167},x_{168},x_{247},x_{248},
\dots ,x_{368}$ & $x_{456}$  \\
\hline
\end{tabular}

\begin{tabular}{|c|c|c|c|c|c|c|c|c|c|c|c|c|c|c|c|c|c|c|c|c|c|c|c|}
\hline
$21,36,46$ &  $52,55,56$  \\
\hline
$x_{178},x_{278},x_{378}$ & $x_{478},x_{578},x_{678}$ \\
\hline
\multicolumn{2}{|c|}{$48,49,50,51,53,54$} \\
\hline
\multicolumn{2}{|c|}{
$x_{457},x_{458},x_{467},x_{468},x_{567},x_{568}$} \\
\hline
\end{tabular}

\\
\hline \hline

\size $63$ & 
\size $M_{[1,7]}\cong \gl_6\times \gl_1^2$
,\quad \size $\Lam^{6,2}_{[2,7]}$ \\

\hline 
\size SP \rule[-14pt]{0cm}{32pt} & 
\size 
\begin{tabular}{|c|c|c|c|c|c|c|c|c|c|c|c|c|c|c|c|c|c|c|c|c|c|c|c|}
\hline
$26,30,33,35,36,40,\dots ,56 $  \\
\hline
$x_{238},x_{248},x_{258},x_{268},x_{278},x_{348},\dots ,x_{678} $\\
\hline
\end{tabular}

\quad 
none

\\
\hline \hline

\size $74$ & 
\size $M_{[2,6,7]}\cong \gl_2\times \times \gl_4 \times \gl_1^2$
,\quad \size 
$\Lam^{2,1}_{[1,2]}\otimes \Lam^{4,1}_{[3,6]} \oplus \Lam^{4,2}_{[3,6]}$ \\

\hline 
\size SP \rule[-24pt]{0cm}{52pt} & 
\size 
\begin{tabular}{|c|c|c|c|c|c|c|c|c|c|c|c|c|c|c|c|c|c|c|c|c|c|c|c|}
\hline
$11,15,18,20,26,30,33,35$ \\
\hline
$x_{138},x_{148},x_{158},x_{168},x_{238},x_{248},x_{258},x_{268}$ \\
\hline
$39,42,44,48,50,53$ \\
\hline
$x_{347},x_{357},x_{367},x_{457},x_{467},x_{567} $\\
\hline
\end{tabular}

\begin{tabular}{|c|c|c|c|c|c|c|c|c|c|c|c|c|c|c|c|c|c|c|c|c|c|c|c|}
\hline
$21,36$ & $46,52,55,56$  \\
\hline
$x_{178},x_{278}$ & $x_{378},x_{478},x_{578},x_{678}$ \\
\hline
\multicolumn{2}{|c|}{
$40,43,45,49,51,54$}  \\
\hline
\multicolumn{2}{|c|}{
$x_{348},x_{358},x_{368},x_{458},x_{468},x_{568}$} \\    
\hline
\end{tabular}

\\
\hline \hline

\size $81$ & 
\size $M_{[1,4,6]}\cong \gl_3\times \gl_2^2\times \gl_1$
,\quad \size $1\oplus \Lam^{3,1}_{[2,4]}\otimes \Lam^{2,1}_{[5,6]}\otimes \Lam^{2,1}_{[7,8]}$ \\

\hline 
\size SP \rule[-14pt]{0cm}{32pt} & 
\size 
\begin{tabular}{|c|c|c|c|c|c|c|c|c|c|c|c|c|c|c|c|c|c|c|c|c|c|c|c|}
\hline
$21$ & $32,33,34,35,42,43,\dots ,51 $  \\
\hline
$x_{178}$ & 
$
x_{257},x_{258},
x_{267},x_{268},
x_{357},x_{358},
\dots ,x_{468} $\\
\hline
\end{tabular}

\begin{tabular}{|c|c|c|c|c|c|c|c|c|c|c|c|c|c|c|c|c|c|c|c|c|c|c|c|}
\hline
$36,46,52$ & $53,54$ & $55,56 $ \\
\hline
$x_{278},x_{378},x_{478}$ & 
$x_{567},x_{568}$ &
$x_{578},x_{678} $\\
\hline
\end{tabular}

\\
\hline \hline

\size $98$ & 
\size $M_{[2]}\cong \gl_6\times \gl_2$
,\quad \size $\Lam^{6,3}_{[3,8]}$ \\

\hline 
\size $\Ex_2(k)$ \rule[-14pt]{0cm}{32pt} & 
\size 
\begin{tabular}{|c|c|c|c|c|c|c|c|c|c|c|c|c|c|c|c|c|c|c|c|c|c|c|c|}
\hline
$37,38,39,40,41,42,43,44,45,\dots ,56 $ \\
\hline
$x_{345},x_{346},x_{347},x_{348},x_{356},x_{357},x_{358},x_{367},
x_{368},\dots ,x_{678} $\\
\hline
\end{tabular}
\quad 
none

\\
\hline \hline

\size $108$ & 
\size $M_{[1,4]}\cong \gl_4\times \gl_3 \times \gl_1$
,\quad \size $\Lam^{3,1}_{[2,4]}\otimes \Lam^{4,2}_{[5,8]}$ \\

\hline 
\size $\mathrm{IQF}_4(k)$ \rule[-24pt]{0cm}{52pt} & 
\size 
\begin{tabular}{|c|c|c|c|c|c|c|c|c|c|c|c|c|c|c|c|c|c|c|c|c|c|c|c|}
\hline
$31,32,33,34,35,36$  \\
\hline
$x_{256},x_{257},x_{258},x_{267},x_{268},x_{278}$ \\
\hline
$41,42,43,44,45,46,47\ccd 52$ \\
\hline
$x_{356},x_{357},x_{358},x_{367},x_{368},x_{378},x_{456}\ccd x_{478}$ \\
\hline
\end{tabular}

\begin{tabular}{|c|c|c|c|c|c|c|c|c|c|c|c|c|c|c|c|c|c|c|c|c|c|c|c|}
\hline
$53,54,55,56 $  \\
\hline
$x_{567},x_{568},x_{578},x_{678}$\\
\hline
\end{tabular}

\\
\hline \hline

\size $152$ & 
\size $M_{[2,4,7]}\cong \gl_2^2 \times \gl_3\times \gl_1$
,\quad \size 
$\Lam^{2,1}_{[1,2]}\otimes \Lam^{3,1}_{[5,7]}
\oplus 1\oplus \Lam^{2,1}_{[3,4]}\otimes \Lam^{3,2}_{[5,7]}$ \\

\hline 
\size SP \rule[-24pt]{0cm}{52pt} & 
\size 
\begin{tabular}{|c|c|c|c|c|c|c|c|c|c|c|c|c|c|c|c|c|c|c|c|c|c|c|c|}
\hline
$18,20,21,33,35,36$ & $40$ \\
\hline
$x_{158},x_{168},x_{178},x_{258},x_{268},x_{278}$ & $x_{348}$ \\
\hline
\multicolumn{2}{|c|}{
$41,42,44,47,48,50$}  \\
\hline
\multicolumn{2}{|c|}{
$x_{356},x_{357},x_{367},x_{456},x_{457},x_{467}$} \\
\hline
\end{tabular}

\begin{tabular}{|c|c|c|c|c|c|c|c|c|c|c|c|c|c|c|c|c|c|c|c|c|c|c|c|}
\hline
\multicolumn{2}{|c|}{$43, 45, 46, 49, 51, 52$} \\ 
\hline
\multicolumn{2}{|c|}{
$x_{358},x_{368},x_{378},x_{458},x_{468},x_{478}$} \\
\hline
$53$ & $54, 55, 56$  \\
\hline
$x_{567}$ & $x_{568},x_{578},x_{678}$ \\
\hline
\end{tabular}

\\
\hline \hline

\size $154$ & 
\size $M_{[1,2,6]}\cong \gl_4\times \gl_2 \times \gl_1^2$
,\quad  \size $1\oplus \Lam^{4,2}_{[3,6]}\otimes \Lam^{2,1}_{[7,8]}$ \\

\hline 
\size $\Ex_2(k)$ \rule[-14pt]{0cm}{32pt} & 
\size 
\begin{tabular}{|c|c|c|c|c|c|c|c|c|c|c|c|c|c|c|c|c|c|c|c|c|c|c|c|}
\hline
$36$ & $39,40,42,43,44,45,48,49,50,51,53,54 $ \\
\hline
$x_{278}$ & $x_{347},x_{348},x_{357},x_{358},x_{367},x_{368},x_{457},x_{458},
x_{467},x_{468},x_{567},x_{568}$\\
\hline
\end{tabular}

\begin{tabular}{|c|c|c|c|c|c|c|c|c|c|c|c|c|c|c|c|c|c|c|c|c|c|c|c|}
\hline 
$46,52,55,56 $  \\
\hline
$x_{378},x_{478},x_{578},x_{678} $\\
\hline
\end{tabular}

\\
\hline 
\end{tabular}

\begin{tabular}{|c|l|}
\hline 

\size $172$ & 
\size $M_{[1,4,7]}\cong \gl_3^2\times \gl_1^2$
,\quad \size $\Lam^{3,1}_{[2,4]}\otimes \Lam^{3,1}_{[5,7]}\oplus 1$ \\

\hline 
\size SP \rule[-14pt]{0cm}{32pt} & 
\size 
\begin{tabular}{|c|c|c|c|c|c|c|c|c|c|c|c|c|c|c|c|c|c|c|c|c|c|c|c|}
\hline
$33,35,36,43,45,46,49,51,52$ & $53 $  \\
\hline
$x_{258},x_{268},x_{278},x_{358},x_{368},x_{378},x_{458},x_{468},x_{478},
$ & $x_{567} $\\
\hline
\end{tabular}

\begin{tabular}{|c|c|c|c|c|c|c|c|c|c|c|c|c|c|c|c|c|c|c|c|c|c|c|c|}
\hline
$54,55,56 $  \\
\hline
$x_{568},x_{578},x_{678} $\\
\hline
\end{tabular}

\\
\hline \hline

\size $173$ & 
\size $M_{[2,5,7]}\cong \gl_3\times \gl_2^2 \times \gl_1$
,\quad \size $\Lam^{2,1}_{[1,2]}\otimes \Lam^{2,1}_{[6,7]}\oplus 
\Lam^{3,1}_{[3,5]}\oplus \Lam^{3,2}_{[3,5]}$ \\

\hline 
\size $\Ex_3(k)$ \rule[-24pt]{0cm}{52pt} & 
\size 
\begin{tabular}{|c|c|c|c|c|c|c|c|c|c|c|c|c|c|c|c|c|c|c|c|c|c|c|c|}
\hline
\multicolumn{2}{|c|}{$20,21,35,36$} \\
\hline
\multicolumn{2}{|c|}{$x_{168},x_{178},x_{268},x_{278}$} \\ 
\hline
$40,43,49$ & $44,50,53$ \\
\hline
$x_{348},x_{358},x_{458}$ & $x_{367},x_{467},x_{567}$ \\
\hline
\end{tabular}

\begin{tabular}{|c|c|c|c|c|c|c|c|c|c|c|c|c|c|c|c|c|c|c|c|c|c|c|c|}
\hline
$45,46,51,52,54,55$ & $56 $  \\
\hline
$x_{368},x_{378},x_{468},x_{478},x_{568},x_{578}$ & $x_{678} $\\
\hline
\end{tabular}

\\
\hline \hline

\size $179$ & 
\size $M_{[2,5]}\cong \gl_2\times \gl_3^2$
,\quad \size $\Lam^{3,1}_{[3,5]}\otimes \Lam^{3,2}_{[6,8]}$ \\

\hline 
\size SP \rule[-14pt]{0cm}{32pt} & 
\size 

\begin{tabular}{|c|c|c|c|c|c|c|c|c|c|c|c|c|c|c|c|c|c|c|c|c|c|c|c|}
\hline
$44, 45, 46, 50, 51, 52, 53, 54, 55$ \\
\hline
$x_{367},x_{368},x_{378},x_{467},x_{468},x_{478},x_{567},x_{568},x_{578} $\\
\hline 
\end{tabular}

\begin{tabular}{|c|c|c|c|c|c|c|c|c|c|c|c|c|c|c|c|c|c|c|c|c|c|c|c|}
\hline
$56$ \\
\hline
$x_{678}$ \\
\hline
\end{tabular}

\\
\hline \hline

\size $180$ & 
\size $M_{[4,7]}\cong \gl_4\times \gl_3 \times \gl_1$,\quad \size $\Lam^{4,2}_{[1,4]}\oplus 1$ \\

\hline 
\size SP  \rule[-25pt]{0cm}{54pt} 
& 
\size

\begin{tabular}{|c|c|c|c|c|c|c|c|c|c|c|c|c|c|c|c|c|c|c|c|c|c|c|c|}
\hline
$6,11,15,26,30,40$ & $53$  \\
\hline
$x_{128},x_{138},x_{148},x_{238},x_{248},x_{348}$ & $x_{567}  $\\
\hline
\end{tabular}

\begin{tabular}{|c|c|c|c|c|c|c|c|c|c|c|c|c|c|c|c|c|c|c|c|c|c|c|c|}
\hline
$18, 20, 21, 33, 35, 36$ & $54, 55, 56$  \\
\hline
$x_{158},x_{168},x_{178},x_{258},x_{268},x_{278}$ 
& $x_{568},x_{578},x_{678}$ \\
\hline
$43, 45, 46, 49, 51, 52$ & \\
\hline
$x_{358},x_{368},x_{378},x_{458},x_{468},x_{478}$ & \\ 
\hline
\end{tabular}\\

\hline  \hline

\size $182$ & 
\size $M_{[3,7]}\cong \gl_3\times \gl_4 \times \gl_1$
,\quad \size $\Lam^{4,2}_{[4,7]}$ \\

\hline 
\size SP \rule[-14pt]{0cm}{32pt} & 
\size 
\begin{tabular}{|c|c|c|c|c|c|c|c|c|c|c|c|c|c|c|c|c|c|c|c|c|c|c|c|}
\hline
$49,51,52,54,55,56 $  \\
\hline
$x_{458},x_{468},x_{478},x_{568},x_{578},x_{678} $\\
\hline
\end{tabular}

\quad 
none

\\
\hline \hline

\size $183$ & 
\size $M_{[5]}\cong \gl_5\times \gl_3$
,\quad  \size $1$ \\

\hline 
\size SP \rule[-14pt]{0cm}{32pt} & 
\size 
\begin{tabular}{|c|c|c|c|c|c|c|c|c|c|c|c|c|c|c|c|c|c|c|c|c|c|c|c|}
\hline
$56 $  \\
\hline
$x_{678}$\\
\hline
\end{tabular}
\quad 
none

\\

\hline 
\end{tabular}

\end{center}

\vskip 10pt

From now on, we use the notation $Y_i,Z_i,W_i$ instead of 
$Y_{\be_i},Z_{\be_i},W_{\be_i}$ but keep using the notation 
$M_{\be_i},P_{\be_i},G_{\st,\be_i},\chi_{\be_i}$ 
to avoid confusion.

\vskip 10pt
\noindent
(1) $\beta_{1}=\tfrac {1} {120} (-5,-5,-5,3,3,3,3,3)$ 

In this case $Z_1\cong \wedge^2 \aff^5\otimes \aff^3$ and 
\begin{equation*}
M_{\be_1} = \gl_1\times M_{[3]}\cong \gl_1\times \gl_3\times \gl_5.
\end{equation*}
We express elements of $M_{\be_1}$ as $g=(t_0,\diag(g_2,g_1))$ 
where $t_0\in\gl_1,g_1\in\gl_5$ and $g_2\in\gl_3$. 
Then 
\begin{equation*}
M^{\st}_{\be_1} = \{\diag(g_2,g_1)\mid g_1\in\gl_5,g_2\in \gl_3,
(\det g_1)(\det g_2)=1\}. 
\end{equation*}
Since 
\begin{equation*}
\chi_{\be_1}(\diag(g_2,g_1))= (\det g_2)^{-5}(\det g_1)^3=(\det g_1)^8
\end{equation*}
on $M_{\be_1}^{\strict}$, 
\begin{align*}
G_{\st,\be_1} & = \left\{\diag(g_2,g_1) \;\;\vrule\;\; 
\begin{matrix}
g_1\in\gl_5,g_2\in \gl_3, \\
(\det g_1)(\det g_2)=1, \; \det g_1=1 
\end{matrix} 
\right\} \\
& = \left\{\diag(g_2,g_1) \;\;\vrule\;\; 
g_1\in\spl_5,g_2\in\spl_3 
\right\}
\cong \spl_5\times \spl_3.   
\end{align*}

For $x\in Z_1$, let
\begin{align*}
A_i(x) & = \begin{pmatrix}
0 & x_{i45} & x_{i46}& x_{i47}& x_{i48} \\
-x_{i45} & 0 & x_{i56} & x_{i57} & x_{i58} \\
-x_{i46} & -x_{i56} & 0 & x_{i67} & x_{i68} \\
-x_{i47} & -x_{i57} & -x_{i67} & 0 & x_{i78} \\
-x_{i48} & -x_{i58} & -x_{i68} & -x_{i78} & 0 
\end{pmatrix} 
\end{align*}
for $i=1,2,3$ and $A(x) = (A_1(x),A_2(x),A_3(x))$. 
We identify $x$ with $A(x)$. 
The action of the above $g\in M_{\be_1}$  
on $x\in Z_1$ is as follows: 
\begin{equation*}
\begin{pmatrix}
A_1(gx) \\
A_2(gx) \\
A_3(gx) 
\end{pmatrix}
= t_0 g_2 
\begin{pmatrix}
g_1 A_1(x) \,{}^t\! g_1 \\
g_1 A_2(x) \,{}^t\! g_1 \\
g_1 A_3(x) \,{}^t\! g_1
\end{pmatrix}.   
\end{equation*}
So $(M_{\be_1}\cap \gl_8,Z_1)$ 
can be identified with the \rep{} in Section III--6.

Let 
\begin{align*}
R_1 & = e_{147} -e_{156} - e_{248} + 2 e_{257} + e_{358} -e_{367}, \\
R_1' & = e_{147} -e_{156} - e_{248} + e_{257} + e_{358} -e_{367}.
\end{align*}
Then $R_1,R_1'$ correspond to $w,w'$ in 
\cite[pp.229,233]{tajima-yukie-GIT3} respectively.  
It is proved in Corollary III--6.13 
that there exists a homogeneous polynomial $P(x)$ of degree $30$ of 
$x\in V$ such that 
$P((g_1,g_2)x)=(\det g_1)^{12}(\det g_2)^{10}P(x)$ and 
$P(w)=4,P(w')=1$. 

Since $(\det g_1)^{12}(\det g_2)^{10}=(\det g_1)^2$ is proportional
to $\chi_{\be_1}(g)$ on $M^{\st}_{\be_1}$, 
$P(x)$ is an invariant polynomial with respect to
$G_{\st,\be_1}$. So $R_1'\in Z_1^{\sst}$. Therefore, 
$S_{\be_1}\not=\emptyset$ regardless of $\ch(k)$. 

\begin{prop}
\label{prop:r-orbit-case1}
Suppose that $\ch(k)\not=2$. Then 
\begin{math}
M_{\be_1\, } \backslash Z_{1\, k}^{\sst} 
\end{math}
is in bijective correspondence with $\Prg_2(k)$ 
by associating $M_{\be_1\, x}/T_0$ to the orbit of 
$x\in Z_{1\, k}^{\sst}$. 
\end{prop}
\begin{proof}
We remind the reader that $T_0$ is defined 
in (\ref{eq:T0-defn}). 
Corollary III--6.7 and Proposition III--6.8 
imply that $(M_{\be_1}\cap \gl_8,Z_1)$ is a regular \pv. 
Proposition III--6.14 says that 
$(M_{\be_1\, k}\cap \gl_8(k)) \backslash Z_{1\, k}^{\sst}$ 
is in bijective correspondence with $\Prg_2(k)$ 
by associating 
\begin{equation*}
(M_{\be_1\, k}\cap \gl_8(k)) \cdot x \mapsto 
M_{\be_1\, x}\cap \gl_8(k)/(T_0\cap \gl_8).  
\end{equation*}
Since the action of 
$t_0\in\gl_1$ of the $\gl_1$-factor of $M_{\be_1}$
can be absorbed to the action of $t_0I_3\in M_{\be_1}\cap \gl_8$, 
the proposition follows.  
\end{proof}

We verify Condition \ref{cond:unipotent-connected} (2), (3). 
In this case, $U_{\be_1}$ consists of all the elements of the form 
\begin{equation*}
n(u) = 
\begin{pmatrix}
1 & 0 & 0 & 0 & 0 & 0 & 0 & 0 \\
0 & 1 & 0 & 0 & 0 & 0 & 0 & 0 \\
0 & 0 & 1 & 0 & 0 & 0 & 0 & 0 \\
u_{41} & u_{42} & u_{43} & 1 & 0 & 0 & 0 & 0 \\  
u_{51} & u_{52} & u_{53} & 0 & 1 & 0 & 0 & 0 \\ 
u_{61} & u_{62} & u_{63} & 0 & 0 & 1 & 0 & 0 \\
u_{71} & u_{72} & u_{73} & 0 & 0 & 0 & 1 & 0 \\
u_{81} & u_{82} & u_{83} & 0 & 0 & 0 & 0 & 1 
\end{pmatrix}.
\end{equation*}
$\dim U_{\be_1}=15$ as an algebraic group and 
$\dim W_1=10$. 

We choose the order of $u_{ij}$ and the coordinate vectors 
of $W_1$ in the following manner: 
\begin{align*}
& u_{42},u_{43},u_{71},u_{81},u_{82}, \quad 
u_{41},u_{62},u_{83},u_{51},u_{52},
u_{61},u_{63},u_{72},u_{73},u_{53}, \\
& e_{456}, e_{468}, e_{678}, e_{457}, e_{458}, 
e_{467}, e_{568}, e_{478}, e_{578}, e_{567}
\end{align*}
($u_{41}\ccd u_{53}$ correspond to 
the coordinates of $W_1$). 

We computed $n(u)R_1-R_1$ using MAPLE using the 
differential forms package as follows:

\begin{verbatim}
 > with(difforms); 
 > defform(x1=0); defform(x2=0); defform(x3=0); defform(x4=0); 
 > defform(x5=0); defform(x6=0); defform(x7=0); defform(x8=0);
 > with(linalg);    
 > defform(u21=0); **** defform(u87=0);  
 > R1:= matrix(8,8,[1 , 0, 0, 0, 0, 0, 0, 0, 
                    0, 1, 0, 0, 0, 0, 0, 0, 
                    0, 0, 1, 0, 0, 0, 0, 0, 
                    u41, u42, u43, 1, 0, 0, 0, 0,  
                    **** 
                    u81, u82, u83, 0, 0, 0, 0, 1 ]);
 > e1:= d(x1); e2:= d(x2); e3:= d(x3); e4:= d(x4); 
 > e5:= d(x5); e6:= d(x6); e7:= d(x7); e8:= d(x8); 
 > m:= matrix(1,8,[e1,e2,e3,e4,e5,e6,e7,e8]);
 > v1:= sum(R1[i,1] * m[1,i],i=1..8); v2:= sum(R1[i,2] * m[1,i],i=1..8);  
   ****
   v7:= sum(R1[i,7] * m[1,i],i=1..8); v8:= sum(R1[i,8] * m[1,i],i=1..8);
 > v1 &^ v4 &^ v7 - v1 &^ v5 &^ v6 - v2 &^ v4 &^ v8 
   + 2* v2 &^ v5 &^ v7 + v3 &^ v5 &^ v8 -v3&^ v6 &^ v7;
 > simpform(%);
\end{verbatim}

We will not show maple commands for other cases.  
The file which contains maple commands for all the cases 
can be seen in the google drive of the second author
as follows. 

\vskip 10pt

\begin{verbatim}
  https://drive.google.com/drive/folders/
    1iHzyw0T2XumpZ4vM688MLIf8a-lH5JtC?usp=drive_link
\end{verbatim}

\vskip 10pt

Explicitly, $n(u)R_1-R_1$ is 
\begin{align*}
& - u_{41} e_{456} 
+ u_{62} e_{468} 
- u_{83} e_{678}
+ (-u_{51}+2u_{42}) e_{457} 
+ (u_{52}+u_{43}) e_{458} \\
& - (u_{61}+u_{43}) e_{467}  
-(u_{63}+u_{81}) e_{568}
+ (u_{72}+u_{81}) e_{478} 
+ (-u_{73}+2u_{82}) e_{578} \\
& - (u_{53} + 2u_{62}+u_{71}) e_{567}.  
\end{align*}
One can observe that Condition \ref{cond:unipotent-connected} (2), (3)
are satisfied and so $\rg\backslash S_{\be_1\, k}$ is also 
in bijective correspondence with $\Prg_2(k)$.

\vskip 5pt
\noindent
(2) $\beta_{4}=\tfrac {1} {56} (-5,-1,-1,-1,-1,3,3,3)$  

In this case, 
$Z_{4}\cong \wedge^2\aff^4\otimes \aff^3\op \wedge^2\aff^3$ 
and 
\begin{equation*}
M_{\be_4}=\gl_1\times M_{[1,5]}\cong \gl_1\times \gl_1\times 
\gl_4 \times \gl_3. 
\end{equation*}
We express elements of $M_{\be_4}$ as 
\begin{math}
g = (t_0,\diag(t_3,g_1,g_2))
\end{math}
where $t_0,t_3\in\gl_1,g_1\in\gl_4$ and $g_2\in\gl_3$. 
Then 
\begin{align*}
M_{\beta_4}^{\st} 
& = \{\diag(t_3,g_1,g_2) \mid 
g_1\in\gl_4,g_2\in\gl_3,t_3\in\gl_1,(\det g_1)(\det g_2)t_3=1\} \\
& = \{\diag((\det g_1)^{-1}(\det g_2)^{-1},g_1,g_2)\mid 
g_1\in\gl_4,g_2\in\gl_3\} \cong \gl_4\times \gl_3.
\end{align*}
Since 
\begin{equation*}
\chi_{\be_4}(g)=((\det g_1)^{-1}(\det g_2)^{-1})^{-5}
(\det g_1)^{-1}(\det g_2)^3
= (\det g_1)^4 (\det g_2)^8
\end{equation*}
on $M^{\st}_{\be_4}$, 
\begin{align*}
G_{\st, \, \be_4}
& = \left\{(\det g_1)^{-1}(\det g_2)^{-1},g_1,g_2)
\;\;\vrule\;\; \begin{matrix}
g_1\in \gl_4,g_2\in\gl_3, \\
(\det g_1)(\det g_2)^2=1 
\end{matrix} 
\right\}.   
\end{align*}

For $x\in Z_4$, let 
\begin{align*}
A_1(x) & = \begin{pmatrix}
0 & x_{236} & x_{246} & x_{256} \\
-x_{236} & 0 & x_{346} & x_{356} \\
-x_{246} & -x_{346} & 0 & x_{456} \\
-x_{256} & -x_{356} & -x_{456} & 0  
\end{pmatrix}, \\ 
B(x) & = x_{178} \,\bbmp_{3,2}\wedge \bbmp_{3,2}
- x_{168} \,\bbmp_{3,1}\wedge \bbmp_{3,3}
+ x_{167} \,\bbmp_{3,1}\wedge \bbmp_{3,2}.
\end{align*}
We define $A_2(x),A_3(x)$ by replacing $6$ in the above indices 
by $7,8$ respectively and put $A(x)=(A_1(x),A_2(x),A_3(x))$. 
We identify $x\in Z_4$ with 
$(A(x),B(x))$. 
The action of the above 
$g\in M_{\be_4}$ on 
$x = (A(x),B(x)) \in Z_4$ is as follows:
\begin{equation}
\label{eq:g-action-A-case4}
\begin{pmatrix}
A_1(gx) \\
A_2(gx) \\
A_3(gx) 
\end{pmatrix}
= t_0 g_2 \begin{pmatrix}
g_1 A_1(x) \,{}^t\! g_1 \\
g_1 A_2(x) \,{}^t\! g_1 \\
g_1 A_3(x) \,{}^t\! g_1
\end{pmatrix}, \;
B(gx) = t_0t_3 \wedge^2\! g_2 \, B(x).  
\end{equation}
So $(M_{\be_4}\cap \gl_8,Z_4)$ 
is the \rep{} considered in Section \ref{sec:orbits2}.

Let $P_1(x),P_2(x)$ be the polynomials on $Z_{4}$ 
which correspond to (\ref{eq:be4-phi1-inv}), (\ref{eq:be4-phi2-inv}) 
respectively and $P(x)=P_1(x)^2P_2(x)$. 
Then 
\begin{equation*}
P((g_1,g_2)x) = t_3^4 (\det g_1)^5 (\det g_2)^6\Phi(x)
= (\det g_1)(\det g_2)^2 P(x)
\end{equation*}
on $M^{\st}_{\be_4}$. 
Since $(\det g_1)(\det g_2)^2$ is proportional to 
$\chi_{\be_4}(g)$ on $M^{\st}_{\be_4}$, 
$P(x)$ is an invariant polynomial with 
respect to $G_{\st,\be_4}$. 

Let 
\begin{equation*}
R_4 = e_{236} + e_{247} + e_{357} + e_{458} -e_{168}.
\end{equation*}
Then $R_4$ corresponds to the element $w$ 
in (\ref{eq:case4-R-defn}). By (\ref{eq:be4-invariants}), 
$P_1(R_4),P_2(R_4),P(R_4)\not=0$. So 
$R_4\in Z_4^{\sst}$. Therefore, 
$S_{\be_4}\not=\emptyset$ regardless of $\ch(k)$.

Suppose that $\ch(k)\not=2$. 
Proposition \ref{prop:be4-stabilizer} implies that 
$(M_{\be_1}\cap \gl_8,Z_1)$ is a regular \pv.  
By Corollaries III--4.7, \ref{cor:case4-single-orbit}, 
$Z_{4\, k^{\sep}}^{\sst} = M_{\be_4\, k^{\sep}} \cdot R_4$.  

\begin{prop}
\label{prop:r-orbit-case4}
$M_{\be_4\, k} \backslash Z_{4\, k}^{\sst}$ 
is in bijective correspondence with $\Gunit_2(k)$.  
\end{prop}
\begin{proof}
Theorem \ref{thm:case4-rational-orbits} implies that 
\begin{math}
(M_{\be_4\, k}\cap \gl_8(k)) \backslash Z_{4\, k}^{\sst} 
\end{math}
is in bijective correspondence with $\Gunit_2(k)$.
Since the action of 
$t_0\in\gl_1$ of the $\gl_1$-factor of $M_{\be_4}$
can be absorbed to the action of 
$\diag(t_0^{-1},I_4,t_0I_3)\in M_{\be_1}\cap \gl_8$, 
the proposition follows.  
\end{proof}

We verify Condition \ref{cond:unipotent-connected} (2), (3). 
In this case, 
$U_{\be_4}$ consists of all the elements of the form 
\begin{equation*}
n(u) = 
\begin{pmatrix}
1 & 0 & 0 & 0 & 0 & 0 & 0 & 0 \\
u_{21} & 1 & 0 & 0 & 0 & 0 & 0 & 0 \\
u_{31} & 0 & 1 & 0 & 0 & 0 & 0 & 0 \\
u_{41} & 0 & 0 & 1 & 0 & 0 & 0 & 0 \\  
u_{51} & 0 & 0 & 0 & 1 & 0 & 0 & 0 \\ 
u_{61} & u_{62} & u_{63} & u_{64} & u_{65} & 1 & 0 & 0 \\
u_{71} & u_{72} & u_{73} & u_{74} & u_{75} & 0 & 1 & 0 \\
u_{81} & u_{82} & u_{83} & u_{84} & u_{85} & 0 & 0 & 1 
\end{pmatrix}.
\end{equation*}
$\dim U_{\be_1}=19$ as an algebraic group and 
$\dim W_4=13$. 

We choose the order of $u_{ij}$ and the coordinate vectors 
of $W_4$ in the following manner: 
\begin{align*}
& u_{21},u_{31},u_{41},u_{51},u_{61},u_{81}, \quad 
u_{62},u_{63},u_{84},u_{85},u_{64},u_{65},u_{72},
u_{73},u_{82},u_{83},u_{74},u_{75},u_{71} \\
& e_{467}, e_{567}, e_{278}, e_{378}, e_{568}, 
e_{468}, e_{367}, e_{267}, e_{368}, e_{268}, 
e_{578}, e_{478}, e_{678}
\end{align*}
($u_{62}\ccd u_{71}$ correspond to 
the coordinates of $W_4$). 
Explicitly, $n(u)R_4-R_4$ is 
\begin{align*}
& - u_{62} e_{467} - u_{63} e_{567} - u_{84} e_{278} - u_{85} e_{378} 
- (u_{64} - u_{51}) e_{568} + (u_{65} + u_{41}) e_{468} \\
& + (u_{72} + u_{65}) e_{367} + (-u_{73} + u_{64}) e_{267} 
+ (u_{82} + u_{31}) e_{368} - (u_{83} - u_{21}) e_{268} \\
& - (u_{74} - u_{83}) e_{578} + (u_{75} + u_{82}) e_{478} 
+ (-u_{71} -u_{62}u_{84} - u_{63}u_{85} 
+ u_{64}u_{75} + u_{64}u_{82} \\
& \quad - u_{65}u_{74} + u_{65}u_{83} + u_{72}u_{83} 
- u_{73}u_{82}) e_{678}. 
\end{align*}
Condition \ref{cond:unipotent-connected} (2), (3)
are satisfied and so $\rg\backslash S_{\be_4\, k}$
is also in bijective correspondence with $\Gunit_2(k)$.  

\vskip 5pt
\noindent
(3) $\beta_{7}=\tfrac {1} {24} (-1,-1,-1,-1,-1,-1,3,3)$ 

In this case, $Z_7\cong \wedge^2 \aff^6 \otimes\aff^2$ and 
\begin{equation*}
M_{\be_7} = \gl_1\times M_{[6]} \cong \gl_1\times \gl_6\times \gl_2.
\end{equation*}
We express elements of $M_{\be_7}$ as 
$g=(t_0,\diag(g_1,g_2))$ where $t_0\in\gl_1$, 
$g_1\in\gl_6$ and $g_2\in\gl_2$.  Then 
\begin{equation*}
M^{\st}_{\be_7} = 
\left\{
\diag(g_1,g_2) \mid g_1\in \gl_6,g_2\in\gl_2, 
(\det g_1)(\det g_2)=1
\right\}.  
\end{equation*}
Since 
\begin{equation*}
\chi_{\be_7}(g) = (\det g_1)^{-1}(\det g_2)^3
= (\det g_2)^4
\end{equation*}
on $M^{\st}_{\be_7}$, 
\begin{equation*}
G_{\st,\be_7} = 
\left\{
\diag(g_1,g_2) \mid g_1\in \spl_6,g_2\in\spl_2
\right\}.  
\end{equation*}

For $x\in Z_7$, let 
\begin{equation*}
A_1(x) = \begin{pmatrix}
0 & x_{127} & x_{137} & x_{147} & x_{157} & x_{167} \\
-x_{127} & 0 & x_{237} & x_{247} & x_{257} & x_{267} \\ 
-x_{137} & -x_{237} & 0 & x_{347} & x_{357} & x_{367} \\ 
-x_{147} & -x_{247} & -x_{347} & 0 & x_{457} & x_{467} \\ 
-x_{157} & -x_{257} & -x_{357} & -x_{457} & 0 & x_{567} \\ 
-x_{167} & -x_{267} & -x_{367} & -x_{467} & -x_{567} & 0  
\end{pmatrix}.
\end{equation*}
We define $A_2(x)$ by replacing $7$ in the above indices by $8$ 
and $A(x)=(A_1(x),A_2(x))$.  We identify $x\in Z_7$ with 
$A(x)$.  
The action of the above $g\in M_{\be_7}$  
on $x\in Z_7$ is as follows:
\begin{equation*}
\begin{pmatrix}
A_1(gx) \\
A_2(gx)
\end{pmatrix} 
= t_0 g_2 \begin{pmatrix}
g_1 A_1(x) \,{}^t\! g_1 \\
g_1 A_2(x) \,{}^t\! g_1 
\end{pmatrix}. 
\end{equation*}
$(M_{\be_7}\cap\gl_8,Z_7)$
is the \rep{} considered in 
Section \ref{sec:orbits1} Case I.  

Let $P(x)$ be the homogeneous polynomial 
of degree $12$ of $x\in Z_7$ 
which corresponds to the polynomial in 
(\ref{eq:62-relative-invariant}).  Then 
$P(gx)=(\det g_1)^4(\det g_2)^6P(x)$ on $M^{\st}_{\be_7}$.  
Since $(\det g_1)^4(\det g_2)^6$ 
$=(\det g_2)^2$ is 
proportional to $\chi_{\be_7}(g)$ on $M^{\st}_{\be_7}$,  
$P(x)$ is an invariant polynomial with respect to
the action of $G_{\st,\be_7}$.  Let 
\begin{equation*}
R_7 = e_{127} - e_{347} + e_{348} - e_{568}.  
\end{equation*}
Then $R_7$ corresponds to the element 
$w$ in (\ref{eq:case7-R-defn}).  
Therefore, 
$Z_{7\, k^{\sep}}^{\sst} = M_{\be_7\, \sep}\cdot R_7$. 

Let $F:Z_7\ni x \mapsto F_x (v) \in \sym^3 \aff^2$ be the map 
in the proof of Proposition \ref{prop:tangent-cubic-h}. 
Since the action of $t_0\in\gl_1$ of the $\gl_1$-factor 
can be absorbed to the action of 
$\diag(I_6,t_0I_2)\in M_{\be_7\, k}\cap \gl_8(k)$, 
Theorem \ref{thm:case7-rational-orbits} implies 
the following proposition. 

\begin{prop}
\label{prop:r-orbit-case7}
$M_{\be_7\, k}\backslash Z_{7\, k}^{\sst}$ 
is in bijective correspondence with $\Ex_3(k)$ 
by associating the extension of $k$ 
generated by a root of $F_x$ to the orbit of 
$x\in Z_{7\, k}^{\sst}$.  
\end{prop}

We verify Condition \ref{cond:unipotent-connected} (2), (3).  
In this case, $U_{\be_7}$ consists of all the elements of the form 
\begin{equation*}
n(u) = 
\begin{pmatrix}
1 & 0 & 0 & 0 & 0 & 0 & 0 & 0 \\
0 & 1 & 0 & 0 & 0 & 0 & 0 & 0 \\
0 & 0 & 1 & 0 & 0 & 0 & 0 & 0 \\
0 & 0 & 0 & 1 & 0 & 0 & 0 & 0 \\  
0 & 0 & 0 & 0 & 1 & 0 & 0 & 0 \\ 
0 & 0 & 0 & 0 & 0 & 1 & 0 & 0 \\
u_{71} & u_{72} & u_{73} & u_{74} & u_{75} & u_{76} & 1 & 0 \\
u_{81} & u_{82} & u_{83} & u_{84} & u_{85} & u_{86} & 0 & 1 
\end{pmatrix}.
\end{equation*}
Then $\dim U_{\be_7}=12$ as an algebraic group and 
$\dim W_7=6$.  

We choose the order of $u_{ij}$ and the 
coordinate vectors of $W_7$ in the following manner:
\begin{align*}
& u_{71},u_{72},u_{83},u_{84},u_{85},u_{86}, \quad
u_{75},u_{76},u_{81},u_{82},u_{73},u_{74} \\
& e_{678},e_{578},e_{278},e_{178},e_{478},e_{378}
\end{align*}
($u_{75}\ccd u_{74}$ correspond to the coordinates of $W_7$). 
Explicitly, $n(u)R_4-R_4$ is 
\begin{align*}
u_{75} e_{678} 
- u_{76} e_{578} 
+ u_{81} e_{278}  
- u_{82} e_{178} 
+ (u_{73} - u_{83}) e_{478} 
+ (- u_{74} + u_{84}) e_{378}.
\end{align*}
Condition \ref{cond:unipotent-connected} (2), (3)
are satisfied and so $\rg\backslash S_{\be_7\, k}$
is also in bijective correspondence with $\Ex_3(k)$.

\vskip 5pt
\noindent
(4) $\beta_{9}=\tfrac {3} {56} (-7,1,1,1,1,1,1,1)$ 

In this case, $Z_9\cong \wedge^3\aff^7$ and 
\begin{equation*}
M_{\be_9} = \gl_1\times M_{[1]} \cong \gl_1\times \gl_1\times \gl_7.
\end{equation*}
We express elements of $M_{\be_9}$ as $g=(t_0,\diag(t_2,g_1))$
where $t_0,t_2\in\gl_1$ and $g_1\in\gl_7$. Then 
\begin{align*}
M^{\st}_{\be_9} 
& = \{\diag(t_2,g_1)\mid g_1\in\gl_7, t_2\in\gl_1, \; (\det g_1)t_2=1\} 
\cong \gl_7.  
\end{align*}
Since 
\begin{equation*}
\chi_{\be_9}(g) = t_2^{-7}(\det g_1) = (\det g_1)^8
\end{equation*}
on $M^{\st}_{\be_9}$, 
\begin{equation*}
G_{\st , \, \beta_9}
= \{\diag(1,g_1)\mid g_1\in\spl_7\} 
\cong \spl_7.  
\end{equation*}

For $x\in Z_9$, let 
\begin{equation*}
F(x) = \sum_{2\leqq i_1<i_2<i_3} 
x_{i_1-1,i_2-1,i_3-1} \bbmp_{7,i_1-1}
\wedge \bbmp_{7,i_2-1} \wedge \bbmp_{7,i_3-1}
\in\wedge^3 \aff^7.  
\end{equation*}
Then the action of $(t_0,\diag(t_1,I_7))\in M_{\be_9}$ 
($t_0,t_1\in\gl_1$) on $x$ do not depend on $t_1$ 
and can be identified with the action of 
Section \ref{sec:orbits1} Case III.  

For $x\in Z_9$, let $Q_x$ be the quadratic form 
defined in \cite[p.1694, Definition (2.13)]{yukiem}
and $P(x)$ the discriminant
of $\tfrac 16 Q_x$ (division by $6$ was justified) in 
Section \ref{sec:orbits1}.  Then $P(x)$ is a homogeneous 
polynomial of degree $21$ of $x\in Z_9$.  
Let $g=\diag((\det g_1)^{-1},g_1)\in M^{\st}_{\be_9}$ where 
$g_1\in \gl_7$. By 
\cite[p.1694, Definition (2.14)]{yukiem}
$Q_{gx} = (\det g_1)g_1 Q_x$. So 
$P(x) = (\det g_1)^9P(x)$.  
Since $(\det g_1)^9$ is proportional to 
$\chi_{\be_9}(g)$ on $M^{\st}_{\be_9}$,  
$P(x)$ is an invariant polynomial
with respect to $G_{\st,\be_9}$. 

Let 
\begin{equation*}
R_9 = e_{345} + e_{678} + e_{236} + e_{247} + e_{258}.  
\end{equation*}
Then $R_9$ corresponds to $w$ in (\ref{eq:w-defn-case9}). 
Since $P(R_9)=P(w)=\pm 1$, $R_9\in Z_{9\, k}^{\sst}$. 
Therefore, $S_{\be_9}\not=\emptyset$. 
The following proposition follows from 
Proposition \ref{prop:r-orbit-G2}.  

\begin{prop}
\label{prop:r-orbit-case9}
$M_{\be_9\, k}\backslash Z_{9\, k}^{\sst}$ 
is in bijective correspondence with $\mathrm{G}_2(k)$ 
by associating $M_{\be_9\, x\, k}/T_0$ to 
$x\in Z_{9\, k}^{\sst}$.  
\end{prop}

Note that there is no restriction on $\ch(k)$. 
If $\ch(k)\not=2$, we explained the outline 
of the argument to show that 
$M_{\be_9\, k}\backslash Z_{9\, k}^{\sst}$ 
is in bijective correspondence with 
the set of $k$-isomorphism classes of $k$-forms 
of the split octonion $\soct$. 

In this case, $W_9=\{0\}$.  Therefore, 
$\rg\backslash S_{\be_9\, k}$ is also in 
bijective correspondence with $\mathrm{G}_2(k)$.  

\vskip 5pt
\noindent
(5) $\beta_{13}=\tfrac {1} {24} (-3,-1,-1,-1,1,1,1,3)$ 

Let $W_1=W_2=\aff^3$. 
In this case, 
\begin{math}
Z_{13}\cong W_1\otimes \wedge^2 W_2 
\oplus \wedge^2 W_1 \oplus W_2
\end{math}
and 
\begin{equation*}
M_{\be_{13}} = \gl_1\times M_{[1,4,7]} 
\cong \gl_1\times \gl_1\times \gl_3^2 \times \gl_1. 
\end{equation*}
We express elements of $M_{\be_{13}}$ as 
$g=(t_0,\diag(t_3,g_1,g_2,t_4))$ where 
$t_0,t_3,t_4\in \gl_1$ and $g_1,g_2\in\gl_3$. 
Then 
\begin{align*}
M^{\st}_{\be_{13}} 
& = \left\{
\diag(t_3,g_1,g_2,t_4) \;\; \vrule \;\; 
\begin{matrix}
g_1,g_2\in\gl_3,t_3,t_4\in\gl_1 \\
(\det g_1)(\det g_2)t_3t_4=1
\end{matrix}
\right\} \\
& = \left\{
\diag((\det g_1)^{-1}(\det g_2)^{-1}t_4^{-1},g_1,g_2,t_4) 
\mid g_1,g_2\in\gl_3,t_4\in\gl_1 
\right\} \\
& \cong \gl_3^2\times \gl_1.  
\end{align*}
Since 
\begin{equation*}
\chi_{\be_{13}} (g) 
= t_3^{-3}(\det g_1)^{-1}(\det g_2) t_4^3
= (\det g_1)^2(\det g_2)^4 t_4^6
\end{equation*}
on $M^{\st}_{\be_[13]}$, 
\begin{equation*}
G_{\st,\be_{13}} = \left\{
\diag((\det g_1)^{-1}(\det g_2)^{-1}t_4^{-1},g_1,g_2,t_4)
\;\; \vrule \;\; 
\begin{matrix}
g_1,g_2\in\gl_3,t_4\in\gl_1 \\
(\det g_1)(\det g_2)^2 t_4^3=1
\end{matrix}
\right\}.  
\end{equation*}

For $x\in Z_{13}$, let
\begin{align*}
A(x) & = \begin{pmatrix}
x_{267} & -x_{257} & x_{256} \\
x_{367} & -x_{357} & x_{356} \\
x_{467} & -x_{457} & x_{456} 
\end{pmatrix}, \;
v_1(x) = 
\begin{pmatrix}
x_{348} \\ -x_{248} \\ x_{238} 
\end{pmatrix}, \; 
v_2(x) = \begin{pmatrix}
x_{158} \\
x_{168} \\
x_{178} 
\end{pmatrix}.
\end{align*}
We identify $x$ with $(A(x),v_1(x),v_2(x))$.  
The action of the above 
\begin{math}
g = \in M_{\be_{13}}
\end{math}
on $x\in Z_{13}$ is as follows:
\begin{align*}
A(gx) & = t_0(\det g_2) g_1 A(x) g_2^{-1}, \; 
v_1(gx) = t_0t_4 (\det g_1) \, {}^t\!g_1^{-1} v_1(x), \\
v_2(gx) & = t_0t_3t_4g_2 v_2(x).  
\end{align*}

Let 
\begin{equation}
\label{eq:case13-invariant}
P_1(x) = \det A(x), \; 
P_2(x) = {}^t v_1(x) A(x) v_2(x), \;
P(x)=P_1(x)P_2(x)^3.     
\end{equation}
Then 
\begin{align*}
& P_1(gx) = (\det g_1)(\det g_2)^2 P_1(x), \;  
P_2(gx) = t_4 P_2(x), \\
& P(gx) = (\det g_1)(\det g_2)^2t_4^3P(x)
\end{align*}
on $M^{\st}_{\be_{13}}$. Since 
$(\det g_1)(\det g_2)^2t_4^3$ is 
proportional to $\chi_{\be_{13}}(g)$ on $M^{\st}_{\be_{13}}$,   
$P(x)$ is an invariant polynomial 
with respect to $G_{\st,\be_{13}}$. 

Let
\begin{equation*}
R_{13} = e_{158} + e_{348} + e_{267} - e_{357} + e_{456}. 
\end{equation*}
Then 
\begin{equation*}
A(R_{13}) = I_3, \; 
v_1(R_{13}) = v_2(R_{13}) = \bbmp_{3,1}=
\begin{pmatrix}
1 \\ 0 \\ 0
\end{pmatrix}.
\end{equation*}
So $P_1(R_{13})=P_2(R_{13})=P(R_{13})=1$, 
which implies that $R_{13}\in Z_{13\, k}^{\sst}$. 
Therefore, $S_{\be_{13}}\not=\emptyset$. 

Since there are three components for $Z_{13}$ 
and there are only two relative invariant polynomials, 
we determine $M_{\be_{13}\, k}\backslash Z_{13\, k}^{\sst}$ 
by explicit computations. 

\begin{prop}
\label{prop:r-orbit-case13}
$Z_{13\, k}^{\sst} = M_{\be_{13}\, k}\cdot R_{13}$.  
\end{prop}
\begin{proof}
Suppose that $x\in Z_{13\, k}^{\sst}$. 
Since $\det A(x)\not=0$, we may assume that 
\begin{math}
A(x) = I_3
\end{math}
by applying an element of $M_{\be_{13}\, k}$ 
if necessary.

Elements of the form 
$\diag(1,g_1,g_1,1)$ where $g_1\in\spl_3$  
do not change $A(x)$. 
By applying an element of this form if necessary, 
we may assume that $v_2(x) = [1,0,0]$.  Then 
\begin{math}
P_2(x) = x_{348} \not=0.  
\end{math}
By applying the element 
$g=\diag(1,I_3,I_3,x_{348}^{-1})$, 
we may assume that $A(x)=I_3,v(x)=[1,0,0]$ and 
$x_{348}=1$.

If 
\begin{equation*}
g_1 = \begin{pmatrix}
1 & -x_{248} & x_{238} \\
0 & 1 & 0 \\
0 & 0 & 1 
\end{pmatrix}, 
\end{equation*}
then $\det g_1=1$ and  
\begin{equation*}
{}^t\!g_1^{-1} = \begin{pmatrix}
1 & 0 & 0 \\
x_{248} & 1 & 0 \\
-x_{238} & 0 & 1 
\end{pmatrix}. 
\end{equation*}
If $g=\diag(1,g_1,g_1,1)$ then $A(gx)=I_3$, 
$v_2(gx)=[1,0,0]$ 
and 
\begin{equation*}
v_2(gx) = \begin{pmatrix}
1 & 0 & 0 \\
x_{248} & 1 & 0 \\
-x_{238} & 0 & 1 
\end{pmatrix}
\begin{pmatrix}
1 \\ -x_{248} \\ x_{238}
\end{pmatrix}
= \begin{pmatrix}
1 \\ 0 \\ 0
\end{pmatrix}.  
\end{equation*}
This implies that $gx=R_{13}$. 
\end{proof}

We verify Condition \ref{cond:unipotent-connected} (2), (3).  
In this case, $U_{\be_{13}}$ consists of all the elements 
of the form 
\begin{equation*}
n(u) = 
\begin{pmatrix}
1 & 0 & 0 & 0 & 0 & 0 & 0 & 0 \\
u_{21} & 1 & 0 & 0 & 0 & 0 & 0 & 0 \\
u_{31} & 0 & 1 & 0 & 0 & 0 & 0 & 0 \\
u_{41} & 0 & 0 & 1 & 0 & 0 & 0 & 0 \\  
u_{51} & u_{52} & u_{53} & u_{54} & 1 & 0 & 0 & 0 \\ 
u_{61} & u_{62} & u_{63} & u_{64} & 0 & 1 & 0 & 0 \\
u_{71} & u_{72} & u_{73} & u_{74} & 0 & 0 & 1 & 0 \\
u_{81} & u_{82} & u_{83} & u_{84} & u_{85} & u_{86} & u_{87} & 1 
\end{pmatrix}.
\end{equation*}
Then $\dim U_{\be_7}=22$ as an algebraic group and 
$\dim W_7=13$.  

We choose the order of $u_{ij}$ and the 
coordinate vectors of $W_7$ in the following manner:
\begin{align*}
& u_{51},u_{53},u_{54},u_{61},u_{62},
u_{71},u_{72},u_{74},u_{81}, \\
& u_{21},u_{64},u_{73},u_{86},u_{87},u_{85},u_{63},
u_{31},u_{41},u_{52},u_{82},u_{83},u_{84} \\ 
& e_{258},e_{368},e_{478},e_{278},e_{268},e_{378}, 
e_{468},e_{358},e_{458},e_{567},e_{678},e_{578},e_{568}
\end{align*}
($u_{21}\ccd u_{84}$ correspond to the coordinates of $W_7$). 
Explicitly, $n(u)R_{13}-R_{13}$ is
\begin{align*}
& u_{21} e_{258} 
+ u_{64} e_{368} 
- u_{73} e_{478}
- u_{86} e_{278}
+ u_{87} e_{268} 
+ (u_{85} + u_{74}) e_{378} 
- (u_{63} + u_{85}) e_{468} \\
& + (u_{31} + u_{54} - u_{87}) e_{358} 
+ (u_{41} - u_{53} + u_{86}) e_{458} 
+ (u_{52} + u_{63}+u_{74}) e_{567} \\ 
& + (u_{82} - u_{62}u_{86} + u_{63}u_{74} + u_{63}u_{85} - u_{64}u_{73} 
-u_{72}u_{87}+u_{74}u_{85}) e_{678} \\
& - (u_{83}+u_{71}+u_{52}u_{86}-u_{53}u_{74}-u_{53}u_{85}
+u_{54}u_{73}-u_{73}u_{87}+u_{74}u_{86}) e_{578} \\
& + (u_{84}-u_{61}+u_{52}u_{87}+u_{53}u_{64}
-u_{54}u_{63}-u_{54}u_{85}+u_{63}u_{87} - u_{64}u_{86}) e_{568}. 
\end{align*}
Condition \ref{cond:unipotent-connected} (2), (3) 
are satisfied and so $\rg\backslash S_{\be_13\, k}$ 
also consists of a single point.

\vskip 5pt
\noindent
(6) $\beta_{32}=\tfrac {1} {8} (-1,-1,-1,-1,1,1,1,1)$ 

In this case, 
$Z_{\beta_{32}}\cong \wedge^2\aff^4\otimes \aff^4$ and 
\begin{equation*}
M_{\be_{32}} = \gl_1\times M_{[4]} \cong \gl_1\times \gl_4^2.
\end{equation*}
We express elements of $M_{\be_{32}}$ as $g=(t_0,\diag(g_1,g_2))$ 
where $t_0\in\gl_1$ and $g_1,g_2\in\gl_4$.  Then 
\begin{equation*}
M^{\st}_{\be_{32}} = \{\diag(g_1,g_2)\mid g_1,g_2\in\gl_4,
(\det g_1)(\det g_2)=1\}. 
\end{equation*}
Since 
\begin{equation*}
\chi_{\be_{32}}(g) = (\det g_1)^{-1}(\det g_2)=(\det g_2)^2
\end{equation*}
on $M^{\st}_{\be_{32}}$, 
\begin{equation*}
G_{\st,\be_{32}} = \{\diag(g_1,g_2)\mid g_1,g_2\in\spl_4\}. 
\end{equation*}

For $x\in Z_{32}$ and $i=1,2,3,4$, let 
\begin{align*}
A_i(x) = \begin{pmatrix}
0 & x_{i56} & x_{i57} & x_{i58} \\
-x_{i56} & 0 & x_{i67} & x_{i568} \\
-x_{i57} & x-_{i67} & 0 & x_{i78} \\
-x_{i58} & x-_{i78} & -x_{i78} & 0 
\end{pmatrix}
\end{align*}
and $A(x)=(A_1(x)\ccd A_4(x))$. We identify $x$ with $A(x)$. 

Since $\dim \wedge^2 \aff^4=6$, 
$(M_{\be_{32}},Z_{32})$ and 
$(\gl_1\times \gl_4\times \gl_2,(\wedge^2 \aff^4)^* \otimes \aff^2)$
are the Castling transforms of each other. 
If $\om_1,\om_2\in \wedge^2 \aff^4$ then 
$\om_1\wedge \om_2\in \wedge^4 \aff^4\cong \aff^1$. 
This is a perfect pairing and we can identify 
$(\wedge^2 \aff^4)^*\cong \wedge^2 \aff^4$. 
Computations similar as in the proofs of 
Lemma III--5.18 and Proposition III--5.20 
show that $(\gl_1\times \gl_4\times \gl_2,\wedge^2 \aff^4\otimes \aff^2)$ 
is a regular \pv{} by simple Lie algebra computations, which implies that 
$(M_{\be_{32}},Z_{32})$ is a regular \pv{} also. Moreover, 
\begin{math}
(\mk\times \gl_4(k)\times \gl_2(k)) \backslash 
(\wedge^2 \aff^4\otimes \aff^2)^{\sst}  
\end{math}
is in bijective correspondence with 
$M_{\be_{32}\, k} \backslash Z_{32\, k}^{\sst}$.  
Since the action of $\gl_1$ in $\gl_1\times \gl_4\times \gl_2$ 
can be absorbed to the action of $\gl_2$, 
\begin{math}
(\mk\times \gl_4(k)\times \gl_2(k)) \backslash 
(\wedge^2 \aff^4\otimes \aff^2)^{\sst}  
\end{math}
is the same as 
\begin{math}
(\gl_4(k)\times \gl_2(k)) \backslash 
(\wedge^2 \aff^4\otimes \aff^2)^{\sst}.   
\end{math}
Therefore, 
Proposition III--5.23 implies the following proposition. 

\begin{prop}
\label{prop:case32-interpretation-orbit}
$M_{\be_{32}\, k} \backslash Z_{32\, k}^{\sst}$
is in bijective correspondence with $\Ex_2(k)$. 
\end{prop}

Let 
\begin{equation*}
w = \left(
\begin{pmatrix}
J & 0 \\
0 & 0 
\end{pmatrix}, 
\begin{pmatrix}
0 & 0 \\
0 & J
\end{pmatrix}
\right). 
\end{equation*}
This is a point in 
$(\wedge^2 \aff^4\otimes \aff^2)^{\sst}$ 
which corresponds to the trivial extension $k/k$. 
To find the point in $Z_{32\, k}^{\sst}$ 
which corresponds to $w$, 
it is enough to compute the Pl\"ucker coordinate of $w$. 
Let 
\begin{align*}
& M_1 = \bbmp_{4,1} \wedge \bbmp_{4,2}, \;
M_2 = \bbmp_{4,1} \wedge \bbmp_{4,3}, \;
M_3 = \bbmp_{4,1} \wedge \bbmp_{4,4}, \\
& M_4 = \bbmp_{4,2} \wedge \bbmp_{4,3}, \;
M_5 = \bbmp_{4,2} \wedge \bbmp_{4,4}, \;
M_6 = \bbmp_{4,3} \wedge \bbmp_{4,4}
\end{align*}
and we choose $S=\{M_1\ccd M_6\}$ as a basis for $\wedge^2 \aff^4$. 
Then $w$ can be identified with 
$([1,0\ccd 0],[0\ccd 0,1])$. 

Let $N_{ij}=M_i\wedge M_j$ for $i<j$.  
The only non-zero coordinate of the 
Plu\"cker coordinate of $([1,0\ccd 0],[0\ccd 0,1])$ is 
\begin{math}
N_{16}. 
\end{math}
Let 
\begin{equation*}
w' = (M_2,M_3,M_4,M_5)\in \wedge^2 \aff^4\otimes \aff^4. 
\end{equation*}
Since the Pl\"ucker coordinate of $w'$
is the dual of $N_{16}$, this is the point
which corresponds to $w$.  
The point in $Z_{32}$ which corresponds to $w'$ is 
\begin{equation*}
R_{32} = e_{157} + e_{258} + e_{367} + e_{468}. 
\end{equation*}

We verify Condition \ref{cond:unipotent-connected} (2), (3).  
In this case, $U_{\be_{32}}$ consists of all the elements 
of the form 
\begin{equation*}
n(u) = 
\begin{pmatrix}
1 & 0 & 0 & 0 & 0 & 0 & 0 & 0 \\
0 & 1 & 0 & 0 & 0 & 0 & 0 & 0 \\
0 & 0 & 1 & 0 & 0 & 0 & 0 & 0 \\
0 & 0 & 0 & 1 & 0 & 0 & 0 & 0 \\  
u_{51} & u_{52} & u_{53} & u_{54} & 1 & 0 & 0 & 0 \\ 
u_{61} & u_{62} & u_{63} & u_{64} & 0 & 1 & 0 & 0 \\
u_{71} & u_{72} & u_{73} & u_{74} & 0 & 0 & 1 & 0 \\
u_{81} & u_{82} & u_{83} & u_{84} & 0 & 0 & 0 & 1 
\end{pmatrix}.
\end{equation*}
Then $\dim U_{\be_7}=16$ as an algebraic group and 
$\dim W_{32}=4$.  

We choose the order of $u_{ij}$ and the 
coordinate vectors of $W_7$ in the following manner:
\begin{align*}
& u_{51},u_{52},u_{53},u_{54},u_{63},u_{64},
u_{71},u_{72},u_{73},u_{74},u_{82},u_{84}, \quad 
u_{61},u_{62},u_{81},u_{83}, \\
& e_{567},e_{568},e_{578},e_{678}
\end{align*}
($u_{61}\ccd u_{83}$ correspond to the coordinates of $W_{32}$). 
Explicitly, $n(u)R_{32}-R_{32}$ is
\begin{equation*}
(- u_{61}+u_{52}) e_{567} 
+ (- u_{62}+u_{54}) e_{568} 
+ (u_{81}-u_{72}) e_{578} 
+ (u_{83}-u_{74}) e_{678}.
\end{equation*}
Condition \ref{cond:unipotent-connected} (2), (3) 
are satisfied and so $\rg\backslash S_{\be_{32}\, k}$ 
is also in bijective correspondence with $\Ex_2(k)$.

\vskip 5pt
\noindent
(7) $\beta_{33}=\tfrac {1} {40} (-7,-3,-3,1,1,1,5,5)$ 

In this case, $Z_{33}\cong \m_2 \otimes \aff^3 \oplus 1 \oplus 1$ and 
\begin{equation*}
M_{\be_{33}} = \gl_1 \times M_{[1,3,6]} 
\cong \gl_1\times \gl_1\times \gl_2\times \gl_3\times \gl_2.
\end{equation*}
We express elements of $M_{\be_{33}}$ as $g=(t_0,\diag(t_4,g_2,g_1,g_3))$ 
where $t_0,t_4 \in\gl_1,g_1\in\gl_3$ and $g_2,g_3\in\gl_2$. Then 
\begin{equation*}
M^{\st}_{\be_{33}} = \left\{
\diag(t_4,g_2,g_1,g_3) \mid (\det g_1)(\det g_2)(\det g_3)t_4=1
\right\}.
\end{equation*}
Since 
\begin{equation*}
\chi_{\be_{33}}(g) = t_4^{-7}(\det g_2)^{-3}(\det g_1)(\det g_3)^5
= (\det g_2)^4(\det g_1)^8(\det g_3)^{12}.
\end{equation*}
on $M^{\st}_{\be_{33}}$, 
\begin{equation*}
G_{\st,\be_{33}} = \left\{
\diag(t_4,g_2,g_1,g_3) \;\; \vrule \;\;
\begin{matrix}
g_1\in\gl_3,g_2,g_3\in\gl_2, \\
t_4 = (\det g_1)^{-1}(\det g_2)^{-1}(\det g_3)^{-1}, \\
(\det g_1)^2(\det g_2)(\det g_3)^3=1.  
\end{matrix}  
\right\}. 
\end{equation*}

For $x\in Z_{33}$, let
\begin{equation*}
A_1(x) = \begin{pmatrix}
x_{247} & x_{248} \\
x_{347} & x_{348}
\end{pmatrix}, \\
A_1(x) = \begin{pmatrix}
x_{257} & x_{258} \\
x_{357} & x_{358}
\end{pmatrix}, 
A_1(x) = \begin{pmatrix}
x_{267} & x_{268} \\
x_{367} & x_{368}
\end{pmatrix}
\end{equation*}
and $A(x)=(A_1(x),A_2(x),A_3(x))$. Then 
$x$ can be identified with $(A(x),x_{178},x_{456})$. 
The action of the above $g\in M_{\be_{33}}$ on $x\in Z_{33}$ 
is as follows:
\begin{equation*}
\begin{pmatrix}
A_1(gx) \\
A_2(gx) \\
A_3(gx)
\end{pmatrix}
= t_0 g_1 \begin{pmatrix}
g_2 A_1(x) \,{}^t\! g_3 \\
g_2 A_2(x) \,{}^t\! g_3 \\
g_2 A_3(x) \,{}^t\! g_3 
\end{pmatrix}
\end{equation*}
and $x_{178},x_{456}$ are multiplied by 
$t_0t_4(\det g_3)$, $t_0\det g_1$ respectively.  

Let $P_1(x)$ be the homogeneous polynomial of degree $6$ 
in (\ref{eq:223-relative-inv-poly}). Then 
\begin{equation*}
P_1(gx) = (\det g_1)^2(\det g_2)^3(\det g_3)^3P_1(x).  
\end{equation*}
We put 
\begin{math}
P(x) = P_1(x) x_{178}^2 x_{456}^2. 
\end{math}
Then 
\begin{equation*}
P(gx)= (\det g_1)^2 (\det g_2)(\det g_3)^3 P(x)
\end{equation*}
on $M^{\st}_{\be_{33}}$. Since 
$(\det g_1)^2 (\det g_2)(\det g_3)^3$ is proportional to 
$\chi_{\be_{33}}(g)$ on $M^{\st}_{\be_{33}}$,  
$P(x)$ is an invariant polynomial 
with respect to $G_{\st,\be_{33}}$.  

Let 
\begin{equation*}
R_{33} = e_{178} -e_{247} + e_{348} + e_{258} + e_{367} + e_{456}. 
\end{equation*}
Then $-e_{247} + e_{348} + e_{258} + e_{367}$ corresponds 
to $w$ in Proposition \ref{prop:322-orbit-SP} (3). 
So $P_1(R_{33})=1$ and the 
$x_{178},x_{456}$-coordinates are $1$. 
This implies that $R_{33}\in Z_{33\, k}^{\sst}$.  
Therefore, $S_{\be_{33}}\not=\emptyset$. 

\begin{prop}
\label{prop:r-orbit-SP-case33}
$Z_{33\, k}^{\sst} = M_{\be_{33}\, k} \cdot R_{33}$. 
\end{prop}
\begin{proof}
Suppose that $x\in Z_{33\, k}^{\sst}$. 
Proposition \ref{prop:322-orbit-SP} (3) implies that 
there exists $g\in M_{\be_{33}\, k}$ such that 
\begin{math}
A(gx) = A(R_{33}). 
\end{math}
So we may assume that $A(x)=A(R_{33})$. 

Let $t_0,t_4,u\in \mk$ and 
$g=(t_0,\diag(t_4,uI_2,I_3,I_2))\in M_{\be_{33}\, k}$.  
Then 
\begin{equation*}
g x = t_0t_4 x_{178} e_{178} + t_0 u (-e_{247} + e_{348} + e_{258} + e_{367})
+ t_0 x_{456} e_{456}.   
\end{equation*}
Choose $t_0=x_{456}^{-1}$, $t_4=x_{178}^{-1}x_{456}$,  $u=x_{456}$. 
Then  
\begin{equation*}
g x = e_{178} -e_{247} + e_{348} + e_{258} + e_{367}
+ e_{456} = R(33).  
\end{equation*}
(we have used the first $\gl_1$-factor). 
\end{proof}

We verify Condition \ref{cond:unipotent-connected} (2), (3).  
In this case, $U_{\be_{33}}$ consists of all the elements 
of the form 
\begin{equation*}
n(u) = 
\begin{pmatrix}
1 & 0 & 0 & 0 & 0 & 0 & 0 & 0 \\
u_{21} & 1 & 0 & 0 & 0 & 0 & 0 & 0 \\
u_{31} & 0 & 1 & 0 & 0 & 0 & 0 & 0 \\
u_{41} & u_{42} & u_{43} & 1 & 0 & 0 & 0 & 0 \\  
u_{51} & u_{52} & u_{53} & 0 & 1 & 0 & 0 & 0 \\ 
u_{61} & u_{62} & u_{63} & 0 & 0 & 1 & 0 & 0 \\
u_{71} & u_{72} & u_{73} & u_{74} & u_{75}& u_{76} & 1 & 0 \\
u_{81} & u_{82} & u_{83} & u_{84} & u_{85} & u_{86} & 0 & 1 
\end{pmatrix}.
\end{equation*}
Then $\dim U_{\be_{33}}=23$ as an algebraic group and 
$\dim W_{33}=11$.  

We choose the order of $u_{ij}$ and the 
coordinate vectors of $W_{33}$ in the following manner:
\begin{align*}
& u_{71},u_{72},u_{73},u_{74},u_{75},u_{76}, 
u_{81},u_{82},u_{83},u_{84},u_{85},u_{86}, \\
& u_{21},u_{31},u_{52},u_{53},u_{62},u_{63},
u_{42},u_{43},u_{41},u_{51},u_{61} \\ 
& e_{278},e_{378},e_{457},e_{567},e_{568},e_{468},  
e_{458},e_{467}, e_{478},e_{578},e_{678}
\end{align*}
($u_{21}\ccd u_{61}$ correspond to the coordinates of $W_{33}$). 
Explicitly, $n(u)R_{33}-R_{33}$ is
\begin{align*}
& (u_{21}+u_{75}+u_{84})e_{278}
+ (u_{31}+u_{74}-u_{86}) e_{378}
+ (u_{52}+u_{76}) e_{457}
+ (u_{53}+u_{74})  e_{567} \\
& + (-u_{62}+u_{84}) e_{568}
+ (u_{63}+u_{85})  e_{468} 
+ (u_{42}-u_{53}+u_{86}) e_{458}
+ (u_{43}+u_{62}+u_{75})  e_{467} \\
& +(u_{41}-u_{73}-u_{82}+u_{42}u_{75}+u_{42}u_{84}
    +u_{43}u_{74}-u_{43}u_{86}+u_{75}u_{86}+u_{76}u_{85}) e_{478} \\
& +(u_{51}-u_{72}+u_{52}u_{75}+u_{52}u_{84}+u_{53}u_{74}
    -u_{53}u_{86}-u_{74}u_{86}+u_{76}u_{84}) e_{578} \\
& +(u_{61}+u_{83} + u_{62}u_{75}+u_{62}u_{84}+u_{63}u_{74}
    -u_{63}u_{86} + u_{74}u_{85}-u_{75}u_{84}) e_{678}
\end{align*}
Condition \ref{cond:unipotent-connected} (2), (3) 
are satisfied and so $\rg\backslash S_{\be_{33}\, k}$ 
also consists of a single point. 

\vskip 5pt
\noindent
(8) $\beta_{37}=\tfrac {1} {104} (-7,-7,-7,1,1,1,9,9)$ 

In this case, $Z_{37}\cong \m_3\otimes \aff^2\oplus 1$ 
and 
\begin{equation*}
M_{\be_{37}} = \gl_1\times M_{[3,6]} 
\cong \gl_1\times \gl_3^2\times \gl_2.
\end{equation*}
We express elements of $M_{\be_{37}}$ as 
$g=(t_0,\diag(g_1,g_2,g_3))$ where $t_0\in\gl_1,g_1,g_2\in \gl_3$ 
and $g_3\in\gl_2$.  Then 
\begin{equation*}
M^{\st}_{\be_{37}} = \{\diag(g_1,g_2,g_3)
\mid (\det g_1)(\det g_2)(\det g_3)=1\}.
\end{equation*}
Since 
\begin{equation*}
\chi_{\be_{37}}(g) = (\det g_1)^{-7}(\det g_2)(\det g_3)^9
= (\det g_2)^8(\det g_3)^{16}
\end{equation*}
on $M^{\st}_{\be_{37}}$, 
\begin{equation*}
G_{\st,\be_{37}} = \left\{
\diag(g_1,g_2,g_3) \;\; \vrule \;\; 
\begin{matrix}
g_1,g_2\in \gl_3, g_3\in\gl_2, \\
(\det g_1) = (\det g_2)^{-1}(\det g_3)^{-1}\}, \\
(\det g_2)(\det g_3)^2=1
\end{matrix}
\right\}.
\end{equation*}

For $x\in Z_{37}$, let 
\begin{equation*}
A_1(x) = \begin{pmatrix}
x_{147} & x_{157} & x_{167} \\
x_{247} & x_{257} & x_{267} \\
x_{347} & x_{357} & x_{367} 
\end{pmatrix}, \;
A_1(x) = \begin{pmatrix}
x_{148} & x_{158} & x_{168} \\
x_{248} & x_{258} & x_{268} \\
x_{348} & x_{358} & x_{368} 
\end{pmatrix}
\end{equation*}
and $A(x)=(A_1(x),A_2(x))$. We identify $x$ with 
$(A(x),x_{456})$.  
The action of the above $g\in M_{\be_{37}}$ 
on $x\in Z_{37}$ is as follows:  
\begin{equation*}
\begin{pmatrix}
A_1(gx) \\
A_2(gx) 
\end{pmatrix}
= t_0g_3 \begin{pmatrix}
g_1 A_1(x) {}^t\! g_2 \\
g_1 A_2(x) {}^t\! g_2 \\
\end{pmatrix}
\end{equation*}
and $x_{456}$ is multiplied by $t_0\det g_2$.  

Let $F_x(v)$ be the binary cubic form of variables $v=(v_1,v_2)$
defined for $A(x)$ before Proposition III--5.15 
and $P_1(x)$ the discriminant of $F_x(v)$. Then 
\begin{equation*}
P_1(gx) = (\det g_1)^4(\det g_2)^4(\det g_3)^6P_1(x)
= (\det g_3)^2P_1(x)
\end{equation*}
on $M^{\st}_{\be_{37}}$. We put 
\begin{math}
P(x) = P_1(x)x_{456}. 
\end{math}
Then 
\begin{equation*}
P(gx) = (\det g_2)(\det g_3)^2P(x)
\end{equation*}
on $M^{\st}_{\be_{37}}$.  
Since $(\det g_2)(\det g_3)^2$ is proportional 
to $\chi_{\be_{37}}(g)$ on $M^{\st}_{\be_{37}}$, 
$P(x)$ is an invariant polynomial 
with respect to $G_{\st,\be_{37}}$. 

Let 
\begin{equation*}
R_{37} = e_{147} - e_{257} + e_{258} - e_{368} + e_{456}. 
\end{equation*}
Then $A(R_{37})$ corresponds to $w$ in (\ref{eq:332-tensor-standard}).
By Proposition III--5.16, $P_1(R_{37})=1$ 
and so $P(R_{37})=1$. This implies that 
$R_{37}\in Z_{37\, k}^{\sst}$.   
Therefore, $S_{\be_{37}}\not=\emptyset$. 

We consider $M_{\be_{37}\, k}\backslash Z_{37\, k}^{\sst}$.  
Suppose that $x\in Z_{37\, k}^{\sst}$.  Then 
$x_{456}\not=0$. The coefficient of $e_{456}$ of  
$(t_0,\diag(g_1,g_2,g_3))x$ is $t_0(\det g_2)x_{456}$. 
So we may assume that $x_{456}=1$ by applying an element 
of $M_{\be_{37}\, k}$. This condition does not change 
by the action of $(t_0,\diag(g_1,g_2,g_3))$ 
if and only if $t_0=(\det g_2)^{-1}$. 
If we put $g_1'=(\det g_2)^{-1}g_1$ then 
\begin{equation*}
((\det g_2)^{-1},\diag(g_1,g_2,g_3))x
= (1,\diag(g_1',g_2,g_3))x. 
\end{equation*}
So we can absorb the factor $(\det g_2)^{-1}$.  
Therefore, 
$M_{\be_{37}\, k}\backslash Z_{37\, k}^{\sst}$
is in bijective correspondence with the set of 
generic rational orbits of the \pv{} (\ref{eq:332-tensor-standard}). 
This implies the following proposition. 
\begin{prop}
\label{prop:r-orbits-case37}
$M_{\be_{37}\, k}\backslash Z_{37\, k}^{\sst}$ is in bijective 
correspondence with $\Ex_3(k)$ by associating the field generated by 
a root of $F_x(v)$.  
\end{prop}

We verify Condition \ref{cond:unipotent-connected} (2), (3).  
In this case, $U_{\be_{37}}$ consists of all the elements 
of the form 
\begin{equation*}
n(u) = 
\begin{pmatrix}
1 & 0 & 0 & 0 & 0 & 0 & 0 & 0 \\
0 & 1 & 0 & 0 & 0 & 0 & 0 & 0 \\
0 & 0 & 1 & 0 & 0 & 0 & 0 & 0 \\
u_{41} & u_{42} & u_{43} & 1 & 0 & 0 & 0 & 0 \\  
u_{51} & u_{52} & u_{53} & 0 & 1 & 0 & 0 & 0 \\ 
u_{61} & u_{62} & u_{63} & 0 & 0 & 1 & 0 & 0 \\
u_{71} & u_{72} & u_{73} & u_{74} & u_{75}& u_{76} & 1 & 0 \\
u_{81} & u_{82} & u_{83} & u_{84} & u_{85} & u_{86} & 0 & 1 
\end{pmatrix}.
\end{equation*}
Then $\dim U_{\be_{33}}=21$ as an algebraic group and 
$\dim W_{33}=12$.  

We choose the order of $u_{ij}$ and the 
coordinate vectors of $W_{37}$ in the following manner:
\begin{align*}
& u_{41},u_{52},u_{62},u_{63},u_{71},u_{72},u_{83},u_{85},u_{86}, \\
& u_{76},u_{84},u_{42},u_{75},u_{61},u_{43},
u_{74},u_{53},u_{51},u_{73},u_{81},u_{82} \\
& e_{378},e_{178},e_{458},e_{278},e_{467},e_{468},
e_{567},e_{568},e_{457},e_{678},e_{478},e_{578}
\end{align*}
($u_{76}\ccd u_{82}$ correspond to the coordinates of $W_{33}$). 
Explicitly, $n(u)R_{37}-R_{37}$ is 
\begin{align*}
& - u_{76} e_{378} 
- u_{84} e_{178} 
+ (u_{42}+u_{86}) e_{458} 
+ (u_{75}+u_{85}) e_{278} 
- (u_{61}+u_{75}) e_{467} \\
& - (u_{43}+u_{85}) e_{468} 
+ (u_{74}+u_{62}) e_{567} 
- (u_{53}+u_{62}-u_{84}) e_{568} 
- (u_{51}+u_{42}-u_{76}) e_{457} \\
& + (u_{73}-u_{61}u_{84}+u_{62}u_{75}+u_{62}u_{85}-u_{63}u_{76} 
  + u_{74}u_{85}-u_{75}u_{84}) e_{678} \\
& + (u_{81}-u_{41}u_{84}+u_{42}u_{75}+u_{42}u_{85}-u_{43}u_{76} 
  +u_{75}u_{86}-u_{76}u_{85}) e_{478} \\
& - (u_{82}+u_{72}+u_{51}u_{84}-u_{52}u_{75}-u_{52}u_{85}+u_{53}u_{76} 
  +u_{74}u_{86}-u_{76}u_{84}) e_{578}
\end{align*}
Condition \ref{cond:unipotent-connected} (2), (3) 
are satisfied and so $\rg\backslash S_{\be_{37}\, k}$ 
is also in bijective correspondence with $\Ex_3(k)$.

\vskip 5pt
\noindent
(9) $\beta_{63}=\tfrac {1} {24} (-9,-1,-1,-1,-1,-1,-1,15)$ 

In this case, $Z_{\beta_{63}}\cong \wedge^2 \aff^6$ and 
\begin{equation*}
M_{\be_{63}} = \gl_1\times M_{[1,7]} \cong 
\gl_1\times \gl_1\times \gl_6\times \gl_1.
\end{equation*}
We express elements of $M_{\be_{63}}$ as 
$g=(t_0,\diag(t_2,g_1,t_3))$ where 
$t_0,t_2,t_3\in\gl_1$ and $g_1\in\gl_6$. 
Then 
\begin{equation*}
M^{\st}_{\be_{63}} = 
\left\{
\diag(t_2,g_1,t_3) \;\;\vrule\;\; 
\begin{matrix}
g_1\in\gl_6,t_2,t_3\in\gl_1, \\
(\det g_1)t_2t_3=1
\end{matrix}
\right\}.
\end{equation*}
Since 
\begin{equation*}
\chi_{\be_{63}}(g)
= t_2^{-9} (\det g_1)^{-1}t_3^{15}
= (\det g_1)^8t_3^{24}
\end{equation*}
on $M^{\st}_{\be_{63}}$, 
\begin{equation*}
G_{\st,\be_{63}} = 
\left\{
\diag((\det g_1)^{-1}t_3^{-1},g_1,t_3) \;\;\vrule\;\; 
\begin{matrix}
g_1\in\gl_6,t_3\in\gl_1, \\
(\det g_1)t_3^3=1
\end{matrix}
\right\}.
\end{equation*}

For $x\in Z_{63}$, let 
\begin{equation*}
A(x) = \begin{pmatrix}
0 & x_{238} & x_{248} & x_{258} & x_{268} & x_{278} \\
-x_{238} & 0 & x_{348} & x_{358} & x_{368} & x_{378} \\
-x_{248} & -x_{348} & 0 & x_{458} & x_{468} & x_{478} \\
-x_{258} & -x_{358} & -x_{458} & 0 & x_{568} & x_{578} \\
-x_{268} & -x_{368} & -x_{468} & -x_{568} & 0 & x_{678} \\
-x_{278} & -x_{378} & -x_{478} & -x_{578} & -x_{678} & 0
\end{pmatrix}.
\end{equation*}
We identify $x$ with $A(x)$.  
The action of the above $g\in M_{\be_{63}}$ 
on $x\in Z_{63}$ is as follows:
\begin{equation*}
A(gx) = t_0t_3 g_1 A(x) \,{}^t\!g_1.
\end{equation*}

Let $P(x)$ be the Pfaffian of $A(x)$. Then 
\begin{equation*}
P(gx) = (\det g_1)t_3^3 P(x)
\end{equation*}
on $M^{\st}_{\be_{63}}$.  
Since $(\det g_1)t_3^3$ is proportional to $\chi_{\be_{63}}(g)$  
on $M^{\st}_{\be_{63}}$, $P(x)$ is an invariant polynomial 
with respect to $G_{\st,\be_{63}}$. 

Let $R_{63}\in Z_{63}$ be the element such that
\begin{equation*}
A(R_{63}) = 
\begin{pmatrix}
J & 0 & 0 \\
0 & J & 0 \\
0 & 0 & J
\end{pmatrix}.  
\end{equation*}
Then $P(R_{63})=1$, which implies that 
$R_{63}\in Z_{63\, k}^{\sst}$. 
Therefore, $S_{\be_{63}}\not=\emptyset$.

If $P(x)\not=0$ then $A(x)$ is non-degenerate. 
By Witt's theorem (Theorem \ref{thm:alternating-matrix}), 
there exists $g\in M_{\be_{63}}$ 
such that $A(gx)=A(R_{63})$. Therefore, 
we obtain the following proposition.
\begin{prop}
\label{prop:r-orbit-case63}
$Z_{63\, k}^{\sst} = M_{\be_{63}\, k}\cdot R_{63}$. 
\end{prop}
Since $W_{63}=\{0\}$, $\rg\backslash S_{\be_{63}\, k}$
also consists of a single point.

\vskip 5pt
\noindent 
(10) $\beta_{74}=\tfrac {1} {88} (-9,-9,-1,-1,-1,-1,7,15)$ 

In this case, $Z_{74}\cong \wedge^2 \aff^4 \oplus \m_{4,2}$ 
and 
\begin{equation*}
M_{\be_{74}} = \gl_1\times M_{[2,6,7]} \cong 
\gl_1\times \gl_2\times \gl_4\times \gl_1^2.  
\end{equation*}
We express elements of $M_{\be_{74}}$ as 
$g=(t_0,\diag(g_2,g_1,t_3,t_4))$ where 
$t_0,t_3,t_4\in\gl_1$, $g_1\in\gl_4$ 
and $g_2\in\gl_2$. Then 
\begin{equation*}
M^{\st}_{\be_{74}} = 
\left\{
\diag(g_2,g_1,t_3,t_4) \;\;\vrule\;\; 
\begin{matrix}
g_1\in\gl_4,g_2\in\gl_2,t_3,t_4\in\gl_1, \\
(\det g_1)(\det g_2)t_3t_4=1
\end{matrix}
\right\}.
\end{equation*}
Since 
\begin{equation*}
\chi_{\be_{74}}(g) = (\det g_2)^{-9}
(\det g_1)^{-1}t_3^7t_4^{15} 
= (\det g_1)^{-16}(\det g_2)^{-24}t_3^{-8}
\end{equation*}
on $M^{\st}_{\be_{74}}$, 
\begin{equation*}
G_{\st,\be_{74}} = 
\left\{
\diag(g_2,g_1,t_3,t_4) 
\;\;\vrule\;\; 
\begin{matrix}
g_1\in\gl_4,g_2\in\gl_2,t_3\in\gl_1, \\
t_4 = (\det g_1)^{-1}(\det g_2)^{-1}t_3^{-1}, \\
(\det g_1)^2(\det g_2)^3t_3=1
\end{matrix}
\right\}.
\end{equation*}

For $x\in Z_{74}$, let
\begin{equation*}
v_1(x) = \begin{pmatrix}
x_{138} \\  
x_{148} \\ 
x_{158} \\ 
x_{168}   
\end{pmatrix}, \;
v_2(x) = \begin{pmatrix}
x_{238} \\  
x_{248} \\ 
x_{258} \\ 
x_{268}  
\end{pmatrix}, \;
A(x) = \begin{pmatrix}
0 & x_{347} & x_{357} & x_{367} \\
-x_{347} & 0 & x_{457} & x_{467} \\
-x_{357} & -x_{457} & 0 & x_{567} \\
-x_{367} & -x_{467} & -x_{567} & 0 
\end{pmatrix}
\end{equation*}
and $B(x)=(v_1(x)\; v_2(x))$. 
We identify $x$ with $(A(x),B(x))$.  
The action of the above 
$g\in M_{\be_{74}}$ on $x\in Z_{74}$ is as follows:
\begin{equation*}
A(gx) = t_0t_3 g_1 A(x) \,{}^t\! g_1, \; 
B(gx) = t_0t_4
g_1 B(x) \,{}^t\! g_2, \;
\end{equation*}
Let $P_1(x) = v_1(x)\wedge v_2(x)\wedge A(x)$ 
regarding $A(x)\in \wedge^2 \aff^4$.  
Let $P_2(x)$ be the Pfaffian of $A(x)$. 
Then 
\begin{equation*}
P_1(gx) = t_4P_1(x), \; 
P_2(gx) = t_3^2(\deg g_1)P_2(x)
\end{equation*}
on $M^{\st}_{\be_{74}}$. Let $P(x)=P_1(x)^3 P_2(x)$. 
Then 
\begin{equation*}
P(gx) = (\det g_1)^{-2}(\det g_2)^{-3}t_3^{-1} P(x)
\end{equation*}
on $M^{\st}_{\be_{74}}$. 
Since $(\det g_1)^{-2}(\det g_2)^{-3}t_3^{-1}$ 
is proportional to $\chi_{\be_{74}}(g)$ on $M^{\st}_{\be_{74}}$,  
$P(x)$ is an invariant polynomial
with respect to $G_{\st,\be_{74}}$. 
So $x\in Z_{74}^{\sst}$ 
if and only if $P_1(x),P_2(x)\not=0$.  

Let
\begin{equation*}
R_{74} = e_{138} + e_{248} + e_{347} + e_{567}. 
\end{equation*}
Then 
\begin{equation*}
v_1(R_{74}) = 
\begin{pmatrix}
1 \\ 0 \\ 0 \\ 0
\end{pmatrix}, \;
v_2(R_{74}) = 
\begin{pmatrix}
0 \\ 1 \\ 0 \\ 0
\end{pmatrix}, \;
A(R_{74})
= \begin{pmatrix}
J & 0 \\
0 & J 
\end{pmatrix}. 
\end{equation*}
Therefore, $P_1(R_{74})=P_2(R_{74})=P(R_{74})=1$, 
which implies that $R_{74}\in Z_{74\, k}^{\sst}$.  
Therefore, $S_{\be_{74}}\not=\emptyset$.  

\begin{prop}
$Z_{74\, k}^{\sst}=M_{\be_{74}\, k}\cdot R_{74}$. 
\end{prop}
\begin{proof}
Suppose that $x\in Z_{74\, k}^{\sst}$. 
Then $A(x)$ is non-degenerate. 
By Witt's theorem 
(Theorem \ref{thm:alternating-matrix}), 
we may assume that $A(x)=A(R_{74})$. 
Let $h_1,h_2\in\spl_2(k)$, $u\in k$, 
$T(h_1,h_2) = \diag(h_1,h_2)$ and 
\begin{align*}
& \tau = \begin{pmatrix}
0 & I_2 \\
I_2 & 0 
\end{pmatrix}, \; 
n_1(u) = \begin{pmatrix}
1 & 0 & 0 & 0 \\
0 & 1 & u & 0 \\
0 & 0 & 1 & 0 \\
u & 0 & 0 & 1 
\end{pmatrix}, \;
n_2(u) = \begin{pmatrix}
1 & 0 & 0 & u \\
0 & 1 & 0 & 0 \\
0 & u & 1 & 0 \\
0 & 0 & 0 & 1 
\end{pmatrix}.
\end{align*}
If $h$ is any of these matrices, 
\begin{equation*}
A(\diag(I_2,h,1,1)R_{74})
= h A(R_{74}) {}^t h = A(R_{74}).   
\end{equation*}

Since $P_1(x)\not=0$, $v_1(x)\not=0$.
By applying $\diag(I_2,\tau,1,1)$ and an element 
of the form $\diag(I_2,T(h_1,h_2),1,1)$  
if necessary, we may assume that $v_1(x)$ is in the form 
$v_1(x)=[1,0,0,*]$.  By applying an element of the form 
$\diag(I_2,n_1(u),1,1)$, 
$v_1(x)$ becomes $[1,0,0,0]$.   
Since $P_1(x)\not=0$, the second entry of $v_2(x)$ 
is non-zero. By applying an element of the form 
$\diag(I_2,T(h_1,h_2),1,1)$ where $h_1$ 
is upper triangular with diagonal entries $1$, 
we may assume that $v_2(x)$ is in the form 
$[0,1,*,0]$. By applying an element of the form 
$\diag(I_2,n_2(u),1,1)$, $v_2(x)$ becomes $[0,1,0,0]$. 
\end{proof}

We verify Condition \ref{cond:unipotent-connected} (2), (3).  
In this case, $U_{\be_{74}}$ consists of all the elements 
of the form 
\begin{equation*}
n(u) = 
\begin{pmatrix}
1 & 0 & 0 & 0 & 0 & 0 & 0 & 0 \\
0 & 1 & 0 & 0 & 0 & 0 & 0 & 0 \\
u_{31} & u_{32} & 1 & 0 & 0 & 0 & 0 & 0 \\
u_{41} & u_{42} & 0 & 1 & 0 & 0 & 0 & 0 \\  
u_{51} & u_{52} & 0 & 0 & 1 & 0 & 0 & 0 \\ 
u_{61} & u_{62} & 0 & 0 & 0 & 1 & 0 & 0 \\
u_{71} & u_{72} & u_{73} & u_{74} & u_{75}& u_{76} & 1 & 0 \\
u_{81} & u_{82} & u_{83} & u_{84} & u_{85} & u_{86} & u_{87} & 1 
\end{pmatrix}.
\end{equation*}
Then $\dim U_{\be_{74}}=21$ as an algebraic group and 
$\dim W_{74}=12$.

We choose the order of $u_{ij}$ and the 
coordinate vectors of $W_{74}$ in the following manner:
\begin{align*}
& u_{31},u_{41},u_{42},u_{71},u_{72},u_{75},u_{76},u_{81},u_{82}, \\
& u_{51},u_{52},u_{61},u_{62},u_{73},u_{74},
u_{87},u_{32},u_{83},u_{84},u_{85},u_{86} \\
& e_{358},e_{458},e_{368},e_{468},e_{178},e_{278} 
e_{568},e_{348},e_{478},e_{378},e_{678},e_{578} 
\end{align*}
($u_{51}\ccd u_{86}$ correspond to the coordinates of $W_{74}$). 
Explicitly, $n(u)R_{74}-R_{74}$ is 
\begin{align*}
& - u_{51} e_{358}
- u_{52} e_{458} 
- u_{61} e_{368} 
- u_{62} e_{468} 
+ u_{73} e_{178} 
+ u_{74} e_{278} 
- u_{87} e_{568} \\
& + (u_{32}+u_{87}-u_{41}) e_{348} 
+ (u_{83}-u_{72}+u_{41}u_{73}+u_{42}u_{74}-u_{73}u_{87}) e_{478} \\
& + (-u_{84}-u_{71}+u_{31}u_{73}+u_{32}u_{74}+u_{74}u_{87}) e_{378} \\
& + (-u_{85} + u_{61}u_{73}+u_{62}u_{74}+u_{75}u_{87}) e_{678} 
+ (u_{86}+u_{51}u_{73}+u_{52}u_{74}-u_{76}u_{87}) e_{578}.   
\end{align*}
Condition \ref{cond:unipotent-connected} (2), (3) 
are satisfied and so $\rg\backslash S_{\be_{37}\, k}$ 
also consists of a single point.

\vskip 5pt
\noindent 
(11) $\beta_{81}=\tfrac {1} {56} (-13,-5,-5,-5,3,3,11,11)$ 

In this case, $Z_{\beta_{81}}\cong \m_2\otimes \aff^3\oplus 1$
and 
\begin{equation*}
M_{\be_{81}} = \gl_1\times M_{[1,4,6]}
\cong \gl_1\times \gl_3\otimes \gl_2^2.  
\end{equation*}
We express elements of $M_{\be_{81}}$ as 
$g=(t_0,\diag(t_4,g_1,g_2,g_3))$ where $t_0,t_4\in \gl_1$, 
$g_1\in \gl_3$ and $g_2,g_3\in\gl_2$.  Then 
\begin{equation*}
M^{\st}_{\be_{81}} = \left\{
\diag(t_4,g_1,g_2,g_3) \;\;\vrule\;\;
\begin{matrix}
g_1\in \gl_3,g_2,g_3\in\gl_2, t_4\in\gl_1\\
(\det g_1)(\det g_2)(\det g_3)t_4=1.
\end{matrix}
\right\}
\end{equation*}
Since 
\begin{align*}
\chi_{\beta_{81}}(g)=t_4^{-13}(\det g_1)^{-5}(\det g_2)^3(\det g_3)^{11}
=(\det g_1)^8(\det g_2)^{16}(\det g_3)^{24}
\end{align*}
on $M^{\st}_{\be_{81}}$, 
\begin{equation*}
G_{\st,\be_{81}} = \left\{
\diag(t_4,g_1,g_2,g_3) \;\;\vrule\;\;
\begin{matrix}
t_4 = (\det g_1)^{-1}(\det g_2)^{-1}(\det g_3)^{-1}, \\
(\det g_1)(\det g_2)^2(\det g_3)^3=1
\end{matrix}
\right\}.  
\end{equation*}

For $x\in Z_{81}$, let 
\begin{equation*}
A_1(x) = \begin{pmatrix}
x_{257} & x_{258} \\
x_{267} & x_{268} 
\end{pmatrix}, \;
A_2(x) = \begin{pmatrix}
x_{357} & x_{358} \\
x_{367} & x_{368} 
\end{pmatrix}, \;
A_3(x) = \begin{pmatrix}
x_{457} & x_{458} \\
x_{467} & x_{468} 
\end{pmatrix} 
\end{equation*}
and $A(x)=(A_1(x),A_2(x),A_3(x))$. 
We identify $x$ with $(A(x),x_{178})$. 
The action of the above 
\begin{math}
g\in M_{\be_{81}}
\end{math}
on $x$ is as follows: 
\begin{equation*}
\begin{pmatrix}
A_1(gx) \\
A_2(gx) \\
A_3(gx) 
\end{pmatrix}
= t_0g_1 \begin{pmatrix}
g_2 A_1(x) \,{}^t g_3 \\
g_2 A_2(x) \,{}^t g_3 \\
g_2 A_3(x) \,{}^t g_3 
\end{pmatrix}
\end{equation*}
and $x_{178}$ is multiplied by $t_0t_4(\det g_3)$. 

Let $P_1(x)$ be the homogeneous polynomial of degree $6$
in (\ref{eq:223-relative-inv-poly}) and $P(x)=P_1(x)x_{178}$. 
Then, $P_1(gx)=(\det g_1)^2(\det g_2)^3(\det g_3)^3P_1(x)$ and
\begin{align*}
P(gx) &= 
(\det g_1)^2(\det g_2)^3(\det g_3)^3P_1(x) 
(\det g_1)^{-1}(\det g_2)^{-1}x_{178} \\
&= (\det g_1)^1(\det g_2)^2(\det g_3)^3P(x)
\end{align*}
on $M^{\st}_{\be_{81}}$.  
Since $(\det g_1)^1(\det g_2)^2(\det g_3)^3$ is 
proportional to $\chi_{\be_{81}}(g)$ on $M^{\st}_{\be_{81}}$, 
$P(x)$ is an invariant polynomial 
with respect to $G_{\st ,\, \be_{81}}$. 

Let 
\begin{equation*}
R_{81} = e_{178} -e_{257} + e_{268} + e_{358} + e_{467} \in Z_{81}. 
\end{equation*}
Then $-e_{257} + e_{268} + e_{358} + e_{467}$ corresponds to $w$ in 
Proposition \ref{prop:322-orbit-SP} (3). So $P(R_{81}=1$, 
which implies that 
$R_{81}\in Z_{81\, k}^{\sst}$.  Therefore, 
$S_{\be_{81}}\not=\emptyset$.

\begin{prop}
\label{prop:r-orbit-case81}
$Z_{\be_{81}\, k}^{\sst}=M_{\be_{81}\, k} \cdot R_{81}$
\end{prop}
\begin{proof}
Suppose that $x\in Z_{81\, k}^{\sst}$. Then $P_1(x),x_{178}\not=0$. 
Proposition \ref{prop:322-orbit-SP} (3) implies that 
there exists $g\in M_{\be_{81}\, k}$ such that 
$A(gx)=A(R_{81})$. Then 
\begin{equation*}
(x_{178}^{-1},\diag(x_{178}I_3,I_2,I_2)) x = R_{81}. 
\end{equation*}
\end{proof}

We verify Condition \ref{cond:unipotent-connected} (2), (3).  
In this case, $U_{\be_{81}}$ consists of all the elements 
of the form 
\begin{equation*}
n(u) = 
\begin{pmatrix}
1 & 0 & 0 & 0 & 0 & 0 & 0 & 0 \\
u_{21} & 1 & 0 & 0 & 0 & 0 & 0 & 0 \\
u_{31} & 0 & 1 & 0 & 0 & 0 & 0 & 0 \\
u_{41} & 0 & 0 & 1 & 0 & 0 & 0 & 0 \\  
u_{51} & u_{52} & u_{53} & u_{54} & 1 & 0 & 0 & 0 \\ 
u_{61} & u_{62} & u_{63} & u_{64} & 0 & 1 & 0 & 0 \\
u_{71} & u_{72} & u_{73} & u_{74} & u_{75}& u_{76} & 1 & 0 \\
u_{81} & u_{82} & u_{83} & u_{84} & u_{85} & u_{86} & 0 & 1 
\end{pmatrix}.
\end{equation*}
Then $\dim U_{\be_{81}}=23$ as an algebraic group and 
$\dim W_{81}=7$.

We choose the order of $u_{ij}$ and the 
coordinate vectors of $W_{81}$ in the following manner:
\begin{align*}
& u_{71}\ccd u_{76},u_{81}\ccd u_{86},u_{52},u_{53},u_{54},u_{64}, \quad
u_{21},u_{31},u_{41},u_{62},u_{63},u_{51},u_{61} \\
& e_{278},e_{378},e_{478},e_{567},e_{568},e_{578},e_{678}
\end{align*}
($u_{21}\ccd u_{61}$ correspond to the coordinates of $W_{81}$). 
Explicitly, $n(u)R_{81}-R_{81}$ is
\begin{align*}
& (u_{21} + u_{76} + u_{85}) e_{278} 
+ (u_{31} + u_{75}) e_{378} 
+ (u_{41} - u_{86}) e_{478} 
+ (u_{62} + u_{54}) e_{567} \\
& + (-u_{63} + u_{52}) e_{568}
+ (u_{51}-u_{73}-u_{82}+u_{52}u_{76}+u_{52}u_{85}
   +u_{53}u_{75}-u_{54}u_{86}) e_{578} \\
& + (u_{61}-u_{72}+u_{84}+u_{62}u_{76}
   +u_{62}u_{85}+u_{63}u_{75}-u_{64}u_{86}) e_{678}
\end{align*}
Condition \ref{cond:unipotent-connected} (2), (3) 
are satisfied and so $\rg\backslash S_{\be_{81}\, k}$ 
also consists of a single point.

\vskip 5pt
\noindent
(12) $\beta_{98}=\tfrac {1} {8} (-3,-3,1,1,1,1,1,1)$ 

In this case, $Z_{\beta_{98}}\cong \wedge^3\aff^6$ and 
\begin{equation*}
M_{\be_{98}} = \gl_1\times M_{[2]} \cong \gl_1\times \gl_2\times \gl_6.
\end{equation*}
We express elements of $M_{\be_{98}}$ as 
$g=(t_0,\diag(g_2,g_1))$ where $t_0\in\gl_1$, $g_2\in\gl_2$ 
and $g_1\in\gl_6$. Then 
\begin{equation*}
M^{\st}_{\be_{98}} = 
\left\{
\diag(g_2,g_1)\mid g_1\in \gl_6,g_2\in \gl_2, \;
(\det g_1)(\det g_2)=1
\right\}.
\end{equation*}
Since
\begin{equation*}
\chi_{98}(g) = (\det g_2)^{-3}(\det g_1)
= (\det g_1)^4
\end{equation*}
on $M^{\st}_{\be_{98}}$, 
\begin{equation*}
G_{\st,\be_{98}} 
= \left\{
\diag(g_2,g_1)\mid g_1\in \spl_6,g_2\in \spl_2
\right\}.
\end{equation*}

Since $\spl_6\times \spl_2$ is semi-simple, 
any relative invariant polynomial on $G_{\st,\be_{98}}$ 
is an invariant polynomial. 
The $\gl_2$-factor acts on $Z_{98}$ trivially and the action of 
$\gl_1\times \gl_6$ on $Z_{98}$ can be identified 
with the action in Case IV of Subsection \ref{regular-wedge36}.  
Therefore, Theorem \ref{thm:case98-interpretation} 
implies the following proposition. 

\begin{prop}
\label{prop:r-orbit-case98}
$M_{\be_{98}\, k} \backslash Z_{98\, k}^{\sst}$ 
is in bijective correspondence with $\Ex_2(k)$. 
\end{prop}

Since $W_{\beta_{98}}=\{0\}$, 
$\rg\backslash S_{\be_{98}\, k}$
is also in bijective correspondence with $\Ex_2(k)$.

\vskip 5pt
\noindent
(13) $\beta_{108}=\tfrac {1} {24} (-9,-1,-1,-1,3,3,3,3)$ 

In this case $Z_{108}\cong \wedge^2 \aff^4\otimes \aff^3$ and 
\begin{equation*}
M_{\be_{108}} = \gl_1\times M_{[1,4]} \cong 
\gl_1\times \gl_1\times \gl_3\times \gl_4.  
\end{equation*}
We  express elements of $M_{\be_{108}}$ as 
$g=(t_0,\diag(t_3,g_2,g_1))$ where $t_0,t_3\in\gl_1$, 
$g_1\in \gl_4$ and $g_2\in\gl_3$. Then 
\begin{equation*}
M^{\st}_{\be_{98}} = \left\{
\diag(t_3,g_2,g_1) \;\; \vrule \;\; 
\begin{matrix}
g_1\in \gl_4,g_2\in\gl_3,t_3\in\gl_1, \\
(\det g_1)(\det g_2)t_3=1
\end{matrix}
\right\}.  
\end{equation*}
Since 
\begin{equation*}
\chi_{\be_{108}}(g) = t_3^{-9}(\det g_2)^{-1}(\det g_1)^3
= (\det g_1)^{12}(\det g_2)^8
\end{equation*}
on $M^{\st}_{\be_{98}}$, 
\begin{equation*}
G_{\st,\be_{108}} = \left\{
(\det g_1)^{-1}(\det g_2)^{-1},g_2,g_1) \mid 
(\det g_1)^3(\det g_2)^2=1
\right\}.
\end{equation*}

For $x\in Z_{108}$, let 
\begin{equation*}
A_1(x) = 
\begin{pmatrix}
0 & x_{256} & x_{257} & x_{258} \\
-x_{256} & 0 & x_{267} & x_{268} \\
-x_{257} & -x_{267} & 0 & x_{278} \\
-x_{258} & -x_{268} & -x_{278} & 0 
\end{pmatrix}. 
\end{equation*}
We define $A_2(x),A_3(x)$ by replacing $2$ in the indices by 
$3,4$ respectively  and $A(x)=(A_1(x),A_2(x),A_3(x))$.  
The action of the above $g\in M_{\be_{108}}$ on 
$x\in Z_{108}$ is the same as the first part of 
(\ref{eq:g-action-A-case4}).

This is the case (a) of Section III--8. 
Let 
\begin{equation*}
R(108) = e_{256} + e_{358} - e_{367} + e_{478}.
\end{equation*}
Then 
\begin{equation*}
A_1(R_{108}) = \begin{pmatrix}
J & 0 \\
0 & 0
\end{pmatrix}, \quad
A_2(R_{108}) = \begin{pmatrix}
0 & J \\
J & 0 
\end{pmatrix}, \quad
A_3(R_{108})  = \begin{pmatrix}
0 & 0 \\
0 & J 
\end{pmatrix}. 
\end{equation*}
So $R_{108}$ corresponds to $w$ in (III--9.1). 
Since $(1,\diag(t_3,I_3,I_4))$ acts on $Z_{108}$ 
trivially and the action of $(t_0,\diag(1,I_3,I_4))$ 
can be absorbed to the action of $(1,\diag(1,t_0I_3,I_4))$, 
the action of $M_{\be_{108}}$ on $Z_{108}$ 
is essentially the same as the case (a) of 
Section III--9.

By Proposition III--9.2 and  
Corollary III--9.10, 
$(M_{\be_{108}},Z_{108})$ is a regular \pv.  
Lemma III--9.13  implies  
that there is a homogeneous polynomial $P(x)$ of degree $12$ 
of $x\in Z_{108}$ such that 
$P(gx)=(\det g_1)^6(\det g_2)^4P(x)$ and $P(R_{108})=1$. 
Therefore, $R_{108}\in Z_{108\, k}^{\sst}$ and 
$M_{\be_{108}}\cdot R_{108}\sub Z_{108}$ is Zariski open.

Let $\mathrm{IQF}_4$ be the set of 
$k$-isomorphism classes of inner forms over $k$ 
of the quadratic form $v_1v_4-v_2v_3$. 
The following proposition follows from 
Proposition III--9.16. 
\begin{prop}
$M_{\be_{108}\, k}\backslash Z_{108\, k}^{\sst}$  
is in bijective correspondence with $\mathrm{IQF}_4(k)$.   
\end{prop}

We verify Condition \ref{cond:unipotent-connected} (2), (3).  
In this case, $U_{\be_{108}}$ consists of all the elements 
of the form 
\begin{equation*}
n(u) = 
\begin{pmatrix}
1 & 0 & 0 & 0 & 0 & 0 & 0 & 0 \\
u_{21} & 1 & 0 & 0 & 0 & 0 & 0 & 0 \\
u_{31} & 0 & 1 & 0 & 0 & 0 & 0 & 0 \\
u_{41} & 0 & 0 & 1 & 0 & 0 & 0 & 0 \\  
u_{51} & u_{52} & u_{53} & u_{54} & 1 & 0 & 0 & 0 \\ 
u_{61} & u_{62} & u_{63} & u_{64} & 0 & 1 & 0 & 0 \\
u_{71} & u_{72} & u_{73} & u_{74} & 0 & 0 & 1 & 0 \\
u_{81} & u_{82} & u_{83} & u_{84} & 0 & 0 & 0 & 1 
\end{pmatrix}.
\end{equation*}
Then $\dim U_{\be_{108}}=19$ as an algebraic group and 
$\dim W_{108}=4$.

We choose the order of $u_{ij}$ and the 
coordinate vectors of $W_{108}$ in the following manner:
\begin{align*}
& u_{7j},u_{8j}\;(j=1,2,3,4),
u_{21},u_{31},u_{41},u_{51},u_{52},u_{61},u_{62}, \quad 
u_{53},u_{54},u_{63},u_{64} \\
& e_{567},e_{578},e_{568},e_{678}
\end{align*}
($u_{53},u_{54},u_{63},u_{64}$ correspond 
to the coordinates of $W_{108}$). 
Explicitly, $n(u)R_{108}-R_{108}$ is 
\begin{align*}
(- u_{53}+u_{72}) e_{567} 
+ (u_{54}-u_{73}) e_{578}
+ (- u_{63}+u_{82}) e_{568} 
+ (u_{64}-u_{83}) e_{678}.  
\end{align*}
Condition \ref{cond:unipotent-connected} (2), (3) 
are satisfied and so $\rg\backslash S_{\be_{108}\, k}$ 
is also in bijective correspondence with 
$\mathrm{IQF}_4(k)$.

\vskip 5pt   
\noindent      
(14) $\be_{152}= \frac 1{72}(-11,-11,-3,-3,5,5,5,13)$ 

In this case, 
\begin{math}
Z_{152} \cong \aff^2\otimes \aff^3 \oplus 
\aff^2\otimes \wedge^2 \aff^3 \oplus 1
\end{math}
and 
\begin{equation*}
M_{\be_{152}} = \gl_1\times M_{[2,4,7]}
\cong \gl_1\times \gl_2^2\times \gl_3\times \gl_1.
\end{equation*}
We express elements of $M_{\be_{152}}$ as 
\begin{math}
g = (t_0,\diag(g_2,g_3,g_1,t_4))
\end{math}
where $t_0,t_4\in\gl_1,g_1\in\gl_3$ and $g_2,g_3\in\gl_2$. 
Then 
\begin{equation*}
M_{\be_{152}} = 
\left\{
\diag(g_2,g_3,g_1,t_4) \;\;\vrule\;\;
\begin{matrix}
g_1\in\gl_3,g_2,g_3\in\gl_2,t_4\in\gl_1, \\
(\det g_1)(\det g_2)(\det g_3)t_4=1
\end{matrix}
\right\}.  
\end{equation*}
Since 
\begin{equation*}
\chi_{\be_{152}}(g) = 
(\det g_2)^{-11}(\det g_3)^{-3}(\det g_1)^5t_4^{13}
= (\det g_1)^{-8}(\det g_2)^{-24}(\det g_3)^{-16}
\end{equation*}
on $M^{\st}_{\be_{152}}$, 
\begin{equation*}
G_{\st,\be_{152}} = 
\left\{
\diag(g_2,g_3,g_1,t_4) \;\;\vrule\;\;
\begin{matrix}
g_1\in\gl_3,g_2,g_3\in\gl_2,t_4\in\gl_1, \\
(\det g_1)(\det g_2)(\det g_3)t_4=1, \\
(\det g_1)(\det g_2)^3(\det g_3)^2
\end{matrix}
\right\}.  
\end{equation*}

For $x\in Z_{152}$, let
\begin{equation*}
A(x) = \begin{pmatrix}
x_{158} & x_{168} & x_{178} \\
x_{258} & x_{268} & x_{278} \\
\end{pmatrix}, \;
B(x) = \begin{pmatrix}
x_{367} & -x_{357} & x_{356} \\
x_{467} & -x_{457} & x_{456}
\end{pmatrix}. 
\end{equation*}
We identify $x$ with $(A(x),B(x),x_{348})$.  
Then the action of the above $g\in M_{\be_{152}}$ 
on $x\in Z_{152}$ is as follows: 
\begin{equation*}
A(gx) = t_0 t_4 g_2 A(x)\, {}^t \!g_1, \;
B(gx) = t_0 (\det g_1)g_3B(x) g_1^{-1}
\end{equation*}
and $x_{348}$ is multiplied by $t_0t_4(\det g_3)$.  

Let
\begin{equation*}
P_1(x) = \det (A(x)\,{}^t\!B(x)), \; P(x) = P_1(x)^2x_{348}.  
\end{equation*}
Then 
\begin{align*}
P_1(gx) & = t_4^2 (\det g_1)^2 (\det g_2)(\det g_3) P_1(x) 
= (\det g_2)^{-1}(\det g_3)^{-1}P_1(x), \\
P(gx) & = t_4 (\det g_3) (\det g_2)^{-2}(\det g_3)^{-2} P(x) 
=  (\det g_1)^{-1}(\det g_2)^{-3}(\det g_3)^{-2}P(x)
\end{align*}
on $M^{\st}_{\be_{152}}$. 
Since $(\det g_1)^{-1}(\det g_2)^{-3}(\det g_3)^{-2}$ 
is proportional to $\chi_{\be_{152}}(g)$ on $M^{\st}_{\be_{152}}$, 
$P(x)$ is an invariant polynomial 
with respect to $G_{\st,\be_{152}}$. This implies that 
$x\in V^{\sst}$ if and only if $P_1(x),x_{348}\not=0$.

Let 
\begin{equation*}
R_{152} = e_{158} + e_{268} + e_{367} - e_{457} + e_{348}. 
\end{equation*}
Then 
\begin{equation*}
A(R_{152}) = \begin{pmatrix}
1 & 0 & 0 \\
0 & 1 & 0 \\
\end{pmatrix}, \;
B(R_{152}) = \begin{pmatrix}
1 & 0 & 0 \\
0 & 1 & 0 
\end{pmatrix}
\end{equation*}
and the coefficient of $e_{348}$ is $1$.
So $P(R_{152})=1$, which implies that 
$R_{152}\in Z_{152\, k}^{\sst}$.  Therefore, 
$S_{\be_{152}}\not=\emptyset$.  

\begin{prop}
\label{prop:r-orbit-case152}
$Z_{152\, k}^{\sst} = M_{\be_{152}\, k} \cdot R_{152}$. 
\end{prop}
\begin{proof}
Since $A(x)\,{}^t\!B(x)$ is non-singular, 
the ranks of $A(x),B(x)$ 
are $2$.  So there exists $g_1\in\gl_3(k)$ such that 
$A(\diag(I_2,I_2,g_1,1)x)=A(R_{152})$. 
Therefore, we may assume that $A(x)=A(R_{152})$.  
Since $A(x)\,{}^t\!B(x)$ is non-singular, 
the first $2\times 2$ block of $B(x)$ is non-singular. 
So there exists $g_3\in \gl_2(k)$ such that 
\begin{equation*}
B(\diag(I_2,g_3,I_3,1)x) = 
\begin{pmatrix}
1 & 0 & * \\
0 & 1 & * 
\end{pmatrix}
\end{equation*}
and $A(x)$ does not change.  

Choose $g_1\in\gl_3(k)$ such that $g_1^{-1}$ 
is in the form 
\begin{equation*}
g_1^{-1} = \begin{pmatrix}
1 & 0 & * \\
0 & 1 & * \\
0 & 0 & 1 
\end{pmatrix}
\end{equation*}
and 
\begin{math}
B(x) g_1^{-1} = B(R_{152}). 
\end{math}
It is easy to see that $A(R_{152})\,{}^t\!g_1 = A(R_{152})$. 
So 
\begin{equation*}
A(\diag(I_2,I_2,g_1,1)x) = A(R_{152}), \;
B(\diag(I_2,I_2,g_1,1)x) = B(R_{152}). 
\end{equation*}
Then 
\begin{equation*}
\diag(x_{348}I_2,I_2,I_3,x_{348}^{-1}) x = R_{152}.  
\end{equation*}
\end{proof}

We verify Condition \ref{cond:unipotent-connected} (2), (3).  
In this case, $U_{\be_{152}}$ consists of all the elements 
of the form 
\begin{equation*}
n(u) = 
\begin{pmatrix}
1 & 0 & 0 & 0 & 0 & 0 & 0 & 0 \\
0 & 1 & 0 & 0 & 0 & 0 & 0 & 0 \\
u_{31} & u_{32} & 1 & 0 & 0 & 0 & 0 & 0 \\
u_{41} & u_{42} & 0 & 1 & 0 & 0 & 0 & 0 \\  
u_{51} & u_{52} & u_{53} & u_{54} & 1 & 0 & 0 & 0 \\ 
u_{61} & u_{62} & u_{63} & u_{64} & 0 & 1 & 0 & 0 \\
u_{71} & u_{72} & u_{73} & u_{74} & 0 & 0 & 1 & 0 \\
u_{81} & u_{82} & u_{83} & u_{84} & u_{85} & u_{86} & u_{87} & 1 
\end{pmatrix}.
\end{equation*}
Then $\dim U_{\be_{152}}=23$ as an algebraic group and 
$\dim W_{152}=10$.  

We choose the order of $u_{ij}$ and the 
coordinate vectors of $W_{152}$ in the following manner:
\begin{align*}
& u_{51},u_{52},u_{53},u_{54},u_{62},u_{63},
u_{8j} \; (j=1\ccd 7), \\
& u_{31},u_{41},u_{42},u_{64},u_{73},
u_{74},u_{32},u_{61},u_{71},u_{72} \\
& e_{358},e_{458},e_{468},e_{567},e_{478},
e_{378},e_{368},e_{568},e_{578},e_{678}
\end{align*}
($u_{31}\ccd u_{72}$ correspond to the coordinates of $W_{152}$). 
Explicitly, $n(u)R_{152}-R_{152}$ is 
\begin{align*}
& (u_{31}+u_{54}) e_{358} 
+ (u_{41}-u_{53}-u_{87}) e_{458} 
+ (u_{42}-u_{63}) e_{468} 
+ (u_{64}+u_{53}) e_{567} \\
& + (- u_{73}+u_{85}) e_{478} 
+ (u_{74}-u_{86}) e_{378} 
+ (u_{32}+ u_{64}+u_{87}) e_{368} \\
& + (- u_{61}+u_{52}+u_{53}u_{64}+u_{53}u_{87}
  -u_{54}u_{63}+u_{64}u_{87})e_{568} \\
& + (-u_{71}-u_{84}+u_{53}u_{74}-u_{53}u_{86}-u_{54}u_{73}
  +u_{54}u_{85}+u_{74}u_{87})e_{578} \\
& + (-u_{72}+u_{83}+u_{63}u_{74}-u_{63}u_{86}-u_{64}u_{73}
   +u_{64}u_{85}-u_{73}u_{87})e_{678}.  
\end{align*}
Condition \ref{cond:unipotent-connected} (2), (3) 
are satisfied and so $\rg\backslash S_{\be_{152}\, k}$ 
is also a single orbit.

\vskip 5pt   
\noindent      
(15) $\beta_{154}=\tfrac {1} {40} (-15,-7,1,1,1,1,9,9)$ 

In this case, $Z_{\beta_{154}}\cong \wedge^2 \aff^4\otimes \aff^2\oplus 1$
and 
\begin{equation*}
M_{\be_{154}} \cong \gl_1\times M_{[1,2,6]}
\cong \gl_1 \times \gl_1^2 \times \gl_4\times \gl_2.
\end{equation*}
We express elements of $M_{\be_{154}}$ as 
$g=(t_0,\diag(t_3,t_4,g_1,g_2)$ where 
$t_0,t_3,t_4\in\gl_1$, $g_1\in\gl_4$ and $g_2\in\gl_2$. 
Then 
\begin{equation*}
M^{\st}_{\be_{154}}
= \left\{
\diag(t_3,t_4,g_1,g_2) \;\;\vrule\;\;
\begin{matrix}
g_1\in\gl_4,g_2\in\gl_2,t_3,t_4\in\gl_1, \\
(\det g_1)(\det g_2)t_3t_4=1
\end{matrix}
\right\}.
\end{equation*}
Since 
\begin{equation*}
\chi_{\be_{154}}(g)=t_3^{-15}t_4^{-7}(\det g_1)(\det g_2)^9
=(\det g_1)^{16}(\det g_2)^{24}t_4^8
\end{equation*}
on $M^{\st}_{\beta_{154}}$, 
\begin{equation*}
G_{\st ,\beta_{154}} = 
\left\{
\diag(t_3,t_4,g_1,g_2) \;\;\vrule\;\;
\begin{matrix}
g_1\in\gl_4,g_2\in\gl_2,t_3,t_4\in\gl_1, \\
(\det g_1)(\det g_2)t_3t_4=1, \\
(\det g_1)^2(\det g_2)^3t_4=1
\end{matrix}
\right\}.
\end{equation*}

For $x\in Z_{154}$, let 
\begin{equation*}
A_1(x) = 
\begin{pmatrix}
0 & x_{347} & x_{357} & x_{367} \\
-x_{347} & 0 & x_{457} & x_{467} \\
-x_{357} & -x_{457} & 0 & x_{567} \\
-x_{367} & -x_{467} & -{567} & 0 
\end{pmatrix}, \;
A_2(x) = 
\begin{pmatrix}
0 & x_{348} & x_{358} & x_{368} \\
-x_{348} & 0 & x_{458} & x_{468} \\
-x_{358} & -x_{458} & 0 & x_{568} \\
-x_{368} & -x_{468} & -x_{568} & 0 
\end{pmatrix}
\end{equation*}
and $A(x)=(A_1(x),A_2(x))$.  
We identify $x$ with $(A(x),x_{278})$. 
The action of the above $g\in M_{\be_{154}}$ on $x\in Z_{154}$ is 
\begin{equation*}
\begin{pmatrix}
A_1(gx) \\
A_2(gx) 
\end{pmatrix}
= t_0g_2 \begin{pmatrix}
g_1 A_1(x) \,{}^t g_1 \\
g_1 A_2(x) \,{}^t g_1 
\end{pmatrix}
\end{equation*}
and $x_{278}$ is multiplied by $t_0t_4(\det g_2)$.  

The representation on the first factor of $Z_{154}$
was considered in \cite{wryu} and the 
assumption on $\ch(k)$ was removed in Subsection III--5.3. 
Let $P_1(x)$ be the homogeneous polynomial of degree $4$  
in Proposition III--5.22 and $P(x)=P_1(x)x_{278}$. Then 
\begin{equation*}
P(gx)=(\det g_1)^2(\det g_2)^2P_1(x) t_4 (\det g_2) x_{278}
=(\det g_1)^2(\det g_2)^3t_4P(x) 
\end{equation*}
on $M^{\st}_{\be_{154}}$.  Since $(\det g_1)^2(\det g_2)^3t_4$ 
is proportional to $\chi_{\be_{154}}(g)$ on $M^{\st}_{\be_{154}}$,  
$P(x)$ is an invariant polynomial with respect to $G_{\st,\be_{154}}$.  

Let 
\begin{equation*}
R_{154} = e_{278} + e_{347} + e_{568} 
\end{equation*}
Then $A(R_{154})$ corresponds to 
$w$ in (III--5.17). So 
$P_1(R_{154})=P(R_{154})=1$, which implies that 
$R_{154}\in Z_{154\, k}^{\sst}$.  
Therefore, $S_{\be_{154}}\not=\emptyset$.  
We would like to prove that 
the consideration of $M_{\be_{154}\, k}\backslash Z_{154\, k}^{\sst}$  
can be reduced to that of the action on the $A(x)$ part. 

Suppose that $x\in Z_{154\, k}^{\sst}$. Then $x_{278}\not=0$. 
By applying $(1,\diag(1,x_{278}^{-1},I_4,I_2))$, 
we may assume that $x_{278}=1$. Then 
$(t_0,\diag(t_3,t_4,g_1,g_2))$ does not change 
$x_{278}=1$ if and only if $t_0t_4(\det g_2)=1$.  
Then the action of 
$(1,\diag(1,(\det g_2)^{-1},g_1,g_2))$ 
on $A(x)$ is the same as the action of 
$(1,\diag(1,1,g_1,g_2))$. Therefore,  
the following proposition follows from  
Proposition III--5.23. 

\begin{prop}
\label{prop:r-orbit-case154}
$M_{\be_{154}\, k}\backslash Z_{154\, k}^{\sst}$  
is in bijective correspondence with $\Ex_2(k)$.  
\end{prop}

We verify Condition \ref{cond:unipotent-connected} (2), (3).  
In this case, $U_{\be_{108}}$ consists of all the elements 
of the form 
\begin{equation*}
n(u) = 
\begin{pmatrix}
1 & 0 & 0 & 0 & 0 & 0 & 0 & 0 \\
u_{21} & 1 & 0 & 0 & 0 & 0 & 0 & 0 \\
u_{31} & u_{32} & 1 & 0 & 0 & 0 & 0 & 0 \\
u_{41} & u_{42} & 0 & 1 & 0 & 0 & 0 & 0 \\  
u_{51} & u_{52} & 0 & 0 & 1 & 0 & 0 & 0 \\ 
u_{61} & u_{62} & 0 & 0 & 0 & 1 & 0 & 0 \\
u_{71} & u_{72} & u_{73} & u_{74} & u_{75} & u_{76} & 1 & 0 \\
u_{81} & u_{82} & u_{83} & u_{84} & u_{85} & u_{86} & 0 & 1 
\end{pmatrix}.
\end{equation*}
Then $\dim U_{\be_{154}}=21$ as an algebraic group and 
$\dim W_{108}=4$.

We choose the order of $u_{ij}$ and the 
coordinate vectors of $W_{108}$ in the following manner:
\begin{align*}
& u_{7j},u_{8j}\;(j=1\ccd 6),
u_{21},u_{31},u_{41},u_{51},u_{61}, \quad  
u_{32},u_{42},u_{52},u_{62} \\
& e_{378},e_{478},e_{578},e_{678}
\end{align*}
($u_{32},u_{42},u_{52},u_{62}$ correspond 
to the coordinates of $W_{154}$). 
Explicitly, $n(u)R_{154}-R_{154}$ is 
\begin{align*}
u_{32} e_{378} + u_{42} e_{478} + u_{52} e_{578} + u_{62} e_{678}.  
\end{align*}
Condition \ref{cond:unipotent-connected} (2), (3) 
are satisfied and so $\rg\backslash S_{\be_{154}\, k}$ 
is also in bijective correspondence with $\Ex_2(k)$.

\vskip 5pt   
\noindent 
(16) $\beta_{172}=\tfrac {1} {8} (-3,-1,-1,-1,1,1,1,3)$ 

In this case, $Z_{172}\cong \m_3\oplus 1$ and 
\begin{equation*}
M_{\be_{172}} = \gl_1\times M_{[1,4,7]}
\cong \gl_1\times \gl_1\times \gl_3^2\times \gl_1.
\end{equation*}
We express elements of $M_{\be_{172}}$ as 
$g=(t_0,\diag(t_3,g_1,g_2,t_4))$ where 
$t_0,t_3,t_4\in\gl_1,g_1$ and $g_2\in \gl_3$. Then 
\begin{equation*}
M^{\st}_{\be_{172}} = 
\left\{
\diag(t_3,g_1,g_2,t_4) \;\;\vrule\;\;
\begin{matrix}
t_0,t_3,t_4\in\gl_1,g_1,g_2\in \gl_3, \\
(\det g_1)(\det g_2)t_3t_4=1
\end{matrix}
\right\}.
\end{equation*}
Since 
\begin{equation*}
\chi_{\be_{172}}(g)=t_3^{-3}(\det g_1)^{-1} (\det g_2)t_4^3
= (\det g_1)^2 (\det g_2)^4t_4^6
\end{equation*}
on $M^{\st}_{\be_{172}}$,
\begin{equation*}
G_{\st,\be_{172}} = 
\left\{
\diag(t_3,g_1,g_2,t_4) \;\;\vrule\;\;
\begin{matrix}
t_0,t_3,t_4\in\gl_1,g_1,g_2\in \gl_3, \\
(\det g_1)(\det g_2)t_3t_4=1, \\
(\det g_1)(\det g_2)^2t_4^3=1
\end{matrix}
\right\}.
\end{equation*}

For $x\in Z_{172}$, let 
\begin{equation*}
A(x) = \begin{pmatrix}
x_{258} & x_{268} & x_{278} \\
x_{358} & x_{368} & x_{378} \\
x_{458} & x_{468} & x_{478} \\
\end{pmatrix}.
\end{equation*}
We identify $x$ with $(A(x),x_{567})$. 
The action of the above $g\in M_{\be_{172}}$ 
on $x\in Z_{172}$ is 
\begin{equation*}
A(gx) = t_0t_4 g_1 A(x) \, {}^t\!g_2
\end{equation*}
and $x_{567}$ is multiplied by $t_0\det g_2$. 
Let $P(x)= (\det A(x)) x_{567}$. 
Then 
\begin{equation*}
P(gx) = (\det g_1)(\det g_2)^2t_4^3P(x)
\end{equation*}
on $M^{\st}_{\be_{172}}$. Since 
$(\det g_1)(\det g_2)^2t_4^3$ is proportional to 
$\chi_{\be_{172}}(g)$ on $M^{\st}_{\be_{172}}$,  
$P(x)$ is an invariant polynomial 
with respect to $G_{\st,\be_{172}}$. 
This implies that $x\in Z_{172}^{\sst}$ 
if and only if $\det A(x),x_{567}\not=0$. 

Let 
\begin{equation*}
R_{172} = e_{258} + e_{368} + e_{478} + e_{567}. 
\end{equation*}
Since $A(R_{172})=I_3$ and the coefficient of $x_{567}$ is $1$, 
$P(R_{172})=1$. So $R_{172}\in Z_{172\, k}^{\sst}$. 
Therefore, $S_{\be_{172}}\not=\emptyset$.  

\begin{prop}
\label{prop:r-orbit-case172}
$Z_{172\, k}^{\sst} = M_{\be_{172}\, k}\cdot R_{172}$. 
\end{prop}
\begin{proof}
Suppose that $x\in Z_{172\, k}^{\sst}$. 
It easy to see that there exists $g\in M_{\be_{172}\, k}$ such that 
$A(gx)=I_3$. So we may assume that $A(x)=I_3$. 
Since $x_{567}\not=0$,  
\begin{equation*}
(x_{567}^{-1},\diag(1,I_3,I_3,x_{567})) x = R_{172}. 
\end{equation*}
\end{proof}

We verify Condition \ref{cond:unipotent-connected} (2), (3).  
In this case, $U_{\be_{172}}$ consists of all the elements 
of the form 
\begin{equation*}
n(u) = 
\begin{pmatrix}
1 & 0 & 0 & 0 & 0 & 0 & 0 & 0 \\
u_{21} & 1 & 0 & 0 & 0 & 0 & 0 & 0 \\
u_{31} & 0 & 1 & 0 & 0 & 0 & 0 & 0 \\
u_{41} & 0 & 0 & 1 & 0 & 0 & 0 & 0 \\  
u_{51} & u_{52} & u_{53} & u_{54} & 1 & 0 & 0 & 0 \\ 
u_{61} & u_{62} & u_{63} & u_{64} & 0 & 1 & 0 & 0 \\
u_{71} & u_{72} & u_{73} & u_{74} & 0 & 0 & 1 & 0 \\
u_{81} & u_{82} & u_{83} & u_{84} & u_{85} & u_{86} & u_{87} & 1 
\end{pmatrix}.
\end{equation*}
Then $\dim U_{\be_{172}}=22$ as an algebraic group and 
$\dim W_{172}=3$.  

We choose the order of $u_{ij}$ and the 
coordinate vectors of $W_{172}$ in the following manner:
\begin{align*}
& u_{i1} \; (i=2\ccd 8),u_{52},u_{62},u_{63},
u_{72},u_{73},u_{74},u_{8j} \; (j=1\ccd 7), \quad 
u_{53},u_{54},u_{64} \\
& e_{568},e_{578},e_{678} 
\end{align*}
($u_{53},u_{54},u_{64} $ correspond 
to the coordinates of $W_{172}$). 
Explicitly, $n(u)R_{172}-R_{172}$ is 
\begin{equation*}
(u_{53}-u_{62}+u_{87})e_{568} 
+ (u_{54}-u_{72}-u_{86}) e_{578} 
+ (u_{64}-u_{73}+u_{85}) e_{678}.
\end{equation*}
Condition \ref{cond:unipotent-connected} (2), (3) 
are satisfied and so $\rg\backslash S_{\be_{172}\, k}$ 
also consists of a single point.

\vskip 5pt   
\noindent 
(17) $\beta_{173}=\tfrac {1} {24} (-5,-5,-1,-1,-1,3,3,7)$ 

In this case, 
$Z_{\beta_{173}}\cong \aff^2\otimes \aff^2\op \wedge^2\aff^3\op \aff^3$ 
and 
\begin{equation*}
M_{\be_{173}} = \gl_1\times M_{[2,5,7]}
\cong \gl_1\times \gl_2\times \gl_3\times \gl_2\times \gl_1.
\end{equation*}
We express elements of $M_{\be_{173}}$ as 
$g=(t_0,\diag(g_2,g_1,g_3,t_4))$ where 
$t_0,t_4\in\gl_1$, $g_1\in\gl_3$ and 
$g_2,g_3\in \gl_2$. Then 
\begin{equation*}
M^{\st}_{\be_{173}} = \left\{
\diag(g_2,g_1,g_3,t_4) \;\;\vrule\;\;
\begin{matrix}
g_1\in\gl_3,g_2,g_3\in \gl_2, t_4\in\gl_1 \\
(\det g_1)(\det g_2)(\det g_3)t_4=1
\end{matrix}
\right\}. 
\end{equation*}
Since 
\begin{equation*}
\chi_{\be_{173}}(g) 
= (\det g_2)^{-5}(\det g_1)^{-1}(\det g_3)^3t_4^7
= (\det g_1)^{-8}(\det g_2)^{-12}(\det g_3)^{-4}
\end{equation*}
on $M^{\st}_{\be_{173}}$, 
\begin{equation*}
G_{\st,\be_{173}} = \left\{
\diag(g_2,g_1,g_3,t_4) \;\;\vrule\;\;
\begin{matrix}
g_1\in\gl_3,g_2,g_3\in \gl_2, \\
t_4 = (\det g_1)^{-1}(\det g_2)^{-1}(\det g_3)^{-1}, \\
(\det g_1)^{-2}(\det g_2)^{-3}(\det g_3)^{-1}=1
\end{matrix}
\right\}. 
\end{equation*}

Let $\{\bbmp_{3,1},\bbmp_{3,2},\bbmp_{3,3}\}$ 
be the standard basis of $\aff^3$, etc.  
For $x\in Z_{173}$, let
\begin{align*}
A(x) & = \begin{pmatrix}
x_{168} & x_{178} \\
x_{268} & x_{278} 
\end{pmatrix}, \;
v_1(x) = \begin{pmatrix}
x_{367} \\
x_{467} \\
x_{567} 
\end{pmatrix}, \\
v_2(x) & = 
x_{458} \,\bbmp_{3,2} \wedge \bbmp_{3,3}
- x_{358} \,\bbmp_{3,1} \wedge \bbmp_{3,3}
+ x_{348} \,\bbmp_{3,1} \wedge \bbmp_{3,2}.  
\end{align*}
We identify $x$ with $(A(x),v_1(x),v_2(x))$. 
The action of the above $g\in M_{\be_{173}}$ 
on $x\in Z_{173}$ is as follows:  
\begin{equation*}
A(gx) = t_0t_4 g_2 A(x) \,{}^t\! g_3, \;
v_1(gx) = t_0(\det g_3) g_1 v_1(x), \;
v_2(gx) = t_0t_4 \wedge^2 g_1 \, v_2(x). 
\end{equation*}

Note that $v_1(x)= x_{367}\bbmp_{3,1}+x_{467}\bbmp_{3,2}+x_{567}\bbmp_{3,3}$. 
We identify $\wedge^3 \aff^3\cong \aff^1$ so that 
$\bbmp_{3,1}\wedge \bbmp_{3,2}\wedge \bbmp_{3,2}$ corresponds to $1\in k$. 
Let 
\begin{equation*}
P_1(x)=\det A(x), \; 
P_2(x) = v_1(x)\wedge v_2(x). 
\end{equation*}
Then 
\begin{align*}
P_1(gx) & = t_4^2 (\det g_2)(\det g_3) P_1(x)
= (\det g_1)^{-2}(\det g_2)^{-1}(\det g_3)^{-1} P_1(x), \\
P_2(gx) & = t_4 (\det g_1)(\det g_3) P_2(x)
= (\det g_2)^{-1} P_2(x) 
\end{align*}
on $M^{\st}_{\be_{173}}$. 
Let $P(x)=P_1(x)P_2(x)^2$.  Then 
\begin{equation*}
P(gx) = (\det g_1)^{-2}(\det g_2)^{-3}(\det g_3)^{-1} P(x) 
\end{equation*}
on $M^{\st}_{\be_{173}}$. 
Since $(\det g_1)^{-2}(\det g_2)^{-3}(\det g_3)^{-1}$ 
is proportional to $\chi_{\be_{173}}(g)$ on $M^{\st}_{\be_{173}}$,  
$P(x)$ is an invariant polynomial with respect to
$G_{\st,\be_{173}}$.  This implies that 
$x\in Z_{173}^{\sst}$ if and only if $P_1(x),P_2(x)\not=0$. 

Let 
\begin{equation*}
R_{173} = e_{168} + e_{278} + e_{367} + e_{458}.
\end{equation*}
Then $A(R_{173})=I_2$, $v_1(R_{173})= [1,0,0]$ and 
$v_2(R_{173})= \bbmp_{3,2} \wedge \bbmp_{3,3}$. 
So $P_1(R_{173})=P_2(R_{173})=1$ and 
$R_{173}\in Z_{173\, k}^{\sst}$.
Therefore, $S_{\be_{173}}\not=\emptyset$.  

\begin{prop}
\label{prop:r-orbit-case173}
$Z_{173\, k}^{\sst} = M_{\be_{173}\, k}\cdot R_{173}$. 
\end{prop}
\begin{proof}
Suppose that $x\in Z_{173\, k}^{\sst}$. 
Since $\det A(x)\not=0$, there exists 
$g\in M_{\be_{173}\, k}$ such that 
$A(gx)=I_2$.   

Elements of the form $\diag(1,g_1,I_2,I_2)$ 
do not change the condition $A(x)=I_2$.  
Since $P_2(x)\not=0$, $v_1(x)\not=[0,0,0]$.  
So there exists $g_1\in \gl_3(k)$ such that 
$v_1(gx)=[1,0,0]$. Since $P_2(x)\not=0$, 
$x_{458}\not=0$. Let $g_1=\diag(1,x_{458}^{-1},1)$. 
By considering $\diag(1,g_1,I_2,I_2)x)$, we may assume that
$x_{458}=1$. Note that the condition $v_1(x)=[1,0,0]$ 
does not change.  

Let $g_1\in \gl_3(k)$ be the element such that
\begin{equation*}
{}^t g_1^{-1} = 
\begin{pmatrix}
1 & 0 & 0 \\
x_{358} & 1 & 0 \\
-x_{348} & 0 & 1
\end{pmatrix}.
\end{equation*}
Then $v_2(\diag(1,g_1,I_2,I_2)x)=\bbmp_{3,2} \wedge \bbmp_{3,3}$. 
Therefore, $x=R_{173}$. 
\end{proof}

We verify Condition \ref{cond:unipotent-connected} (2), (3).  
In this case, $U_{\be_{173}}$ consists of all the elements 
of the form 
\begin{equation*}
n(u) = 
\begin{pmatrix}
1 & 0 & 0 & 0 & 0 & 0 & 0 & 0 \\
0 & 1 & 0 & 0 & 0 & 0 & 0 & 0 \\
u_{31} & u_{32} & 1 & 0 & 0 & 0 & 0 & 0 \\
u_{41} & u_{42} & 0 & 1 & 0 & 0 & 0 & 0 \\  
u_{51} & u_{52} & 0 & 0 & 1 & 0 & 0 & 0 \\ 
u_{61} & u_{62} & u_{63} & u_{64} & u_{65} & 1 & 0 & 0 \\
u_{71} & u_{72} & u_{73} & u_{74} & u_{75} & 0 & 1 & 0 \\
u_{81} & u_{82} & u_{83} & u_{84} & u_{85} & u_{86} & u_{87} & 1 
\end{pmatrix}.
\end{equation*}
Then $\dim U_{\be_{173}}=23$ as an algebraic group and 
$\dim W_{173}=7$.  

We choose the order of $u_{ij}$ and the 
coordinate vectors of $W_{173}$ in the following manner:
\begin{align*}
& u_{61},u_{63},u_{64},u_{65},u_{7j}\; (j=1..5),u_{8j}\;(j=1..7),
\quad u_{31},u_{32},u_{41},u_{42},u_{51},u_{52},u_{62} \\
& e_{368},e_{378},e_{468},e_{478},e_{568},e_{578},e_{678} 
\end{align*}
($u_{31}\ccd u_{62} $ correspond 
to the coordinates of $W_{173}$). 
Explicitly, $n(u)R_{173}-R_{173}$ is 
\begin{align*}
& (u_{31} + u_{87}) e_{368}
+ (u_{32} - u_{86}) e_{378} 
+ (u_{41} + u_{65}) e_{468} \\
& + (u_{42} + u_{75}) e_{478} 
+ (u_{51} - u_{64}) e_{568} 
+ (u_{52} - u_{74}) e_{578} \\
& - (u_{62}-u_{71}+u_{83}-u_{63}u_{86}+u_{64}u_{75}
    -u_{65}u_{74}-u_{73}u_{87}) e_{678}. 
\end{align*}
Condition \ref{cond:unipotent-connected} (2), (3) 
are satisfied and so $\rg\backslash S_{\be_{173}\, k}$ 
also consists of a single point.

\vskip 5pt   
\noindent 
(18) $\be_{179} = \frac 1{24}(-9,-9,-1,-1,-1,7,7,7)$. 

In this case, $Z_{179}\cong \m_3$ and 
\begin{equation*}
M_{\be_{179}} = \gl_1\times M_{[2,5]} \cong \gl_1\times 
\gl_2\times \gl_3^2.
\end{equation*}
We express elements of $M_{\be_{179}}$ as 
$g=(t_0,\diag(g_3,g_1,g_2))$ where $t_0\in\gl_1$, 
$g_1,g_2\in\gl_3$ and $g_3\in\gl_2$. Then 
\begin{equation*}
M^{\st}_{\be_{179}} 
= \left\{
\diag(g_3,g_1,g_2) \mid g_1,g_2\in\gl_3,g_3\in\gl_2, \;
(\det g_1)(\det g_2)(\det g_3)=1
\right\}
\end{equation*}
Since 
\begin{equation*}
\chi_{\be_{179}}(g) = (\det g_3)^{-9}(\det g_1)^{-1}(\det g_2)^7
= (\det g_1)^8(\det g_2)^{16}
\end{equation*}
on $M^{\st}_{\be_{179}}$, 
\begin{equation*}
G_{\st,\be_{179}} = \left\{
\diag(g_3,g_1,g_2) \in M_{\be_{179}} \;\;\vrule\;\;
\begin{matrix}
(\det g_1)(\det g_2)(\det g_3)=1, \\
(\det g_1)(\det g_2)^2=1
\end{matrix}
\right\}.
\end{equation*}

For $x\in Z_{179}$, let 
\begin{equation*}
A(x) = \begin{pmatrix}
x_{378} & -x_{368} & x_{367} \\
x_{478} & -x_{468} & x_{467} \\
x_{578} & -x_{568} & x_{567}  
\end{pmatrix}.
\end{equation*}
We identify $x$ with $A(x)$. The action of 
the above $g\in M_{\be_{179}}$ on $x\in Z_{179}$ 
is as follows:
\begin{equation*}
A(gx) = t_0(\det g_2) g_1 A(x) g_2^{-1}. 
\end{equation*}
Let $P(x) = \det A(x)$. Then 
\begin{equation*}
P(gx)=(\det g_1)(\det g_2)^2P(x)
\end{equation*}
on $M^{\st}_{\be_{179}}$. 
Since $(\det g_1)(\det g_2)^2$ is proportional to 
$\chi_{\be_{179}}(g)$ on $M^{\st}_{\be_{179}}$, 
$P(x)$ is an invariant polynomial
with respect to $G_{\st,\be_{179}}$. 

Let 
\begin{equation*}
R_{179} = e_{378} - e_{468} +  e_{567}. 
\end{equation*}
Then $A(R_{179})=I_3$. Therefore, $P(R_{179})=1$, 
which implies that $R_{179}\in Z_{179\, k}^{\sst}$.
Therefore, $S_{\be_{179}}\not=\emptyset$.    

The following proposition is obvious. 
\begin{prop}
$Z^{\sst}_{179\, k}=M_{\be_{179}\, k} \cdot R_{179}$. 
\end{prop}

We verify Condition \ref{cond:unipotent-connected} (2), (3).  
In this case, $U_{\be_{179}}$ consists of all the elements 
of the form 
\begin{equation*}
n(u) = 
\begin{pmatrix}
1 & 0 & 0 & 0 & 0 & 0 & 0 & 0 \\
0 & 1 & 0 & 0 & 0 & 0 & 0 & 0 \\
u_{31} & u_{32} & 1 & 0 & 0 & 0 & 0 & 0 \\
u_{41} & u_{42} & 0 & 1 & 0 & 0 & 0 & 0 \\  
u_{51} & u_{52} & 0 & 0 & 1 & 0 & 0 & 0 \\ 
u_{61} & u_{62} & u_{63} & u_{64} & u_{65} & 1 & 0 & 0 \\
u_{71} & u_{72} & u_{73} & u_{74} & u_{75} & 0 & 1 & 0 \\
u_{81} & u_{82} & u_{83} & u_{84} & u_{85} & 0 & 0 & 1 
\end{pmatrix}.
\end{equation*}
Then $\dim U_{\be_{179}}=21$ as an algebraic group and 
$\dim W_{179}=1$.  

Since 
\begin{equation*}
n(u)R_{179}-R_{179} = (u_{63} + u_{74} + u_{85}) e_{678},  
\end{equation*}
Condition \ref{cond:unipotent-connected} (2), (3) 
are satisfied and so $\rg\backslash S_{\be_{179}\, k}$ 
also consists of a single point.

\vskip 5pt   
\noindent  
(19) $\beta_{180}=\tfrac {1} {40} (-3,-3,-3,-3,1,1,1,9)$ 

In this case, $Z_{180}\cong \wedge^2 \aff^4 \oplus 1$ and 
\begin{equation*}
M_{\be_{180}} = \gl_1\times M_{[4,7]}\cong 
\gl_1\times \gl_4\times \gl_3\times \gl_1.
\end{equation*}
We express elements of $M_{\be_{180}}$ as 
$g=(t_0,\diag(g_1,g_2,t_3))$ where $t_0,t_3\in\gl_1$, 
$g_1\in\gl_4$ and $g_2\in\gl_3$.  Then 
\begin{equation*}
M^{\st}_{\be_{180}} = \left\{
\diag(g_1,g_2,t_3) 
\mid g_1\in\gl_4,g_2\in\gl_3, \; (\det g_1)(\det g_2)=1
\right\}. 
\end{equation*}
Since 
\begin{equation*}
\chi_{\be_{180}}(g) = 
(\det g_1)^{-3}(\det g_2)t_4^9
= (\det g_1)^{-12}(\det g_2)^{-8}
\end{equation*}
on $M^{\st}_{\be_{180}}$, 
\begin{equation*}
G_{\st,\be_{180}} = 
\left\{
\diag(g_1,g_2,(\det g_1)^{-1}(\det g_2)^{-1}) 
\in M^{\st}_{\be_{180}} 
\mid (\det g_1)^{-3}(\det g_2)^{-2}=1
\right\}.  
\end{equation*}

For $x\in Z_{\be_{180}}$, let 
\begin{align*}
A(x)=\begin{pmatrix} 0 & x_{128} & x_{138} & x_{148} \\
-x_{128} & 0 & x_{238} & x_{248} \\
-x_{138} & -x_{238} & 0 & x_{348} \\
-x_{148} & -x_{248} & -x_{348} & 0
\end{pmatrix}.
\end{align*}
We identify $x\in Z_{\be_{180}}$ with $(A(x),x_{567})$. 
The action of the above $g\in M_{\be_{180}}$ on $x\in Z_{180}$ is 
as follows:  
\begin{equation*}
A(gx)=t_0t_3 g_1A(x)\,{}^t\!g_1
\end{equation*}
and $x_{567}$ is multiplied by $t_0(\det g_2)$. 
Let $P_1(x)$ be the Pfaffian of $A(x)$ and 
$P(x)=P_1(x)^3x_{567}^4$. Then 
\begin{align*}
P_1(gx) & = t_3^2 (\det g_1) P_1(x)
= (\det g_1)^{-1}(\det g_2)^{-2}P_1(x), \\ 
P(x) & = (\det g_1)^{-3}(\det g_2)^{-2}P(x)
\end{align*}
on $M^{\st}_{\be_{180}}$. Since 
$(\det g_1)^{-3}(\det g_2)^{-2}$ is proportional 
to $\chi_{\be_{180}}(g)$ on $M^{\st}_{\be_{180}}$,  
$P(x)$ is an invariant polynomial with respect to 
$G_{\st,\be_{180}}$.  This implies that 
$x\in Z_{180}^{\sst}$ if and only if 
$\det A(x),x_{567}\not=0$. 

Let 
\begin{equation*}
R_{180} = e_{128} + x_{348} +  e_{567}. 
\end{equation*}
Then
\begin{equation*}
A(R_{180}) = 
\begin{pmatrix}
J & 0 \\ 0 & J
\end{pmatrix}
\end{equation*}
and the $x_{567}$-coordinate of $R_{180}$ is $1$. 
So $P(R_{180})=1$, which implies that 
$R_{180}\in Z_{180\, k}^{\sst}$. 
Therefore, $S_{\be_{180}}\not=\emptyset$.   

\begin{prop}
$Z^{\sst}_{180\, k}=M_{\be_{180}\, k} R_{180}$. 
\end{prop}
\begin{proof}
By Witt's theorem (Theorem \ref{thm:alternating-matrix}), 
there exists $g\in M_{\be_{180}\, k}$ such that 
$A(gx)=A(R_{180})$. So we may assume that 
$A(x)=A(R_{180})$. Since $x_{567}\not=0$, 
\begin{equation*}
(x_{567}^{-1},\diag(I_4,I_3,x_{567}) x = R_{180}. 
\end{equation*}
\end{proof}

We verify Condition \ref{cond:unipotent-connected} (2), (3).  
In this case, $U_{\be_{180}}$ consists of all the elements 
of the form 
\begin{equation*}
n(u) = 
\begin{pmatrix}
1 & 0 & 0 & 0 & 0 & 0 & 0 & 0 \\
0 & 1 & 0 & 0 & 0 & 0 & 0 & 0 \\
0 & 0 & 1 & 0 & 0 & 0 & 0 & 0 \\
0 & 0 & 0 & 1 & 0 & 0 & 0 & 0 \\  
u_{51} & u_{52} & u_{53} & u_{54} & 1 & 0 & 0 & 0 \\ 
u_{61} & u_{62} & u_{63} & u_{64} & 0 & 1 & 0 & 0 \\
u_{71} & u_{72} & u_{73} & u_{74} & 0 & 0 & 1 & 0 \\
u_{81} & u_{82} & u_{83} & u_{84} & u_{85} & u_{86} & u_{87} & 1 
\end{pmatrix}.
\end{equation*}
Then $\dim U_{\be_{180}}=19$ as an algebraic group and 
$\dim W_{180}=15$.  

We choose the order of $u_{ij}$ and the 
coordinate vectors of $W_{180}$ in the following manner:
\begin{align*}
& u_{81},u_{82},u_{83},u_{84}, \quad 
u_{51},u_{52},u_{53},u_{54},u_{61},u_{62},u_{63},u_{64},
u_{71},u_{72},u_{73},u_{74},u_{85},u_{86},u_{87} \\
& e_{258},e_{158},e_{458},e_{358},e_{268},e_{468},e_{368},
e_{278},e_{178},e_{478},e_{378},e_{678},e_{578},e_{568}
\end{align*}
($u_{51}\ccd u_{87}$ correspond to the coordinates of $W_{180}$). 
Explicitly, $n(u)R_{180}-R_{180}$ is 
\begin{align*}
& - u_{51} e_{258}
+ u_{52} e_{158} 
- u_{53} e_{458}
+ u_{54} e_{358}
- u_{61} e_{268}
+ u_{62} e_{168} \\
& - u_{63} e_{468} + u_{64} e_{368}
- u_{71} e_{278}
+ u_{72} e_{178}
- u_{73} e_{478}
+ u_{74} e_{378} \\
& + (u_{85} + u_{61}u_{72} - u_{62}u_{71} 
+ u_{63}u_{74} - u_{64}u_{73}) e_{678} \\
& + (-u_{86} + u_{51}u_{72} - u_{52}u_{71} 
+ u_{53}u_{74} - u_{54}u_{73}) e_{578} \\
& + (u_{87} + u_{51}u_{62} - u_{52}u_{61} 
+ u_{53}u_{64} - u_{54}u_{63}) e_{568}
\end{align*}
Condition \ref{cond:unipotent-connected} (2), (3) 
and so $\rg\backslash S_{\be_{180}\, k}$ 
also consists of a single point.

\vskip 5pt   
\noindent 
(20) $\beta_{182}=\tfrac {1} {8} (-3,-3,-3,1,1,1,1,5)$ 

In this case $Z_{182}\cong \wedge^2 \aff^4$ and 
\begin{equation*}
M_{\be_{182}} = \gl_1\times M_{[3,7]} = \gl_1\times M_{[3,7]}
\cong \gl_1\times \gl_3\times \gl_4\times \gl_1.
\end{equation*}
We express elements of 
as $g=(t_0,\diag(g_2,g_1,t_3))$ where $t_0,t_3\in\gl_1$, 
$g_1\in\gl_4$ and $g_2\in\gl_3$.  Then 
\begin{equation*}
M^{\st}_{\be_{182}} = 
\left\{
\diag(g_2,g_1,t_3)
\mid g_1\in\gl_4,g_2\in\gl_3, \;
(\det g_1)(\det g_2)=1 
\right\}.
\end{equation*}
Since 
\begin{equation*}
\chi_{\be_{182}}(g) = (\det g_2)^{-3}(\det g_1)t_3^5
= (\det g_1)^{-4} (\det g_2)^{-8}
\end{equation*}
on $M^{\st}_{\be_{182}}$, 
\begin{equation*}
G_{\st,\be_{182}}
= \left\{
\diag(g_2,g_1,(\det g_1)^{-1}(\det g_2)^{-1}) 
\in M^{\st}_{\be_{182}}, \; (\det g_1)^{-1}(\det g_2)^{-2}=1
\right\}.  
\end{equation*}

For $x\in Z_{182}$, let 
\begin{align*}
A(x)=\begin{pmatrix} 0 & x_{458} & x_{468} & x_{478} \\
-x_{458} & 0 & x_{568} & x_{578} \\
-x_{468} & -x_{568} & 0 & x_{678} \\
-x_{478} & -x_{578} & -x_{678} & 0
\end{pmatrix}.
\end{align*}
We identify $x\in Z_{\be_{182}}$ with $A(x)$. 
The action of the above $g\in M_{\be_{182}}$ 
on $x\in Z_{182}$ is as follows:  
\begin{equation*}
A(gx)=t_0t_3g_1A(x)\,{}^t\!g_1. 
\end{equation*}
Let $P(x)$ be the Pfaffian of $A(x)$. Then, 
\begin{equation*}
P(gx)= t_3^2\det g_1 P(x)= (\det g_1)^{-1}(\det g_2)^{-2} P(x)
\end{equation*}
on $M^{\st}_{\be_{182}}$. Since $(\det g_1)^{-1}(\det g_2)^{-2}$ 
is proportional to $\chi_{\be_{182}}(g)$ on $M^{\st}_{\be_{182}}$, 
$P(x)$ is an invariant polynomial with respect to $G_{\st,\be_{182}}$. 

Let 
\begin{equation*}
R_{182} = e_{458}+e_{678}. 
\end{equation*}
Then 
\begin{equation*}
A(R_{182})=
\begin{pmatrix}
J & 0 \\ 0 & J 
\end{pmatrix}. 
\end{equation*}
So $P(R_{182})=1$, which implies that 
$R_{182}\in Z_{182\, k}^{\sst}$. 
Therefore, $S_{\be_{182}}\not=\emptyset$.   

The following proposition 
follows from  Lemma II--4.6. 
\begin{prop}
$Z_{\be_{182}\, k}^{\sst}=M_{\be_{182}\, k} \cdot R_{182}$. 
\end{prop}
Since $W_{\beta_{182}}=\{0\}$, 
$\rg\backslash S_{\be_{182}\, k}$ 
also consists of a single point. 

\vskip 5pt   
\noindent 
(21) $\beta_{183}=\tfrac {1} {8} (-3,-3,-3,-3,-3,5,5,5)$ 

In this case, $Z_{183} \cong 1$ (the trivial \rep). So 
by the construction of $Z_{\be_{183}}$, 
$G_{\st,\be_{183}}$ acts trivially on $Z_{183}$.  
Therefore, $P(x)=x_{678}$ is an invariant polynomial
with respect to $G_{\st,\be_{183}}$. This implies that 
\begin{equation*}
Z_{183}^{\sst}=\{x_{678}e_{678}\in Z_{183}\mid x_{678}\not=0\}
=M_{\beta_{183}\, k}\cdot R_{183}
\end{equation*}
where $R_{183}=e_{678}$. Since  Since $W_{\beta_{183}}=\{0\}$, 
$\rg\backslash S_{\beta_{183}\, k}$ also consists of a single point.

\section{Empty strata}
\label{sec:empty-strata}

\newcommand{\tallstrut}{\rule[-12pt]{-0.2cm}{28pt}}
\newcommand{\tinsert}{\quad\tallstrut}

In this section we prove that $S_{\be_i}=\emptyset$ 
for $i$ not in (\ref{eq:list-non-empty-tri8}) for 
the \pv{} in this paper. 

We used some lemmas from Part II in Part III. 
We make slight improvement for those lemmas 
to begin with as follows.  

For some lemmas, the only improvement is 
to include the contragriedient \rep{} 
in the statement.  

The following lemma is an analogue of Lemma II--4.1. 

\begin{lem}
\label{lem:eliminate-standard}
Let $G=\spl_n$ and $V=\aff^n$ (the standard \rep) and 
$V^*$ its contragriedient \rep{} of $V$ with dual basis. 
Then for any $x\in V_k$ or $x\in V^*_k$, 
there exists $g\in G_k$ such that 
there is at most one non-zero entry for $gx$.   
\end{lem}
\begin{proof}
Note that $g\in G$ acts on $V^*\cong \aff^n$ 
by $x\mapsto {}^t\! g^{-1}x$.  If the \rep{}
is the standard \rep{} then there is $g\in G_k$ 
such that $gx$ is in the form $[0\ccd 0,*]$. 
Since there are elements in $G_k$ which 
exchanges coordinates up to sign, 
the lemma follows.
Put $h={}^t\! g^{-1}$. Then ${}^t\!h^{-1}=g$. 
Let $\rho,\rho^*$ be the standard \rep{} 
and its contragriedient \rep{} respectively. 
Then $\rho^*(h)(x)=\rho(g)(x)$.  
\end{proof}

We improve Lemma II--4.2 as follows.

\begin{lem}
\label{lem:eliminate-nm-matrix}
Let $n\geqq m>0$, $G=\spl_n$ and $V=\m_{n,m}$.  
Then for any $A\in V_k$ there exists $g,h\in G_k$ such that 
if $B=(b_{ij})=gA,C=(c_{ij})={}^t\!h^{-1}A$ 
then $b_{ij}=c_{ij}=0$ for 
\begin{equation*}
i=1\ccd n-m,j=1\ccd m, \; 
i=m+1,j=1\ccd m-1\ccd 
i=n-1,j=1.
\end{equation*}
It is possible to change the columns for which 
entries are eliminated.  
\end{lem}
\begin{proof}
Since the action of $G$ on the first column of $A$ is 
the standard \rep{} or the contragriedient \rep, 
we can make the first column into the form 
$[0\ccd 0,*]$ by Lemma \ref{lem:eliminate-standard}. 
Then the subgroup $H$ of $G$ consisting of elements
of the form 
\begin{equation*}
\begin{pmatrix}
h &  \\
 & 1
\end{pmatrix}
\end{equation*}
where $h\in \gl_{n-1}$ acts on $A$ without changing the condition
that the first column is in the form 
$[0\ccd 0,*]$. Then the lemma follows by induction on $n$. 
\end{proof}

For example, if $n=6,m=4$ then $B,C$ are in the following form 
\begin{equation*}
\begin{pmatrix}
0 & 0 & 0 & 0 \\
0 & 0 & 0 & 0 \\
0 & 0 & 0 & * \\
0 & 0 & * & * \\
0 & * & * & * \\
* & * & * & * \\
\end{pmatrix}. 
\end{equation*}
It is possible to change the columns as follows also.
\begin{equation*}
\begin{pmatrix}
0 & 0 & 0 & * \\
0 & * & * & * \\
* & * & * & * \\
0 & 0 & * & * \\
0 & 0 & 0 & 0 \\
0 & 0 & 0 & 0 
\end{pmatrix}. 
\end{equation*}

The following lemma is an improvement of 
Lemma II--4.3.

 \begin{lem}
\label{lem:eliminate-2m(2)}
Let $G=\spl_2\times \spl_2$ and $V=\aff^2\oplus \m_2$
where $g=(g_1,g_2)\in G$, $g_1$ acts on $\aff^2$ by 
the standard \rep{} or its contragriedient \rep, 
and $g$ acts on $\m_2$ by the multiplication 
of $g_1$ or ${}^t\! g_1^{-1}$ from the left 
and the multiplication of 
${}^t\!g_2$ or $g_2^{-1}$ from the right. 
Then for any element $(v,A)\in V_k$, 
$(i,j)=(1,1),(1,2),(2,1)$ or $(2,2)$ and $k=1$ or $2$,  
there exists $g\in G_k$ such that 
if $g(v,A)=(w,B)$ then the $k$-th entry of $w$ and the $(i,j)$-entry of $B$ 
are $0$.  
 \end{lem}
\begin{proof}
Note that there are eight possibilities of the action 
of $G$ on $V$. What the statement says is that we can make any 
one entry of $v$ and any one entry of $A$ zero. 

We can eliminate either entry of $v$ by Lemma \ref{lem:eliminate-standard}. 
Then elements of the form $(I_2,g_2)$ do not change $v$.  
Since the action of $g_2$ on rows of $A$ can be identified 
with the standard \rep{} or its contragriedient \rep{} 
(which is isomorphic to the standard \rep{} since the group is $\spl_2$), 
we can eliminate any one entry of $A$. 
\end{proof}

The following lemma is an improvement of 
Lemma II--4.4.

\begin{lem}
\label{lem:eliminate-2m(32)}
Let $G=\spl_3\times \spl_2$ and $V=\aff^2\oplus \m_{3,2}$
where $(g_1,g_2)\in G$ acts on $V$ by $V\ni (v,A)\mapsto (g_2v,g_1A{}^tg_2)$.
Then for any element $(v,A)\in V_k$, there exists $g\in G_k$ such that 
if $g(v,A)=(w,B)$ then either entry of $w$, the first row of $B$ 
and the $(2,1)$-entry are $0$. The actions of $g_1,g_2$ can be changed to 
the contragriedient \rep s.  
\end{lem}
\begin{proof}
We can eliminate either entry of $v$ 
by applying Lemma \ref{lem:eliminate-standard} to $g_2$. 
Then we can apply Lemma \ref{lem:eliminate-nm-matrix} to $g_1$
\end{proof}

The following lemma is an improvement of 
Lemma II--4.5 and the proof is similar. 

\begin{lem}
\label{lem:eliminate-2m(3-32)}
Let $G=\spl_3\times \spl_2$ and $V=\aff^3\oplus \m_{3,2}$
where $(g_1,g_2)\in G$ acts on $V$ by $V\ni (v,A)\mapsto (g_1v,g_1A{}^t\!g_2)$.
Then for any element $(v,A)\in V$, there exists $g\in G$ such that 
if $g(v,A)=(w,B)$ then any two entries of $w$ and any 
one entry of $B$ are $0$. The actions of $g_1,g_2$ can be changed to 
the contragriedient \rep s.   
\end{lem}
\begin{proof}
We can eliminate any two entries of $w$ by the action of 
$\spl_3$.  Then $\spl_2$ does not change $w$ and acts on any 
row of $B$. Therefore, we can make any one entry of $B$ $0$ 
by Lemma \ref{lem:eliminate-standard}.  
\end{proof}
\begin{lem}
\label{lem:empty-criterion}
Let $\be\in\gB$. Suppose that $x\in Z_{\be}$, 
$M_{\be}x\sub Z_{\be}$ is Zariski open and that $x$ is 
unstable with respect to the action of {\upshape $G_{\text{st},\be}$}. 
Then $Z^{\sst}_{\be}=\emptyset$. 
\end{lem}
\begin{proof}
Let $T_0$ be as in (\ref{eq:T0-defn}).
Then $M_{\be}=M^{\st}_{\be}\cdot T_0$. 
Let $T_1\sub G$ be the subgroup  generated by 
$T_0$ and $\{\lam_{\be}(t)\mid t\in \gl_1\}$ 
(see Section I--2). 
Then $M_{\be} = G_{\text{st},\be}\cdot T_1$ 
and $T_1$ acts on $Z_{\be}$ by scalar multiplication.

Suppose that $Z^{\sst}_{\be}\not=\emptyset$. 
Since $Z^{\sst}_{\be},M_{\be}x\sub Z_{\be}$ are Zariski open, 
there exists $g\in G_{\text{st},\be}$, $t\in T_1$ 
such that $gtx \in Z^{\sst}_{\be}$. Since 
$Z^{\sst}_{\be}$ is $G_{\text{st},\be}$-invariant, 
$tx\in Z^{\sst}_{\be}$. Since $tx$ is a scalar 
multiple of $x$, $x\in Z^{\sst}_{\be}$, 
which is a contradiction. 
\end{proof}

\begin{rem}
\label{remark:not-all-zero}
When we apply the above lemmas, 
we do not necessarily have to make 
all the claimed entries zero.  
For example, when we apply Lemma \ref{lem:eliminate-nm-matrix}, 
we can just make the first ($n-m$) rows zero and
use $\spl_{m}$ for other \rep s to eliminate other 
entries.  
\end{rem}

In the following table, for each $i$, we list 
which coordinates of $x$ we can eliminate and 
the 1PS with the property that weights of non-zero
coordinates are all positive. 

\vskip 10pt

\begin{center}


\end{center}

We change the format of the table to save space. 

\begin{center}
%
\\
\hline
\end{tabular}

\end{center}

\vskip 10pt

We assume that 
$x=\sum_j y_j \mathbbm a_j\in Z_{\be_i}$ 
unless otherwise stated.  The reason why we 
we use this coordinate instead of $x_{ijk}$ is that 
it was convenient that way to check the computation by 
MAPLE. 

We verify the informations (which coordinates we can eliminate) 
in the above table in the following order. 

\vskip 10pt
\begin{itemize}
\item[(a)] 
\setlength{\itemsep=2pt}
The easiest cases are those where 
each factor of $M_{\be_i}$ 
except for $\gl_1$ have at most one standard \rep{} 
or its contragriedient \rep{} in 
$Z_{\be_i}$ and it is enough to apply Lemma \ref{lem:eliminate-standard}
to all such standard \rep s to show that $S_{\be_i}=\emptyset$.  
\item[(b)] 
The next easy cases are those where we can apply 
Lemma \ref{lem:eliminate-standard}, 
but either we only use some of the standard \rep s 
or the contragriedient \rep s or 
we have to choose which standard \rep{} of a particular 
component of $M_{\be_i}$. Note that there are cases 
where a $\gl_2$ component of $M_{\be_i}$ 
has two standard \rep s in $Z_{\be_i}$ and 
eliminating a coordinate for only one \rep{} 
works. 
\item[(c)] 
There are cases where it is enough to apply 
Lemma \ref{lem:eliminate-2m(2)} once. 
For that purpose there has to be a factor 
$\gl_2\times \gl_2$ in $M_{\be_i}$ 
and a factor $\m_2\oplus \aff^2$ in $Z_{\be_i}$. 
If this pair is unique then we do not have to 
worry about to which factor we apply Lemma \ref{lem:eliminate-2m(2)}. 
\item[(d)] 
There are cases where it is enough to apply 
Lemma \ref{lem:eliminate-nm-matrix} once.  
For that purpose there has to be a factor 
$\gl_n$ in $M_{\be_i}$
and a factor $\m_{n,m}$ with $n>m$ in $Z_i$. 
\item[(e)] 
The remaining cases require individual considerations. 
\end{itemize}

\vskip 10pt
\noindent
{\bf Case (a)} 
We list the vector $\be$ of the applicable cases 
as follows.  

The numbering is according to the order in the table. 

\vskip 10pt

$\beta_{15}=\tfrac {9} {632} (-13,-5,-5,-5,3,3,3,19)$ 

$\beta_{16}=\tfrac {11} {312} (-7,-7,-7,1,1,1,9,9)$ 

$\beta_{17}=\tfrac {3} {1208} (-23,-15,1,1,1,9,9,17)$ 

$\beta_{19}=\tfrac {5} {568} (-9,-9,-9,-1,-1,7,7,15)$ 

$\beta_{26}=\tfrac {3} {440} (-7,-7,-7,1,1,1,1,17)$ 

$\beta_{31}=\tfrac {11} {952} (-15,-15,-7,1,1,9,9,17)$ 

$\beta_{39}=\tfrac {7} {888} (-19,-11,-3,-3,5,5,13,13)$ 

$\beta_{40}=\tfrac {1} {40} (-3,-3,-1,-1,1,1,3,3)$ 

$\beta_{41}=\tfrac {1} {1912} (-21,-13,-13,-5,3,3,19,27)$ 
 
$\beta_{43}=\tfrac {1} {216} (-17,-9,-9,-1,-1,7,7,23)$ 

$\beta_{46}=\tfrac {5} {24} (-1,-1,-1,-1,-1,-1,3,3)$ 

$\beta_{48}=\tfrac {7} {88} (-3,-3,-3,-3,1,1,5,5)$ 
 
$\beta_{50}=\tfrac {1} {16} (-1,-1,-1,-1,-1,1,1,3)$ 

$\beta_{59}=\tfrac {1} {24} (-5,-5,-5,-1,-1,-1,3,15)$ 
 
$\beta_{60}=\tfrac {1} {184} (-45,-45,-45,-5,11,11,51,67)$  

$\beta_{64}=\tfrac {1} {232} (-47,-15,-15,-15,-3,17,17,61)$  

$\beta_{67}=\tfrac {1} {120} (-29,-5,-5,-5,-5,3,19,27)$ 

$\beta_{69}=\tfrac {3} {152} (-7,-7,1,1,1,1,5,5)$ 

$\beta_{79}=\tfrac {1} {56} (-11,-11,-9,-9,-9,15,17,17)$  

$\beta_{80}=\tfrac {1} {104} (-19,-15,-15,3,3,7,7,29)$

$\beta_{86}=\tfrac {1} {104} (-39,3,3,5,5,7,7,9)$

$\beta_{89}=\tfrac {1} {56} (-13,-13,-9,-1,-1,7,11,19)$ 
 
$\beta_{90}=\tfrac {1} {376} (-69,-69,-29,11,11,27,51,67)$  
 
$\beta_{91}=\tfrac {1} {112} (-6,-5,-3,2,2,3,3,4)$ 

$\beta_{95}=\tfrac {1} {4} (-1,-1,0,0,0,0,1,1)$ 
  
$\beta_{96}=\tfrac {1} {88} (-21,-9,-5,-1,-1,3,15,19)$ 
 
$\beta_{97}=\tfrac {1} {296} (-39,-15,-7,1,1,9,17,33)$  
 
$\beta_{99}=\tfrac {1} {112} (-5,-3,-3,-1,-1,3,3,7)$ 
 
$\beta_{100}=\tfrac {1} {56} (-9,-9,-1,-1,3,3,7,7)$ 
 
$\beta_{101}=\tfrac {1} {88} (-13,-9,-1,-1,3,3,7,11)$ 

$\beta_{103}=\tfrac {1} {312} (-13,-5,-1,-1,3,3,7,7)$

$\beta_{105}=\tfrac {1} {168} (-23,-15,-7,-7,1,9,17,25)$ 
 
$\beta_{110}=\tfrac {1} {40} (-15,-7,-7,-7,-7,-7,25,25)$  
 
$\beta_{111}=\tfrac {1} {104} (-23,-23,-23,-23,1,1,25,65)$  
 
$\beta_{114}=\tfrac {1} {136} (-27,-27,-27,-19,-19,13,21,85)$  
 
$\beta_{115}=\tfrac {1} {136} (-51,-27,-27,-27,13,13,53,53)$ 

$\beta_{116}=\tfrac {1} {56} (-13,-13,-13,-5,-5,11,19,19)$
  
$\beta_{117}=\tfrac {1} {88} (-33,-9,-1,-1,7,7,7,23)$ 
 
$\beta_{118}=\tfrac {1} {168} (-15,-15,-7,-7,-7,9,9,33)$ 

$\beta_{119}=\tfrac {1} {136} (-51,-19,-19,-3,-3,-3,13,85)$ 

$\beta_{120}=\tfrac {1} {296} (-71,-15,-15,-15,1,1,41,73)$ 

$\beta_{122}=\tfrac {1} {104} (-39,1,1,1,1,9,9,17)$

$\beta_{124}=\tfrac {1} {232} (-31,-31,-7,-7,-7,25,25,33)$ 
 
$\beta_{136}=\tfrac {1} {88} (-33,-17,-17,-1,7,7,23,31)$ 
 
$\beta_{137}=\tfrac {1} {328} (-59,-51,-51,-3,21,29,29,85)$ 

$\beta_{138}=\tfrac {1} {104} (-23,-23,-15,-15,9,9,17,41)$
 
$\beta_{139}=\tfrac {1} {152} (-33,-25,-9,-9,7,7,15,47)$  
 
$\beta_{140}=\tfrac {1} {88} (-33,-33,-1,7,7,15,15,23)$  

$\beta_{141}=\tfrac {1} {328} (-123,-3,13,13,21,21,29,29)$ 
 
$\beta_{142}=\tfrac {1} {328} (-59,-51,-51,-3,-3,53,53,61)$ 

$\beta_{143}=\tfrac {1} {104} (-39,-15,-15,1,1,17,17,33)$

$\beta_{145}=\tfrac {1} {136} (-15,-15,-3,1,1,5,5,21)$

$\beta_{146}=\tfrac {1} {120} (-21,-21,-13,-5,11,11,19,19)$

$\beta_{147}=\tfrac {1} {136} (-15,-7,-7,-3,-3,1,13,21)$ 
 
$\beta_{148}=\tfrac {1} {584} (-35,-27,-11,-11,13,21,21,29)$  

$\beta_{155}=\tfrac {1} {8} (-3,-3,-1,-1,-1,-1,5,5)$  

$\beta_{156}=\tfrac {1} {40} (-15,-7,-7,-7,1,1,9,25)$  

$\beta_{157}=\tfrac {1} {16} (-3,-3,-3,1,1,1,3,3)$ 

$\beta_{158}=\tfrac {3} {8} (-1,-1,-1,-1,1,1,1,1)$ 

$\beta_{159}=\tfrac {1} {16} (-6,0,0,0,1,1,1,3)$ 

$\beta_{164}=\tfrac {1} {40} (-7,-7,-3,-3,-3,-3,1,25)$

$\beta_{165}=\tfrac {1} {8} (-3,-3,-1,-1,1,1,3,3)$ 

$\beta_{166}=\tfrac {1} {40} (-15,-7,-7,-3,-3,9,13,13)$, 

$\beta_{167}=\tfrac {1} {40} (-15,-15,-3,1,5,5,9,13)$  

$\beta_{169}=\tfrac {1} {104} (-23,-15,-15,1,1,9,17,25)$ 

$\beta_{170}=\tfrac {1} {136} (-15,-3,-3,1,1,1,9,9)$ 

$\beta_{171}=\tfrac {1} {104} (-11,-7,-7,-3,-3,5,5,21)$ 
 
$\beta_{174}=\tfrac {1} {24} (-9,-9,-9,-1,-1,-1,15,15)$  

$\beta_{175}=\tfrac {1} {24} (-9,-9,-9,-9,7,7,7,15)$ 

$\beta_{176}=\tfrac {1} {56} (-21,-21,-5,-5,3,3,11,35)$  

$\beta_{177}=\tfrac {1} {56} (-21,-21,-21,3,11,11,19,19)$ 

$\beta_{178}=\tfrac {1} {88} (-33,-9,-9,7,7,7,15,15)$  

$\beta_{181}=\tfrac {1} {8} (-3,-3,-3,-3,1,1,5,5)$

\vskip 10pt
\noindent
{\bf Case (b)}
We list the vector $\be$ of the applicable cases 
as follows.  We also list the subspaces to which we apply 
Lemma \ref{lem:eliminate-standard} 
to show that $S_{\be_i}=\emptyset$. 

\vskip 10pt

$\beta_{18}=\tfrac {1} {184} (-9,-5,-1,-1,3,3,3,7)$
\begin{math}
\langle \mathbbm a_{18},\mathbbm a_{20},\mathbbm a_{21}\rangle,
\langle \mathbbm a_{26},\mathbbm a_{30}\rangle 
\end{math}

$\beta_{42}=\tfrac {7} {1528} (-27,-11,-3,-3,5,5,13,21)$, 
$\langle \mathbbm a_{32}, \mathbbm a_{34}\rangle$

$\beta_{44}=\tfrac {1} {1272} (-21,-13,-5,-5,3,11,11,19)$, 
$\langle \mathbbm a_{26},\mathbbm a_{30}\rangle$

$\beta_{45}=\tfrac {1} {616} (-23,-15,-7,-7,1,9,17,25)$, 
$\langle \mathbbm a_{26},\mathbbm a_{30}\rangle$

$\beta_{65}=\tfrac {3} {152} (-7,-3,-3,1,1,1,1,9)$ 
$\langle \mathbbm a_{47},\mathbbm a_{48},\mathbbm a_{50}\rangle$

$\beta_{68}=\tfrac {1} {280} (-25,-17,-9,-9,-9,-1,31,39)$ 
$\lan \mathbbm a_{39},\mathbbm a_{42},\mathbbm a_{48}\ran$

$\beta_{75}=\tfrac {1} {216} (-5,-5,-5,-1,-1,3,3,11)$ 
$\langle \mathbbm a_{6},\mathbbm a_{11},\mathbbm a_{26}\rangle$, 
$\langle \mathbbm a_{47},\mathbbm a_{48}\rangle$

$\beta_{92}=\tfrac {1} {44} (-3,-3,-2,-2,-1,0,5,6)$, 
$\langle \mathbbm a_{42},\mathbbm a_{48}\rangle $

$\beta_{93}=\tfrac {1} {88} (-17,-9,-9,-9,3,11,11,19)$, 
$\langle \mathbbm a_{20},\mathbbm a_{21}\rangle$

$\beta_{94}=\tfrac {1} {22} (-3,-1,-1,0,0,1,1,3)$, 
$\langle \mathbbm a_{34},\mathbbm a_{44}\rangle$

$\beta_{102}=\tfrac {1} {88} (-9,-6,-6,-2,1,5,5,12)$, 
$\langle \mathbbm a_{47},\mathbbm a_{48}\rangle$

$\beta_{107}=\tfrac {1} {248} (-29,-17,-5,-1,-1,11,15,27)$, 
$\langle \mathbbm a_{39},\mathbbm a_{42}\rangle$

$\beta_{109}=\tfrac {1} {72} (-5,-3,-1,-1,1,1,3,5)$,
$\langle \mathbbm a_{18},\mathbbm a_{20}\rangle$

$\beta_{123}=\tfrac {1} {16} (-2,-1,-1,0,0,0,1,3)$,  
$\langle \mathbbm a_{48},\mathbbm a_{50},\mathbbm a_{53}\rangle$

$\beta_{144}=\tfrac {1} {152} (-57,-9,-1,7,7,15,15,23)$, 
$\langle \mathbbm a_{40},\mathbbm a_{43}\rangle$

$\beta_{150}=\tfrac {1} {152} (-25,-17,-9,-1,-1,15,15,23)$, 
$\langle \mathbbm a_{40},\mathbbm a_{43}\rangle$

$\beta_{151}=\tfrac {1} {456} (-43,-27,-27,-19,-3,5,45,69)$, 
$\langle \mathbbm a_{34},\mathbbm a_{44}\rangle$

$\beta_{153}=\tfrac {1} {456} (-59,-43,-3,5,5,13,13,69)$, 
$\langle \mathbbm a_{47},\mathbbm a_{48}\rangle$

$\beta_{168}=\tfrac {1} {8} (-3,-1,0,0,0,1,1,2)$, 
$\langle \mathbbm a_{35}, \mathbbm a_{36}\rangle$

\vskip10pt
\noindent
{\bf Case (c)} 
We list the vector $\be$ of the applicable cases 
with the coordinates we can eliminate 
to show that $S_{\be_i}=\emptyset$ as follows.  

\vskip 10pt

$\beta_{34}=\tfrac {3} {424} (-21,-13,-5,-5,3,11,11,19)$, 
$y_{20},y_{47}=0$

$\beta_{35}=\tfrac {1} {312} (-9,-5,-5,-1,-1,3,7,11)$, 
$y_{15},y_{29}=0$

$\beta_{82}=\tfrac {1} {152} (-37,-9,-5,-5,-1,-1,27,31)$, 
$y_{33},y_{42}=0$

$\beta_{83}=\tfrac {1} {152} (-25,-25,-17,-9,11,19,19,27)$, 
$y_{36},y_{48}=0$

$\beta_{84}=\tfrac {1} {12} (-2,-1,-1,0,0,1,1,2)$, 
$y_{40},y_{44}=0$

$\beta_{87}=\tfrac {1} {24} (-4,-3,-1,0,0,2,3,3)$, 
$y_{34},y_{39}=0$

$\beta_{104}=\tfrac {1} {88} (-5,-3,-3,-1,1,3,3,5)$, 
$y_{40},y_{42}=0$

$\beta_{106}=\tfrac {3} {248} (-7,-7,-3,-3,1,5,5,9)$ 
$y_{19},y_{48}=0$

$\beta_{129}=\tfrac {1} {264} (-59,-35,-35,-3,13,13,37,69)$, 
$y_{33},y_{48}=0$

$\beta_{132}=\tfrac {1} {232} (-15,-15,1,1,1,9,9,9)$, 
$y_{19},y_{34}=0$

$\beta_{133}=\tfrac {1} {264} (-99,-3,5,5,13,13,29,37)$, 
$y_{35},y_{50}=0$

$\beta_{134}=\tfrac {1} {264} (-59,-35,-35,-3,-3,29,53,53)$, 
$y_{35},y_{49}=0$

$\beta_{149}=\tfrac {1} {328} (-35,-11,-3,-3,5,5,21,21)$, 
$y_{25},y_{41}=0$

$\beta_{161}=\tfrac {1} {24} (-9,-9,-9,3,3,7,7,7)$, 
$y_{50},y_{53}=0$

$\beta_{162}=\tfrac {1} {56} (-21,-5,-5,3,3,7,7,11)$, 
$y_{35},y_{47}=0$

\vskip 10pt
\noindent
{\bf Case (d)}
We list the vector $\be$ of the applicable cases 
If necessary, we write the subspace to which 
Lemma \ref{lem:eliminate-nm-matrix} is applied.  

\vskip 10pt

$\beta_{6}=\tfrac {3} {56} (-3,-3,-3,1,1,1,1,5)$ 

$\beta_{8}=\tfrac {9} {248} (-5,-5,-5,-5,3,3,3,11)$ 

$\beta_{10}=\tfrac {7} {120} (-3,-3,-3,-3,-3,5,5,5)$ 

$\beta_{14}=\tfrac {1} {120} (-5,-5,-1,-1,-1,3,3,7)$, 
$\langle \mathbbm a_{11}\ccd \mathbbm a_{33}\rangle $

$\beta_{28}=\tfrac {3} {248} (-11,-3,-3,1,1,1,5,9)$ 

$\beta_{30}=\tfrac {9} {632} (-13,-13,-5,3,3,3,11,11)$
 
$\beta_{47}=\tfrac {3} {40} (-5,-1,-1,1,1,1,1,3)$ 

$\beta_{52}=\tfrac {1} {24} (-3,-3,-3,-3,-1,-1,-1,15)$ 

$\beta_{53}=\tfrac {1} {40} (-7,-7,-7,-7,-7,5,5,25)$ 

$\beta_{54}=\tfrac {1} {40} (-9,-5,-5,1,1,1,5,11)$ 

$\beta_{58}=\tfrac {1} {20} (-2,-1,-1,-1,0,0,2,3)$ 

$\beta_{70}=\tfrac {1} {24} (-9,-3,-3,-3,-3,7,7,7)$ 

$\beta_{71}=\tfrac {1} {40} (-5,-4,0,0,1,1,1,6)$ 

$\beta_{72}=\tfrac {1} {136} (-19,-19,-3,-3,-3,13,17,17)$ 

$\beta_{76}=\tfrac {1} {16} (-3,-3,-1,-1,1,1,1,5)$ 

$\beta_{77}=\tfrac {1} {88} (-17,-17,-17,-13,7,11,11,35)$ 

$\beta_{78}=\tfrac {1} {16} (-6,-1,0,1,1,1,2,2)$ 

$\beta_{85}=\tfrac {1} {56} (-21,-11,-11,3,3,3,17,17)$ 

$\beta_{112}=\tfrac {1} {136} (-51,-11,-11,5,13,13,13,29)$ 

$\beta_{113}=\tfrac {1} {40} (-7,-7,-7,1,1,5,5,9)$ 

$\beta_{126}=\tfrac {1} {136} (-19,-19,-19,-11,-11,-3,-3,85)$ 

$\beta_{127}=\tfrac {3} {136} (-17,-17,-1,-1,7,7,7,15)$ 

$\beta_{128}=\tfrac {1} {136} (-51,-19,-19,-19,-11,37,37,45)$ 

$\beta_{131}=\tfrac {1} {136} (-51,-51,-11,13,13,13,37,37)$

$\beta_{135}=\tfrac {1} {88} (-17,-17,-17,-17,-9,23,23,31)$ 

$\beta_{160}=\tfrac {1} {24} (-9,-9,-1,-1,-1,3,3,15)$ 

$\beta_{163}=\tfrac {1} {8} (-3,-1,-1,-1,-1,1,1,5)$

\vskip 10pt
\noindent
{\bf Case (e)}
We consider the remaining cases one by one. 
We sometimes have to mention the semi-simple parts 
of subgroups of $M_{\be}$ for the cases below.  
We use the notation such as $\spl_m^{(i,j)}$. 
This means that it is a group which is isomorphic to 
$\spl_m$ from the $(i,i)$-coordinate to the 
$(j,j)$-coordinate in $\gl_8$. In terms of 
the notion such as $M_{[3,5]}$, 
$\spl_m^{(i,j)}$ is the semi-simple part 
of $M_{[1,2\ccd i-1,j,j+1\ccd 7]}$. 
For example, 
\begin{equation*}
\spl_3^{(3,5)}=M_{[1,2,5,6,7]}
\end{equation*}
consists of elements of the form 
$\diag(I_2,g,I_3)$ where $g\in \spl_3$.  

\vskip 5pt
\noindent
(1) $\be_2=\tfrac {1} {376} (-13,-5,-5,3,3,3,3,11)$ 

In this case, 
$Z_2\cong \aff^4\oplus \wedge^2 \aff^4\otimes \aff^2$ 
and 
\begin{equation*}
M_{\be_2} = \gl_1\times M_{[1,3,7]}
\cong \gl_1\times \gl_4\times \gl_2\times \gl_1.
\end{equation*}
We express elements of $M_{\be_2}$ as 
$g=(t_0,\diag(g_2,g_1,t_3))$ 
where $t_0,t_3\in\gl_1$, $g_1\in\gl_4,g_2\in\gl_2$. 

For $x\in Z_2$, let 
\begin{equation*}
v(x) = [y_{15},y_{18},y_{20},y_{21}]
= [x_{148},x_{158},x_{168},x_{178}]
\end{equation*}
and $A(x)=(A_1(x),A_2(x))$ where 
\begin{align*}
A_1(x) & = \begin{pmatrix}
0 & y_{27} & y_{28} & y_{29} \\
-y_{27} & 0 & y_{31} & y_{32} \\
-y_{28} & -y_{31} & 0 & y_{34} \\
-y_{29} & -y_{32} & -y_{34} & 0 
\end{pmatrix} 
= \begin{pmatrix}
0 & x_{245} & x_{246} & x_{247} \\
-x_{245} & 0 & x_{256} & x_{257} \\
-x_{246} & -x_{256} & 0 & x_{267} \\
-x_{247} & -x_{257} & -x_{267} & 0
\end{pmatrix}, 
\end{align*}
\begin{align*}
A_2(x) & = \begin{pmatrix}
0 & y_{37} & y_{38} & y_{39} \\
-y_{37} & 0 & y_{41} & y_{42} \\
-y_{38} & -y_{41} & 0 & y_{44} \\
-y_{39} & -y_{42} & -y_{44} & 0 
\end{pmatrix}
= \begin{pmatrix}
0 & x_{345} & x_{346} & x_{347} \\
-x_{345} & 0 & x_{356} & x_{357} \\
-x_{346} & -x_{356} & 0 & x_{367} \\
-x_{347} & -x_{357} & -x_{367} & 0
\end{pmatrix}.
\end{align*}
We identify $x$ with $(v(x),A(x))$. 

Let 
\begin{equation*}
R_2 = e_{158}+e_{178}+x_{267}+x_{345}. 
\end{equation*}
Then 
\begin{equation}
\label{eq:vR2etc}
v(R_2) = \begin{pmatrix}
0 \\ 1 \\ 0 \\ 1 \end{pmatrix}, \;
A_1(R_2) = \begin{pmatrix}
0 & 0 & 0 & 0 \\
0 & 0 & 0 & 0 \\
0 & 0 & 0 & 1 \\
0 & 0 & -1 & 0 
\end{pmatrix}, \;
A_2(R_2) = \begin{pmatrix}
0 & 1 & 0 & 0 \\
-1 & 0 & 0 & 0 \\
0 & 0 & 0 & 0 \\
0 & 0 & 0 & 0 
\end{pmatrix}. 
\end{equation}

We show that $M_{\be_2} \cdot R_2\sub Z_2$ is Zariski open. 
Then the following coordinates of $R_2$ 
\begin{equation*}
y_{15},y_{20},y_{27},y_{28},y_{29},y_{31},y_{32},
y_{38},y_{39},y_{41},y_{42},y_{44}
\end{equation*}
are $0$. Since all the weights of non-zero coordinates 
of $R_2$ are positive with respect to the 
1PS in the table, $S_{\be_2}=\emptyset$ 
by Lemmas \ref{lem:empty-criterion}.

Let $H\sub M_{\be_2}$ be the following subgroup
\begin{equation*}
H = \{g=\diag(g_2,g_1,1)\mid g_1\in\gl_4,g_2\in\gl_2\}
\cong \gl_4\times \gl_2\sub M_{\be_2}. 
\end{equation*}
Let $g=\diag(g_2,g_1,1)\in H$. 
The action of $g$ on $Z_2$ is 
\begin{equation*}
v(gx) = g_1 v(x), \; 
\begin{pmatrix}
A_1(gx) \\
A_2(gx) 
\end{pmatrix}
= g_2 
\begin{pmatrix}
g_1 A_1(x)\,{}^t\!g_1 \\
g_1 A_2(x)\,{}^t\!g_1
\end{pmatrix}. 
\end{equation*}
If $H\cdot R_2\sub Z_2$ is Zariski open, 
$M_{\be_2}\cdot R_2\sub Z_2$ is Zariski open. 
Therefore, it is enough to show the following proposition. 

\begin{prop}
\label{prop:4241}
$H\cdot R_2\sub Z_2$ is Zariski open.  
\end{prop}
\begin{proof}
Let $H_{R_2}$ be the stabilizer of $R_2$ as a group scheme and 
$\mathrm{T}_e(H_{R_2})$ its tangent space at the unit element 
$e=(I_4,I_2)$. Since $\dim H = 20$ and $\dim Z_2=16$, 
it is enough to prove that $\dim \mathrm{T}_e(H_{R_2})=4$ by 
Proposition \ref{prop:dimensions-compatible}. 

Let $v=v(R_2)$ and $X_1=A_1(R_{2}),X_2=A_2(R_2)$. 
We use these notations since we would like to use the letter
$A$ to describe elements of the tangent space. 

Let $k[\ep]/(\ep^2)$ be the ring of dual numbers. 
Elements of $\mathrm{T}_e(H)$ can be expressed as 
$(I_4+\ep A,I_2+\ep B)$ where 
\begin{equation*}
A=(a_{ij})\in \m_4(k),\quad  B = (b_{ij})\in \m_2(k). 
\end{equation*}
Then $(I_4+\ep A,I_2+\ep B)\in \mathrm{T}_e(H_{R_2})$ if and only if 
the following conditions  
\begin{align*}
& (I+\ep A)v  = v, \\
& (1+\ep b_{11})[(I+\ep A)X_1(1+\ep {}^t A)] 
+ \ep b_{12}[(I+\ep A)X_2(1+\ep {}^t A)] = X_1, \\
& \ep b_{21}[(I+\ep A)X_1(1+\ep {}^t A)] 
+(1+\ep b_{22})[(I+\ep A)X_2(1+\ep {}^t A)] = X_2
\end{align*}
are satisfied. 

These are equivalent to the following conditions
\begin{align*}
& a_{i2}=a_{i4}=0 \; (i=1,2,3,4), \\
& b_{11}X_1 + AX_1 + X_1{}^t A + b_{12}X_2 = 0, \\
& b_{21}X_1 + AX_2 + X_2{}^t A + b_{22}X_2 = 0.
\end{align*}

Explicitly, 
\begin{align*}
& b_{11}\begin{pmatrix}
0 & 0 & 0 & 0 \\
0 & 0 & 0 & 0 \\
0 & 0 & 0 & 1 \\
0 & 0 & -1 & 0 
\end{pmatrix} + 
\begin{pmatrix}
a_{11} & 0 & a_{13} & 0 \\
a_{21} & 0 & a_{23} & 0 \\
a_{31} & 0 & a_{33} & 0 \\
a_{41} & 0 & a_{43} & 0
\end{pmatrix} 
\begin{pmatrix}
0 & 0 & 0 & 0 \\
0 & 0 & 0 & 0 \\
0 & 0 & 0 & 1 \\
0 & 0 & -1 & 0 
\end{pmatrix} \\
& + \begin{pmatrix}
0 & 0 & 0 & 0 \\
0 & 0 & 0 & 0 \\
0 & 0 & 0 & 1 \\
0 & 0 & -1 & 0 
\end{pmatrix}
\begin{pmatrix}
a_{11} & a_{21} & a_{31} & a_{41} \\
0      &    0   &   0    &  0     \\
a_{13} & a_{23} & a_{33} & a_{43} \\
0      &    0   &   0    &  0
\end{pmatrix}
+ b_{12}
\begin{pmatrix}
0 & 1 & 0 & 0 \\
-1 & 0 & 0 & 0 \\
0 & 0 & 0 & 0 \\
0 & 0 & 0 & 0 
\end{pmatrix} = 0, \\
& b_{21} 
\begin{pmatrix}
0 & 0 & 0 & 0 \\
0 & 0 & 0 & 0 \\
0 & 0 & 0 & 1 \\
0 & 0 & -1 & 0 
\end{pmatrix}
+ \begin{pmatrix}
a_{11} & 0 & a_{13} & 0 \\
a_{21} & 0 & a_{23} & 0 \\
a_{31} & 0 & a_{33} & 0 \\
a_{41} & 0 & a_{43} & 0
\end{pmatrix} 
\begin{pmatrix}
0 & 1 & 0 & 0 \\
-1 & 0 & 0 & 0 \\
0 & 0 & 0 & 0 \\
0 & 0 & 0 & 0 
\end{pmatrix} \\
& + \begin{pmatrix}
0 & 1 & 0 & 0 \\
-1 & 0 & 0 & 0 \\
0 & 0 & 0 & 0 \\
0 & 0 & 0 & 0 
\end{pmatrix}
\begin{pmatrix}
a_{11} & a_{21} & a_{31} & a_{41} \\
0      &    0   &   0    &  0     \\
a_{13} & a_{23} & a_{33} & a_{43} \\
0      &    0   &   0    &  0
\end{pmatrix}
+ b_{22}
\begin{pmatrix}
0 & 1 & 0 & 0 \\
-1 & 0 & 0 & 0 \\
0 & 0 & 0 & 0 \\
0 & 0 & 0 & 0 
\end{pmatrix} = 0 
\end{align*}
%

%

By long but straightforward computations, 
\begin{equation*}
A = 
\begin{pmatrix}
a_{11} & 0 & 0 & 0 \\
a_{21} & 0 & 0 & 0 \\
0 & 0 & a_{33} & 0 \\
0 & 0 & a_{43} & 0
\end{pmatrix}, \;
\begin{pmatrix}
b_{11} & 0 \\
0 & b_{22}
\end{pmatrix}, \;
b_{11}+a_{33}=0,\; a_{11}+b_{22}=0.
\end{equation*}
Therefore, $\dim \mathrm{T}_e(H_{R_2})=4$. 
\end{proof}

\vskip 5pt
\noindent
(2) $\beta_{3}=\tfrac {3} {184} (-7,-7,1,1,1,1,1,9)$ 

In this case, $Z_3=\m_{5,2} \oplus \wedge^2 (\aff^5)^*$  
($(\aff^5)^*$ is the dual space) and 
\begin{equation*}
M_{\be_3} = \gl_1\times M_{[2,7]}
\cong \gl_1\times \gl_5\times \gl_2.
\end{equation*}
We express elements of $M_{\be_2}$ as 
$g=(t_0,\diag(g_2,g_1,t_3))$ 
where $t_0,t_3\in\gl_1$, $g_1\in\gl_5,g_2\in\gl_2$. 

For $x\in Z_3$, let 
\begin{equation*}
A(x) = 
\begin{pmatrix}
y_{11} & y_{26} \\ 
y_{15} & y_{30} \\ 
y_{18} & y_{33} \\ 
y_{20} & y_{35} \\ 
y_{21} & y_{36} 
\end{pmatrix}
= \begin{pmatrix}
x_{138} & x_{238} \\
x_{148} & x_{248} \\
x_{158} & x_{258} \\
x_{168} & x_{268} \\
x_{178} & x_{278} 
\end{pmatrix}
\end{equation*}
and 
\begin{align*}
B(x) & = \begin{pmatrix}
0 & y_{53} & -y_{50} & y_{48} & -y_{47} \\
& 0 & y_{44} & -y_{42} & y_{41} \\
&& 0 & y_{39} & -y_{38} \\
&&& 0 & y_{37} \\
&&&& 0
\end{pmatrix} 
= \begin{pmatrix}
0 & x_{567} & -x_{467} & x_{457} & -x_{456} \\
& 0 & x_{367} & -x_{357} & x_{356} \\
&& 0 & x_{347} & -x_{346}  \\
&&& 0 & x_{345} \\
&&&& 0  
\end{pmatrix}
\end{align*}
($B(x)$ is an alternating matrix and 
the lower triangular part is suppressed). 
Elements of $Z_{3}$ can be identified with $(A(x),B(x))$. 
The action of $g=\diag(g_2,g_1,1)$ 
($g_1\in \spl_5,g_2\in\spl_2$) on $x$ is 
\begin{equation*}
A(gx) = g_1 A(x)\, {}^t\!g_2, \;
B(gx) = {}^t\!g_1^{-1} B(x) g_1^{-1}. 
\end{equation*}

We first apply Witt's theorem (Theorem \ref{thm:alternating-matrix}) 
and can make 
\begin{equation*}
y_{38},y_{39},y_{41},y_{42},y_{47},y_{48},y_{50},y_{53}=0
\end{equation*}
(note that $y_{37},y_{44}\not=0$ if the rank of $B(x)$ is two).  
Then $\spl_2^{(1,2)}$ does not change 
these conditions. This group acts on the subspace 
$\lan \mathbbm a_{11},\mathbbm a_{26}\ran$ and so 
we can make $\mathbbm a_{11}=0$.  
Then $\spl_2^{(4,5)}\times \spl_2^{(6,7)}$ does not change 
these conditions. This group acts on the subspaces 
$\lan \mathbbm a_{15}, \mathbbm a_{18}\ran$,  
$\lan \mathbbm a_{20}, \mathbbm a_{21}\ran$ and
so we can make $y_{15},y_{20}=0$.  
Therefore, we may assume that
\begin{equation*}
y_{11},y_{15},y_{20},y_{38},y_{39},y_{41},
y_{42},y_{47},y_{48},y_{50},y_{53}=0.
\end{equation*}
The remaining coordinates are 
\begin{equation*}
y_{18},y_{21},y_{26},y_{30},y_{33},y_{35},
y_{36},y_{37},y_{44}. 
\end{equation*}
Then all the weights of these coordinates 
are positive with respect to the 
1PS in the table and $S_{\be_3}=\emptyset$ 

\vskip 5pt
\noindent
(3) $\beta_{5}=\tfrac {1} {248} (-5,-5,-5,-5,3,3,3,11)$ 

In this case, 
$Z_5=\wedge^2 \aff^4\oplus \wedge^2 \aff^3\otimes \aff^4$ 
and 
\begin{equation*}
M_{\be_5} = \gl_1\times M_{[4,7]}
\cong \gl_1\times \gl_4\times \gl_3\times \gl_1.
\end{equation*}
We express elements of $M_{\be_5}$ as 
$g=(t_0,\diag(g_1,g_2,t_3))$ 
where $t_0,t_3\in\gl_1$, $g_1\in\gl_4,g_3\in\gl_3$. 

For $x\in Z_5$, let 
\begin{align*}
& A_1(x) 
= \begin{pmatrix}
y_{19} \\
-y_{17} \\
y_{16}
\end{pmatrix}
= \begin{pmatrix}
x_{167} \\
-x_{157} \\
x_{156}
\end{pmatrix}, \;
A_2(x) = \begin{pmatrix}
y_{34} \\
-y_{32} \\
y_{31}
\end{pmatrix}, \;
= \begin{pmatrix}
x_{267} \\
-x_{257} \\
x_{256}
\end{pmatrix}, \\
& A_3(x) = \begin{pmatrix}
y_{44} \\
-y_{42} \\
y_{41}
\end{pmatrix} 
= \begin{pmatrix}
x_{367} \\
-x_{357} \\
x_{356}
\end{pmatrix}, \;
A_4(x) = \begin{pmatrix}
y_{50} \\
-y_{48} \\
y_{47}
\end{pmatrix} 
= \begin{pmatrix}
x_{467} \\
-x_{457} \\
x_{456}
\end{pmatrix}, \\
& A(x) = (A_1(x)\; A_2(x)\; A_3(x)\; A_4(x)), \\
& B(x) = \begin{pmatrix}
0 & y_6 & y_{11} & y_{15} \\
-y_6 & 0 & y_{26} & y_{30} \\
-y_{11} & -y_{26} & 0 & y_{40} \\
-y_{15} & -y_{30} & -y_{40} & 0
\end{pmatrix} 
= \begin{pmatrix}
0 & x_{128} & x_{138} & x_{148} \\
-x_{128} & 0 & x_{238} & x_{248} \\
-x_{138} & -x_{238} & 0 & x_{348} \\
-x_{148} & -x_[248] & -x_{348} & 0
\end{pmatrix}.  
\end{align*}
The action of $g=\diag(g_1,g_2,1)$ 
($g_1\in\spl_4,g_2\in\spl_3$) is 
\begin{equation*}
\begin{pmatrix}
A_1(gx) \\
\vdots \\
A_4(gx)
\end{pmatrix}
= g_1 \begin{pmatrix}
{}^t\!g_2^{-1} A_1(x) \\
\vdots \\
{}^t\!g_2^{-1} A_4(x)
\end{pmatrix} \, {}^t\!g_1, \;
B(gx) = g_1 B(x) \, {}^t\!g_1
\end{equation*}

By Lemma \ref{lem:eliminate-nm-matrix}, 
we may assume that $y_{16},y_{17},y_{19}=0$. 
Then $\spl_3^{(2,4)}\times \spl_3^{(5,7)}\sub M_{\be_5}$
does not change this condition. 
Since $\lan \mathbbm a_6,\mathbbm a_{11},\mathbbm a_{15}\ran$
is the standard \rep{} of $\spl_3$, 
we can make $y_6,y_{11}=0$. 
Then $\spl_2^{(2,3)}$
does not change these conditions.  
So we may assume that $y_{30}=0$ since 
$\lan \mathbbm a_{30},\mathbbm a_{40}\ran$ is the standard 
\rep{} of $\spl_2$.  Then 
$\spl_2^{(5,7)}$ does not change these conditions. 
So we may assume that
$y_{31},y_{32},y_{41}=0$ by Lemma \ref{lem:eliminate-nm-matrix}
($n=m=3$).  Therefore, we may assume that 
\begin{equation*}
y_6,y_{11},y_{16},y_{17},y_{19},y_{30},
y_{31},y_{32},y_{41}=0.   
\end{equation*}
Then the 1PS in the table shows that $S_{\be_5}=\emptyset$.

\vskip 10pt
\noindent
(4) $\beta_{11}=\tfrac {5} {184} (-9,-1,-1,-1,-1,-1,7,7)$ 

In this case, 
$Z_{11} \cong 1\oplus \wedge^2 \aff^5\otimes \aff^2$ 
and 
\begin{equation*}
M_{\be_{11}} = \gl_1\times M_{[1,6]}
\cong \gl_1\times \gl_1\times \gl_5\times \times \gl_2.  
\end{equation*}
We express elements of $M_{\be_{11}}\cap \gl_8$ as 
$g=\diag(t_3,g_1,g_2)$ where $g_1\in\gl_5,g_2\in\gl_2$. 

For $x\in Z_{11}$, let
\begin{align*}
A_1(x) = 
\begin{pmatrix}
0 & y_{25} & y_{29} & y_{32} & y_{34} \\
-y_{25} & 0 & y_{39} & y_{42} & y_{44} \\
-y_{29} & -y_{39} & 0 & y_{48} & y_{50} \\
-y_{32} & -y_{42} & -y_{48} & 0 & y_{53} \\
-y_{34} & -y_{44} & -y_{50} & -y_{53} & 0 
\end{pmatrix}
= \begin{pmatrix}
0 & x_{237} & x_{247} & x_{257} & x_{267} \\
-x_{237} & 0 & x_{347} & x_{357} & x_{367} \\
-x_{247} & -x_{347} & 0 & x_{457} & x_{467} \\
-x_{257} & -x_{357} & -x_{457} & 0 & x_{567} \\
-x_{267} & -x_{367} & x_{467} & -x_{567} & 0 
\end{pmatrix}, \\
A_2(x) = 
\begin{pmatrix}
0 & y_{26} & y_{30} & y_{33} & y_{35} \\
-y_{26} & 0 & y_{40} & y_{43} & y_{45} \\
-y_{30} & -y_{40} & 0 & y_{49} & y_{51} \\
-y_{33} & -y_{43} & -y_{49} & 0 & y_{54} \\
-y_{35} & -y_{45} & -y_{51} & -y_{54} & 0 
\end{pmatrix}
= \begin{pmatrix}
0 & x_{238} & x_{248} & x_{258} & x_{268} \\
-x_{238} & 0 & x_{348} & x_{358} & x_{368} \\
-x_{248} & -x_{348} & 0 & x_{458} & x_{468} \\
-x_{258} & -x_{358} & -x_{458} & 0 & x_{568} \\
-x_{268} & -x_{368} & x_{468} & -x_{568} & 0 
\end{pmatrix}
\end{align*}
and $A(x)=(A_1(x),A_2(x))$. We identify $x$ with 
$(y_{21}=x_{178},A(x))$. The action of 
$g=diag(1,g_1,g_2)$ ($g_1\in\spl_5,g_2\in\spl_2$) is  
\begin{equation*}
\begin{pmatrix}
A_1(gx) \\
A_2(gx) 
\end{pmatrix}
= g_2 \begin{pmatrix}
g_1 A_1(x) \,{}^t\!g_1 \\
g_1 A_2(x) \,{}^t\!g_1 
\end{pmatrix}
\end{equation*}
and $x_{178}$ is multiplied by $t_3(\det g_2)$.

Let 
\begin{equation}
\label{eq:X1X2-defn}
\begin{aligned}
R_{11} & = e_{178} + e_{347} + e_{567} + e_{238} + e_{458}, \\
X_1 & = \begin{pmatrix}
0 & 0 & 0 & 0 & 0 \\
0 & 0 & 1 & 0 & 0 \\
0 & -1 & 0 & 0 & 0 \\
0 & 0 & 0 & 0 & 1 \\
0 & 0 & 0 & -1 & 0
\end{pmatrix}, \quad 
X_2 = \begin{pmatrix}
0 & 1 & 0 & 0 & 0 \\
-1 & 0 & 0 & 0 & 0 \\
0 & 0 & 0 & 1 & 0 \\
0 & 0 & -1 & 0 & 0 \\
0 & 0 & 0 & 0 & 0
\end{pmatrix}.
\end{aligned}
\end{equation}
Then $A(R_{11})=(X_1,X_2)$ and the coefficient of 
$e_{178}$ is $1$.  

\begin{prop}
\label{prop:522}
$M_{\be_{11}}\cdot R_{11} \sub Z_{11}$ is Zariski open.  
\end{prop}
\begin{proof}
It is enough to show that 
$(M_{\be_{11}}\cap \gl_8)\cdot R_{11} \sub Z_{11}$ 
is Zariski open.  

We determine the dimension of the stabilizer.  
Let $H=M_{\be_{11}}\cap \gl_8$ and 
$H_{R_{11}}$ be the stabilizer of $R_{11}$ 
as a group scheme and 
$\mathrm{T}_e(H_{R_{11}})$ its tangent space at the unit element 
$e=(1,I_4,I_2)$. Since $\dim G = 30$ and $\dim V=21$, 
it is enough to prove that $\dim \mathrm{T}_e(H_{R_{11}})=9$ by 
Proposition \ref{prop:dimensions-compatible}. 

Let $k[\ep]/(\ep^2)$ be the ring of dual numbers. 
Elements of $\mathrm{T}_e(G)$ can be expressed as 
$(1+\ep c,I_4+\ep A,I_2+\ep B)$ where 
\begin{equation*}
c\in k, \quad A=(a_{ij})\in \m_5(k),\quad  
B = (b_{ij})\in \m_2(k). 
\end{equation*}
Then $(1+\ep c,I_5+\ep A,I_2+\ep B)\in \mathrm{T}_e(H_{R_{11}})$ 
if and only if the following conditions  
\begin{align*}
& 1+\ep(c+\trace(B)) = 1, \\
& (1+\ep b_{11})[(I+\ep A)X_1(1+\ep {}^t A)] 
+ \ep b_{12}[(I+\ep A)X_2(1+\ep {}^t A)] = X_1, \\
& \ep b_{21}[(I+\ep A)X_1(1+\ep {}^t A)] 
+(1+\ep b_{22})[(I+\ep A)X_2(1+\ep {}^t A)] = X_2
\end{align*}
are satisfied. 
These conditions are equivalent to the 
following conditions
\begin{align*}
& c=-\trace(B), \\
& b_{11}X_1 + AX_1 + X_1{}^t A + b_{12}X_2 = 0, \\
& b_{21}X_1 + AX_2 + X_2{}^t A + b_{22}X_2 = 0.
\end{align*}

Explicitly, the second and the third equations are 
\begin{align*}
& b_{11}
\begin{pmatrix}
0 & 0 & 0 & 0 & 0 \\
0 & 0 & 1 & 0 & 0 \\
0 & -1 & 0 & 0 & 0 \\
0 & 0 & 0 & 0 & 1 \\
0 & 0 & 0 & -1 & 0
\end{pmatrix} 
+ \begin{pmatrix}
a_{11} & a_{12} & a_{13} & a_{14} & a_{15} \\
a_{21} & a_{22} & a_{23} & a_{24} & a_{25} \\
a_{31} & a_{32} & a_{33} & a_{34} & a_{35} \\
a_{41} & a_{42} & a_{43} & a_{44} & a_{45} \\
a_{51} & a_{52} & a_{53} & a_{54} & a_{55} 
\end{pmatrix} 
\begin{pmatrix}
0 & 0 & 0 & 0 & 0 \\
0 & 0 & 1 & 0 & 0 \\
0 & -1 & 0 & 0 & 0 \\
0 & 0 & 0 & 0 & 1 \\
0 & 0 & 0 & -1 & 0
\end{pmatrix} \\
& + \begin{pmatrix}
0 & 0 & 0 & 0 & 0 \\
0 & 0 & 1 & 0 & 0 \\
0 & -1 & 0 & 0 & 0 \\
0 & 0 & 0 & 0 & 1 \\
0 & 0 & 0 & -1 & 0
\end{pmatrix}
\begin{pmatrix}
a_{11} & a_{21} & a_{31} & a_{41} & a_{51} \\
a_{12} & a_{22} & a_{32} & a_{42} & a_{52} \\
a_{13} & a_{23} & a_{33} & a_{43} & a_{53} \\
a_{14} & a_{24} & a_{34} & a_{44} & a_{54} \\
a_{15} & a_{25} & a_{35} & a_{45} & a_{55} 
\end{pmatrix}
+ b_{12}
\begin{pmatrix}
0 & 1 & 0 & 0 & 0 \\
-1 & 0 & 0 & 0 & 0 \\
0 & 0 & 0 & 1 & 0 \\
0 & 0 & -1 & 0 & 0 \\
0 & 0 & 0 & 0 & 0
\end{pmatrix}  = 0, \\
& b_{21} 
\begin{pmatrix}
0 & 0 & 0 & 0 & 0 \\
0 & 0 & 1 & 0 & 0 \\
0 & -1 & 0 & 0 & 0 \\
0 & 0 & 0 & 0 & 1 \\
0 & 0 & 0 & -1 & 0
\end{pmatrix}
+ \begin{pmatrix}
a_{11} & a_{12} & a_{13} & a_{14} & a_{15} \\
a_{21} & a_{22} & a_{23} & a_{24} & a_{25} \\
a_{31} & a_{32} & a_{33} & a_{34} & a_{35} \\
a_{41} & a_{42} & a_{43} & a_{44} & a_{45} \\
a_{51} & a_{52} & a_{53} & a_{54} & a_{55} 
\end{pmatrix} 
\begin{pmatrix}
0 & 1 & 0 & 0 & 0 \\
-1 & 0 & 0 & 0 & 0 \\
0 & 0 & 0 & 1 & 0 \\
0 & 0 & -1 & 0 & 0 \\
0 & 0 & 0 & 0 & 0
\end{pmatrix} \\
& + \begin{pmatrix}
0 & 1 & 0 & 0 & 0 \\
-1 & 0 & 0 & 0 & 0 \\
0 & 0 & 0 & 1 & 0 \\
0 & 0 & -1 & 0 & 0 \\
0 & 0 & 0 & 0 & 0
\end{pmatrix}
\begin{pmatrix}
a_{11} & a_{21} & a_{31} & a_{41} & a_{51} \\
a_{12} & a_{22} & a_{32} & a_{42} & a_{52} \\
a_{13} & a_{23} & a_{33} & a_{43} & a_{53} \\
a_{14} & a_{24} & a_{34} & a_{44} & a_{54} \\
a_{15} & a_{25} & a_{35} & a_{45} & a_{55} 
\end{pmatrix}
+ b_{22}
\begin{pmatrix}
0 & 1 & 0 & 0 & 0 \\
-1 & 0 & 0 & 0 & 0 \\
0 & 0 & 0 & 1 & 0 \\
0 & 0 & -1 & 0 & 0 \\
0 & 0 & 0 & 0 & 0
\end{pmatrix} = 0 
\end{align*}
%

%

By long but straightforward computations, 
\begin{align*}
A & = \begin{pmatrix}
-a_{22} & 0 & -a_{42} & 0 & 0 \\
a_{21} & a_{22} & a_{23} & -a_{53} & a_{43} \\
2a_{53} & 0 & -a_{44} & 0 & -2a_{42} \\
a_{41} & a_{42} & a_{43} & a_{44} & a_{45} \\
0 & 0 & a_{53} & 0 & a_{55} 
\end{pmatrix}, \quad \\
B & = \begin{pmatrix}
-a_{22}+a_{44}=-a_{44}-a_{55} & -a_{42} \\
a_{53} & 0
\end{pmatrix} 
\end{align*}
Since $c$ is determined by $B$, $\dim T_e(H_{R_{11}})=9$.  
\end{proof}

The above lemma and Lemma \ref{lem:empty-criterion} 
imply that we may assume that 
\begin{equation*}
y_{25},y_{29},y_{32},y_{34},
y_{42},y_{44},
y_{48},y_{50},
y_{30},y_{33},y_{35},
y_{40},y_{43},y_{45},
y_{51},y_{54}=0.  
\end{equation*}
Then the 1PS in the table shows that $S_{\be_{11}}=\emptyset$. 

\vskip 10pt
\noindent
(5) $\beta_{12}=\tfrac {5} {312} (-9,-9,-1,-1,-1,7,7,7)$ 

For $x\in Z_{12}$, let
\begin{align*}
A(x) & = 
\begin{pmatrix}
\bbma_{21} & \bbma_{36} \\
-\bbma_{20} & -\bbma_{35} \\ 
\bbma_{19} & \bbma_{34}
\end{pmatrix}
= \begin{pmatrix}
e_{178} & e_{278} \\
-e_{168} & -e_{268} \\
e_{167} & e_{267} 
\end{pmatrix}, \\
B(x) & = 
\begin{pmatrix}
\bbma_{47} & \bbma_{48} & \bbma_{49} \\
-\bbma_{41} & -\bbma_{42} & -\bbma_{43} \\
\bbma_{38} & \bbma_{39} & \bbma_{40},
\end{pmatrix}
=\begin{pmatrix}
e_{456} & e_{457} & e_{458} \\
-e_{356} & e_{357} & e_{358} \\
e_{346} & e_{347} & e_{348}
\end{pmatrix}. 
\end{align*}

Let $g_1,g_2\in\spl_3$. Then 
the action of $g=\diag(I_2,g_1,g_2)\in  M_{\be_{12}}$ 
on $x$ is 
\begin{equation*}
A(gx) = {}^tg_2^{-1} A(x), \;
B(gx) = {}^tg_1^{-1} B(x) \,{}^t g_2.  
\end{equation*}
By Lemma \ref{lem:eliminate-nm-matrix}, 
we may assume that $y_{19},y_{34}=0$.  
Then $\spl_3^{(3,5)}\sub M_{\be_{12}}$
does not change this condition.  
We apply Lemma \ref{lem:eliminate-nm-matrix} 
to the action of this $\spl_3$ on $B(x)$ 
and can make $y_{38},y_{39},y_{41}=0$. 
Then the 1PS in the table shows that 
$S_{\be_{12}}=\emptyset$. 

\vskip 5pt
\noindent
(6) $\beta_{20}=\tfrac {1} {232} (-15,-7,-7,1,1,1,9,17)$  

For $x\in Z_{20}$, let
\begin{equation*}
A(x) = \begin{pmatrix}
y_{15} \\
y_{18} \\
y_{20} 
\end{pmatrix}
= \begin{pmatrix}
x_{148} \\
x_{158} \\
x_{168} 
\end{pmatrix}, \;
B(x) = \begin{pmatrix}
y_{29} & y_{39} \\
y_{32} & y_{42} \\
y_{34} & y_{44}
\end{pmatrix}
= \begin{pmatrix}
x_{247} & x_{347} \\
x_{257} & x_{357} \\
x_{267} & x_{367}
\end{pmatrix}.   
\end{equation*}

Let $g_1\in\spl_3$, $g_2\in\spl_2$ and 
$g=\diag(1,g_2,g_1,1,1)\in M_{\be_{20}}$.  
Then the action of $g$ is 
\begin{equation*}
A(gx) = g_1 A(x), \;
B(gx) = g_1 B(x) \,{}^t\! g_2.
\end{equation*}
By Lemma \ref{lem:eliminate-2m(3-32)}, we may assume that 
$y_{15},y_{18},y_{29}=0$. 

\vskip 5pt
\noindent
(7) $\beta_{21}=\tfrac {3} {440} (-15,-7,1,1,1,1,9,9)$

For $x\in Z_{21}$, let
\begin{equation*}
A(x) = \begin{pmatrix}
y_{25} & y_{26} \\
y_{29} & y_{30} \\
y_{32} & y_{33} \\
y_{34} & y_{35}
\end{pmatrix}
=\begin{pmatrix}
x_{237} & x_{238} \\
x_{247} & x_{248} \\
x_{257} & x_{258} \\
x_{267} & x_{268} 
\end{pmatrix}, \;
B(x) = \begin{pmatrix}
y_{47} \\
-y_{41} \\
y_{38} \\
-y_{37}
\end{pmatrix}
= \begin{pmatrix}
x_{456} \\
- x_{356} \\
x_{346} \\
- x_{345}
\end{pmatrix}. 
\end{equation*}

Let $g_1\in\spl_4$, $g_2\in\spl_2$ and 
$g=\diag(1,1,g_1,g_2)\in M_{\be_{21}}$.  
Then the action of $g$ is 
\begin{equation*}
A(gx) = g_1 A(x) \, {}^t g_2, \;
B(gx) = {}^t\! g_1^{-1} B(x).
\end{equation*}
We apply Lemma \ref{lem:eliminate-nm-matrix} and make 
$y_{25},y_{26},y_{29},y_{30}=0$, but we do not 
eliminate $y_{32}$.  
Then $\spl_2^{(3,4)},\spl_2^{(5,6)}\sub M_{\be_{21}}$
do not change this condition and 
$\lan \mathbbm a_{41},\mathbbm a_{47}\ran$
$\lan \mathbbm a_{37},\mathbbm a_{38}\ran$, 
are the standard \rep s of 
$\spl_2^{(3,4)},\spl_2^{(5,6)}\sub M_{\be_{21}}$ 
respectively. we can make $y_{37},y_{41}=0$ 
by Lemma \ref{lem:eliminate-standard}. 
Therefore, we may assume that 
\begin{equation*}
y_{25},y_{26},y_{29},y_{30},y_{37},y_{41}=0.
\end{equation*}

\vskip 5pt
\noindent
(8) $\beta_{22}=\tfrac {7} {376} (-11,-3,-3,-3,-3,5,5,13)$ 

For $x\in Z_{22}$, let $A(x)$ be the following alternating matrix:
\begin{equation*}
A(x) = \begin{pmatrix}
y_{20} \\ y_{21}
\end{pmatrix}
= \begin{pmatrix}
x_{168} \\
x_{178}
\end{pmatrix}, \;
B(x) = \begin{pmatrix}
0 & y_{26} & y_{30} & y_{33} \\
& 0 & y_{40} & y_{43} \\
&& 0 & y_{49} \\
&&& 0
\end{pmatrix}
= \begin{pmatrix}
0 & x_{238} & x_{248} & x_{258} \\
& 0 & x_{348} & x_{258} \\
&& 0 & x_{458} \\
&&& 0
\end{pmatrix}
\end{equation*}
(the lower triangular part is suppressed). 

Let $g_1\in\spl_4,g_2\in\spl_2$ and 
$g=\diag(1,g_1,g_2,1)\in M_{\be_{22}}$. 
Then the action of $g$ is 
\begin{equation*}
A(gx) = g_2 A(x), \; 
B(gx) = g_1 A(x) \, {}^t g_1.
\end{equation*}
By applying Lemma \ref{lem:eliminate-standard} 
to $A(x)$, we may assume that $y_{21}=0$. Then 
$\spl_4^{(2,5)}\sub M_{\be_{22}}$ 
does not change this condition. 
By Witt's theorem (Theorem \ref{thm:alternating-matrix}), 
we may assume that $y_{30},y_{33},y_{40},y_{43}=0$. 
Therefore, we may assume that 
\begin{equation*}
y_{21},y_{30},y_{33},y_{40},y_{43}=0.
\end{equation*}

\vskip 5pt
\noindent
(9) $\beta_{23}=\tfrac {5} {696} (-17,-9,-1,-1,-1,7,7,15)$ 

For $x\in Z_{23}$, let
\begin{equation*}
A(x) = \begin{pmatrix}
y_{20} \\ y_{21}
\end{pmatrix}
= \begin{pmatrix}
x_{167} \\ x_{178}
\end{pmatrix}, \;
B(x) = 
\begin{pmatrix}
y_{47} & y_{48} \\
-y_{41} & -y_{42} \\
y_{38} & y_{39} 
\end{pmatrix}
= \begin{pmatrix}
x_{456} & x_{457} \\
-x_{356} & -x_{357} \\
x_{346} & x_{347} 
\end{pmatrix}. 
\end{equation*}

Let $g_1\in \spl_3,g_2\in\spl_2$ 
and $g=\diag(1,1,g_1,g_2)\in M_{\be_{23}}$. 
Then 
\begin{equation*}
A(gx) = g_2 A(x), \; 
B(gx) = {}^t\!g_1^{-1} B(x) \,{}^t\! g_2.
\end{equation*}
We first apply Lemma \ref{lem:eliminate-standard} 
to $A(x)$ and make $y_{20}=0$.  Then 
$\spl_3^{(3,5)}$ does not change this condition. 
By applying Lemma \ref{lem:eliminate-nm-matrix}
to $B(x)$, we can make $y_{38},y_{39},y_{41}=0$. 
Therefore, we may assume that 
\begin{equation*}
y_{20},y_{38},y_{39},y_{41}=0. 
\end{equation*}

\vskip 5pt
\noindent
(10) $\beta_{24}=\tfrac {1} {504} (-13,-5,-5,-5,3,3,11,11)$ 

For $x\in Z_{24}$, let
\begin{align*}
A(x) & = \begin{pmatrix}
y_{39} & y_{40} \\
-y_{29} & -y_{30} \\
y_{25} & y_{26} 
\end{pmatrix}
= \begin{pmatrix}
x_{347} & x_{348} \\
-x_{247} & -x_{248} \\
x_{237} & x_{238} 
\end{pmatrix}, \\
B(x) & = \begin{pmatrix}
y_{31} \\
y_{41} \\
y_{47}
\end{pmatrix}
= \begin{pmatrix}
x_{256} \\
x_{356} \\
x_{456}
\end{pmatrix}, \;
C(x) = 
\begin{pmatrix}
y_{17} & y_{18} \\
y_{19} & y_{20}
\end{pmatrix}
= \begin{pmatrix}
x_{157} & x_{158} \\
x_{167} & x_{168}
\end{pmatrix}, \;
\end{align*}

Let $g_1\in \spl_3,g_2,g_3\in\spl_2$ and 
$g=\diag(1,g_1,g_2,g_3)\in M_{\be_{24}}$. 
Then 
\begin{equation*}
A(gx) = {}^t g_1^{-1} A(x) \, {}^t\!g_3, \;
B(gx) = g_1 B(x)
C(gx) = g_2 \, C(x) \, {}^t\! g_3. 
\end{equation*}
We first apply Lemma \ref{lem:eliminate-nm-matrix}
to $A(x)$ and make $y_{25},y_{26}=0$. Then 
$\spl_2^{(2,3)}\sub M_{\be_{24}}$ 
does not change this condition. 
So by applying Lemma \ref{lem:eliminate-standard} 
to the subspace $\lan \mathbbm a_{31},\mathbbm a_{41}\ran$, 
we can make $y_{31}=0$. 
Then $\spl_2^{(5,6)}\sub M_{\be_{24}}$ 
does not change the above conditions. So we can apply 
Lemma \ref{lem:eliminate-standard} to the subspace 
$\lan \mathbbm a_{17},\mathbbm a_{19}\ran$ and make 
$y_{17}=0$. Therefore, we may assume that 
\begin{equation*}
y_{17},y_{25},y_{26},y_{31}=0. 
\end{equation*}

\vskip 5pt
\noindent
(11) $\beta_{25}=\tfrac {3} {824} (-15,-7,-7,-7,1,9,9,17)$ 

For $x\in Z_{25}$, let
\begin{equation*}
A(x) = 
\begin{pmatrix}
y_{40} \\
-y_{30} \\
y_{26}
\end{pmatrix}
= \begin{pmatrix}
x_{348} \\
-x_{248} \\
x_{238}
\end{pmatrix}, \;
B(x) = 
\begin{pmatrix}
y_{31} & y_{32} \\
y_{41} & y_{42} \\
y_{47} & y_{48}
\end{pmatrix}
= \begin{pmatrix}
x_{256} & x_{257} \\
x_{356} & x_{357} \\
x_{456} & x_{457}
\end{pmatrix}.  
\end{equation*}

Let $g_1\in\spl_3,g_2\in\spl_2$ and 
$g=\diag(1,g_1,1,g_2,1)\in M_{\be_{25}}$. 
Then 
\begin{equation*}
A(gx) = {}^tg_1^{-1} A(x), \; 
B(gx) = g_1 B(x) \, {}^t g_2. 
\end{equation*}
Lemma \ref{lem:eliminate-2m(3-32)} 
for $(A(x),B(x))$ implies that 
we can make  and make $y_{26},y_{30},y_{42}=0$.

\vskip 5pt
\noindent
(12) $\beta_{27}=\tfrac {5} {56} (-1,-1,-1,-1,-1,-1,-1,7)$ 

In this case, 
$Z_{27}$ can be identified with $\wedge^2 \aff^7$.  
Let $g_1\in\gl_7$ and $g=\diag(g_1,1)\in M_{\be_{27}}$.
If we regard $x$ as a $7\times 7$ alternating matrix
with diagonal entries $0$ then the action of 
$g$ is $x\mapsto g_1 x \, {}^t\! g_1$. 
By Theorem \ref{thm:alternating-matrix}, 
we can make the first row
and the first column $0$. 
Therefore, we may assume that $x_{128}\ccd x_{178}=0$, i.e., 
\begin{equation*}
y_6,y_{11},y_{15},y_{18},y_{20},y_{21}=0.
\end{equation*}

\vskip 5pt
\noindent
(13) $\beta_{29}=\tfrac {13} {568} (-9,-9,-9,-1,-1,7,7,15)$ 

For $x\in Z_{29}$, let 
\begin{equation*}
A(x) = 
\begin{pmatrix}
y_{20} & y_{21} \\
y_{35} & y_{36} \\
y_{45} & y_{46}
\end{pmatrix}
=\begin{pmatrix}
x_{168} & x_{178} \\
x_{268} & x_{278} \\
x_{368} & x_{378} 
\end{pmatrix}, \;
B(x) = \begin{pmatrix}
y_{50} \\
y_{53}
\end{pmatrix}
= \begin{pmatrix}
x_{467} \\
x_{567}
\end{pmatrix}. 
\end{equation*}

Let $g_1\in\spl_3,g_2,g_3\in\spl_2$ and 
$g=\diag(g_1,g_2,g_3,1)\in M_{\be_{29}}$. 
Then 
\begin{equation*}
A(gx) = g_1 A(x) \, {}^t\!g_3, \;
B(gx) = g_2 B(x). 
\end{equation*}
We first apply Lemma \ref{lem:eliminate-nm-matrix} 
to $A(x)$ and make $y_{45},y_{46}=0$. Then 
$\spl_2^{(4,5)}\sub M_{\be_{29}}$ 
does not change this conditions. So we apply 
Lemma \ref{lem:eliminate-standard} to 
$B(x)$ and make $y_{50}=0$. Therefore, we may assume that 
\begin{equation*}
y_{45},y_{46},y_{50}=0. 
\end{equation*}

\vskip 5pt
\noindent
(14) $\beta_{36}=\tfrac {1} {888} (-13,-13,-5,-5,3,3,11,19)$ 

For $x\in Z_{36}$, let 
\begin{align*}
A(x) & = \begin{pmatrix}
y_{11} & y_{15} \\
y_{26} & y_{30}
\end{pmatrix}
= \begin{pmatrix}
x_{138} & x_{148} \\
x_{238} & x_{248} 
\end{pmatrix}, \;
B(x) = \begin{pmatrix}
y_{41} \\ y_{47}
\end{pmatrix}
= \begin{pmatrix}
x_{356} \\ x_{456}
\end{pmatrix}, \\
C(x) & = 
\begin{pmatrix}
y_{17} & y_{19} \\
y_{32} & y_{34}
\end{pmatrix}
=\begin{pmatrix}
x_{157} & x_{167} \\
x_{257} & x_{267}
\end{pmatrix}. 
\end{align*}

Let $g_1,g_2,g_3\in\gl_2$ and 
$g=\diag(g_1,g_2,g_3,1,1)\in M_{\be_{36}}$. 
Then 
\begin{equation*}
A(gx) = g_1 A(x) \,{}^t\! g_2,\;
B(gx) = g_2 B(x), \; 
C(gx) = g_1 C(x) \,{}^t\! g_3.  
\end{equation*}
We first apply Lemma \ref{lem:eliminate-2m(2)}  
to $(A(x),B(x))$ and make $y_{15},y_{47}=0$. 
Then $\spl_2^{(5,6)}\sub M_{\be_{36}}$ 
does not change this conditions. 
So we apply Lemma \ref{lem:eliminate-standard} 
to the action of $\spl_2^{(5,6)}$ to 
$C(x)$ and make $y_{19}=0$. 
Therefore, we may assume that 
\begin{equation*}
y_{15},y_{19},y_{47}=0.
\end{equation*}

\vskip 5pt
\noindent
(15) $\beta_{38}=\tfrac {1} {760} (-21,-5,-5,3,3,3,11,11)$ 

For $x\in Z_{38}$, let 
\begin{equation*}
A(x) = \begin{pmatrix}
y_{25} \\ y_{26}
\end{pmatrix}
= \begin{pmatrix}
x_{237} \\ x_{238}
\end{pmatrix}, \;
B(x) = \begin{pmatrix}
y_{31} & -y_{28} & y_{27} \\
y_{41} & -y_{38} & y_{37}
\end{pmatrix}
= \begin{pmatrix}
x_{256} & -x_{246} & x_{245} \\ 
x_{356} & -x_{346} & x_{345}  
\end{pmatrix}.
\end{equation*}

Let $g_1\in\spl_3,g_21\in\spl_2$ and 
$g=\diag(1,I_2,g_1,g_2)\in M_{\be_{38}}$. 
Then 
\begin{equation*}
A(gx) = g_2 A(x), \;
B(gx) = B(x) \,g_1^{-1}.
\end{equation*}
We apply Lemma \ref{lem:eliminate-nm-matrix} 
to B(x) and make $y_{27},y_{37}=0$. 
Then $\spl_2^{(7,8)}\sub M_{\be_{38}}$ 
does not change this conditions. 
So we Lemma \ref{lem:eliminate-standard} 
to $A(x)$ and make $y_{25}=0$. Therefore, 
we may assume that
\begin{equation*}
y_{25}=y_{27},y_{37}=0. 
\end{equation*}

\vskip 5pt
\noindent
(16) $\beta_{49}=\tfrac {1} {88} (-5,-5,-1,-1,3,3,3,3)$ 

For $x\in Z_{49}$, let 
\begin{equation*}
v(x) = [37,38,39,40]
= [x_{345},x_{346},x_{347},x_{348}]
\end{equation*}
and 
\begin{align*}
& A_1(x) = \begin{pmatrix}
0 & y_{16} & y_{17} & y_{18} \\
-y_{16} & 0 & y_{19} & y_{20} \\
-y_{17} & -y_{19} & 0 & y_{21} \\
-y_{18} & -y_{20} & -y_{21} & 0
\end{pmatrix}
= \begin{pmatrix}
0 & x_{156} & x_{157} & x_{158} \\
-x_{156} & 0 & x_{167} & x_{168} \\
-x_{157} & -x_{167} & 0 & x_{178} \\
-x_{158} & -x_{168} & -x_{178} & 0 
\end{pmatrix}, \\ 
& A_2(x) = \begin{pmatrix}
0 & y_{31} & y_{32} & y_{33} \\
-y_{31} & 0 & y_{34} & y_{35} \\
-y_{32} & -y_{34} & 0 & y_{36} \\
-y_{33} & -y_{35} & -y_{36} & 0
\end{pmatrix}
= \begin{pmatrix}
0 & x_{256} & x_{257} & x_{258} \\
-x_{256} & 0 & x_{267} & x_{268} \\
-x_{257} & -x_{267} & 0 & x_{278} \\
-x_{258} & -x_{268} & -x_{278} & 0 
\end{pmatrix}. 
\end{align*}
Let 
\begin{equation*}
R_{49} = e_{178} + e_{256} + e_{346} + e_{348}. 
\end{equation*}
Then $v(R_{49}),A_1(R_{49}),A_2(R_{49})$ 
are equal to $v(R_2),A_1(R_2),A_2(R_2)$ 
in (\ref{eq:vR2etc}).  

By almost the same computations as in 
the proof of Proposition \ref{prop:4241}, 
$M_{\be_{49}}\cdot R_{49}\sub Z_{49}$ is Zariski 
open. Therefore, we may assume that 
\begin{equation*}
y_{16},y_{17},y_{18},y_{19},y_{20},
y_{32},y_{33},y_{34},y_{35},y_{36},
y_{37},y_{39}=0 
\end{equation*}
and the 1PS in the table shows that 
$S_{\be_{49}}=\emptyset$.  

\vskip 5pt
\noindent
(17) $\beta_{51}=\tfrac {1} {12} (-2,-2,-2,0,1,1,1,3)$

\vskip 5pt
\noindent
(18) $\beta_{55}=\tfrac {1} {24} (-9,0,1,1,1,2,2,2)$

For $x\in Z_{55}$, let 
\begin{align*}
A(x) & = \begin{pmatrix}
y_{36} \\ -y_{35} \\ y_{34}
\end{pmatrix}
= \begin{pmatrix}
x_{278} \\ -x_{268} \\ x_{267}
\end{pmatrix}, \\
B(x) & = \begin{pmatrix}
y_{47} & y_{48} & y_{49} \\
-y_{41} & -y_{42} & -y_{43} \\
y_{38} & y_{39} & y_{40}
\end{pmatrix}
= \begin{pmatrix}
x_{456} & x_{457} & x_{458} \\
-x_{356} & -x_{357} & -x_{358} \\ 
x_{346} & x_{347} & x_{348} 
\end{pmatrix}.
\end{align*}

Let $g_1,g_2\in\spl_3$ and 
$g=\diag(1,1,g_1,g_2)\in M_{\be_{55}}$. 
Then 
\begin{equation*}
A(gx) = {}^t g_2^{-1} A(x), \; 
B(gx) = {}^t g_1^{-1} B(x) \, {}^t\! g_2.  
\end{equation*}
We apply Lemma \ref{lem:eliminate-standard} 
to $A(x)$ and make $y_{34},y_{35}=0$. 
Then $\spl_3^{(3,5)}\sub M_{\be_{55}}$ 
does not change this condition. 
Then we apply Lemma \ref{lem:eliminate-nm-matrix} 
($n=m=3$) to $B(x)$ and make 
$y_{38},y_{39},y_{41}=0$.
Therefore, we may assume that 
\begin{equation*}
y_{34},y_{35},y_{38},y_{39},y_{41}=0. 
\end{equation*}

\vskip 5pt
\noindent
(19) $\beta_{56}=\tfrac {1} {6} (-1,-1,-1,0,0,1,1,1)$ 

For $x\in Z_{56}$, let
\begin{equation*}
A(x) = \begin{pmatrix}
y_{21} & -y_{20} & y_{19} \\
y_{36} & -y_{35} & y_{34} \\
y_{46} & -y_{45} & y_{44}
\end{pmatrix}
= \begin{pmatrix}
x_{178} & -x_{168} & x_{167} \\
x_{278} & -x_{268} & x_{267} \\
x_{378} & -x_{368} & x_{367} 
\end{pmatrix},\;
v(x) = \begin{pmatrix}
y_{47} \\
y_{48} \\
y_{49}
\end{pmatrix}
= \begin{pmatrix}
x_{456} \\
x_{457} \\
x_{458} 
\end{pmatrix}. 
\end{equation*}

Let $g_1,g_2\in\spl_3$ and 
$g=\diag(g_1,1,1,g_2)\in M_{\be_{56}}$.  
Then 
\begin{equation*}
A(gx) = g_1 A(x) g_2^{-1}, \;
v(gx) = g_2 v(x). 
\end{equation*}
We first apply Lemma \ref{lem:eliminate-standard} to 
$v(x)$ and make $y_{47},y_{48}=0$. 
Then the action of $\spl_3^{(1,3)}\sub M_{\be_{56}}$ 
does not change this condition.  So 
we apply Lemma \ref{lem:eliminate-nm-matrix} 
to $A(x)$ and make $y_{19},y_{20},y_{34}=0$.  
Therefore, we may assume that 
\begin{equation*}
y_{19},y_{20},y_{34},y_{47},y_{48}=0. 
\end{equation*}

\vskip 5pt
\noindent
(20) $\beta_{57}=\tfrac {1} {40} (-15,1,1,1,1,1,5,5)$ 

For $x\in Z_{57}$, let $A(x)$ be the following alternating matrix:
\begin{align*}
& A_1(x) = \begin{pmatrix}
0 & y_{25} & y_{26} & y_{29} & y_{30} \\
& 0 & y_{32} & y_{33} & y_{34} \\
&& 0 & y_{35} & y_{39} \\
&&& 0 & y_{40} \\
&&&& 0  
\end{pmatrix}
= \begin{pmatrix}
0 & x_{237} & x_{247} & x_{257} & x_{267} \\
& 0 & x_{347} & x_{357} & x_{367} \\
&& 0 & x_{457} & x_{467} \\
&&& 0 & x_{567} \\
&&&& 0  
\end{pmatrix}, \\ 
& A_2(x) = \begin{pmatrix}
0 & y_{42} & y_{43} & y_{44} & y_{45} \\
& 0 & y_{48} & y_{49} & y_{50} \\
&& 0 & y_{51} & y_{53} \\
&&& 0 & y_{54} \\
&&&& 0  
\end{pmatrix}
= \begin{pmatrix}
0 & x_{238} & x_{248} & x_{258} & x_{268} \\
& 0 & x_{348} & x_{358} & x_{368} \\
&& 0 & x_{458} & x_{468} \\
&&& 0 & x_{568} \\
&&&& 0  
\end{pmatrix} 
\end{align*}
(the lower triangular part is suppressed). 
Let 
\begin{equation*}
R_{57} = e_{347}+e_{567}+e_{238}+e_{458}. 
\end{equation*}
Then $A_1(R_{57}),A_2(R_{57})$ are equal to 
$X_1,X_2$ in (\ref{eq:X1X2-defn}).  
By almost the same computations as in the proof 
of Proposition \ref{prop:522}, 
$M_{\be_{57}}\cdot R_{57}\sub Z_{57}$ is Zariski open. 
Therefore, Lemma \ref{lem:empty-criterion} 
implies that we may assume that 
\begin{equation*}
y_{25},y_{29},y_{32},y_{34},
y_{42},y_{44},
y_{48},y_{50},
y_{30},y_{33},y_{35},
y_{40},y_{43},y_{45},
y_{51},y_{54}=0.  
\end{equation*}
Then the 1PS in the table shows that 
$S_{\be_{57}}=\emptyset$. 

\vskip 5pt
\noindent
(21) $\beta_{61}=\tfrac {1} {184} (-13,-13,-13,-5,-5,3,19,27)$ 

For $x\in Z_{61}$, let
\begin{equation*}
A(x) = \begin{pmatrix}
\bbma_{15} & \bbma_{18} \\
\bbma_{30} & \bbma_{33} \\
\bbma_{40} & \bbma_{43}
\end{pmatrix}
= \begin{pmatrix}
x_{148} & x_{158} \\
x_{248} & x_{258} \\
x_{348} & x_{358} 
\end{pmatrix}, \;
v(x) = \begin{pmatrix}
x_{167} \\
x_{267} \\
x_{367}
\end{pmatrix}
= \begin{pmatrix}
\bbma_{19} \\ 
\bbma_{34} \\ 
\bbma_{44}
\end{pmatrix}. 
\end{equation*}
Let $g_1\in \spl_3,g_2\in\spl_2$ and 
$g=\diag(g_1,g_2,1,1,1)\in M_{\be_{61}}$. Then 
\begin{equation*}
A(gx) = g_1 A(x) \, {}^t\! g_2, \; 
v(gx) = g_1 v(x). 
\end{equation*}
Lemma \ref{lem:eliminate-2m(3-32)} implies that 
we may assume that 
\begin{equation*}
y_{19},y_{33},y_{34}=0. 
\end{equation*}

\vskip 5pt
\noindent
(22) $\beta_{62}=\tfrac {1} {48} (-5,-1,-1,-1,1,1,3,3)$ 

Similarly as in the case (21), 
Lemma \ref{lem:eliminate-2m(3-32)} implies that we may assume that 
\begin{equation*}
y_{29},y_{41},y_{47}=0. 
\end{equation*}

\vskip 5pt
\noindent
(23)  $\beta_{66}=\tfrac {1} {40} (-15,-1,1,1,3,3,3,5)$ 

For $x\in Z_{66}$, let
\begin{equation*}
v(x) = \begin{pmatrix}
y_{33} \\ y_{35} \\ y_{36}
\end{pmatrix}
= \begin{pmatrix}
x_{258} \\ x_{268} \\ x_{278}
\end{pmatrix}, \; 
A(x) = \begin{pmatrix}
y_{44} & -y_{42} & y_{41} \\
y_{50} & -y_{48} & y_{47}
\end{pmatrix}
= \begin{pmatrix}
x_{367} & -x_{357} & x_{356} \\
x_{467} & -x_{457} & x_{456}
\end{pmatrix}. 
\end{equation*}

Let $g_1,\in\spl_3,g_2\in\spl_2$ and 
$g=\diag(1,1,g_2,g_1,1)\in M_{\be_{66}}$. 
Then 
\begin{equation*}
v(gx) = g_1 v(x), \; 
A(gx) = g_2 B(x) g_1^{-1}. 
\end{equation*}
We apply Lemma \ref{lem:eliminate-nm-matrix} 
to $A(x)$ and make $y_{41},y_{47}=0$. 
Then $\spl_2^{(5,6)}\sub M_{\be_{66}}$ 
does not change this condition.  
By applying Lemma \ref{lem:eliminate-standard} 
to the subspace $\lan y_{33},y_{35}\ran$ 
and we can make $y_{33}=0$. Therefore, we may assume that 
\begin{equation*}
y_{33},y_{41},y_{47}=0. 
\end{equation*}

\vskip 5pt
\noindent
(24) $\beta_{73}=\tfrac {1} {48} (-5,-2,0,0,0,1,3,3)$ 

For $x\in Z_{73}$, let
\begin{equation*}
A(x) = \begin{pmatrix}
y_{25} & y_{26} \\
y_{29} & y_{30} \\
y_{32} & y_{33}
\end{pmatrix}
= \begin{pmatrix}
x_{237} & x_{238} \\
x_{247} & x_{248} \\
x_{257} & x_{258}
\end{pmatrix}, \;
v(x) = \begin{pmatrix}
y_{47} \\
-y_{41} \\
y_{38}
\end{pmatrix}
= \begin{pmatrix}
x_{456} \\
-x_{356} \\
x_{346} 
\end{pmatrix}.
\end{equation*}
Let $g_1\in\spl_3,g_2\in\spl_2$ and 
$g=\diag(1,1,g_1,1,g_2)\in M_{\be_{73}}$. 
Then 
\begin{equation*}
A(gx) = g_1 A(x) \, {}^t\!g_2, \; v(gx) = {}^t g_1^{-1} v(x).
\end{equation*}
We first apply 
Lemma \ref{lem:eliminate-nm-matrix} 
to $A(x)$ and make $y_{25},y_{26}=0$. Then 
the $\spl_2^{(4,5)}\sub M_{\be_{73}}$ 
does not change this condition.  
By applying Lemma \ref{lem:eliminate-standard} 
to the subspace $\lan \bbma_{38},\bbma_{41}\ran$, 
we can make $y_{38}=0$.  
Then $\spl_2^{(7,8)}\sub M_{\be_{73}}$ 
does not change this condition. 
Since this $\spl_2$ acts on the subspace 
$\lan\mathbbm a_{29},\mathbbm a_{30}\ran$, 
we can make $y_{29}=0$. Therefore, we may assume that 
\begin{equation*}
y_{25},y_{26},y_{29},y_{38}=0.  
\end{equation*}

\vskip 5pt
\noindent
(25) $\beta_{88}=\tfrac {1} {104} (-5,-5,-3,1,1,3,3,5)$ 

By applying Lemma \ref{lem:eliminate-2m(2)} 
to the subspace 
\begin{math}
\lan\mathbbm a_{15},\mathbbm a_{33}\ran\oplus 
\lan\mathbbm a_{19},\mathbbm a_{34}\ran, 
\end{math}
we can make $y_{18},y_{19}=0$.
Then $\spl_2^{(6,7)}\sub M_{[2,3,5,7]}=M_{\be_{88}}$ 
does not change this condition and acts on the 
subspace $\lan \mathbbm a_{41},\mathbbm a_{42}\ran$. 
So we can make $y_{41}=0$ by Lemma \ref{lem:eliminate-standard}. 
Therefore, we may assume that 
\begin{equation*}
y_{18},y_{19},y_{41}=0.
\end{equation*}

\vskip 5pt
\noindent
(26)  $\beta_{121}=\tfrac {1} {136} (-51,-11,-11,5,5,21,21,21)$

Similarly as in the case (21), 
Lemma \ref{lem:eliminate-2m(3-32)} implies that we may assume that 
\begin{equation*}
y_{47},y_{48},y_{34}=0.
\end{equation*}

\vskip 5pt
\noindent
(27) $\beta_{125}=\tfrac {1} {40} (-15,-15,1,1,1,1,1,25)$ 

For $x\in Z_{125}$, let $A(x)$ be the following 
alternating matrix:
\begin{align*}
A(x) = \begin{pmatrix}
0 & y_{40} & y_{43} & y_{45} & y_{46} \\
& 0 & y_{49} & y_{51} & y_{52} \\
&& 0 & y_{54} & y_{55} \\
&&& 0 & y_{56} \\
&&&& 0 
\end{pmatrix}
= \begin{pmatrix} 0 & x_{348} & x_{358} & x_{368}  & x_{378} \\
& 0 & x_{458} & x_{468}  &  x_{478} \\
&& 0 & x_{568} & x_{578} \\ 
&&&  0 & x_{678} \\
&&&& 0
\end{pmatrix}
\end{align*}
(the lower triangular part is suppressed). 
We identify $x$ with $A(x)$. 

Let $g_1\in \spl_5$ and $g=\diag(1,1,g_1,1)$. Then 
\begin{math}
A(gx) = g_1 A(x) \,{}^t\! g_1.
\end{math}
By Witt's theorem (Theorem \ref{thm:alternating-matrix}), 
we can make the first row and the first column, 
$(2,4)$, $(2,5)$, $(3,4)$, $(3,5)$-entries $0$. 
Therefore, we may assume that 
\begin{equation*}
y_{40},y_{43},y_{45},y_{46},y_{51},y_{52},y_{54},y_{55}=0. 
\end{equation*}

\vskip 5pt
\noindent
(28)  $\beta_{130}=\tfrac {3} {40} (-5,-5,-5,3,3,3,3,3)$ 

For $x\in Z_{130}$, let $A(x)$ be the following alternating matrix:
\begin{align*}
A(x) = 
\begin{pmatrix}
0 & y_{56} & -y_{55} & y_{54} & -y_{53} \\
& 0 & y_{52} & -y_{51} & y_{50} \\
&& 0 & -y_{49} & y_{48} \\
&&& 0 & -y_{47} \\
&&&& 0
\end{pmatrix}
= \begin{pmatrix} 
0 & x_{678} & -x_{578} & x_{568}  & - x_{567} \\
& 0 & x_{478} & -x_{468} &  x_{467} \\
&& 0  & -x_{458} & x_{457} \\
&&& 0 & -x_{456} \\
&&&& 0
\end{pmatrix}
\end{align*}
(the lower triangular part is suppressed).  
Since $\wedge^3 \aff^5$ can be identified with the 
dual space of $\wedge^2 \aff^5$, we identify 
$x\in Z_{130}$ with $A(x)$. 

Let $g_1\in\spl_5$ and $g=\diag(I_3,g_1)\in M_{\be_{130}}$. 
Then 
\begin{math}
A(gx) = {}^t g_1^{-1} A(x) g_1^{-1}. 
\end{math}
%
By Witt's theorem (Theorem \ref{thm:alternating-matrix}), 
we may assume that 
%
\begin{equation*}
y_{48},y_{49},y_{50},y_{51},y_{53},y_{54},y_{55},y_{56}=0.
\end{equation*}

\bibliographystyle{plain} 
\bibliography{ref4} 

\end{document}